\documentclass[11pt,reqno]{amsart}
\usepackage[margin=1.18in]{geometry}
\usepackage{amsmath,amssymb,amsfonts,mathtools,mathrsfs}
\usepackage{enumitem}
\usepackage{comment}
\usepackage{microtype}
\usepackage{booktabs,array,graphicx}
\usepackage[dvipsnames]{xcolor}
\usepackage[colorlinks=true,citecolor=MidnightBlue,linkcolor=BrickRed,
  urlcolor=RoyalBlue]{hyperref}
\newcommand{\scalarshelean}{%
  \href{https://github.com/hanyi162013-Yihan/scalar-she-lean}{scalar-she-lean}}
\numberwithin{equation}{section}
\newtheorem{theorem}{Theorem}[section]
\newtheorem{lemma}[theorem]{Lemma}
\newtheorem{proposition}[theorem]{Proposition}
\newtheorem{corollary}[theorem]{Corollary}
\theoremstyle{definition}
\newtheorem{candidategeneralization}[theorem]{Candidate generalization}
\newtheorem{definition}[theorem]{Definition}

\theoremstyle{remark}\newtheorem{remark}[theorem]{Remark}
\title[Weak uniqueness with nondegenerate diffusion]
{Weak uniqueness for the stochastic heat equation with nondegenerate diffusion:
the full H\"older regime}
\author{Yi Han}
\thanks{Department of Mathematics, Massachusetts Institute of Technology,
Cambridge, MA.\newline Email: \texttt{hanyi16@mit.edu}.}
\thanks{The accompanying Lean code is available in the \scalarshelean{} repository.
See \hyperref[sec:ai-lean]{Section~\getrefnumber{sec:ai-lean}}.}
\date{}
\hypersetup{
  pdftitle={Weak uniqueness for the stochastic heat equation with nondegenerate diffusion: the full Hölder regime},
  pdfauthor={Yi Han}
}
\begin{document}
\begin{abstract}
We prove weak well-posedness for the scalar stochastic heat equation
on the one-dimensional torus with uniformly nondegenerate diffusion
that is H\"older continuous of any positive order.
This resolves the weak uniqueness problem in the full positive
H\"older range, extending our earlier result for exponents above
two thirds~\cite{Han2026Threshold}.
The proof couples a weak solution to a smooth companion and uses spatial cancellation
to estimate the signed covariance difference in a negative Sobolev
norm. A finite Gaussian expansion on each Littlewood--Paley block
gives the decay needed to remove the coupling drift by a change of
measure. A localized Girsanov argument extends the result to
measurable drifts of linear growth.
\end{abstract}
\maketitle
\enlargethispage{2pt}
\tableofcontents
\clearpage

\section{Introduction}
\label{sec:main-results-and-mechanism}
Let $\mathbb T=\mathbb R/\mathbb Z$ be the one-dimensional torus with
normalized Lebesgue measure. We consider the stochastic heat equation
\[
 du=\Delta u\,dt+g(u)\,dW,\qquad u(0)=u_0,
\]
where $\Delta=\partial_x^2$, $W$ is space-time white noise, and
$g:\mathbb R\to\mathbb R$ acts pointwise on the solution. We write
$H=L^2(\mathbb T;\mathbb R)$ and $S_t=e^{t\Delta}$. By weak uniqueness
we mean uniqueness of the law of the solution on $C([0,T];H)$ for
each finite $T$.

For a drift-free diffusion in finite dimensions, continuity and
uniform ellipticity of the covariance suffice for uniqueness in law by the
martingale-problem theory of Stroock and Varadhan~\cite{StroockVaradhan2006}.
Existing infinite-dimensional analogues impose additional covariance
structure and leave open the general pointwise model. For
one-dimensional space-time white noise, Mytnik and
Perkins~\cite{MytnikPerkins2011} proved pathwise uniqueness for
H\"older exponents above $3/4$ on the real line, allowing the diffusion
to vanish.
The counterexamples below $3/4$ of Mueller, Mytnik and
Perkins~\cite{MuellerMytnikPerkins2014} use coefficients that do vanish.
This leaves the question of what regularity is required for weak
uniqueness when the diffusion is uniformly nondegenerate.

In~\cite{Han2026Threshold} we proved weak uniqueness for exponents
above $2/3$, giving the first general result below $3/4$ for the
uniformly nondegenerate pointwise model. A natural conjecture was
that weak uniqueness should hold for every positive H\"older exponent.
In that work we suggested that space-time cancellation could help
reach smaller exponents. The following theorem confirms the conjecture
for bounded diffusion coefficients.

\begin{theorem}[Weak existence and uniqueness at every positive H\"older exponent]
\label{thm:main-weak-uniqueness}
Let $g:\mathbb R\to\mathbb R$ satisfy, for some $\beta_0>0$,
$L_g <\infty$ and $0<c_0\le C_0<\infty$,
\[
 c_0\le |g(z)|\le C_0,\qquad
 |g(z)-g(z')|\le L_g |z-z'|^{\beta_0},
 \qquad z,z'\in\mathbb R.
\]
For every deterministic $u_0\in L^2(\mathbb T)$ and every $T>0$,
a weak mild solution of the stochastic heat equation above with
continuous $L^2(\mathbb T)$ paths exists, and its law on
$C([0,T];L^2(\mathbb T))$ is unique.
\end{theorem}

The next theorem allows measurable drifts of linear growth and
also treats unbounded diffusion for bounded initial data.

\begin{theorem}[Scalar unbounded diffusion and measurable drift]
\label{thm:main-extensions}
Let $g:\mathbb R\to\mathbb R$ satisfy $|g(z)|\ge c_0>0$ for every $z$.
Let $b:[0,\infty)\times\mathbb T\times\mathbb R\to\mathbb R$ be Borel
measurable and, for every $T<\infty$, suppose that
$|b(t,x,z)|\le L_{b,T}(1+|z|)$ for $0\le t\le T$.
Write $b(t,u)(x)=b(t,x,u(x))$. The scalar equation
\begin{equation}
 du=[\Delta u+b(t,u)]\,dt+g(u)\,dW,\qquad u(0)=u_0,
 \label{eq:main-drifted-system}
\end{equation}
has a weak mild solution and a unique law on $C([0,T];H)$ for every
$T>0$ and every deterministic $u_0$ in either of the following cases:
\begin{enumerate}[label=(\roman*),leftmargin=*]
\item $g$ is bounded and globally $\beta_0$-H\"older for some
$\beta_0>0$, and $u_0\in L^2(\mathbb T)$.
\item $g$ is globally $\beta$-H\"older for some $0<\beta\le1$,
with no uniform upper bound required, and $u_0\in L^\infty(\mathbb T)$.
\end{enumerate}
\end{theorem}

The bounded initial data in Theorem~\ref{thm:main-extensions}(ii) are
required by our localization proof of weak uniqueness for unbounded
diffusion. They ensure that the solution stays below a fixed spatial
supremum-norm cutoff for a positive time. See
Section~\ref{sec:unbounded-coefficient-linf-localization}. A finer
argument near time zero may allow us to replace this hypothesis by
$u_0\in L^2(\mathbb T)$. We leave this extension for future work.

Theorems~\ref{thm:main-weak-uniqueness} and~\ref{thm:main-extensions}
use classical weak existence for the drift-free equation with continuous
diffusion of at most linear growth
\cite{GatarekGoldys1994,DaPratoZabczyk2014}. In
Theorem~\ref{thm:main-extensions}, a localized Girsanov argument
then gives existence with the measurable drift. See
Section~\ref{sec:scalar-extensions}. The main contribution of this
paper is weak uniqueness for the drift-free equation with uniformly
nondegenerate diffusion coefficients of any positive H\"older order.

\subsection{Weak uniqueness and nondegeneracy}
\label{sec:literature-review}
Our results identify the law of the solution path. Pathwise uniqueness
asks a different question: whether solutions with the same initial
state and the same noise on a common probability space agree almost
surely. Strong well-posedness, which requires strong existence and
pathwise uniqueness, is outside the scope of this paper. For these
notions and their relation to uniqueness of the joint law of the
solution and the noise, see
\cite{Ondrejat2004,CriensRitter2022}.

We begin with pathwise uniqueness results for H\"older diffusion coefficients. In one-dimensional SDEs, the
Yamada--Watanabe criterion~\cite{YamadaWatanabe1971} gives pathwise
uniqueness for $1/2$-H\"older diffusion coefficients. For the
stochastic heat equation driven by space-time white noise, Mytnik
and Perkins~\cite{MytnikPerkins2011} proved pathwise uniqueness
above $3/4$ on the real line. Below $3/4$, the example
$g(z)=|z|^\gamma$ with $0<\gamma<3/4$ admits both a zero and a
nonzero solution from zero initial data
\cite{MuellerMytnikPerkins2014}. Thus allowing the coefficient
to vanish can obstruct even uniqueness in law. For spatially colored
noise, Mytnik, Perkins and Sturm~\cite{MytnikPerkinsSturm2006} relate
the regularity required for pathwise uniqueness to the spatial
covariance.

The square-root equation $du=\Delta u\,dt+\sqrt{u}\,dW$, with
$u\ge0$, is one of the few classical examples of weak uniqueness
below $3/4$ for the stochastic heat equation
with pointwise diffusion and space-time white noise. Its coefficient
$g(z)=\sqrt z$ is $1/2$-H\"older on $[0,\infty)$. Among nonnegative
solutions, uniqueness in law comes from the branching structure of
one-dimensional super-Brownian motion, whose density solves this
equation~\cite{Dawson1993,Perkins2002,BassPerkins2001}.

Uniform ellipticity opens up a different set of methods for weak
uniqueness. In finite dimensions, these methods can handle rougher
coefficients than the pathwise criteria above. The martingale problem
uses the generator to determine the law, and continuity of the
covariance suffices under uniform ellipticity
\cite{StroockVaradhan2006}. For H\"older covariance matrices, the perturbation argument of
Bass and Perkins~\cite{BassPerkins2009} uses smoothing by a diffusion
with frozen coefficients.

A weak uniqueness theory of comparable generality has been missing
for the stochastic heat equation with pointwise diffusion. Our aim
is to establish such a theory for uniformly nondegenerate coefficients
throughout the positive H\"older range. A central difficulty is
that the pointwise map $u\mapsto g(u(\cdot))$ couples infinitely
many spatial modes.

Earlier infinite-dimensional martingale-problem results leave
the general pointwise H\"older class open. Zambotti's Schauder
theory treats covariance perturbations of an Ornstein--Uhlenbeck
operator~\cite{Zambotti1999,Zambotti2000}. Athreya, Bass, Gordina
and Perkins~\cite{AthreyaBassGordinaPerkins2006} impose conditions
on the covariance matrix and coordinate damping. In Bass and
Perkins~\cite[Section~1]{BassPerkins2012}, the exponent above $1/2$
refers to H\"older dependence on the solution in a weak Wasserstein
norm. Their result also requires spatial regularity of the noise
coefficient. These hypotheses allow dependence on suitable smooth
averages of the solution. A scalar H\"older bound on $g$ gives no
such control of $u\mapsto g(u(x))$ in terms of spatial averages.
The $1/2$ threshold in that work therefore does not apply to the
general pointwise equation studied here.

The generalized coupling method used here grew out of work on
asymptotic coupling for ergodicity by Hairer, Mattingly and
Scheutzow~\cite{HairerMattinglyScheutzow2011}. Kulik and Scheutzow
developed generalized coupling criteria for Markov processes
and applied the method to weak well-posedness of stochastic delay
equations~\cite{KulikScheutzow2018,KulikScheutzow2020}. Butkovsky,
Kulik and Scheutzow used it to obtain quantitative ergodic
estimates for SPDEs~\cite{ButkovskyKulikScheutzow2020}.
In the version used here, a control drift keeps a companion
process close to a prescribed solution. A Girsanov change of measure
removes the control, and an estimate on the control energy bounds
the resulting change in law. We used this approach to prove weak
uniqueness above $3/4$ for the white-noise heat equation and to
treat wave equations with H\"older noise~\cite{Han2025,Han2024Wave}.

Lindvall and Rogers developed reflection coupling for
multidimensional diffusions~\cite{LindvallRogers1986}. Da Prato, Debussche and
Tubaro~\cite{DaPratoDebusscheTubaro2005} and Wang~\cite{Wang2015}
adapted reflection methods to SPDEs. In our earlier heat equation
proof~\cite{Han2026Threshold}, uniform nondegeneracy allowed us
to split $g^2$ into a positive constant and a state-dependent
remainder and to realize the two parts using independent white
noises. We then reflected a component of the additive noise
along the difference of the solutions. This produced a favorable term in
It\^o's formula for a concave function of their negative Sobolev
distance and gave weak uniqueness above $2/3$.

Here we use the same additive-noise split and drive both solutions
with the same additive white noise. It cancels in the difference
equation and provides the noise through which the change of measure
removes the control. The remaining noise is coupled through its
spatial covariance. We measure the signed difference of these
covariances in a negative Sobolev norm, so that heat smoothing
can average the coefficient difference across space before it
enters the estimate on the solutions.

\subsection{Beyond the scalar periodic equation}
There are two potential directions for extending the scalar result of
Theorem~\ref{thm:main-extensions}. The first is to treat vector-valued
solutions on the torus. We expect the coupling argument to extend to
this setting because a uniformly positive covariance matrix still allows
us to split off a common additive white noise. The second is to consider
other spatial domains, boundary conditions, or generators, which would
require heat-kernel estimates adapted to the new setting.

\begin{candidategeneralization}[Vector generalization]
\label{ext:vector-extension-statement}
Fix $1\le d<\infty$ and put $H_d=L^2(\mathbb T;\mathbb R^d)$.
Write $G^*$ for transpose and $I_d$ for the $d\times d$ identity matrix.
Let $G:\mathbb R^d\to\mathbb R^{d\times d}$ satisfy
$G(z)G(z)^*\ge c_0^2I_d$ for some $c_0>0$ and every $z$.
Let $b:[0,\infty)\times\mathbb T\times\mathbb R^d\to\mathbb R^d$
be Borel measurable and suppose that, for every $T<\infty$, there is
$L_{b,T}<\infty$ such that
\[
 |b(t,x,z)|\le L_{b,T}(1+|z|),
 \qquad 0\le t\le T,\quad x\in\mathbb T,\quad z\in\mathbb R^d.
\]
Write $b(t,u)(x)=b(t,x,u(x))$. Consider
\begin{equation}
 du=[\Delta u+b(t,u)]\,dt+G(u)\,dW,\qquad u(0)=u_0,
 \label{eq:vector-drifted-system}
\end{equation}
where the $d$ components of $W$ are independent space-time white noises.
H\"older continuity of $G$ is measured using the Euclidean norm on
$\mathbb R^d$ and the Hilbert--Schmidt norm on matrices.
For every $T>0$, a weak mild solution with continuous $H_d$ paths
exists and its law on $C([0,T];H_d)$ is unique in either of the
following cases:
\begin{enumerate}[label=(\roman*),leftmargin=*]
\item $G$ is bounded and globally $\beta_0$-H\"older for some
$\beta_0>0$, and $u_0\in L^2(\mathbb T;\mathbb R^d)$.
\item $G$ is globally $\beta$-H\"older for some $0<\beta\le1$,
with no uniform upper bound required, and
$u_0\in L^\infty(\mathbb T;\mathbb R^d)$.
\end{enumerate}
We expect this vector extension to hold and leave its proof for future work.
\end{candidategeneralization}

Remark~\ref{rem:boundary-fractional-discussion} discusses other
generalizations, including boundary conditions on an interval and
fractional powers of the Laplacian.

\begin{remark}[Boundary conditions and fractional powers]
\label{rem:boundary-fractional-discussion}
{\interlinepenalty=10000
The proof on the torus uses the heat kernel to calculate the variance of
the Gaussian reference fields \eqref{eq:gaussian-profile-variance-map}
and to average the signed covariance difference across space.\par}

\emph{Boundary conditions on an interval.}
With Neumann conditions, short-time Gaussian fluctuations have the same
order near the endpoints as in the interior. The spatial cancellation
estimate \eqref{eq:haar-heat-signed} would need a version for the reflected
heat kernel. With Dirichlet conditions, the solution is pinned to zero at
the endpoints and the fluctuations weaken near the boundary. The
comparison would need a variance bound that tracks the distance to the
endpoints~\cite{GyryaSaloffCoste2011}.

\emph{Fractional diffusion.} Replacing $\Delta$ by $-(-\Delta)^\alpha$
changes the size and spatial dependence of the Gaussian fluctuations
even on the torus. An extension would need a counterpart to the
two-point bound \eqref{eq:sharp-two-point-annular}, which controls the
expected product of two frequency-localized tests in terms of their
spatial separation. For each fixed $\beta>0$, the proof's parameter
inequalities are strict at $\alpha=1$. They can persist for nearby powers
if the required estimates have exponents that vary continuously with
$\alpha$.

For comparison, our earlier work gives weak uniqueness for the spectral
fractional Dirichlet Laplacian on an interval at H\"older exponents
$\beta>1/2+1/(4\alpha)$ under the remaining hypotheses
of~\cite[Theorem~1.2]{Han2025}.
\end{remark}

\subsection{Main ideas of the proof}
\label{sec:comparison-setup}
To prove Theorem~\ref{thm:main-weak-uniqueness}, we approximate solutions
by those of equations with smooth diffusion coefficients. For each smooth
coefficient, the equation has a unique law from the prescribed
initial state. To compare directly with an arbitrary weak solution
$U$, we construct a smooth companion $V$ on an enlargement of its
probability space, start it from the same state, and add the drift
$\Lambda(U-V)$, with $\Lambda\ge1$, to keep it close to $U$.
This follows the generalized coupling strategy
of~\cite{KulikScheutzow2020,ButkovskyKulikScheutzow2020,Han2025}.

On a fixed horizon $0<T\le1$, the integrated mean-square distance
between the solutions is
$E_\Lambda=\mathbb E\int_0^T\|U_t-V_t\|_{L^2(\mathbb T)}^2\,dt$.
Removing the added drift by Girsanov costs, up to a fixed noise
factor, its integrated mean square $\Lambda^2E_\Lambda$.
Thus $E_\Lambda$ must decay faster than $\Lambda^{-2}$
(Proposition~\ref{prop:law-comparison}). For each $\Lambda$, we
smooth the coefficient finely enough that the error due to
smoothing also vanishes at the scale of the control energy.

\subsubsection{Averaging the signed covariance difference}
Pointwise H\"older continuity by itself gives too little decay at
small exponents. For finite-dimensional diffusions damping contributes
a factor $\Lambda^{-1}$ to the mean-square estimate.
For the heat equation it has factor $\Lambda^{-1/2}$ because
the squared heat-kernel norm behaves like $t^{-1/2}$ at short times.
With the pointwise H\"older bound, the respective rates for
$E_\Lambda$ are $\Lambda^{-1/(1-\beta)}$ and
$\Lambda^{-1/[2(1-\beta)]}$ at exponent $0<\beta<1$, with
sufficiently small smoothing errors. To beat $\Lambda^{-2}$, these
rates require $\beta>1/2$ and $\beta>3/4$, respectively.

The spatial estimate exploits a feature that a pointwise bound
misses. If the covariance error is positive on one half of a
short interval and negative on the other, its signed spatial
average can be small even when its pointwise magnitude is large.
To make this cancellation available in the stochastic isometry,
we build a noise factor adapted to spatial averages.

First reserve a common additive white noise of amplitude $m_0$,
where $0<m_0<\inf_{z\in\mathbb R}|g(z)|$. It cancels in the difference equation and
provides the noise through which Girsanov removes the control.
Write $a=g^2-m_0^2$ for the remaining variance and
$a_{n_{\mathrm{sm}}}$ for its smooth approximation, where
$n_{\mathrm{sm}}$ indexes the smoothing. The smooth equation
has diffusion $\sqrt{m_0^2+a_{n_{\mathrm{sm}}}}$.
For the remaining noise we use a weighted Haar factor
(Lemma~\ref{lem:analytic-covariance-factor}). Its weights depend
on the variance masses of the two halves of each dyadic interval
and make the corresponding column integrate to zero. Heat propagation
can then exploit cancellation both within each column and in the
differences of the local masses.

Our earlier proof~\cite{Han2026Threshold} used a negative Sobolev
norm on the solution difference. Here the norm is on the signed
covariance error $a(U)-a_{n_{\mathrm{sm}}}(V)$. Write
$\mathsf d_{\mathrm{coef}}(U,V)$ for its integrated size in the
spatial norm $\|(\Lambda-\Delta)^{-\rho/2}\,\cdot\,\|_{L^2(\mathbb T)}$,
where $0<\rho<1/2$, including the smoothing allowance specified in
\eqref{eq:coefficient-discrepancy-norm}. The signed heat estimate and
stochastic isometry give
\begin{equation}
 E_\Lambda^{1/2}
 \le C\Lambda^{\rho/2-1/4}\mathsf d_{\mathrm{coef}}(U,V).
 \label{eq:introduction-signed-energy}
\end{equation}
To close this inequality, we need to estimate
$\mathsf d_{\mathrm{coef}}(U,V)$ in terms of $E_\Lambda^{1/2}$.

\subsubsection{How state oscillations give spatial decay}
A Littlewood--Paley block $a_\lambda(z)$ of the coefficient
oscillates on scale $\lambda^{-1}$ in the scalar state $z$,
where $\lambda$ ranges over dyadic frequencies. A Gaussian shift
of this size begins to average away the oscillation. For a concrete
calculation, let $G_h$ be the Gaussian heat field driven by the
common additive white noise over a short window of length $h$, starting
from zero. At spatial points $x,y\in\mathbb T$,
\[
 \mathbb E\cos\!\bigl(\lambda[G_h(x)-G_h(y)]\bigr)
 =\exp\!\left[-\tfrac12\lambda^2
           \operatorname{Var}(G_h(x)-G_h(y))\right].
\]
Up to separation $\sqrt h$, the variance of the difference is
at least a constant times the torus distance between $x$ and $y$.
Thus the state scale $\lambda^{-1}$ leads to a spatial separation
scale $\lambda^{-2}$. The two-point estimate
\eqref{eq:sharp-two-point-annular} makes this decay precise for
general frequency-localized tests. After integration in the
negative Sobolev norm, it gives a factor $\lambda^{-2\rho}$ in the
norm, up to an exponentially small term.

The same scaling explains the time window. At one point the
Gaussian fluctuation after time $h$ is of order $h^{1/4}$, so
smoothing a state oscillation on scale $\lambda^{-1}$ begins
at time of order $\lambda^{-4}$. Below the cutoff of order
$\Lambda^{1/4}$, the blocks can be estimated by Bernstein's
inequality and absorbed using \eqref{eq:introduction-signed-energy}.
Above it we use $h_\lambda=\lambda^{-3\beta}$, where $\beta$ is the
effective H\"older exponent fixed in Section~\ref{sec:framework}.
These windows shrink with frequency, making the expansion
remainder smaller, while $h_\lambda$ is long enough for Gaussian
smoothing to act strongly on the band.

To use this calculation for the nonlinear equation, we condition
on the past at the start of the window and construct a Gaussian
reference for each solution. Its residual variance profile is
chosen to equal the conditional mean of $a$ evaluated on that very
reference field. One can picture the construction at each point
as averaging $a$ over a Gaussian cloud and then using the average
to set the width of the cloud through the heat equation.
Proposition~\ref{prop:self-profile} constructs this profile by
contraction. The defining property means that the residual
coefficient minus the profile has zero conditional mean when
evaluated on the reference field.

On a single band, the main term of the Gaussian comparison
\eqref{eq:sharp-main-Gaussian-band} is
\[
 \lambda^{-\beta}\,\lambda^{-2\rho}
          \min(1,\lambda E_\Lambda^{1/2})
 \qquad\text{in the integrated negative Sobolev norm.}
\]
The factors are the amplitude of the block, the spatial gain, and
the size of the change when the two reference fields are compared.
The last factor grows linearly with the solution distance until
that distance reaches one oscillation scale. Summing the bands
on either side of $\lambda E_\Lambda^{1/2}=1$ gives the power
$E_\Lambda^{(\beta+2\rho)/2}$ when $0<\beta+2\rho<1$
(Theorem~\ref{thm:coupled-gaussian-band-sum}).

\subsubsection{What one step of the expansion adds}
Fix a band and a window $[t_0,t_0+h]$ with $h=h_\lambda$.
At a terminal point $x$, put $Y_s=(S_{h-s}U_{t_0+s})(x)$ for
$0\le s\le h$, so $Y_h=U_{t_0+h}(x)$. Let $F_s$ be
$a_\lambda$ averaged against the future noise of the Gaussian
reference from elapsed time $s$ to $h$. In It\^o's formula for
$F_s(Y_s)$, the backward time derivative cancels the reference
part of the covariance. Alongside the initial Gaussian average
and a stochastic integral, the formula leaves the integral in
$s$ and $y\in\mathbb T$ of
\[
 \tfrac12 p_{h-s}(x-y)^2
 \bigl[a(U_{t_0+s}(y))-\bar a_s(y)\bigr] F_s''(Y_s).
\]
Here $p_t$ is the heat kernel and $\bar a_s$ is the reference
profile. The second derivative still comes from the original
terminal test at $x$. The bracket introduces a new coefficient
observation at $(t_0+s,y)$. Applying the same step to the enlarged
list of observations produces the finite expansion
(Lemma~\ref{lem:finite-expansion}). In the comparison of $U$ and
$V$, the first Gaussian average and stochastic integral are grouped
before estimation so that their cancellation is preserved
(Proposition~\ref{prop:joint-first-pair}).

The danger in repeating the step is that two observation times
can nearly coincide. A second Gaussian derivative across a time
gap $\delta$ gives a weight of order $\delta^{-1}$, but that
weight has Gaussian mean zero. We can therefore subtract the
expression with the corresponding Gaussian increment suppressed.
The singular weight is now paired with a difference that vanishes
as the gap shrinks. A terminal
block gives a factor $\min(1,\lambda\delta^{1/4})$ relative to
its amplitude, while an inserted H\"older coefficient gives
$\delta^{\beta/4}$. If several increments act on the same rough
coefficient, the smallest gap controls the shared gain
(Lemma~\ref{lem:centered-word-time-kernel}). The time integrals
at successive orders give enough decay in the remainder at depth
$N_\beta=\lceil6/\beta^2\rceil$, fixed once $\beta$ is chosen.

For the companion, the backward observation in
\eqref{eq:joint-scalar-backward-coordinate} still ends at
$V_{t_0+h}(x)$ and cancels the control drift. Its covariance has
an extra term with a factor $1-e^{-\Lambda\tau}\le\Lambda\tau$,
where $\tau$ is the time remaining to the observation. The factor
reduces the time singularity near that observation. We bound the extra
covariance term when it first appears and continue expanding the
ordinary insertions (Theorem~\ref{thm:joint-finite-comparison}).

\subsubsection{From the energy bound to uniqueness in law}
Increasing $\rho$ improves the power of the solution distance in the
Gaussian comparison and weakens the damping factor in
\eqref{eq:introduction-signed-energy}. With $\rho=1/2-3\beta/5$, the comparison
exponent becomes $\beta+2\rho=1-\beta/5$. Substituting this
term in \eqref{eq:introduction-signed-energy} gives
$C\Lambda^{-3\beta/10}E_\Lambda^{(1-\beta/5)/2}$ in the bound
for $E_\Lambda^{1/2}$. Young's inequality converts it to a
squared-energy contribution of order $\Lambda^{-3}$. The other
terms in the finite expansion allow the common rate
$E_\Lambda=O(\Lambda^{-5/2})$ after the smoothing is chosen
(Lemmas~\ref{lem:parameter-energy-closure}
and~\ref{lem:common-approximation-control}).

Consequently the control energy tends to zero. For two prescribed
weak solutions we use the same smooth coefficient at each
$\Lambda$, so the two comparisons give the same limiting law.
This proves local uniqueness in law. Conditional restart extends
it to all finite horizons (Section~\ref{sec:conditional-global}).

\subsection{AI assistance and Lean formalization}
\label{sec:ai-lean}
This work grew out of the author's proof of weak uniqueness above
the H\"older exponent $2/3$~\cite{Han2026Threshold}. The author developed
the comparison argument through sustained interaction with OpenAI's
GPT-6.1, which helped uncover spatial cancellation in the signed
covariance difference and the tree structure of the Gaussian expansion.
The trees organize the centering of successive coefficient insertions
through a depth that is finite for each H\"older exponent and can be
large for small exponents. GPT-6.1 also assisted with the estimates
and exposition. The author takes responsibility for the manuscript.

The author developed the end-to-end Lean~4 formalization of
Theorems~\ref{thm:main-weak-uniqueness} and~\ref{thm:main-extensions}
in the \scalarshelean{} repository. The project builds on Mathlib
and the bundled BrownianMotion, KolmogorovExtension4 and MathFin
libraries. The author's development includes the following tools,
which can be imported from the \emph{SHEStochastic} modules for use
in other stochastic models.

\begin{itemize}[leftmargin=*,nosep]
\item \emph{Cylindrical It\^o integration.} The modules
\emph{FiniteCylindricalIto} and \emph{HilbertTerminalIto} construct
stochastic integrals for predictable square-integrable inputs with
real or Hilbert-valued outputs. They include the isometry and
compatibility with continuous linear functionals.
\item \emph{Burkholder--Davis--Gundy moment bounds.}
\emph{NativeScalarBDGAll} and \emph{NativeHilbertBDG} bound terminal
moments by moments of the integrated squared integrand norm. The scalar
bounds cover all even integer orders, and the Hilbert-valued bounds
cover orders $2,4,8,\ldots$.
\item \emph{Cylindrical Girsanov.} \emph{NativeTerminalGirsanovTransport}
gives the joint laws of the shifted Brownian coordinates and the formula
for the corresponding change in the stochastic integral.
\item \emph{Weak limits and pasting of laws.}
\emph{WeakLimitCylinderMartingale} passes martingale identities tested
against bounded continuous functions of the past to weak limits on
path space. \emph{KernelMartingalePasting} gives the martingale identities
for a stopped path continued using a probability kernel.
\end{itemize}

The formalization also proves weak existence for the scalar heat
equation. For bounded diffusion, \emph{GlobalNativePositiveWeakExistence}
constructs a limit of the joint laws of smooth solutions and their
Brownian drivers and identifies the limiting stochastic integral.
\emph{ScalarUnboundedNativeExistence} supplies the corresponding
construction for unbounded diffusion. The file \emph{ScalarSHE.lean}
imports the main theorem declarations. The audit in \emph{Audit.lean}
reports their dependencies on \emph{propext}, \emph{Classical.choice},
and \emph{Quot.sound}.

\section{Solution framework and the route from energy to laws}
\label{sec:framework}
Unless another domain is displayed, $\|\cdot\|_p$ is the
spatial $L^p(\mathbb T)$ norm. We write the spatial square norm
explicitly as $\|f\|_{L^2(\mathbb T)}=(\int_{\mathbb T}|f(x)|^2\,dx)^{1/2}$. A star denotes transpose or adjoint;
$I$ denotes the identity operator, and $I_d$ the identity matrix in
$\mathbb R^d$. We write $\mathsf M_f$ for multiplication by a scalar
or matrix field $f$, $\|\cdot\|_{\mathrm{op}}$ for the operator norm,
$\|\cdot\|_{\mathrm{HS}}$ for the Hilbert--Schmidt norm, and
$\|A\|_{\mathrm F}=(\operatorname{tr}(A^*A))^{1/2}$ for the Frobenius
norm of a matrix. Here $\operatorname{tr}$ denotes the trace.
The symbols $c,C$ denote positive constants whose values may change.
Unless stated otherwise, they may depend on the fixed exponent,
coefficient bounds, finite order, dimension and time horizon, and are uniform in the damping parameter and approximation indices.
$A\lesssim B$ means $A\le CB$ and $A\asymp B$ means $cB\le A\le CB$.
Expectations are taken under the probability measure of the given
solution, denoted by $\mathbb P$, with filtration $(\mathcal F_t)$
on sample space $\Omega$.

\subsection{State space, noise and observable mild equations}
Put $H=L^2(\mathbb T;\mathbb R)$, with normalized Lebesgue measure
on the torus. The heat kernel is denoted by $p_t$; thus
$S_tf=p_t*f$. All path spaces below carry their Borel sigma fields.
For each horizon $h>0$, we use the path spaces
\[
 \mathcal X_h=C([0,h];H),\qquad
 \mathcal Y_h=L^2([0,h];H).
\]
Their norms are
$\|u\|_{\mathcal X_h}=\sup_{0\le t\le h}\|u(t)\|_{L^2(\mathbb T)}$ and
$\|u\|_{\mathcal Y_h}=(\int_0^h\|u(t)\|_{L^2(\mathbb T)}^2dt)^{1/2}$.
The inclusion $\iota:\mathcal X_h\to\mathcal Y_h$ sends a continuous
path to its equivalence class under equality at almost every time.
Both are separable complete metric spaces. Their roles differ:
$\mathcal Y_h$ is convenient for coupling estimates, whereas the
conclusion of the theorem concerns $\mathcal X_h$.

\begin{definition}[Integrated norms for windows of length $h$]
\label{def:window-integrated-norms}
Fix a finite horizon $T>0$ and a probability space
$(\Omega,\mathcal F,\mathbb P)$. For a Hilbert space $E$ and
$0\le h\le T$, set
\[
 \begin{aligned}
 \mathscr L_h^2(E)
 &:=L^2([0,T-h]\times\Omega,\,dt_0\otimes\mathbb P;E),\\
 \|F\|_{\mathscr L_h^2(E)}^2
 &:=\int_0^{T-h}\mathbb E\|F(t_0)\|_E^2\,dt_0.
 \end{aligned}
\]
The horizon $T$ and probability space are fixed in each application.
For $h>0$, $F(t_0)$ is a random output from the window
$[t_0,t_0+h]$; $t_0$ ranges continuously over its starting times,
while the subscript $h$ records the common duration. In particular,
$\mathscr L_0^2(E)$ is the full-horizon space on $[0,T]$.
\end{definition}

\begin{definition}[Weak mild solution]
\label{def:weak-mild-solution}
Write $\mathsf M_f$ for multiplication by $f$.
A weak solution consists of a filtered probability space, an adapted
continuous $H$-valued process $U$, and a cylindrical Brownian motion
$W$ on $H$, such that
\begin{equation}
 U_t=S_tu_0+\int_0^t S_{t-r}\mathsf M_{g(U_r)}\,dW_r.
 \label{eq:framework-mild}
\end{equation}
\end{definition}

The integrand is
Hilbert--Schmidt after heat smoothing. Indeed boundedness of $g$ gives
\[
 \|S_t\mathsf M_{g(u)}\|_{\mathrm{HS}}^2
 =\int_{\mathbb T^2}p_t(x-y)^2g(u(y))^2\,dx\,dy
 \le C\|p_t\|_{L^2(\mathbb T)}^2\le C(1+t^{-1/2}).
\]
The singularity is integrable at zero, so the mild formulation is
well defined for $L^2$ initial states.

An equivalent and useful formulation tests \eqref{eq:framework-mild}
against positive-time heat rows. For $q\ge0$, $\varepsilon>0$ and
$x\in\mathbb T$, it reads
\begin{equation}
 S_\varepsilon U_q(x)=S_{q+\varepsilon}u_0(x)
 +\int_0^q\langle g(U_r)p_{q+\varepsilon-r}(x-\cdot),dW_r\rangle_H.
 \label{eq:framework-observed}
\end{equation}
This scalar formulation will be used when conditioning. It avoids
an implicit choice of point evaluations of an $L^2$ state. The family
of positive-time heat tests separates $H$: if all such tests of $f$
vanish, then $S_\varepsilon f=0$ for every $\varepsilon>0$, and
strong continuity gives $f=0$.

We use cylindrical It\^o integration, finite-dimensional It\^o's formula,
Burkholder--Davis--Gundy (BDG) inequalities, stochastic Fubini,
L\'evy's characterization of Brownian motion, Novikov's criterion,
and cylindrical Girsanov. For finite-dimensional stochastic calculus
we use \cite[Chapters~IV--V and VIII]{RevuzYor1999}; for cylindrical
integration, stochastic Fubini and Girsanov we use
\cite[Chapters~4 and 10]{DaPratoZabczyk2014}. The equivalent
martingale-measure formulation of the heat stochastic
integrals is developed in \cite{Walsh1986}.
The conditional restart argument is given
in Section~\ref{sec:conditional-global}.

\subsection{Coefficient reduction and smooth approximations}
For a scalar function $f$ let
$[f]_{C^\beta}=\sup_{z\ne z'}|f(z)-f(z')|/|z-z'|^\beta$.
Fix an effective exponent
\[
 0<\beta\le\min(\beta_0,4/25).
\]
Boundedness of $g$ allows its original H\"older estimate to be
lowered to this exponent: use the original estimate for distances
at most one and the uniform bound for larger distances. Increase $L_g$ to this effective H\"older bound if necessary. Choose
$0<m_0<c_0$. Set $\underline a=c_0^2-m_0^2$, $\overline a=C_0^2-m_0^2$,
and choose $L_a\ge[g^2]_{C^\beta}$. Then
\begin{equation}
 a(z)=g(z)^2-m_0^2,\qquad
 0<\underline a\le a(z)\le\overline a,
 \qquad [a]_{C^\beta}\le L_a.
 \label{eq:framework-residual}
\end{equation}
Only $g^2$ enters the covariance estimates. The noise rotation uses the sign of the original coefficient;
the state-law argument is unchanged.

The smoothing index $n_{\mathrm{sm}}\ge1$ is an integer.
Fix a nonnegative smooth compactly supported function $\eta$ on
$\mathbb R$ with integral one, put
$\eta_{n_{\mathrm{sm}}}(z)=n_{\mathrm{sm}}\eta(n_{\mathrm{sm}}z)$,
and let $a_{n_{\mathrm{sm}}}=\eta_{n_{\mathrm{sm}}}*a$ be the state mollification.
Then
\begin{equation}
 \underline a\le a_{n_{\mathrm{sm}}}\le\overline a,\quad
 [a_{n_{\mathrm{sm}}}]_{C^\beta}\le L_a,\quad
 \epsilon_{n_{\mathrm{sm}}}:=\|a_{n_{\mathrm{sm}}}-a\|_\infty\longrightarrow0.
 \label{eq:framework-smoothing}
\end{equation}
Each $a_{n_{\mathrm{sm}}}$ is globally Lipschitz, although its derivative bounds
need not be uniform in $n_{\mathrm{sm}}$. The canonical smooth equation has
diffusion $\sqrt{m_0^2+a_{n_{\mathrm{sm}}}}$. All quantitative comparison constants
depend on the three common bounds in
\eqref{eq:framework-smoothing}.

Write $\bar a_t(x)$ for a frozen residual profile, whose
self-consistent construction is given in Section~\ref{sec:self-consistent-profile}.
We distinguish residual and total variances:
\[
 a,\ \bar a\quad\hbox{are residual profiles},\qquad
 \mathfrak a=m_0^2+a,\quad \bar{\mathfrak a}=m_0^2+\bar a\quad\hbox{are total profiles}.
\]
The insertion difference is $\mathfrak a-\bar{\mathfrak a}=a-\bar a$, but every physical Gaussian
bracket includes the additive contribution $m_0^2$, which supplies
the uniform lower bound used by the Gaussian derivative estimates.
The positive residual bound also gives the invertible covariance
factor and its derivative estimates in
Lemma~\ref{lem:analytic-covariance-factor}. The positive total variance
provides the Gaussian regularization used to construct the reference
profiles in Proposition~\ref{prop:self-profile} and to control the
Gaussian derivative weights in Lemmas~\ref{lem:gaussian-one-input-scores}
and~\ref{lem:physical-heat-score}, which are used throughout Sections~5--6.
The separate additive channel has two further roles: it preserves
Gaussian smoothing after conditioning on the residual noise in
Corollary~\ref{cor:coupled-gaussian-band}, and it allows the Girsanov
change of measure that removes the control in
Proposition~\ref{prop:law-comparison}.

\subsection{What the quantitative argument must prove}
For a prescribed weak solution $U$, damping $\Lambda\ge1$ and
smooth index $n_{\mathrm{sm}}$, we construct a companion $V^{n_{\mathrm{sm}},\Lambda}$ from
the same initial state. Its extra drift is $\Lambda(U-V^{n_{\mathrm{sm}},\Lambda})$.
For a fixed horizon $T>0$ and difference $D=U-V$, measure the
solution difference by its integrated norm and define the coefficient discrepancy
$\mathsf d_{\mathrm{coef}}(U,V)$ by
\begin{equation}
 \begin{aligned}
 \|U-V\|_{\mathscr L_0^2(H)}^2&:=\mathbb E\int_0^T\|D_t\|_{L^2(\mathbb T)}^2\,dt,\\
 \mathsf d_{\mathrm{coef}}(U,V)^2&:=\mathbb E\int_0^T
 \|a(U_t)-a_{n_{\mathrm{sm}}}(V_t)\|_{H_\Lambda^{-\rho}}^2\,dt
                                      +\varepsilon_{\mathrm{app}},
 \end{aligned}
 \label{eq:coefficient-discrepancy-norm}
\end{equation}
where
$\|f\|_{H_\Lambda^{-\rho}}=\|(\Lambda-\Delta)^{-\rho/2}f\|_{L^2(\mathbb T)}$,
$0<\rho<1/2$, and $\varepsilon_{\mathrm{app}}\ge0$ is a deterministic approximation source.
Here $\mathscr L_0^2(H)=L^2([0,T]\times\Omega,dt\otimes\mathbb P;H)$;
its window variants are defined in Definition~\ref{def:window-integrated-norms}.
We write $\mathsf d_{\mathrm{coef}}(U,V;T)$ when the horizon
needs to be explicit.
We denote the integral in the definition of $\mathsf d_{\mathrm{coef}}(U,V)^2$ by $J_\Lambda$,
so $\mathsf d_{\mathrm{coef}}(U,V)^2=J_\Lambda+\varepsilon_{\mathrm{app}}$.
The operator $Q_\Lambda=(\Lambda-\Delta)^{-\rho}$ will be used in the spatial estimates.

Set
\[
 \|U-V^{n_{\mathrm{sm}},\Lambda}\|_{L^2(\Omega;\mathcal Y_h)}^2
 =\mathbb E\int_0^h\|U_t-V_t^{n_{\mathrm{sm}},\Lambda}\|_{L^2(\mathbb T)}^2\,dt.
\]
The target is a deterministic choice of indices such that
\begin{equation}
 \Lambda^2\|U-V^{n_{\mathrm{sm}}(\Lambda),\Lambda}\|_{L^2(\Omega;\mathcal Y_h)}^2\longrightarrow0.
 \label{eq:framework-target}
\end{equation}
This is stronger than convergence of $V$ to $U$. The extra factor
$\Lambda^2$ is the energy of the drift that must be removed by
Girsanov. Preliminary damping alone gives only $\|U-V\|_{L^2(\Omega;\mathcal Y_h)}^2\lesssim
\Lambda^{-1/2}$, so most of the argument is spent improving this rate.

In this local comparison, the coefficient discrepancy is also taken over
$[0,h]$. The signed covariance bound yields
$\|U-V\|_{L^2(\Omega;\mathcal Y_h)}^2\lesssim\Lambda^{\rho-1/2}\mathsf d_{\mathrm{coef}}(U,V)^2$. The finite comparison yields a
bound on $\mathsf d_{\mathrm{coef}}(U,V)$ in terms of absorbable linear terms, three sublinear
powers of $\|U-V\|_{L^2(\Omega;\mathcal Y_h)}$, and approximation errors. After choosing
parameters, the resulting estimate is
\[
 \|U-V^{n_{\mathrm{sm}}(\Lambda),\Lambda}\|_{L^2(\Omega;\mathcal Y_h)}^2\le C\Lambda^{-5/2}.
\]
The exact exponent $5/2$ is dispensable; any exponent larger than
two would suffice. Its positive margin accommodates the finitely
many logarithmic and moment losses.

\begin{lemma}[Integrated laws determine continuous path laws]
\label{lem:integrated-laws}
The inclusion $\iota:\mathcal X_h\to\mathcal Y_h$ is continuous
and injective. If two probability measures on $\mathcal X_h$ have
the same pushforward by $\iota$, they are equal.
\end{lemma}
\begin{proof}
The bound $\|f\|_{\mathcal Y_h}\le\sqrt h\|f\|_{\mathcal X_h}$
gives continuity, and continuity of paths gives injectivity.
Since both spaces are Polish, $\iota$ has a Borel inverse on its
image~\cite[Corollary~15.2]{Kechris1995}. Applying this inverse to
the common pushforward proves equality of the measures.
\end{proof}

For later reference equip $\mathcal Y_h$ with its norm metric
truncated at one. If $F$ is bounded by one and is one-Lipschitz,
\begin{equation}
 |\mathbb EF(U)-\mathbb EF(V)|
 \le\mathbb E\min(1,\|U-V\|_{\mathcal Y_h})\le\|U-V\|_{L^2(\Omega;\mathcal Y_h)}.
 \label{eq:framework-bl}
\end{equation}
Thus the integrated coupling norm controls the distance between path laws.

\subsection{Order of choices and limits}
The effective exponent $\beta$ is fixed first. Next fix the finite
depth $N_\beta=\lceil6/\beta^2\rceil$ and all moment exponents required
at that depth. Small positive losses are then chosen below the
strict margins displayed in the parameter calculation. Only after
these choices does $\Lambda$ tend to infinity.

At fixed $\Lambda$ and expansion depth, we first regularize the
finite words away from zero time gaps, then remove these restrictions
using the uniform estimates of
Lemma~\ref{lem:exact-bracket-cutoff-limits}. The coefficient-insertion
estimate in Appendix~\ref{sec:uniform-coefficient-diagonal} gives a
deterministic smoothing index $n_{\mathrm{sm}}(\Lambda)$, chosen to
work for both prescribed weak solutions.

The comparison is local in its Gaussian windows. We first use it
to compare laws over one fixed short horizon from a deterministic
$H$ state. Standard conditional restart and finite induction then extend this
uniqueness to every finite horizon; see Section~\ref{sec:conditional-global}.

\raggedbottom
\subsection{Notation guide}
\label{sec:notation-guide}
The following table fixes the meanings of the symbols.
Local integration variables and finite lists are specified where they
are introduced. Here $a$ denotes the residual coefficient, $c,C$ are generic positive
constants, and $D=U-V$ is the difference of the physical solutions.
\begingroup
\renewcommand{\arraystretch}{1.12}
\begin{center}
\begin{tabular}{@{}>{\raggedright\arraybackslash}p{.27\linewidth}
>{\raggedright\arraybackslash}p{.67\linewidth}@{}}
\toprule
Symbol & Meaning\\
\midrule
$a,\ a_{n_{\mathrm{sm}}}$ & Residual coefficient and its smooth approximation;
the smoothing index is $n_{\mathrm{sm}}$.\\
$\bar a_t[w](x),\ \bar a_t^\phi(x)$ & Gaussian profile constructed from starting state $w$, or from a prescribed mean $\phi$.\\
$\mathfrak a,\ \bar{\mathfrak a}_t$ &
Total variances $m_0^2+a$ and $m_0^2+\bar a_t$.\\
$\kappa,\widetilde\kappa$ & Generic positive spatial densities used as arguments of $L$.\\
$D=U-V$ & Physical solution difference.\\
$\|U-V\|_{\mathscr L_0^2(H)}$ & Full-horizon $L^2$ norm in time, probability and space.\\
$\mathsf d_{\mathrm{coef}}(U,V)$ & Integrated signed coefficient discrepancy from
\eqref{eq:coefficient-discrepancy-norm}.\\
$\Sigma,\mathsf G$ & Finite observation covariance and constant-profile Gram matrix.\\
$\mathsf P_j,N_j,R_j$ & Gaussian average term, stochastic integral, and finite remainder.\\
\bottomrule
\end{tabular}
\end{center}
\begin{center}
\begin{tabular}{@{}>{\raggedright\arraybackslash}p{.27\linewidth}
>{\raggedright\arraybackslash}p{.67\linewidth}@{}}
\toprule
Symbol & Meaning\\
\midrule
$\lambda_*,\lambda_{\mathrm{cut}}$ & Lower high-frequency threshold and chosen upper
state-frequency cutoff.\\
$\gamma_{\mathrm{time}}$ & Real exponent in
$h_\lambda=\lambda^{-\gamma_{\mathrm{time}}}$; chosen as $3\beta$
in the main statements.\\
$N_\beta$ & Finite expansion depth, chosen once for the fixed exponent $\beta$.\\
$n_{\mathrm r},n_{\mathrm b}$ & Root-interval count and nonempty rough-block count.\\
$\partial_\kappa^j L,\partial_z^j a_\lambda$ &
Heat-regularized derivatives in a spatial density and ordinary derivatives
in the scalar state.\\
$\mathsf k_t,\mathscr K$ & Heat-square kernel $p_t^2$ and a native integration kernel.\\
\bottomrule
\end{tabular}
\end{center}
\endgroup
Gaussian convolution and its time-inhomogeneous versions are written $\mathcal P$; orthogonal projections other than $\Pi_l$ are written $\mathsf{Proj}$. The Haar functions $e_I$ and averaging projections $\Pi_l$ are defined
in the proof of Lemma~\ref{lem:analytic-covariance-factor}. The Bessel
kernel and its properties are given in
Section~\ref{sec:dyadic-bessel-preliminaries}.
\flushbottom

\clearpage
\section{Covariance coupling and spatial transport}

This section provides the coupling and its spatial estimates. The four
main conclusions are stated first, with their notation introduced locally.
Their proofs and the supporting spatial lemmas are collected at the end.
\begin{enumerate}[label=(\roman*),leftmargin=*]
\item \textbf{Signed covariance factor.} Lemma~\ref{lem:analytic-covariance-factor} turns a signed coefficient difference into a heat-smoothed noise estimate.
\item \textbf{Controlled companion.} Proposition~\ref{prop:construct-companion} constructs the smooth comparison process on the same noise space.
\item \textbf{Removal of the control.} Proposition~\ref{prop:law-comparison} converts small control energy into equality of state-path laws.
\item \textbf{Spatial trace transport.} Proposition~\ref{prop:unified-dyadic-column-transport} supplies the one-input estimates used by the finite expansion.
\end{enumerate}

\subsection{Main results}

\subsubsection{Analytic properties of the covariance factor}
\label{sec:haar-covariance-coupling}
\label{sec:haar-exact-coupling-law}

Identify $\mathbb T$ with $[0,1)$, with normalized Lebesgue measure,
and put $H=L^2(\mathbb T;\mathbb R)$. Our Fourier and heat conventions are
\[
 \widehat f(k)=\int_0^1 f(x)e^{-2\pi ikx}\,dx,\qquad
 p_t(x)=\sum_{k\in\mathbb Z}e^{-4\pi^2k^2t}e^{2\pi ikx},
 \qquad S_tf=p_t*f.
\]
For $\Lambda\ge1$ and $0\le\rho<1/2$, write
\[
 Q_\Lambda=(\Lambda-\Delta)^{-\rho},\qquad
 \|f\|_{H_\Lambda^{-\rho}}^2
   =\sum_{k\in\mathbb Z}(\Lambda+4\pi^2k^2)^{-\rho}
                                      |\widehat f(k)|^2.
\]
The positive space $H_\Lambda^\rho$ has the reciprocal Fourier weight
and is the dual space under the $L^2$ pairing. In particular,
$H_\Lambda^0=H$. For fixed $0<\kappa_-\le\kappa_+<\infty$ let
\[
 \mathcal C=\{\kappa\in L^\infty(\mathbb T):
               \kappa_-\le\kappa(x)\le\kappa_+
                       \text{ for almost every }x\}.
\]
We write $\mathsf M_\kappa$ for multiplication by $\kappa$.
Throughout this section pointwise heat estimates use $0<t\le1$;
damped integrals use the stated range $[0,\infty)$.
Constants may depend on $\kappa_-,\kappa_+,\rho$ and the fixed
derivative order, but are uniform in the profiles and $\Lambda$.
The following estimates are the properties of the factor used in the
comparison argument.

The next lemma constructs the common residual-noise factor. Its
signed estimates convert a coefficient difference in $H^{-\rho}_\Lambda$
into a heat-smoothed noise difference, which is the input to the
coupling energy estimate.

\begin{lemma}[Signed covariance factor and heat estimates]
\label{lem:analytic-covariance-factor}
Let $\mathcal L(H)$ denote the bounded linear operators on $H$.
There is a single map $L:\mathcal C\to\mathcal L(H)$, strongly Borel
measurable for the $L^2$ topology on $\mathcal C$, such that
\begin{equation}
 L(\kappa)L(\kappa)^*=\mathsf M_\kappa,\qquad
 \|L(\kappa)\|_{\mathrm{op}}\le\sqrt{\kappa_+},\qquad
 \|L(\kappa)^{-1}\|_{\mathrm{op}}\le\kappa_-^{-1/2}.
 \label{eq:covariance-factor-identity}
\end{equation}
For every $0\le\rho<1/2$, uniformly in $\kappa,\widetilde\kappa\in\mathcal C$,
$\Lambda\ge1$, and $0<t\le1$, it satisfies the signed estimates
\begin{align}
 \|S_t[L(\kappa)-L(\widetilde\kappa)]\|_{\mathrm{HS}}
 &\le C\bigl(\Lambda^{\rho/2}t^{-1/4}+t^{-1/4-\rho/2}\bigr)
                   \|\kappa-\widetilde\kappa\|_{H_\Lambda^{-\rho}},
 \label{eq:haar-heat-signed}\\
 \int_0^\infty e^{-2\Lambda t}
    \|S_t[L(\kappa)-L(\widetilde\kappa)]\|_{\mathrm{HS}}^2\,dt
 &\le C\Lambda^{\rho-1/2}
                   \|\kappa-\widetilde\kappa\|_{H_\Lambda^{-\rho}}^2.
 \label{eq:haar-damped-signed}
\end{align}
Here $C$ depends only on $\kappa_-,\kappa_+$ and $\rho$.
Thus the right-hand sides preserve the cancellation in the signed
profile difference. At $\rho=0$, these bounds give an $L^2$ Lipschitz
estimate with heat factor $t^{-1/4}$ and squared integral factor
$\Lambda^{-1/2}$.

For each $t>0$, the heat-regularized map $\kappa\mapsto S_tL(\kappa)$
extends smoothly, with values in the Hilbert--Schmidt operators, to
$\{\kappa\in L^\infty:\operatorname*{ess\,inf}\kappa>0\}$.
Its derivatives are the heat-regularized derivatives described in the proof.
\end{lemma}

\subsubsection{The canonical controlled companion}
\label{sec:canonical-companion}

Fix a prescribed weak solution $U$ from $u_0$. On an enlarged
space take orthogonal cylindrical Brownian motions $(W^0,B)$ such that
\[
 dU=\Delta U\,dt+m_0dW^0+L(a(U))dB.
\]
The noise enlargement establishing this representation is proved in
Lemma~\ref{lem:noise-enlargement}. The companion uses the same noises.

\begin{proposition}[Construction with prescribed forcing]
\label{prop:construct-companion}
For every $n_{\mathrm{sm}}$ and $\Lambda\ge1$ there is a unique adapted continuous
$H$-valued solution, on the enlarged space, of
\begin{equation}
 dV=(\Delta V+\Lambda(U-V))dt+m_0dW^0+L(a_{n_{\mathrm{sm}}}(V))dB,
 \qquad V_0=u_0.
 \label{eq:companion-controlled}
\end{equation}
The same construction without the control gives a canonical smooth
state-path law $\nu_{n_{\mathrm{sm}}}$ on $\mathcal X_h$, independent of the space
carrying the Brownian pair.
\end{proposition}

\subsubsection{Difference energy and removal of the control}
Fix a horizon $h>0$. For this comparison
$\|U-V\|_{L^2(\Omega;\mathcal Y_h)}^2=\mathbb E\int_0^h\|D_t\|_{L^2(\mathbb T)}^2dt$ and
$J_\Lambda=\mathbb E\int_0^h
\|a(U_t)-a_{n_{\mathrm{sm}}}(V_t)\|_{H^{-\rho}_\Lambda}^2dt$.
\begin{proposition}[Local law comparison from small control energy]
\label{prop:law-comparison}
Fix $h>0$ and the coupled pair $(U,V^{n_{\mathrm{sm}},\Lambda})$
from Proposition~\ref{prop:construct-companion}. Abbreviate
$V=V^{n_{\mathrm{sm}},\Lambda}$ and put
$\|U-V\|_{L^2(\Omega;\mathcal Y_h)}^2=\mathbb E\int_0^h\|U_t-V_t^{n_{\mathrm{sm}},\Lambda}\|_{L^2(\mathbb T)}^2dt$ and
$\mathcal C=m_0^{-2}\Lambda^2\|U-V\|_{L^2(\Omega;\mathcal Y_h)}^2$.
Recall that $\iota:C([0,h];H)\to L^2([0,h];H)$ is the inclusion
of continuous paths into time-square-integrable paths, and
$\nu_{n_{\mathrm{sm}}}$ is the law, from the same initial state $u_0$,
of the uncontrolled smooth equation with diffusion
$\sqrt{m_0^2+a_{n_{\mathrm{sm}}}}$.
For a bounded measurable test $G$ on $\mathcal X_h$,
$\nu_{n_{\mathrm{sm}}}(G)=\int_{\mathcal X_h}G\,d\nu_{n_{\mathrm{sm}}}$. For every test $F$ bounded by
one and one-Lipschitz on $\mathcal Y_h$ with metric
$\min(1,\|u-w\|_{\mathcal Y_h})$,
\begin{equation}
 |\mathbb EF(U)-\nu_{n_{\mathrm{sm}}}(F\circ\iota)|
 \le\|U-V\|_{L^2(\Omega;\mathcal Y_h)}+2\sqrt{\mathcal C}+2\mathcal C.
 \label{eq:companion-law-bound}
\end{equation}
In particular, if the same deterministic indices $n_{\mathrm{sm}}(\Lambda)$
are chosen for two given weak solutions and, for each coupling,
\[
 \Lambda^2\mathbb E\int_0^h
       \|U_t-V_t^{n_{\mathrm{sm}}(\Lambda),\Lambda}\|_{L^2(\mathbb T)}^2dt
       \longrightarrow0\qquad(\Lambda\to\infty),
\]
then the two weak solutions have the same law on $\mathcal X_h$.
\end{proposition}

\subsubsection{From localized columns to spatial trace estimates}
\label{sec:haar-one-l2-trace}
\label{sec:haar-negative-to-negative-weighted-trace}
\label{sec:haar-scale-separated-negative-input}

Fix a profile $\kappa\in\mathcal C$, an integer $j\ge1$, and
bounded scalar directions $f_1=f,f_2,\ldots,f_j$ on $\mathbb T$.
For each vector $e$ of the fixed orthonormal input basis used to
construct $L$, define the columnwise derivative $T_j=T_j(\kappa)$ by
\[
 \bigl(T_j[f_1,\ldots,f_j]\bigr)e
 =\left.\frac{\partial^j}{\partial\varepsilon_1\cdots\partial\varepsilon_j}
       \left(L\!\left(\kappa+\sum_{i=1}^j\varepsilon_i f_i\right)e\right)
   \right|_{\varepsilon_1=\cdots=\varepsilon_j=0}.
\]
The derivative is taken in $H$, keeping $e$ fixed and perturbing
$\kappa$ within strictly positive profiles. Extend $T_j$ linearly to
finite input sums. For $t>0$, its heat regularization is precisely the
$j$th Fr\'echet derivative of $\kappa\mapsto S_tL(\kappa)$ from
Lemma~\ref{lem:analytic-covariance-factor}.
Write $\mathsf A_{\mathrm{col}}=\prod\nolimits_{i=2}^j\|f_i\|_\infty$, with $\mathsf A_{\mathrm{col}}=1$ when $j=1$; the distinguished
input $f_1=f$ is measured separately in the estimates below.

\begin{proposition}[Unified trace transport]
\label{prop:unified-dyadic-column-transport}
Fix $0<\rho<1/2$, $\Lambda\ge1$, and $0<t\le1$.
With $T_j$ and $\mathsf A_{\mathrm{col}}$ as above and a bounded
scalar multiplier $b$, the following estimates hold pathwise, also
for jointly random profiles and inputs:
\begin{align}
 \|Q_\Lambda^{1/2}\mathsf M_bS_tT_j[f,f_2,\ldots]\|_{\rm HS}
 &\le C_j\|b\|_\infty t^{\rho/2-1/4}\|f\|_{L^2(\mathbb T)}\mathsf A_{\mathrm{col}} ,
                                                        \label{eq:haar-one-l2-trace}\\
 \|Q_\Lambda^{1/2}\mathsf M_bS_tT_j[f,f_2,\ldots]\|_{\rm HS}
 &\le C_j\|b\|_\infty t^{-1/4}\|f\|_{H_\Lambda^{-\rho}}\mathsf A_{\mathrm{col}} .
                                                        \label{eq:haar-negative-negative-trace}
\end{align}
For $0<\delta\le w\le1$, let $\mathcal A$ be an integral operator
with a possibly adapted kernel, and let $\mathsf A_{\mathcal A}\ge0$
be an envelope constant such that
$|\mathcal A(x,y)|\le \mathsf A_{\mathcal A}p_{Cw}(x-y)$. Then
\begin{equation}
 \|Q_\Lambda^{1/2}\mathcal AS_\delta
                   T_j[f,f_2,\ldots]\|_{\rm HS}
 \le C_j\mathsf A_{\mathcal A}w^{\rho/2-1/4}\delta^{-\rho/2}
 \sqrt{1+\log(w/\delta)}\,\|f\|_{H_\Lambda^{-\rho}}\mathsf A_{\mathrm{col}} .
 \label{eq:haar-scale-separated-negative-input}
\end{equation}
\end{proposition}

\subsection{Proofs and supporting spatial estimates}
\label{sec:spatial-proof-details}
\subsubsection{Construction and estimates for the covariance factor}
\nopagebreak[4]
\begin{proof}[Proof of Lemma~\ref{lem:analytic-covariance-factor}]

\emph{Construction of the factor.}
Fix a dyadic system independent of $\kappa,\rho,\Lambda$.

For each integer $l\ge0$ let
\[
 \mathcal I_l=\{I_{l,k}=[k2^{-l},(k+1)2^{-l}):0\le k<2^l\}.
\]
If $I\in\mathcal I_l$, its left and right halves are denoted by
$I_{\mathrm L}$ and $I_{\mathrm R}$. For a measurable set $A\subset\mathbb T$, $1_A$ denotes its
indicator and $|A|$ its Lebesgue measure. The real Haar basis of
$H=L^2(\mathbb T)$ is
\begin{equation}
 e_\star(x)=1,\qquad
 e_I(x)=|I|^{-1/2}
            \bigl(1_{I_{\mathrm L}}(x)-1_{I_{\mathrm R}}(x)\bigr),
 \qquad I\in\bigcup_{l=0}^{\infty}\mathcal I_l .
 \label{eq:ordinary-haar-basis}
\end{equation}
Its order is the constant first, then successive generations, with
intervals ordered from left to right within a generation. A generation-$l$ prefix spans the $2^l$ dyadic step functions,
whose union is dense in $H$. The $e_I$ are the orthonormal input
basis; the weighted functions $\chi_I$ below form a separate family.

For $f\in L^1(\mathbb T)$ the averaging operator is
\begin{equation}
 \Pi_l f=\sum_{I\in\mathcal I_l}
             \left(\frac1{|I|}\int_I f(y)\,dy\right)1_I.
 \label{eq:dyadic-averaging-formula}
\end{equation}
On $L^2$ it is the orthogonal projection onto the stated prefix
space, and hence
\begin{align}
 \Pi_l f&=\langle f,e_\star\rangle e_\star
       +\sum_{j=0}^{l-1}\sum_{I\in\mathcal I_j}
                               \langle f,e_I\rangle e_I,\nonumber\\
 \mathsf H_l f:=(\Pi_{l+1}-\Pi_l)f
       &=\sum_{I\in\mathcal I_l}\langle f,e_I\rangle e_I .
 \label{eq:dyadic-haar-projection-formula}
\end{align}
For $l=0$ the first sum is empty. These projections are
self-adjoint contractions, $\Pi_l\Pi_j=\Pi_{\min(l,j)}$, and
$\Pi_l f\to f$ in $L^2$. Interval Poincar\'e and the heat
derivative bound give
\[
 \|S_t\mathsf H_l\|_{\mathrm{op}}
   \le C\min(1,2^{-l}t^{-1/2}).
\]
For a distribution $f\in H_\Lambda^{-\rho}$, $0<\rho<1/2$,
\eqref{eq:dyadic-averaging-formula} is defined by replacing
$\int_I f(y)\,dy$ with the dual pairing $\langle f,1_I\rangle$.
The indicator bound in Lemma~\ref{lem:bessel-properties} makes this
definition meaningful.

For a dyadic interval $I=I_{\mathrm L}\cup I_{\mathrm R}$ put
$\mathfrak m_J=\int_J\kappa(x)\,dx$ and
$\mathfrak m=\mathfrak m_{I_{\mathrm L}}+\mathfrak m_{I_{\mathrm R}}$.
To shorten the next formulas write
$\mathfrak m_{\mathrm L}=\mathfrak m_{I_{\mathrm L}}$ and
$\mathfrak m_{\mathrm R}=\mathfrak m_{I_{\mathrm R}}$. Define
\[
 \chi_I=
   \sqrt{\frac{\mathfrak m_{\mathrm R}}
                    {\mathfrak m_{\mathrm L}\mathfrak m}}\,1_{I_{\mathrm L}}
  -\sqrt{\frac{\mathfrak m_{\mathrm L}}
                    {\mathfrak m_{\mathrm R}\mathfrak m}}\,1_{I_{\mathrm R}},
 \qquad
 \chi_\star=\left(\int_{\mathbb T}\kappa(x)\,dx\right)^{-1/2},
\]
and prescribe the columns by
\begin{equation}
 L(\kappa)e_I=\kappa\chi_I,\qquad
 L(\kappa)e_\star=\kappa\chi_\star.
 \label{eq:haar-local-column}
\end{equation}
The $\chi_I,\chi_\star$ are an orthonormal basis of $L^2(\kappa\,dx)$:
each nonconstant function has weighted mean zero, coarser functions
are constant on finer supports, and the generation spans are the
dense dyadic step spaces.
The map sending $e_I$ to $\sqrt\kappa\chi_I$ is therefore unitary on
$H$. Multiplication by $\sqrt\kappa$ then gives $L(\kappa)$, proving
\eqref{eq:covariance-factor-identity}, including invertibility.
Each column depends continuously in $H$ on $\kappa\in\mathcal C$
with its $L^2$ topology. The uniform operator bound and density of
finite linear combinations imply strong Borel measurability.

\emph{Signed finite-dimensional estimate.}
In the chosen Haar order, $L(\kappa)$ is lower triangular with
positive diagonal: a later column has integral zero against every
earlier Haar function. For an integer prefix size $n_{\mathrm H}\ge1$, its first
$n_{\mathrm H}$ rows contain
no later columns. Let $L_{n_{\mathrm H}}$ be this square compression
and let $\Sigma_{\mathrm H}$ be the corresponding compression of
$\mathsf M_\kappa$. Then
$L_{n_{\mathrm H}}L_{n_{\mathrm H}}^*=\Sigma_{\mathrm H}$.
Differentiating along a bounded direction $f$ gives
\[
 L_{n_{\mathrm H}}^{-1}\partial_\kappa L_{n_{\mathrm H}}[f]
   =\mathcal T\bigl(L_{n_{\mathrm H}}^{-1}
                  (\partial_\kappa\Sigma_{\mathrm H}[f])
                                     L_{n_{\mathrm H}}^{-*}\bigr),
\]
where $\mathcal T$ keeps the strict lower triangle and half the
diagonal. Since $\mathcal T$ is a Frobenius contraction and the
factor and its inverse have uniform operator bounds, integration
along the segment between $\kappa$ and $\widetilde\kappa$ yields,
for $f=\kappa-\widetilde\kappa$,
\begin{equation}
 \|\Pi_l[L(\kappa)-L(\widetilde\kappa)]\|_{\mathrm{HS}}^2
 \le C\|\Pi_l\mathsf M_f\Pi_l\|_{\mathrm{HS}}^2
 =C2^l\|\Pi_lf\|_{L^2(\mathbb T)}^2.
 \label{eq:haar-prefix-signed}
\end{equation}
The last equality holds because the compressed multiplication
operator is diagonal in normalized interval indicators.

\emph{Passage to the heat estimates.}
For $0<\rho<1/2$, duality and the step-function inverse inequality give
\[
 \|\Pi_lf\|_{L^2(\mathbb T)}
   \le C(\Lambda+2^{2l})^{\rho/2}
                                      \|f\|_{H_\Lambda^{-\rho}}.
\]
For a generation-$l$ step function $v$,
$\|v(\cdot+z)-v\|_{L^2(\mathbb T)}^2\le C\min(1,2^l|z|)\|v\|_{L^2(\mathbb T)}^2$.
Integrating against $|z|^{-1-2\rho}$ and adding the mass term
$\Lambda^\rho\|v\|_{L^2(\mathbb T)}^2$ proves the inverse inequality.
For $\rho=0$, use $L^2$ contraction.

Write $A=L(\kappa)-L(\widetilde\kappa)$ and use
$A=\Pi_0A+\sum_{l=0}^{\infty}\mathsf H_lA$. Applying
\eqref{eq:haar-prefix-signed} at level $l+1$ and the cancellation
bound for $S_t\mathsf H_l$ gives
\begin{align*}
 \|S_tA\|_{\mathrm{HS}}
 &\le C\left[\Lambda^{\rho/2}
       +\sum_{l=0}^{\infty}\min(1,2^{-l}t^{-1/2})\,
                2^{l/2}(\Lambda^{\rho/2}+2^{l\rho})\right]
                                      \|f\|_{H_\Lambda^{-\rho}}.
\end{align*}
Splitting the sum at $2^l\asymp t^{-1/2}$ gives
\eqref{eq:haar-heat-signed}. The high-generation sum converges
because $\rho<1/2$. Squaring and integrating the two resulting
powers of $t$ gives
\[
 \int_0^1e^{-2\Lambda t}
        (\Lambda^\rho t^{-1/2}+t^{-1/2-\rho})\,dt
                         \le C\Lambda^{\rho-1/2}.
\]
For $t>1$, factor $S_t=S_{t-1/2}S_{1/2}$ and use contractivity;
the exponentially damped tail obeys the same bound.
This proves \eqref{eq:haar-damped-signed} on the stated full time range.

\emph{Heat-regularized derivatives.}
For every fixed integer $j\ge1$, bounded directions
$f_1,\ldots,f_j\in L^\infty(\mathbb T)$, and $0<t\le1$,
\begin{equation}
 \|S_t\partial_\kappa^jL(\kappa)[f_1,\ldots,f_j]\|_{\mathrm{HS}}
 \le C_j t^{-1/4}(1+|\log t|)^{2j-1}
                                      \prod_{i=1}^j\|f_i\|_\infty.
 \label{eq:factor-heat-derivatives}
\end{equation}
The derivatives in this formula are taken after heat regularization.

The column formulas are smooth functions of the positive interval
masses and of $\kappa$. For example, if $F_J=\int_J f(y)\,dy$,
$F=F_{I_{\mathrm L}}+F_{I_{\mathrm R}}$, and
$\epsilon_I=1_{I_{\mathrm L}}-1_{I_{\mathrm R}}$, then
\begin{align*}
 \partial_\kappa L(\kappa)[f]e_I
   &=f\chi_I+\kappa\,\partial_\kappa\chi_I[f],\\
 \partial_\kappa\chi_I[f]
   &=\frac{\chi_I}{2}
       \left[\epsilon_I
          \left(\frac{F_{I_{\mathrm R}}}{\mathfrak m_{\mathrm R}}
                    -\frac{F_{I_{\mathrm L}}}{\mathfrak m_{\mathrm L}}\right)
                       -\frac{F}{\mathfrak m}\right].
\end{align*}
The triangular map $\mathcal T$ has operator norm at most
$C(1+\log(n_{\mathrm H}+1))$: split the strict lower triangle into
that many dyadic levels. At each level the rectangles have disjoint
row and column supports, so their joint norm is bounded by the
original matrix norm.
Repeated differentiation of
$L_{n_{\mathrm H}}L_{n_{\mathrm H}}^*=\Sigma_{\mathrm H}$, induction,
and polarization now give
\[
 \|\partial_\kappa^jL_{n_{\mathrm H}}(\kappa)
                        [f_1,\ldots,f_j]\|_{\mathrm{op}}
   \le C_j(1+\log(n_{\mathrm H}+1))^{2j-1}
                                      \prod_{i=1}^j\|f_i\|_\infty.
\]
At generation $l$, the Hilbert--Schmidt norm adds at most
$C2^{l/2}$. The preceding heat decomposition bounds the derivative by
\[
 C_j\left[1+\sum_{l=0}^{\infty}
   \min(1,2^{-l}t^{-1/2})\,2^{l/2}(1+l)^{2j-1}\right]
                                      \prod_{i=1}^j\|f_i\|_\infty.
\]
Splitting at $2^l\asymp t^{-1/2}$ proves
\eqref{eq:factor-heat-derivatives}. For fixed $t>0$ the derivative
tails converge locally uniformly on the positive cone. Thus they
are the Fr\'echet derivatives of $S_tL$, and the integral Taylor
formula passes to the limit. Its remainder after degree $j$ has
heat bound $C_jt^{-1/4}(1+|\log t|)^{2j+1}\|f\|_\infty^{j+1}$
along segments contained in $\mathcal C$.

Finally, the spatial regularity used later follows directly from
\eqref{eq:haar-heat-signed}. Choose $0<\varepsilon<1/2-\rho$ and set
$s_{\mathrm{reg}}=1/2-\rho-\varepsilon$. Factor half a heat time
and use the heat derivative bound to obtain
\begin{equation}
 \int_0^\infty e^{-2\Lambda t}
   \|(-\Delta)^{s_{\mathrm{reg}}/2}S_tA\|_{\mathrm{HS}}^2\,dt
 \le C\Lambda^{-\varepsilon}\|f\|_{H_\Lambda^{-\rho}}^2.
 \label{eq:factor-safe-spatial-energy}
\end{equation}
The small-time powers are $t^{-1+\varepsilon}$ and
$\Lambda^\rho t^{-1+\rho+\varepsilon}$; the large-time tail
uses the same semigroup contraction.
\end{proof}
\subsubsection{Bessel kernel properties used in the comparison}
\label{sec:dyadic-bessel-preliminaries}

The heat kernel is the periodization of
$(4\pi t)^{-1/2}e^{-x^2/(4t)}$. In particular $p_t\ge0$,
$\int_{\mathbb T}p_t(x)\,dx=1$, and $p_t*p_s=p_{t+s}$.
Let $\Gamma(\rho)=\int_0^\infty s^{\rho-1}e^{-s}\,ds$ be the
Euler gamma function. The heat representation below is the
Bessel-potential construction; see \cite[Chapter~V]{Stein1970}
for its classical Euclidean form. We give the periodic estimates
with their dependence on $\Lambda$ explicitly.
For $\Lambda\ge1$ and $0<\rho<1/2$ define
\begin{equation}
 \mathsf b_{\Lambda,\rho}(x)
  =\frac1{\Gamma(\rho)}
        \int_0^\infty s^{\rho-1}e^{-\Lambda s}p_s(x)\,ds.
 \label{eq:bessel-kernel-definition}
\end{equation}
The integral defines a nonnegative $L^1$ kernel; its value at
the origin may be infinite. Its Fourier coefficients and mass are
\begin{equation}
 \widehat{\mathsf b}_{\Lambda,\rho}(k)
       =(\Lambda+4\pi^2k^2)^{-\rho},\qquad
 \int_{\mathbb T}\mathsf b_{\Lambda,\rho}(x)\,dx=\Lambda^{-\rho}.
 \label{eq:bessel-symbol-mass}
\end{equation}
Thus $Q_\Lambda f=\mathsf b_{\Lambda,\rho}*f$ is exactly
$(\Lambda-\Delta)^{-\rho}$ and
$\|f\|_{H_\Lambda^{-\rho}}=\|Q_\Lambda^{1/2}f\|_{L^2(\mathbb T)}$.

The Bessel kernel expresses the negative Sobolev norm as a positive
quadratic form. Its interval and heat-row masses give the local
spatial gains used in the transport estimates.

\begin{lemma}[Bessel kernel properties used in the comparison]
\label{lem:bessel-properties}
Constants below may depend on $\rho$, but are uniform for
$\Lambda\ge1$, $0<t\le1$, and intervals $I$ of length $\ell\le1$.
Writing $\operatorname{dist}(x,0)$ for torus distance, one has
\begin{align}
 0\le\mathsf b_{\Lambda,\rho}(x)
   &\le\mathsf b_{1,\rho}(x)
    \le C_\rho\operatorname{dist}(x,0)^{2\rho-1},
                  &&x\ne0,\label{eq:bessel-pointwise}\\
 \sup_y\int_I\mathsf b_{\Lambda,\rho}(x-y)\,dx
   &\le C_\rho\min(\Lambda^{-\rho},\ell^{2\rho}),
       \label{eq:bessel-interval-mass}\\
 \|\mathsf b_{\Lambda,\rho}*p_t\|_\infty
   &\le C_\rho\min(\Lambda^{-\rho}t^{-1/2},
                                      t^{\rho-1/2}).
       \label{eq:bessel-heat-row}
\end{align}
For real $f\in L^2$ the negative norm is the positive-kernel
quadratic form
\begin{equation}
 \|f\|_{H_\Lambda^{-\rho}}^2
  =\int_{\mathbb T^2}\mathsf b_{\Lambda,\rho}(x-y)
                                      f(x)f(y)\,dx\,dy.
 \label{eq:bessel-quadratic-form}
\end{equation}
For $f\in L^2(\mathbb T)$ and bounded scalar $u:\mathbb T\to\mathbb R$,
one also has
\begin{align}
 \|f\|_{H_\Lambda^{-\rho}}
      &\le\Lambda^{-\rho/2}\|f\|_{L^2(\mathbb T)},\nonumber\\
 \|1_I\|_{H_\Lambda^\rho}
      &\le C_\rho\ell^{1/2}
                         (\Lambda^{\rho/2}+\ell^{-\rho}),\nonumber\\
 \|u1_I\|_{H_\Lambda^{-\rho}}
      &\le C_\rho\|u\|_\infty\ell^{1/2}
                   \min(\Lambda^{-\rho/2},\ell^\rho).
 \label{eq:bessel-indicator-bounds}
\end{align}

\end{lemma}
\begin{proof}[Proof of Lemma~\ref{lem:bessel-properties}]
The heat representation and the Gamma integral give
\eqref{eq:bessel-symbol-mass}; Plancherel gives
\eqref{eq:bessel-quadratic-form} and the $L^2$ embedding.
For $d=\operatorname{dist}(x,0)>0$, the periodic Gaussian bound gives
\[
 \mathsf b_{1,\rho}(x)
 \le C_\rho\left(1+\int_0^1s^{\rho-3/2}e^{-d^2/(Cs)}\,ds\right)
 \le C_\rho d^{2\rho-1}.
\]
Together with $\mathsf b_{\Lambda,\rho}\le\mathsf b_{1,\rho}$
and total mass $\Lambda^{-\rho}$, this proves
\eqref{eq:bessel-pointwise}--\eqref{eq:bessel-interval-mass}.
Convolving the heat representation with $p_t$ replaces $p_s$
by $p_{s+t}$. Its supremum is bounded both by
$C\Lambda^{-\rho}t^{-1/2}$ and by
\[
 C_\rho\int_0^\infty s^{\rho-1}e^{-s}
                    [1+(s+t)^{-1/2}]\,ds
 \le C_\rho t^{\rho-1/2},\qquad 0<t\le1,
\]
which proves \eqref{eq:bessel-heat-row}.

For the positive indicator norm, use
$|\widehat{1_I}(k)|\le C\min(\ell,|k|^{-1})$ for $k\ne0$ and
$(\Lambda+4\pi^2k^2)^\rho\le C_\rho(\Lambda^\rho+|k|^{2\rho})$.
Splitting the Fourier sum at $|k|=\ell^{-1}$ gives
$\|1_I\|_{H_\Lambda^\rho}^2
\le C_\rho(\Lambda^\rho\ell+\ell^{1-2\rho})$;
the high-frequency sum converges because $\rho<1/2$.
For the negative norm, \eqref{eq:bessel-quadratic-form} and the
interval mass give
$\|u1_I\|_{H_\Lambda^{-\rho}}^2
\le C_\rho\|u\|_\infty^2\ell\min(\Lambda^{-\rho},\ell^{2\rho})$.
Taking square roots proves \eqref{eq:bessel-indicator-bounds}.
\end{proof}

\begin{samepage}
\begin{corollary}[Spatial trace consequences of the Bessel bounds]
\label{cor:bessel-spatial-traces}
Let $\Lambda\ge1$, $0<\rho<1/2$ and $0<t\le1$, and let
$Q_\Lambda=(\Lambda-\Delta)^{-\rho}$ as above.
\begin{enumerate}[label=(\roman*),leftmargin=*]
\item If $v\in L^2(\mathbb T)$ is supported in an interval $I$
of length $\ell\le1$ and $b$ is bounded, then
\begin{equation}
 \|Q_\Lambda^{1/2}\mathsf M_bS_tv\|_{L^2(\mathbb T)}^2
 \le C_\rho\|b\|_\infty^2
         \min(\ell t^{\rho-1/2},\ell^{2\rho})\|v\|_{L^2(\mathbb T)}^2.
 \label{eq:bessel-supported-trace}
\end{equation}
\item If $L_{\mathrm{bd}}:H\to H$ satisfies
$L_{\mathrm{bd}}L_{\mathrm{bd}}^*\le C I$ and $b\in L^2(\mathbb T)$, then
\begin{equation}
 \|Q_\Lambda^{1/2}\mathsf M_bS_tL_{\mathrm{bd}}\|_{\mathrm{HS}}^2
             \le C_\rho t^{\rho-1/2}\|b\|_{L^2(\mathbb T)}^2 .
 \label{eq:bessel-output-trace}
\end{equation}
\end{enumerate}
Both estimates are pathwise and allow jointly random factors
satisfying the displayed hypotheses.
\end{corollary}
\end{samepage}
\begin{proof}[Proof of Corollary~\ref{cor:bessel-spatial-traces}]
Indeed positivity allows $|S_tv|\le S_t|v|$ in the quadratic
form. The resulting kernel on $I\times I$ is
$\mathsf b_{\Lambda,\rho}*p_{2t}$. Its row mass is bounded
both by $C\ell t^{\rho-1/2}$ and by $C\ell^{2\rho}$, by
\eqref{eq:bessel-heat-row} and the interval mass bound averaged
against the positive heat kernel. Schur's inequality proves
\eqref{eq:bessel-supported-trace}.

To verify this, first remove $L_{\mathrm{bd}}$ by its covariance
bound. The squared trace for $\mathsf M_bS_t$ is bounded by
$\int_{\mathbb T^2}\mathsf b_{\Lambda,\rho}(x-y)|b(x)b(y)|
p_{2t}(x-y)\,dx\,dy$. Apply
$2|b(x)b(y)|\le|b(x)|^2+|b(y)|^2$ and
\eqref{eq:bessel-heat-row} at zero.
All rough multiplications occur inside the positive quadratic
form after the indicated heat smoothing.
\end{proof}
At $\rho=0$, $Q_\Lambda=I$.

\subsubsection{Noise enlargement and the controlled companion}
\begin{lemma}[Exact noise enlargement]
\label{lem:noise-enlargement}
Every prescribed weak solution can be enlarged, without changing its
state law, so that
\begin{equation}
 dU=\Delta U\,dt+m_0dW^0+L(a(U))dB,
 \label{eq:companion-original-rotated}
\end{equation}
where $(W^0,B)$ is a pair of orthogonal cylindrical Brownian motions
in the enlarged filtration.
\end{lemma}
\begin{proof}[Proof of Lemma~\ref{lem:noise-enlargement}]
On a product extension carrying an independent cylindrical Brownian
motion $\overline W$, set
\[
 dW^0=\frac{m_0}{g(U)}dW+
             \frac{\sqrt{a(U)}}{g(U)}d\overline W,
 \qquad
 dW^1=\frac{\sqrt{a(U)}}{g(U)}dW-
             \frac{m_0}{g(U)}d\overline W.
\]
Since $g^2=m_0^2+a$, this is a predictable orthogonal rotation.
L\'evy's characterization applied to each finite set of coordinates
gives joint cylindrical Brownian motion $(W^0,W^1)$, and the
original noise becomes $m_0dW^0+\mathsf M_{\sqrt{a(U)}}dW^1$.
The exact covariance factor is invertible and satisfies
$L(\kappa)L(\kappa)^*=\mathsf M_\kappa$, so
\[
 J(\kappa):=L(\kappa)^{-1}\mathsf M_{\sqrt\kappa},\qquad
 J(\kappa)J(\kappa)^*=J(\kappa)^*J(\kappa)=I.
\]
The rotation $dB=J(a(U))dW^1$ gives joint cylindrical Brownian
motion $(W^0,B)$, and $L(a(U))dB=\mathsf M_{\sqrt{a(U)}}dW^1$ gives
\eqref{eq:companion-original-rotated}. Measurability follows from
the local Haar formulas. The product extension preserves the state law.
\end{proof}

\begin{proof}[Proof of Proposition~\ref{prop:construct-companion}]
The signed heat estimate at $\rho=0$ and Lipschitz continuity of
$a_{n_{\mathrm{sm}}}$ give
\begin{equation}
 \|S_t[L(a_{n_{\mathrm{sm}}}(v))-L(a_{n_{\mathrm{sm}}}(w))]\|_{\mathrm{HS}}
 \le C_{n_{\mathrm{sm}}}  t^{-1/4}\|v-w\|_{L^2(\mathbb T)},
 \qquad 0<t\le1.
 \label{eq:companion-picard-row}
\end{equation}
Use the semigroup $e^{-\Lambda t}S_t$ and prescribed forcing
$\Lambda U$. The squared Lipschitz kernel is $t^{-1/2}$, so the
standard mild Picard argument applies~\cite[Chapter~7]{DaPratoZabczyk2014}:
on an interval of length $b$ its contraction factor is at most
$2C_{n_{\mathrm{sm}}}\sqrt b$. Localization of the forcing and bounded
noise covariance give a continuous global solution by successive
intervals, with pathwise and stopped uniqueness.
The construction interval may depend on $n_{\mathrm{sm}}$.

Without the control, the strong solution is a measurable function
of the initial state and Brownian coordinates, so its state-path law
$\nu_{n_{\mathrm{sm}}}$ is independent of the noise space.
Rotating the residual channel by $J(a_{n_{\mathrm{sm}}}(X))^{-1}$
and combining it with $m_0W^0$ gives covariance
$m_0^2+a_{n_{\mathrm{sm}}}(X)$. L\'evy's characterization therefore
identifies $\nu_{n_{\mathrm{sm}}}$ with the ordinary smooth equation
of diffusion $\sqrt{m_0^2+a_{n_{\mathrm{sm}}}}$.
\end{proof}
\subsubsection{The small-control law comparison}
\nopagebreak[4]
\begin{proof}[Proof of Proposition~\ref{prop:law-comparison}]
Subtracting \eqref{eq:companion-controlled} from
\eqref{eq:companion-original-rotated} gives the exact equation
\begin{equation}
 D_t=\int_0^t e^{-\Lambda(t-r)}S_{t-r}
       [L(a(U_r))-L(a_{n_{\mathrm{sm}}}(V_r))]dB_r,\qquad D_0=0.
 \label{eq:companion-difference}
\end{equation}
The additive channel cancels. Stochastic isometry, Tonelli and
\eqref{eq:haar-damped-signed} imply
\begin{equation}
 \|U-V\|_{\mathscr L_0^2(H)}^2\le C\Lambda^{\rho-1/2}
       \mathbb E\int_0^h
        \|a(U_r)-a_{n_{\mathrm{sm}}}(V_r)\|_{H_\Lambda^{-\rho}}^2dr.
 \label{eq:companion-signed-energy}
\end{equation}
Including an additional nonnegative approximation source $\varepsilon_{\mathrm{app}}$
on the right gives the convention $\mathsf d_{\mathrm{coef}}(U,V)^2=J_\Lambda+\varepsilon_{\mathrm{app}}$ used below.
Bounded covariance in \eqref{eq:companion-difference} also gives
$\|U-V\|_{\mathscr L_0^2(H)}^2\le C_h\Lambda^{-1/2}$. This preliminary estimate supplies the bounded moments used below.

Use $\vartheta_t=\Lambda D_t/m_0$ and the stopping time
\[
 \tau_{\mathrm{cap}}=\inf\left\{t:\int_0^t\|\vartheta_r\|_{L^2(\mathbb T)}^2dr\ge1\right\},
 \qquad\inf\varnothing=\infty.
\]
The coupling still satisfies \eqref{eq:companion-controlled};
only the measure-change integrand is stopped.
Define
\[
 Z=\exp\left(-\int_0^{h\wedge\tau_{\mathrm{cap}}}\langle\vartheta_r,dW_r^0\rangle
       -\frac12\int_0^{h\wedge\tau_{\mathrm{cap}}}\|\vartheta_r\|_{L^2(\mathbb T)}^2dr\right),
 \qquad d\mathbb Q=Z\,d\mathbb P.
\]
The bracket is at most one. Novikov's criterion makes $Z$ mean one,
and the cylindrical Girsanov theorem
\cite[Section~10.2]{DaPratoZabczyk2014} makes
$\widetilde W_t^0=W_t^0+\int_0^{t\wedge\tau_{\mathrm{cap}}}\vartheta_rdr$
Brownian under $\mathbb Q$. The residual $B$ remains Brownian
jointly with it, since the measure-change direction is in the
orthogonal additive channel.

Under $\mathbb Q$ the drift of $V$ is removed up to $\tau_{\mathrm{cap}}$.
Let $X$ be the uncontrolled smooth solution driven by
$(\widetilde W^0,B)$ on the whole interval. Stopped uniqueness in
Proposition~\ref{prop:construct-companion} gives $V=X$ up to
$\tau_{\mathrm{cap}}$. The law of $X$ is $\nu_{n_{\mathrm{sm}}}$.

Take total variation to mean $\sup_E|\mathbb P(E)-\mathbb Q(E)|$.
Writing $H(\mathbb P\mid\mathbb Q)
=\mathbb E_{\mathbb P}\log(d\mathbb P/d\mathbb Q)$ for relative
entropy, the mean-zero stopped stochastic integral gives
\[
 H(\mathbb P\mid\mathbb Q)
 =\frac12\mathbb E_{\mathbb P}\int_0^{h\wedge\tau_{\mathrm{cap}}}
                    \|\vartheta_r\|_{L^2(\mathbb T)}^2dr\le\frac{\mathcal C}{2},
 \qquad \|\mathbb P-\mathbb Q\|_{\mathrm{TV}}
                    \le\frac12\sqrt{\mathcal C}.
\]
Also $\mathbb P(\tau_{\mathrm{cap}}\le h)\le\mathcal C$ by Markov's inequality,
so $\mathbb Q(\tau_{\mathrm{cap}}\le h)\le\mathcal C+\sqrt{\mathcal C}/2$.
For $|F|\le1$, compare $V$ under the two measures and then use
$V=X$ off the exit event. The total error is at most
$2\sqrt{\mathcal C}+2\mathcal C$. Add
\eqref{eq:framework-bl} to obtain \eqref{eq:companion-law-bound}.

For two solutions, a common deterministic smoothing sequence
compares both with $\nu_{n_{\mathrm{sm}}(\Lambda)}$.
The triangle inequality and Lemma~\ref{lem:integrated-laws}
give equality of their continuous path laws.
\end{proof}
\subsubsection{Localized columns and trace transport}
At generation $l\ge0$, put
$\ell=2^{-l}$ and
$v_I=T_j[f,f_2,\ldots,f_j]e_I$ for $I\in\mathcal I_l$.
Thus sums over $I\in\mathcal I_l$ range over the $2^l$ spatial
dyadic intervals $[k2^{-l},(k+1)2^{-l})$, $k=0,\ldots,2^l-1$,
defined in \eqref{eq:ordinary-haar-basis}. Let $x_I$ be the midpoint
of $I$, and write $1_I$ for its indicator function. Write
\[
 \mathscr S_{\mathrm x}=\{2^{-l}:l=0,1,2,\ldots\}
\]
for the set of spatial dyadic lengths used in sums over $\ell$.
The constant column $e_\star$ is treated separately as an additional
length-one column and is estimated without cancellation.

\begin{lemma}[Localization and cancellation of differentiated columns]
\label{lem:localized-haar-heat-packets}
Differentiating the covariance factor preserves the spatial support
and zero integral of each nonconstant Haar column. More precisely,
$v_I$ is supported in $I$, $\int_Iv_I=0$, and
\begin{equation}
 |v_I(x)|\le C_j\mathsf A_{\mathrm{col}} \ell^{-1/2}1_I(x)
                   \left(|f(x)|+\ell^{-1}\!\int_I|f(y)|\,dy\right).       \label{eq:haar-one-slot-column}
\end{equation}
When the interval is smaller than the heat scale, $\ell\le\sqrt t$,
its zero integral gives the additional decay
\[
 |S_tv_I(x)|\le C\ell t^{-1/2}p_{Ct}(x-x_I)\|v_I\|_1.
\]
Thus each column is a localized, mean-zero input whose heat evolution
is concentrated near $I$. These are the two features used to sum
columns in the following trace estimates.
\end{lemma}
\begin{proof}[Proof of Lemma~\ref{lem:localized-haar-heat-packets}]

The local formula is linear in $\kappa (x)$ and otherwise depends only on
two child masses. Each differentiated term therefore contains $f$
either locally times a bounded child step function, or in one child
integral times a bounded local numerator. Thus each term has at most
one unaveraged slot. Ellipticity bounds normalized mass derivatives.
This proves the pointwise and $L^p$ column bounds; differentiation
of $\int_{\mathbb T}\kappa(x)\chi_I(x)\,dx=0$ proves cancellation.
For $\ell\le\sqrt t$, write
$S_tv_I(x)=\int_I[p_t(x-y)-p_t(x-x_I)]v_I(y)\,dy$.
The heat-kernel derivative bound and $|y-x_I|\le\ell$ give the
stated decay.
\end{proof}

\begin{proof}[Proof of Proposition~\ref{prop:unified-dyadic-column-transport}]

\emph{Localized estimates used in the proof.}
Use $0<\rho<1/2$, $\Lambda\ge1$, $0<t\le1$.
Below $u,b$ are bounded scalar functions and $v$ is supported in $I$.
The column and localization bounds needed for the traces are
\begin{align}
 \sum_{I\in\mathcal I_l}\|v_I\|_{L^2(\mathbb T)}^2&\le C_j\mathsf A_{\mathrm{col}} ^2\ell^{-1}\|f\|_{L^2(\mathbb T)}^2,
 \qquad \sum_{I\in\mathcal I_l}\|v_I\|_1^2\le C_j\mathsf A_{\mathrm{col}} ^2\|f\|_{L^2(\mathbb T)}^2,
                                                        \label{eq:haar-one-slot-sums}\\
 \sum_{I\in\mathcal I_l}\|1_If\|_{H_\Lambda^{-\rho}}^2
   &\le C\|f\|_{H_\Lambda^{-\rho}}^2,                      \label{eq:haar-negative-disjoint-localization}\\
 \|1_I\|_{H_\Lambda^\rho}
   &\le C\ell^{1/2}(\Lambda^{\rho/2}+\ell^{-\rho}),        \label{eq:haar-indicator-positive}\\
 \|u1_I\|_{H_\Lambda^{-\rho}}
   &\le C\|u\|_\infty\ell^{1/2}
                          \min(\Lambda^{-\rho/2},\ell^\rho),\label{eq:haar-indicator-negative}\\
 \|v_I\|_{H_\Lambda^{-\rho}}
   &\le C_j\mathsf A_{\mathrm{col}} \ell^{-1/2}\|1_If\|_{H_\Lambda^{-\rho}}.      \label{eq:haar-negative-column}
\end{align}
For $v$ supported in $I$, positivity of the Bessel kernel gives
\begin{equation}
 \|Q_\Lambda^{1/2}\mathsf M_bS_tv\|_{L^2(\mathbb T)}^2
 \le C\|b\|_\infty^2
       \min(\ell t^{\rho-1/2},\ell^{2\rho})\|v\|_{L^2(\mathbb T)}^2.
 \label{eq:haar-one-slot-coarse}
\end{equation}
For $\ell\ge\sqrt t$ its negative-input version is
\begin{equation}
 \|Q_\Lambda^{1/2}\mathsf M_bS_tv\|_{L^2(\mathbb T)}
 \le C\|b\|_\infty\min(\Lambda^{-\rho/2},\ell^\rho)
       (\Lambda^{\rho/2}+t^{-\rho/2})\|v\|_{H_\Lambda^{-\rho}}.
 \label{eq:haar-negative-coarse-local}
\end{equation}
For $\ell\le\sqrt t$, cancellation instead gives
\begin{align*}
 |S_tv_I(x)|\le C_j\mathsf A_{\mathrm{col}} \ell(\Lambda^{\rho/2}+\ell^{-\rho})
 t^{-1/2}p_{Ct}(x-x_I)\|1_If\|_{H_\Lambda^{-\rho}},
 \\
 |S_tv_I(x)|\le C\ell t^{-1/2}p_{Ct}(x-x_I)\|v_I\|_1.
\end{align*}

The column sums follow from Lemma~\ref{lem:localized-haar-heat-packets};
the remaining bounds follow as below.

Indicators act boundedly on $H_\Lambda^\rho$, uniformly in interval
and mass. In the Gagliardo form the extra boundary term is the
fractional Hardy weight $\operatorname{dist}(x,\partial I)^{-2\rho}$.
To prove its bound, near an endpoint average
$|u(x)|^2\le(1+\eta)|u(y)|^2+(1+\eta^{-1})|u(x)-u(y)|^2$
over $2x\le y\le3x$, then integrate against $x^{-2\rho}dx$.
The first coefficient is
$(1+\eta)(3^{2\rho}-2^{2\rho})/(2\rho)<1$ for small $\eta$;
absorb it; the part away from the endpoint is bounded in $L^2$.
The second term is bounded
by the Gagliardo form since $|x-y|\asymp x$.
For functions $u_I$ supported in the disjoint intervals
$I\in\mathcal I_l$, $|u_I(x)-u_J(y)|^2\le
2|u_I(x)|^2+2|u_J(y)|^2$ gives
$\|\sum_{I\in\mathcal I_l}1_Iu_I\|_{H_\Lambda^\rho}^2\le
C\sum_{I\in\mathcal I_l}\|u_I\|_{H_\Lambda^\rho}^2$; duality proves disjoint
negative localization. Indicator scaling proves the positive bound;
the Bessel row mass on $I$ is
$C\min(\Lambda^{-\rho},\ell^{2\rho})$, proving the negative bound.
For a mass term with child interval $J\in\{I_{\mathrm L},I_{\mathrm R}\}$, use
$\ell^{-1}|\int_J f(y)\,dy|\le
C\ell^{-1/2}(\Lambda^{\rho/2}+\ell^{-\rho})\|1_If\|_{H_\Lambda^{-\rho}}$.
Multiply by the bounded numerator estimate, noting
$(\Lambda^{\rho/2}+\ell^{-\rho})
\min(\Lambda^{-\rho/2},\ell^\rho)\le2$.
The indicator multiplier gives the negative column bound for the local terms.

By Lemma~\ref{lem:bessel-properties}, the positive-kernel quadratic form, dominated by $|v|$, has Schur
row mass
\[
                         C\min(\ell t^{\rho-1/2},\ell^{2\rho})
\]
after heat regularization, proving \eqref{eq:haar-one-slot-coarse} even with $b$.
Heat smoothing of a supported negative input costs
$\mathsf s_t=\Lambda^{\rho/2}+t^{-\rho/2}$ and is localized to its
$\sqrt t$ enlargement with Gaussian tails. For $n\ge0$, let $E_n$ be the shell at distance comparable to
$n\sqrt t$ from $I$, with $E_0$ including $I$. For a dual test
$g\in L^2(\mathbb T)$ choose a smooth $\psi=1$ on $I$, supported
in its $\sqrt t/4$ enlargement, with $|\partial\psi|\le Ct^{-1/2}$.
For $\ell\ge\sqrt t$, heat-kernel derivative bounds and Schur's
test give
\begin{align*}
 \|\psi S_t(1_{E_n}g)\|_{L^2(\mathbb T)}&\le Ce^{-c n^2}\|g\|_{L^2(\mathbb T)},\\
 \|\partial_x[\psi S_t(1_{E_n}g)]\|_{L^2(\mathbb T)}
 &\le Ct^{-1/2}e^{-c n^2}\|g\|_{L^2(\mathbb T)}.
\end{align*}
Enlarge constants for $n\le2$.
Interpolation and the massive $H^\rho$ norm therefore bound their
$H_\Lambda^\rho$ norm by $Ce^{-c n^2}\mathsf s_t\|g\|_{L^2(\mathbb T)}$.
Duality for $v$ supported in $I$ yields
$\|1_{E_n}S_tv\|_{L^2(\mathbb T)}\le Ce^{-c n^2}\mathsf s_t\|v\|_{H_\Lambda^{-\rho}}$.
Since $|E_n|\le C(\ell+\sqrt t)$, these bounds also sum in
shell-weighted $L^1$, as used for the middle group below.
For $\ell\ge\sqrt t$, apply
the Bessel row bound on each enlarged shell and sum its Gaussian
tail, proving \eqref{eq:haar-negative-coarse-local}.
If $\ell\le\sqrt t$, subtract $p_t(x-x_I)$; the indicator
and Gagliardo estimates give
\[
 \|1_I[p_t(x-\cdot)-p_t(x-x_I)]\|_{H_\Lambda^\rho}
 \le C\ell^{3/2}(\Lambda^{\rho/2}+\ell^{-\rho})
 t^{-1/2}p_{Ct}(x-x_I).
\]
Duality and the negative column bound give the first fine packet;
the mean-value bound gives its $L^1$ version. The constant column
uses only the coarse estimates. All these arguments are pathwise.

\emph{Summing columns over the spatial scales.}

For the $L^2$ input, the coarse squares sum to
\[
 C\|b\|_\infty^2\mathsf A_{\mathrm{col}} ^2\|f\|_{L^2(\mathbb T)}^2
 \sum_{\substack{\ell\in\mathscr S_{\mathrm x}\\\ell\ge\sqrt t}}\ell^{2\rho-1}.
\]
For fine columns the $L^1$ packet and the Bessel form of $p_{Ct}$ give
\[
 C\|b\|_\infty^2\mathsf A_{\mathrm{col}} ^2\|f\|_{L^2(\mathbb T)}^2
 t^{\rho-3/2}\sum_{\substack{\ell\in\mathscr S_{\mathrm x}\\\ell\le\sqrt t}}\ell^2.
\]
Both sums are $O(t^{\rho-1/2})$, proving the first trace.

For the two-scale estimate put
$H_\rho(w)=\min(\Lambda^{-\rho}w^{-1/2},w^{\rho-1/2})$
and $\mathsf s_\delta=\Lambda^{\rho/2}+\delta^{-\rho/2}$.
Divide lengths into $\ell\ge\sqrt w$,
$\sqrt\delta\le\ell\le\sqrt w$, and $\ell\le\sqrt\delta$.
The coarse localized smoothing packet, followed by the positive
envelope of $\mathcal A$, yields a square sum bounded by
\[
 C\mathsf A_{\mathcal A}^2\mathsf A_{\mathrm{col}} ^2\|f\|_{H_\Lambda^{-\rho}}^2\mathsf s_\delta^2
 \sum_{\substack{\ell\in\mathscr S_{\mathrm x}\\\ell\ge\sqrt w}}\ell^{-1}
                        \min(\Lambda^{-\rho},\ell^{2\rho})
 \le C\mathsf A_{\mathcal A}^2\mathsf A_{\mathrm{col}} ^2\|f\|_{H_\Lambda^{-\rho}}^2
                  \mathsf s_\delta^2H_\rho(w).
\]
The sum is split at $\ell=\Lambda^{-1/2}$; $2\rho-1<0$ makes
both tails geometric. For each middle generation,
$\|S_\delta v_I\|_1\le
C\ell^{1/2}\mathsf s_\delta\|v_I\|_{H_\Lambda^{-\rho}}$,
including Gaussian-weighted shells, so the output is at most
$C\mathsf A_{\mathcal A}\mathsf A_{\mathrm{col}} \mathsf s_\delta\|1_If\|_{H_\Lambda^{-\rho}}
p_{C'w}(x-x_I)$. Its Bessel square costs $H_\rho(w)$;
disjoint localization and the $O(1+\log(w/\delta))$ middle
generations give that extra factor. For fine generations cancellation
occurs before $\mathcal A$; the fine packet followed by its positive
envelope gives the same prefactor times
\[
 H_\rho(w)\delta^{-1}
 \sum_{\substack{\ell\in\mathscr S_{\mathrm x}\\\ell\le\sqrt\delta}}\ell^2(\Lambda^\rho+\ell^{-2\rho})
 \le CH_\rho(w)(\Lambda^\rho+\delta^{-\rho}).
\]
All three groups are therefore bounded by the square of the input
times $C \mathsf A_{\mathcal A}^2\mathsf A_{\mathrm{col}} ^2(1+\log(w/\delta))
(\Lambda^\rho+\delta^{-\rho})H_\rho(w)$.
Splitting at $\Lambda w=1$ and using $\delta\le w$ bounds the last
product by $2w^{\rho-1/2}\delta^{-\rho}$.
This proves the two-scale bound. Taking
$\mathcal A=\mathsf M_bS_{t/2}$ and $w=\delta=t/2$ gives the
matched trace; summing columns removes the finite-prefix restriction.
\end{proof}

The localization lemma and the transport proposition also hold for fixed
finite Euclidean input and output dimensions. Use the scalar Bessel
kernel times identity, vector norms for columns, Frobenius norms for
the distinguished input, and an operator-norm envelope for a matrix
left kernel. Their proofs use only Cauchy--Schwarz, scalar positive
kernels, and the displayed column hypotheses, so summing components
proves these versions with dimension-dependent constants. The scalar comparison below uses the one-dimensional version.

\begin{corollary}[Adapted stochastic integrals with one spatial input]
\label{cor:adapted-one-input-ito}
Fix $T>0$, $0<h\le\min(1,T)$, $\Lambda\ge1$,
$0<\rho<1/2$ and an integer $j\ge1$. Use the column derivative
$T_j=\partial_\kappa^j L(\kappa)$ from
Proposition~\ref{prop:unified-dyadic-column-transport}, with all
inputs except the displayed $f_s$ uniformly bounded. Assume that
the stochastic integrands below are predictable. Let $\lambda\ge1$
be a state frequency and $b_{\lambda,t}$ a scalar multiplier.
Use the window norms of Definition~\ref{def:window-integrated-norms}.
If
$\|b_{\lambda,t}\|_\infty\le
C\lambda^{1-\beta}e^{-c\lambda^2\sqrt t}$ and the one-slot
input $f_s$ satisfies
$\sup_{0\le s\le h}\|f_s\|_{\mathscr L_h^2(H)}\le C\|U-V\|_{\mathscr L_0^2(H)}^u$
for an exponent $u\ge0$,
the first trace and isometry give
\begin{equation}
 \left\|\int_0^h Q_\Lambda^{1/2}
 \mathsf M_{b_{\lambda,h-s}}S_{h-s}T_j[f_s,\ldots]dB_s\right\|_{\mathscr L_h^2(H)}
 \le C_j\lambda^{-\beta-2\rho}\|U-V\|_{\mathscr L_0^2(H)}^u.
 \label{eq:haar-one-slot-nonlinear-ito}
\end{equation}
For a frequency $\nu\ge1$ and exponent $\gamma>0$, the matched
negative trace similarly gives $C\nu^{-\gamma}
\sup_{0\le s\le h}\|f_s\|_{\mathscr L_h^2(H_\Lambda^{-\rho})}$ for coefficient
amplitude $\nu^{1-\gamma}e^{-c\nu^2\sqrt t}$, since its squared
time integral uses $t^{-1/2}$. The same bounds hold after integration over the window start $t_0$.
\end{corollary}
\begin{proof}[Proof of Corollary~\ref{cor:adapted-one-input-ito}]
Apply the corresponding trace bounds of
Proposition~\ref{prop:unified-dyadic-column-transport} and stochastic isometry.
Indeed the squared scalar integral is
$\lambda^{2-2\beta}\int_0^h
e^{-c\lambda^2\sqrt t}t^{\rho-1/2}dt\lesssim
\lambda^{-2\beta-4\rho}$.
For the matched negative trace, the time density is $t^{-1/2}$;
the same change of variables gives $\nu^{-2\gamma}$ before taking
square roots. Moving-past integrated bounds follow
by Tonelli and deterministic time translation.
\end{proof}

\subsubsection{The final Gaussian derivative and actual noise}

\begin{corollary}[The last Gaussian derivative with a signed noise input]
\label{cor:last-derivative-signed-input}
Use $\Lambda\ge1$, $0<\rho<1/2$ and the column derivative $T_j$
of Proposition~\ref{prop:unified-dyadic-column-transport}. Write
$\mathsf A_{\mathrm{col}}=\prod_{i=2}^j\|f_i\|_\infty$
for the bounded other inputs, with empty product one.
A native word is one of the stochastic terms in the finite expansion
of Section~\ref{sec:recursive-gaussian-expansion}. Let $h$ be its
terminal time, $q<h$ its last stochastic integration time, and
$0<\delta\le w=h-q$ its gap of the last Gaussian derivative.
Suppose that, after removing the deterministic word factors, its
spatial operator has the factorization
$\mathscr K=\widetilde{\mathscr K}S_{\delta/2}$ and envelope
$|\widetilde{\mathscr K}(x,y)|\le C_{\rm word}(\omega)
\delta^{-1/2}w^{1/4}p_{Cw}(x-y)$.
Here $\omega\in\Omega$ is the sample point and $C_{\rm word}$ is
a nonnegative random envelope.
For every such factorization, the two-scale trace gives
\begin{equation}
 \|Q_\Lambda^{1/2}\widetilde{\mathscr K}S_{\delta/2}
 T_j[f,f_2,\ldots]\|_{\rm HS}
 \le CC_{\rm word}\delta^{-(1+\rho)/2}w^{\rho/2}
 \sqrt{1+\log(w/\delta)}\,\|f\|_{H_\Lambda^{-\rho}}\mathsf A_{\mathrm{col}}.
 \label{eq:haar-native-negative-input-score}
\end{equation}
Its last-gap Dirichlet degree is $(1-\rho)/2>0$; the additional
loss is $\delta^{-\rho/2}$.
\end{corollary}
\begin{proof}[Proof of Corollary~\ref{cor:last-derivative-signed-input}]
Apply the separated-scale bound of
Proposition~\ref{prop:unified-dyadic-column-transport} with heat
scales $w$ and $\delta/2$ and envelope
$C_{\rm word}\delta^{-1/2}w^{1/4}$.
The factors $w^{1/4}$ and $w^{\rho/2-1/4}$ combine to
$w^{\rho/2}$; the powers of $\delta$ combine to
$\delta^{-(1+\rho)/2}$. This proves
\eqref{eq:haar-native-negative-input-score}.
\end{proof}

The actual-noise comparison below combines this trace with clipped
coefficient differences and the normalized-kernel moment truncation of
Lemma~\ref{lem:normalized-kernel-truncation}. Its complete
square-time argument is given in the native comparison table.

\clearpage
\section{Gaussian references and finite expansions}
\label{sec:gaussian-references}
The Gaussian comparison has three main outputs. We first define the
state bands, then state these outputs with their input profiles,
reference fields and word operators. Throughout the window comparisons,
$t_0$ is the physical starting time, $h$ is the duration, and a local
time $s$ corresponds to physical time $t_0+s$.
The fixed point construction and proofs of the three outputs
are in Sections~\ref{sec:gaussian-proof-details}--\ref{sec:finite-expansion-proof}.
Section~\ref{sec:gaussian-input-estimates} records the Gaussian derivative and
kernel estimates used later; their technical proofs are in
Appendices~\ref{app:gaussian-geometry}--\ref{app:time-kernels}.
\begin{enumerate}[label=(\roman*),leftmargin=*]
\item \textbf{Compatible Gaussian references.} Proposition~\ref{prop:self-profile-compatible-gaussian} constructs a common-noise Gaussian pair with exactly centered coefficient errors and bounds its difference.
\item \textbf{The total Gaussian comparison error.} Theorem~\ref{thm:coupled-gaussian-band-sum} bounds the sum of all high state-frequency bands by the physical solution difference, using the annular spatial gain.
\item \textbf{An exact finite expansion.} Lemma~\ref{lem:finite-expansion} decomposes a coefficient observation into Gaussian average terms, martingale words and a finite remainder.
\end{enumerate}
The finite expansion supplies the exact decomposition of physical
coefficient observations. The Gaussian band estimates control the
leading reference terms when the expansions for $U$ and $V$ are compared;
Proposition~\ref{prop:joint-first-pair} makes this connection explicit.

\subsection{Main results}

We use the heat semigroup $S_t=e^{t\Delta}$ and its kernel $p_t$
on the normalized torus $\mathbb T$ from the preceding sections.
Unless a different domain is displayed, $\|\cdot\|_p$ denotes the
spatial $L^p(\mathbb T)$ norm; for a function of the state variable,
$\|\cdot\|_\infty$ denotes its supremum over $\mathbb R$.
We use the negative Sobolev norm
$\|f\|_{H^{-\rho}_\Lambda}=\|(\Lambda-\Delta)^{-\rho/2}f\|_{L^2(\mathbb T)}$,
with $0<\rho<1/2$ and damping parameter $\Lambda\ge1$.
The residual coefficient $a$, its H\"older exponent $\beta\in(0,1)$,
the bound $L_a\ge[a]_{C^\beta}$, and the smooth approximations
$a_{n_{\mathrm{sm}}}$ with errors
$\epsilon_{n_{\mathrm{sm}}}=\|a_{n_{\mathrm{sm}}}-a\|_\infty$
are as in Section~\ref{sec:framework}.
Here $[a]_{C^\beta}=\sup_{z\ne z'}|a(z)-a(z')|/|z-z'|^\beta$.
Constants and adjoints follow the conventions of Section~1.
For symmetric matrices, $A\succeq B$ means that $A-B$ is positive semidefinite.

\subsubsection{State-frequency notation}
\label{sec:state-bands}
We write $z\in\mathbb R$ for the state variable and $x\in\mathbb T$
for the spatial variable. We use the Littlewood--Paley decomposition
and Bernstein estimates of \cite[Chapter~2]{BahouriCheminDanchin2011}.
Choose a smooth inhomogeneous Fourier
partition of unity, with smooth compactly supported multipliers
$\chi,\psi$, Fourier variable $\xi\in\mathbb R$, and fixed
annular radii $0<c<C<\infty$:
\[
 \chi(\xi)+\sum_{j=0}^{\infty}\psi(2^{-j}\xi)=1,
 \qquad \operatorname{supp}\psi\subset\{c\le|\xi|\le C\}.
\]
Write $a_{\mathrm{low}}=\chi(-i\partial_z)a$ and
$a_\lambda=\psi(-i\partial_z/\lambda)a$, where
\[
 \lambda\in\mathscr S_{\mathrm z}:=\{2^j:j=0,1,2,\ldots\}.
\]
Thus $\mathscr S_{\mathrm z}$ indexes state-frequency bands, whereas
$\mathscr S_{\mathrm x}=\{2^{-l}:l=0,1,\ldots\}$ indexes spatial interval lengths. The notation $\chi(-i\partial_z)$ denotes the Fourier multiplier
with symbol $\chi$, and similarly for $\psi(-i\partial_z/\lambda)$.
The annular kernels have integral zero. Convolution against
$a(z-y)-a(z)$ consequently gives,
for every nonnegative integer derivative order $r$ and dyadic
upper cutoff $\lambda_{\mathrm{cut}}\in\mathscr S_{\mathrm z}$,
\begin{equation}
 \|\partial_z^r a_\lambda\|_\infty
       \le C_r L_a\lambda^{r-\beta},\qquad
 \Bigl\|a-a_{\mathrm{low}}-
              \sum_{\substack{\lambda\in\mathscr S_{\mathrm z}\\\lambda\le\lambda_{\mathrm{cut}}}}a_\lambda\Bigr\|_\infty
       \le C L_a \lambda_{\mathrm{cut}} ^{-\beta}.
 \label{eq:bands-bernstein-tail}
\end{equation}
The series converges uniformly. The low block is smooth and bounded.
The same estimates hold uniformly for $a_{n_{\mathrm{sm}}}$. For differences,
$\|\partial_z^r(a_{n_{\mathrm{sm}}}-a)_\lambda\|_\infty
\le C_r\lambda^r\epsilon_{n_{\mathrm{sm}}}$; this last estimate is used only
on finitely many bands.

\subsubsection{Compatible Gaussian references}

The comparison uses a Gaussian pair with exactly centered coefficient
errors and controlled dependence on the starting states. We describe
its law first, then define its profile and state the comparison bounds.

For the scalar profiles fix $0<\beta<1$, additive amplitude $m_0>0$, and
$a\in C_b^\beta(\mathbb R;\mathbb R)$ with
$0<\underline a\le a\le \overline a$.

\emph{The reference as a Gaussian cloud.}
Fix a deterministic starting state $w\in H$, and let $W^{\mathrm{ref}}$
be a cylindrical Brownian motion on $H=L^2(\mathbb T)$.
The desired reference satisfies, in the mild sense,
\begin{equation}
 dX_t=\Delta X_t\,dt+
 \sqrt{m_0^2+\mathbb E[a(X_t(\cdot))]}\,dW_t^{\mathrm{ref}},
 \qquad X_0=w.
 \label{eq:reference-law-spde}
\end{equation}
The square root acts by pointwise multiplication in space.
The expectation averages over the law of $X$ at each spatial point,
so the noise coefficient is deterministic. Any such mild solution
is therefore a Gaussian heat field with mean $S_tw$; the profile
construction below proves short-time existence and uniqueness.

To picture this equation, imagine many independent realizations
starting from $w$. At a fixed point $(t,x)$, their values form a
Gaussian cloud centered at $S_tw(x)$. Evaluate $a$ across this cloud
and average: the resulting number
$\bar a_t(x)=\mathbb E[a(X_t(x))]$ is the residual noise variance
used by every realization at $(t,x)$. The noise determines the width
of the cloud, while that width determines the average of $a$.
The fixed point below makes these two descriptions consistent and
gives the exact centering
$\mathbb E[a(X_t(x))-\bar a_t(x)]=0$.

For example, if $w\equiv c$ for a constant $c\in\mathbb R$ and
$a(z)=2+\cos z$, spatial homogeneity makes the variance
$\upsilon_t=\operatorname{Var}(X_t(x))$ independent of $x$, and
$\bar a_t=2+e^{-\upsilon_t/2}\cos c$.
A narrow cloud samples $a$ near $c$; a wider cloud averages away
more of its oscillation.
This is closely related to Gaussian moment closure for
finite-dimensional diffusions, where the diffusion covariance is
averaged against the Gaussian law being determined; see
\cite[Section~9.1, Algorithm~9.1]{SarkkaSolin2019}.
Equation~\eqref{eq:reference-law-spde} describes the law of one
reference. The Haar representation in
Proposition~\ref{prop:self-profile-compatible-gaussian} also
specifies how two references share noise for the comparison.

For scalar variance $v\ge0$ write
$(\mathcal P_v a)(z)=\mathbb E[a(z+\sqrt v\,\zeta_{\mathrm G})]$, where $\zeta_{\mathrm G}\sim N(0,1)$.

\begin{definition}[Gaussian variance map]
\label{def:gaussian-variance-map}
For any $h\in(0,1]$ and deterministic bounded measurable candidate residual profile
$\kappa:(0,h]\times\mathbb T\to[\underline a,\overline a]$, define
$\mathsf k_t(x)=p_t(x)^2$ and
\begin{equation}
 \mathcal V_r[\kappa](x)
 =\int_0^r\int_{\mathbb T}
       p_{r-s}(x-y)^2\bigl(m_0^2+\kappa_s(y)\bigr)\,dy\,ds,
 \qquad 0<r\le h.
 \label{eq:gaussian-profile-variance-map}
\end{equation}
The input is the entire profile $\kappa$ up to time $r$; the output
is its scalar pointwise variance at $(r,x)$.
\end{definition}

\begin{definition}[Self-consistent reference profile]
\label{def:self-consistent-reference-profile}
For $w\in H$ and $0<h\le1$, a self-consistent reference profile
is a bounded measurable function
$\kappa:(0,h]\times\mathbb T\to[\underline a,\overline a]$
satisfying, with the variance map of
Definition~\ref{def:gaussian-variance-map},
\begin{equation}
 \kappa_r(x)=\bigl(\mathcal P_{\mathcal V_r[\kappa](x)}a\bigr)(S_rw(x)),
 \qquad 0<r\le h.
 \label{eq:reference-profile-self-consistency}
\end{equation}
When this profile exists and is unique, we denote it by $\bar a[w]$.
For a prescribed deterministic mean $\phi_r(x)$, replace $S_rw(x)$
by $\phi_r(x)$ in \eqref{eq:reference-profile-self-consistency};
the unique profile, when it exists, is denoted by $\bar a^\phi$.
Thus $\bar a[w]=\bar a^{S_\cdot w}$.
The next proposition asserts existence for the references used here
and states the properties of the associated Gaussian fields.
\end{definition}

\begin{proposition}[Compatible Gaussian references]
\label{prop:self-profile-compatible-gaussian}
There are $h_0\in(0,1]$ and $C<\infty$, depending only on
$\beta,m_0,\underline a,\overline a$ and $[a]_{C^\beta}$, such that
the following holds for $0<h\le h_0$.
For every $w,\widetilde w\in H$, there are unique self-consistent
reference profiles $\bar a=\bar a[w]$ and
$\widetilde{\bar a}=\bar a[\widetilde w]$ in the sense of
Definition~\ref{def:self-consistent-reference-profile}. On the same pair
of independent cylindrical Brownian motions $(W^0,B)$ on
$L^2(\mathbb T)$, define
\[
 X_r=S_rw+m_0\int_0^r S_{r-z}dW_z^0+\int_0^r S_{r-z}L_{\bar a_z}dB_z,
\]
and $\widetilde X$ with tilded data, where $L_{\bar a} $ is the exact Haar
factor of the multiplication operator $\mathsf M_{\bar a}$.
Here $L_{\bar a}=L(\bar a)$, and $\widetilde X$ uses the same two
noises and replaces $w,\bar a$ by $\widetilde w,\widetilde{\bar a}$.
The pair is jointly Gaussian and, with $\zeta=w-\widetilde w$,
\begin{align}
 \mathbb E a(X_r(x))=\bar a_r(x),\quad
 \mathbb E a(\widetilde X_r(x))&=\widetilde{\bar a}_r(x),
 \label{eq:self-profile-exact-centering}\\
 \mathbb E\|X_h-\widetilde X_h-S_h\zeta\|_{L^2(\mathbb T)}^2
          &\le Ch^{\beta/2}\|\zeta\|_{L^2(\mathbb T)}^2,
 \label{eq:self-profile-noise-stability}\\
 \mathbb E\|X_h-\widetilde X_h\|_{L^2(\mathbb T)}^2&\le C\|\zeta\|_{L^2(\mathbb T)}^2.
 \label{eq:self-profile-total-stability}
\end{align}
The same bounds hold for mean difference
$e^{-\Lambda r}S_r\zeta$ when both fixed-point equations use those means,
with $S_h\zeta$ replaced by $e^{-\Lambda h}S_h\zeta$ in the centered bound.
\end{proposition}

\emph{Why the coefficient error is centered.}
The stochastic integrals defining $X$ have deterministic coefficients.
Thus $X_r(x)$ is Gaussian with mean $S_rw(x)$ and variance
$\mathcal V_r[\bar a](x)$. Consequently,
\[
 \mathbb E a(X_r(x))
 =\bigl(\mathcal P_{\mathcal V_r[\bar a](x)}a\bigr)(S_rw(x))
 =\bar a_r(x),
\]
where the last equality is the self-consistency relation in
Definition~\ref{def:self-consistent-reference-profile}, namely
\eqref{eq:reference-profile-self-consistency}. Thus
$\mathbb E[a(X_r(x))-\bar a_r(x)]=0$; the same calculation applies
to $\widetilde X$.
At a window starting from $U_0$, we take $w=U_0$ after conditioning
on the initial past, so the reference profile is fixed throughout
the window.

\subsubsection{Coupled references and their total Gaussian comparison error}
\label{sec:sharp-annular-main-Gaussian}
The next definition fixes the two reference fields. The theorem then
bounds their combined coefficient error over all high state-frequency
bands. Its spatial gain comes from cancellation on scale
$\lambda^{-2}$; the one-band estimate and covariance calculation
are given in Section~\ref{sec:gaussian-band-estimates-proofs}.

\begin{definition}[Compatible Gaussian references for the coupled solution pair]
\label{def:coupled-gaussian-references}
Fix a horizon $T>0$, a prescribed weak solution $U$, and its controlled
companion $V=V^{n_{\mathrm{sm}},\Lambda}$ from
Proposition~\ref{prop:construct-companion}, with $\Lambda\ge1$.
They use the common Brownian pair $(W^0,B)$ from Section 3.
Put $D=U-V$. On a window $[t_0,t_0+h]\subset[0,T]$, $0<h\le h_0$,
condition on $\mathcal F_{t_0}$ and prescribe the means
\begin{equation}
 \begin{aligned}
 \phi_s^{U,t_0}&:=S_sU_{t_0},\\
 \phi_s^{V,t_0}&:=S_sU_{t_0}-e^{-\Lambda s}S_sD_{t_0}
              =S_sV_{t_0}+(1-e^{-\Lambda s})S_sD_{t_0},
 \qquad 0\le s\le h.
 \end{aligned}
 \label{eq:coupled-reference-means}
\end{equation}
For $i\in\{U,V\}$, let $\bar a^{i,t_0}:=\bar a^{\phi^{i,t_0}}$
be the self-consistent profile of
Definition~\ref{def:self-consistent-reference-profile}, whose existence
and uniqueness are proved in Proposition~\ref{prop:self-profile}.
Both profiles use the residual map $a$; the actual companion $V$
uses $a_{n_{\mathrm{sm}}}$ in its equation.
Write $W_s^{0,t_0}:=W^0_{t_0+s}-W^0_{t_0}$ and $B_s^{t_0}:=B_{t_0+s}-B_{t_0}$.
Define the Gaussian reference fields on these same future noises by
\begin{equation}
 X_s^{i,t_0}:=\phi_s^{i,t_0}
       +m_0\int_0^sS_{s-r}\,dW_r^{0,t_0}
       +\int_0^sS_{s-r}L(\bar a_r^{i,t_0})\,dB_r^{t_0},
 \qquad i\in\{U,V\}.
 \label{eq:coupled-gaussian-reference-fields}
\end{equation}
These fields start at $X_0^{U,t_0}=U_{t_0}$ and $X_0^{V,t_0}=V_{t_0}$.
At a fixed window start $t_0$, abbreviate
$\phi^i=\phi^{i,t_0}$, $\bar a^i=\bar a^{i,t_0}$ and $X^i=X^{i,t_0}$.
Thus the superscripts $U,V$ specify the two prescribed means and
their respective covariance profiles. Conditional on $\mathcal F_{t_0}$,
the pair is jointly Gaussian and has exactly centered errors
$\mathbb E[a(X_s^{i,t_0}(x))-\bar a_s^{i,t_0}(x)\mid\mathcal F_{t_0}]=0$,
as in Proposition~\ref{prop:self-profile-compatible-gaussian}.
\end{definition}

\begin{theorem}[Total Gaussian comparison error above the cutoff]
\label{thm:coupled-gaussian-band-sum}
Fix $0<\beta<1$, a horizon $T>0$, a weak solution $U$ and its
controlled companion $V=V^{n_{\mathrm{sm}},\Lambda}$, with
$n_{\mathrm{sm}}\ge1$ and $\Lambda\ge1$, from
Proposition~\ref{prop:construct-companion}, with their common noise
pair. Write $D=U-V$ and $H=L^2(\mathbb T)$. Choose
$0<\rho<(1-\beta)/2$ and a tail power $p_{\mathrm{tail}}>0$.
For the residual coefficient $a$ of \eqref{eq:framework-residual},
with $[a]_{C^\beta}\le L_a$, let
$a_\lambda=\psi(-i\partial_z/\lambda)a$ be its state-frequency
blocks from Section~\ref{sec:state-bands}, indexed by
$\lambda\in\mathscr S_{\mathrm z}=\{2^j:j=0,1,\ldots\}$.
Set $h_\lambda=\lambda^{-3\beta}$ and
\begin{equation}
 \lambda_*(\Lambda):=2^{\lceil\frac14\log_2\Lambda\rceil},\qquad
 \Lambda^{1/4}\le\lambda_*(\Lambda)<2\Lambda^{1/4}.
 \label{eq:gaussian-band-dyadic-cutoff}
\end{equation}
For every window $[t_0,t_0+h_\lambda]\subset[0,T]$, let
$X^{U,t_0},X^{V,t_0}$ be the fields
\eqref{eq:coupled-gaussian-reference-fields}, with respective means
$S_sU_{t_0}$ and $S_sU_{t_0}-e^{-\Lambda s}S_sD_{t_0}$ and their
self-consistent profiles for the common coefficient $a$, as in
Definition~\ref{def:coupled-gaussian-references}.
Write $X^i=X^{i,t_0}$ when the starting time is suppressed.
Use the spatial norm
$\|f\|_{H_\Lambda^{-\rho}}=\|(\Lambda-\Delta)^{-\rho/2}f\|_H$
and the window norm
$\|Q\|_{\mathscr L_h^2(H_\Lambda^{-\rho})}^2
=\int_0^{T-h}\mathbb E\|Q(t_0)\|_{H_\Lambda^{-\rho}}^2\,dt_0$.
The full-horizon input is
$\|U-V\|_{\mathscr L_0^2(H)}^2
=\int_0^T\mathbb E\|D_t\|_H^2\,dt$.

There are $C<\infty$ and $\Lambda_0\ge1$, depending only on
$T,\beta,\rho,p_{\mathrm{tail}}$ and the common coefficient bounds,
such that for $\Lambda\ge\Lambda_0$ all chosen durations satisfy
$h_\lambda\le\min(h_0,T)$, where $h_0$ is the reference-existence
threshold of Proposition~\ref{prop:self-profile-compatible-gaussian},
and
\begin{equation}
 \begin{aligned}
 &\sum_{\substack{\lambda\in\mathscr S_{\mathrm z}\\\lambda\ge\lambda_*}}
 \|a_\lambda(X_{h_\lambda}^U)-a_\lambda(X_{h_\lambda}^V)
                         \|_{\mathscr L_{h_\lambda}^2(H_\Lambda^{-\rho})}\\
 &\quad=\sum_{\substack{\lambda\in\mathscr S_{\mathrm z}\\\lambda\ge\lambda_*}}
 \left(\int_{h_\lambda}^T\mathbb E
 \left\|a_\lambda(X_{h_\lambda}^{U,t-h_\lambda})
       -a_\lambda(X_{h_\lambda}^{V,t-h_\lambda})
 \right\|_{H_\Lambda^{-\rho}}^2\,dt\right)^{1/2}\\
 &\quad\le C\|U-V\|_{\mathscr L_0^2(H)}^{\beta+2\rho}
                                    +\Lambda^{-p_{\mathrm{tail}}}.
 \end{aligned}
 \label{eq:sharp-main-Gaussian-sum}
\end{equation}
The same bound holds when, in each window, $X_{h_\lambda}^{V,t_0}$
is replaced by
$X_{h_\lambda}^{U,t_0}-e^{-\Lambda h_\lambda}S_{h_\lambda}D_{t_0}$.
\end{theorem}

At physical time $t$, band $\lambda$ uses the reference window
$[t-h_\lambda,t]$. The equality in
\eqref{eq:sharp-main-Gaussian-sum} is the change of variables
$t=t_0+h_\lambda$; the sum is over state-frequency bands, with each
duration already fixed by its band. Higher frequencies use shorter
windows, while their Gaussian damping improves because
$\lambda^2\sqrt{h_\lambda}=\lambda^{2-3\beta/2}\to\infty$.
Thus closeness of the physical solutions controls the total
high-band Gaussian error, up to a tail smaller than any prescribed
power of $\Lambda^{-1}$. The low bands are estimated separately.
Section~\ref{sec:gaussian-band-comparison-proof} proves the theorem
by summing the one-band bound of
Corollary~\ref{cor:coupled-gaussian-band}.

\subsubsection{The finite-word identity}
\label{sec:recursive-gaussian-expansion}
The purpose of this identity is to expand the coefficient observation
$a_\lambda(U_h(x))$ into a Gaussian average, a martingale term and
successive covariance corrections. The following definitions fix its
initial test and all subsequent words. The vector $\mathbf Y_t^Z$
defined below collects observations of the physical solution at the
space--time targets in $Z$.

The process $U$ is the weak solution of \eqref{eq:framework-mild},
equivalently \eqref{eq:framework-observed}. On a window $[t_0,t_0+h]$,
we relabel $U_{t_0+t}$ as $U_t$, $0\le t\le h$, and condition on
$\mathcal F_{t_0}$. Thus $U_0$ is the state at the start of this window.

\begin{definition}[Data and initial test for the coefficient-band expansion]
\label{def:coefficient-band-data}
Fix $0<h\le h_0$, $x\in\mathbb T$ and
$\lambda\in\mathscr S_{\mathrm z}$. The residual covariance map is
$a:\mathbb R\to[\underline a,\overline a]$,
$a(z)=g(z)^2-m_0^2$, as in \eqref{eq:framework-residual}.
Its state-frequency band is the smooth function
$a_\lambda=\psi(-i\partial_z/\lambda)a:\mathbb R\to\mathbb R$
from Section~\ref{sec:state-bands}, with bounds
\eqref{eq:bands-bernstein-tail}. Define the initial observation list
and initial test by
\begin{equation}
 \boxed{\quad Z_0:=((h,x)),\qquad
 F_0:\mathbb R\to\mathbb R,\quad F_0(z):=a_\lambda(z).\quad}
 \label{eq:recursive-initial-test}
\end{equation}
Thus the expansion starts from $F_0(U_h(x))=a_\lambda(U_h(x))$.

For its Gaussian reference, use the unique profile
$\bar a:=\bar a[U_0]$ of
Definition~\ref{def:self-consistent-reference-profile}, whose existence
is supplied by Proposition~\ref{prop:self-profile-compatible-gaussian},
for the heat mean $S_tU_0$. It satisfies the fixed point equation
\eqref{eq:reference-profile-self-consistency}, with mean $S_tU_0$ and
variance \eqref{eq:gaussian-profile-variance-map}. After conditioning on
$\mathcal F_{t_0}$, this profile is deterministic.
Write $\bar{\mathfrak a}_t=m_0^2+\bar a_t$ for the total reference
variance. The actual total variance at $(t,x)$ is
$m_0^2+a(U_t(x))$, so the covariance discrepancy is
$a(U_t(x))-\bar a_t(x)$.
\end{definition}

\begin{definition}[Observation vectors, Gaussian smoothing and covariance insertions]
\label{def:gaussian-observation-operators}
Use the reference profile and total variance from
Definition~\ref{def:coefficient-band-data}. Let
$Z=((T_i,x_i))_{i=1}^d$ be an ordered list of observation targets,
with $d\ge1$, $0<T_i\le h$ and $x_i\in\mathbb T$.
For $0\le t<\min_iT_i$, define its scalar coordinates and observation
vector by
\begin{equation}
 \begin{aligned}
 Y_t^{Z,i}&:=(S_{T_i-t}U_t)(x_i),\qquad 1\le i\le d,\\
 \mathbf Y_t^Z&:=(Y_t^{Z,1},\ldots,Y_t^{Z,d})\in\mathbb R^d.
 \end{aligned}
 \label{eq:physical-observation-vector}
\end{equation}
Thus every coordinate is a backward heat observation of the same
physical solution $U$; the list $Z$ specifies its target time and
spatial point. For $Z_0=((h,x))$, there is one coordinate,
$Y_t^{Z_0,1}=(S_{h-t}U_t)(x)$, whose terminal value is $U_h(x)$.
The associated heat-row coordinates and vector are
\[
 q_t^{Z,i}(\xi):=p_{T_i-t}(x_i-\xi),\qquad
 \mathbf q_t^Z(\xi):=(q_t^{Z,1}(\xi),\ldots,q_t^{Z,d}(\xi))\in\mathbb R^d,
 \qquad \xi\in\mathbb T.
\]
For $0\le t<s\le\min_iT_i$, let $\mathcal P^Z_{t,s}$ be Gaussian
convolution with covariance
\begin{equation}
 \Sigma^Z_{t,s}=\int_t^s\int_{\mathbb T}
 \bar{\mathfrak a}_r(\xi)\mathbf q_r^Z(\xi)\mathbf q_r^Z(\xi)^*\,d\xi\,dr.
 \label{eq:observation-increment-covariance}
\end{equation}
Equivalently, $\mathcal P^Z_{t,s}$ is convolution by the centered
Gaussian measure $N(0,\Sigma^Z_{t,s})$, supported on the covariance
range when this matrix is singular.

\begin{samepage}
\noindent Define the covariance-insertion operator by
\begin{equation}
 \mathscr I^Z_{s,\xi}
 :=\tfrac12[a(w)-\bar a_s(\xi)]
                   (\mathbf q_s^Z(\xi)\cdot\nabla_z)^2,
 \qquad (z,w)\in\mathbb R^d\times\mathbb R.
 \label{eq:recursive-covariance-insertion}
\end{equation}
\end{samepage}
The derivatives act on the existing $d$ coordinates $z$; multiplication
by $a(w)-\bar a_s(\xi)$ appends the new state coordinate $w$.
Thus the operator sends a function of $d$ coordinates to one of
$d+1$ coordinates. In the recursion below, it acts on the Gaussian
smoothing of the preceding test $F_{j-1}$.
\end{definition}

\begin{definition}[Finite coefficient-band words]
\label{def:finite-coefficient-band-words}
Use the initial data of Definition~\ref{def:coefficient-band-data}
and the operators of Definition~\ref{def:gaussian-observation-operators}.
The integer $n\ge0$ counts covariance insertions. At depth $n$,
a history consists of ordered elapsed times
$h=s_0>s_1>\cdots>s_n>0$ and points
$\xi_1,\ldots,\xi_n\in\mathbb T$. In the physical window
$[t_0,t_0+h]$, insertion $j$ occurs at time $t_0+s_j$;
$s_0=h$ labels the terminal observation.
These times range continuously over the ordered simplex; their gaps
are arbitrary positive numbers, not a prescribed time mesh.
Starting from $Z_0=((h,x))$ and $F_0(z)=a_\lambda(z)$, define
\begin{equation}
 Z_j:=Z_{j-1}\cup((s_j,\xi_j)),\qquad
 F_j:=\mathscr I^{Z_{j-1}}_{s_j,\xi_j}
                  \mathcal P^{Z_{j-1}}_{s_j,s_{j-1}}F_{j-1},
 \qquad 1\le j\le n.
 \label{eq:recursive-operator-word}
\end{equation}
Here $\cup$ appends the new observation, including coincident entries.
At a fixed history, $F_j:\mathbb R^{j+1}\to\mathbb R$ is the
resulting test on the enlarged observation list. In particular,
\[
 F_1(z,w)=\tfrac12[a(w)-\bar a_{s_1}(\xi_1)]
 p_{h-s_1}(x-\xi_1)^2
 \partial_z^2\bigl(\mathcal P^{Z_0}_{s_1,h}a_\lambda\bigr)(z).
\]
Write $\mathbf s=(s_1,\ldots,s_n)$ and
$\boldsymbol\xi=(\xi_1,\ldots,\xi_n)$. Define
\[
 \mathcal H_n(h):=\bigl\{(\mathbf s,\boldsymbol\xi)\in
 (0,h)^n\times\mathbb T^n:h>s_1>\cdots>s_n>0\bigr\},
 \qquad d\mathbf s\,d\boldsymbol\xi:=\prod_{j=1}^n ds_j\,d\xi_j.
\]
For $n=0$, this is the empty history with mass one and $s_0=h$.
The remainder, initial-past term and martingale word are, respectively,
\begin{align*}
 R_n&:=\int_{\mathcal H_n(h)}F_n(\mathbf Y_{s_n}^{Z_n})
                                  \,d\mathbf s\,d\boldsymbol\xi,\\
 \mathsf P_{n+1}&:=\int_{\mathcal H_n(h)}
       (\mathcal P^{Z_n}_{0,s_n}F_n)(\mathbf Y_0^{Z_n})
                                  \,d\mathbf s\,d\boldsymbol\xi,\\
 N_{n+1}&:=\int_{\mathcal H_n(h)}
       \left[\int_0^{s_n}\nabla\mathcal P^{Z_n}_{q,s_n}
                   F_n(\mathbf Y_q^{Z_n})\cdot d\mathbf Y_q^{Z_n}\right]
                                  \,d\mathbf s\,d\boldsymbol\xi.
\end{align*}
The subscript of $R_n$ counts its $n$ insertions. Applying one more
It\^o decomposition to this remainder produces $\mathsf P_{n+1}$
and $N_{n+1}$, both still containing $n$ insertions, and the next
remainder $R_{n+1}$. Thus $\mathsf P_j$ and $N_j$ contain $j-1$
insertions. In particular, $R_0=F_0(U_h(x))=a_\lambda(U_h(x))$;
$\mathsf P_1$ and $N_1$ are the terms with no insertion.
\end{definition}

Each step smooths the previous test and inserts one covariance
correction. The approximation and positive-gap limits for these
words are justified by Lemma~\ref{lem:exact-bracket-cutoff-limits}.
For the smooth companion, the covariance insertion uses
$a_{n_{\mathrm{sm}}}$ from \eqref{eq:framework-smoothing} and its
chosen frozen reference, while the terminal test may remain $a_\lambda$.
The control contributes the first-derivative term described in the proof.

\begin{lemma}[Exact finite expansion]
\label{lem:finite-expansion}
For the words $R_n,\mathsf P_{n+1},N_{n+1}$ of
Definition~\ref{def:finite-coefficient-band-words}, the terms
$\mathsf P_{n+1}$ are measurable at the initial past, and
$N_{n+1}$ are conditionally centered. For the uncontrolled physical
marginal $U$, the finite recursion and expansion are
\begin{align}
 R_n&=\mathsf P_{n+1}+N_{n+1}+R_{n+1},
 \label{eq:recursive-master-identity}\\
 a_\lambda(U_h)&=\sum_{j=1}^n(\mathsf P_j+N_j)+R_n.
 \label{eq:recursive-finite-expansion}
\end{align}
\end{lemma}
Lemma~\ref{lem:finite-expansion} follows by iterating It\^o's formula
for the observation vector $\mathbf Y^Z$;
Corollary~\ref{cor:coupled-gaussian-band} supplies
a separate Gaussian estimate. Both enter the later comparison:
Lemma~\ref{lem:finite-expansion} organizes the physical coefficient
observations into words, while Corollary~\ref{cor:coupled-gaussian-band} and
Theorem~\ref{thm:coupled-gaussian-band-sum} bound the Gaussian difference
that appears after regrouping the first predictable and martingale
terms in Proposition~\ref{prop:joint-first-pair}. Sections 5 and 6
estimate the remaining martingale, predictable and remainder differences.

\subsection{Gaussian damping of one state-frequency block}
\label{sec:gaussian-proof-details}
\label{sec:annular-damping-proof}
\raggedbottom
Fix one dyadic state-frequency block $a_\lambda$ from
Section~\ref{sec:state-bands}. The goal here is the exponential
damping of this block under scalar Gaussian convolution; the
variance and derivative order are the only additional inputs.

\begin{lemma}[Annular Gaussian damping]
\label{lem:bands-gaussian-damping}
If $\mathcal P _\upsilon$ is convolution with the scalar Gaussian of variance $\upsilon$,
with $\mathcal P_0$ the identity, then for nonnegative integers
$r$ and variances $\upsilon\ge0$,
\begin{equation}
 \|\partial_z^r\mathcal P _\upsilon a_\lambda\|_\infty
 \le C_r\lambda^{r-\beta}e^{-c \upsilon\lambda^2}.
 \label{eq:bands-gaussian-damping}
\end{equation}
In particular, for a time gap $\tau>0$, a reserved variance
comparable to $\sqrt\tau$ gives
the factor $e^{-c\lambda^2\sqrt\tau}$ used below.
\end{lemma}
\begin{proof}[Proof of Lemma~\ref{lem:bands-gaussian-damping}]
The frequency-localized heat estimate
\cite[Section~2.1.2]{BahouriCheminDanchin2011}, with heat time
$\upsilon/2$ in the state variable, gives
$\|\mathcal P_\upsilon a_\lambda\|_\infty
\le Ce^{-c\upsilon\lambda^2}\|a_\lambda\|_\infty$.
Its Fourier support remains in the same annulus. Bernstein's
inequality~\cite[Lemma~2.1]{BahouriCheminDanchin2011} therefore
adds $C_r\lambda^r$ for $r$ derivatives, and
\eqref{eq:bands-bernstein-tail} supplies $\|a_\lambda\|_\infty
\le C\lambda^{-\beta}$. The case $\upsilon=0$ is the same
Bernstein bound without damping.
\end{proof}

The annular support supplies the Gaussian damping. In a multivariate
observation vector, reserve an independent Gaussian component in
the common state-translation direction. Applying the scalar lemma
to that component extracts the annular factor before estimating
the remaining derivative weights and leaves by fixed moments.

\begin{samepage}
\subsection{Self-consistent profiles: construction and stability}
\label{sec:reference-profile-construction-stability}

\subsubsection{Constructing the self-consistent reference profile}
\label{sec:self-consistent-profile}
This subsection constructs the variance profile used in
Proposition~\ref{prop:self-profile-compatible-gaussian}.
The input is a prescribed mean $\phi_t(x)$: for the ordinary
reference it is $S_tw$, while the controlled reference uses the
mean in \eqref{eq:coupled-reference-means}. The residual coefficient
$a$ takes values in $[\underline a,\overline a]$, and $m_0>0$
is the common additive-noise amplitude. We use the scalar Gaussian
convolution $\mathcal P_v$ and variance map $\mathcal V$ of
Definition~\ref{def:gaussian-variance-map}.

\end{samepage}
\begin{definition}[Gaussian averaging map]
\label{def:gaussian-averaging-map}
Fix $0<h\le1$. Let $\mathscr C_h$ be the closed subset of
$L^\infty((0,h]\times\mathbb T)$ consisting of profiles valued in
$[\underline a,\overline a]$. For a prescribed deterministic
measurable mean $\phi:(0,h]\times\mathbb T\to\mathbb R$, use the variance map of
Definition~\ref{def:gaussian-variance-map} to define
\[
 \begin{aligned}
 (\mathsf T_\phi\kappa)_t(x)
 &:=\bigl(\mathcal P_{\mathcal V_t[\kappa](x)}a\bigr)(\phi_t(x))\\
 &=\frac1{\sqrt{2\pi}}\int_{\mathbb R}
 a\!\left(\phi_t(x)+\sqrt{\mathcal V_t[\kappa](x)}\,z\right)
                    e^{-z^2/2}\,dz,
 \qquad\kappa\in\mathscr C_h.
 \end{aligned}
\]
For a starting state $w\in H$, the prescribed mean is $\phi_t=S_tw$.
\end{definition}

The next proposition gives a sufficient short-time condition for
this map to have a unique fixed point. To state it explicitly, write
\[
 \Theta(h):=\int_0^h\|p_s\|_{L^2(\mathbb T)}^2\,ds,\qquad
 b_-:=m_0^2+\underline a,\qquad
 c_\beta:=\tfrac12\mathbb E\bigl[|\zeta_{\mathrm G}^2-1|\,|\zeta_{\mathrm G}|^\beta\bigr],
 \quad \zeta_{\mathrm G}\sim N(0,1).
\]
Here $\Theta(h)$ is the variance mass of a unit-amplitude heat
convolution. On $0<h\le1$, it is comparable to $\sqrt h$.

\begin{proposition}[Self-consistent Gaussian profiles: existence and stability with respect to the mean]
\label{prop:self-profile}
For any $0<h\le1$ satisfying
\begin{equation}
 q(h):=c_\beta[a]_{C^\beta}\,
 b_-^{\beta/2-1}\Theta(h)^{\beta/2}\le\tfrac12,
 \label{eq:self-profile-contraction-condition}
\end{equation}
and every prescribed deterministic measurable mean
$\phi:(0,h]\times\mathbb T\to\mathbb R$, the map $\mathsf T_\phi$
on $\mathscr C_h$ from Definition~\ref{def:gaussian-averaging-map}
is a contraction with factor at most $q(h)$.
It has a unique fixed point $\bar a^\phi\in\mathscr C_h$, satisfying
\begin{equation}
 \bar a_r^\phi(x)
   =\bigl(\mathcal P_{\mathcal V_r[\bar a^\phi](x)}a\bigr)(\phi_r(x)).
 \label{eq:self-profile-fixed-point}
\end{equation}
For the heat mean $\phi_t=S_tw$, this fixed point is the reference
profile $\bar a[w]=\bar a^{S_\cdot w}$ of
Definition~\ref{def:self-consistent-reference-profile}.

There are $h_0\in(0,1]$ and $C<\infty$, depending only on
$\beta,m_0,\underline a,\overline a$ and $[a]_{C^\beta}$, such
that \eqref{eq:self-profile-contraction-condition} holds for every
$h\le h_0$ and the following stability estimates hold for
$0<r\le h\le h_0$.
The fixed point stays close to the coefficient at its prescribed mean:
\begin{equation}
 |\bar a_r^\phi(x)-a(\phi_r(x))|\le Cr^{\beta/4}.
 \label{eq:self-profile-displacement}
\end{equation}
For two input means with
$\mathsf d_2=\operatorname*{ess\,sup}_{0<s\le h}
\|\phi_s-\widetilde\phi_s\|_{L^2(\mathbb T)}<\infty$, the corresponding fixed points satisfy
\begin{equation}
 \|\bar a_r^\phi-\bar a_r^{\widetilde\phi}\|_{L^2(\mathbb T)}
 \le C\min\{r^{(\beta-1)/4}\mathsf d_2,\mathsf d_2^\beta\}.
 \label{eq:self-profile-lp}
\end{equation}
In particular, if $\phi_s-\widetilde\phi_s=S_s\zeta$ or
$e^{-\Lambda s}S_s\zeta$ with $\Lambda\ge0$, one may take
$\mathsf d_2=\|\zeta\|_{L^2(\mathbb T)}$, where $\zeta\in L^2(\mathbb T)$ is
the initial mean difference.
\end{proposition}

\begin{proof}[Proof of Proposition~\ref{prop:self-profile} (existence)]
\emph{Gaussian bounds and contraction.}
The heat-kernel estimates give
\begin{equation}
 \|\mathsf k_t\|_1\le Ct^{-1/2},\qquad
 c\sqrt r\le\mathcal V_r[\kappa](x)\le C\sqrt r,
 \quad 0<r,t\le1,\quad \underline a\le\kappa\le\overline a.
 \label{eq:self-profile-kernel}
\end{equation}
Gaussian differentiation and the H\"older bound give
$\|\partial_z\mathcal P_va\|_\infty\le C[a]_{C^\beta}v^{(\beta-1)/2}$ and
\[
 \partial_v\mathcal P_va(z)=\frac1{2v}\mathbb E\!\left[
 (\zeta_{\mathrm G}^2-1)
 (a(z+\sqrt v\,\zeta_{\mathrm G})-a(z))\right],\qquad
 \|\partial_v\mathcal P_va\|_\infty
 \le c_\beta[a]_{C^\beta}v^{\beta/2-1}.
\]
Gaussian averaging maps the complete closed set $\mathscr C_h$
into itself. The variance map satisfies
$\mathcal V_t[\kappa]\ge b_-\Theta(t)$ and
$\|\mathcal V_t[\kappa]-\mathcal V_t[\widetilde\kappa]\|_\infty
\le\Theta(t)\|\kappa-\widetilde\kappa\|_\infty$.
The variance derivative bound therefore yields
\[
 \|\mathsf T_\phi\kappa-\mathsf T_\phi\widetilde\kappa\|_\infty
 \le q(h)\|\kappa-\widetilde\kappa\|_\infty.
\]
Banach's theorem gives the unique fixed point under
\eqref{eq:self-profile-contraction-condition}, as the limit of
\begin{equation}
 \kappa_t^{(0)}(x):=a(\phi_t(x)),\qquad
 \kappa^{(m+1)}:=\mathsf T_\phi\kappa^{(m)},\qquad m\ge0.
 \label{eq:self-profile-explicit-iteration}
\end{equation}
Since $\Theta(h)\asymp\sqrt h$, the condition holds for all
sufficiently small $h_0$. The fixed point equation and
\eqref{eq:self-profile-kernel} also give
$|\bar a_r^\phi-a(\phi_r)|\le C\mathcal V_r[\bar a^\phi]^{\beta/2}
\le Cr^{\beta/4}$, proving \eqref{eq:self-profile-displacement}.

\end{proof}

\subsubsection{Stability of the profile under a change of mean}
\label{sec:profile-stability-proof}
Work on $(0,h]$ with $h\le h_0$. The goal is to compare the fixed points
$\bar a^\phi$ and $\bar a^{\widetilde\phi}$ constructed in
Section~\ref{sec:self-consistent-profile}, keeping the coefficient
$a$ fixed. Here $\phi,\widetilde\phi$ are prescribed mean paths,
$\mathcal V_r[\bar a^\phi]$ is the pointwise reference variance,
and $\|\cdot\|_p$ denotes $L^p(\mathbb T)$.
We prove both the linear and the H\"older bounds in
\eqref{eq:self-profile-lp}, as well as the variance bound used later.

\begin{proof}[Proof of Proposition~\ref{prop:self-profile} (stability)]
\emph{Linear stability.}
For $1\le p\le\infty$, put
$\mathsf d_p:=\operatorname*{ess\,sup}_{0<s\le h}
\|\phi_s-\widetilde\phi_s\|_p$ and
$f(r):=\|\bar a_r^\phi-\bar a_r^{\widetilde\phi}\|_p$.
Split the Gaussian average into changes of mean and variance.
The Gaussian mean and variance derivative bounds from
Section~\ref{sec:self-consistent-profile}, followed by Young's inequality, yield
\begin{equation}
 f(r)\le Cr^{(\beta-1)/4}\mathsf d_p+
 Cr^{(\beta-2)/4}\int_0^r(r-s)^{-1/2}f(s)\,ds.
 \label{eq:self-profile-volterra}
\end{equation}
In the finite weighted supremum
$\operatorname*{ess\,sup}_{0<r\le h}r^{(1-\beta)/4}f(r)$,
the substitution $s=rt$ gives
\[
 r^{(\beta-2)/4}\int_0^r(r-s)^{-1/2}s^{(\beta-1)/4}\,ds
 =r^{(2\beta-1)/4}B(1/2,(3+\beta)/4),
\]
where $B(u,v):=\int_0^1t^{u-1}(1-t)^{v-1}dt$ for $u,v>0$.
Multiplication by $r^{(1-\beta)/4}$ leaves the small factor
$r^{\beta/4}\le h_0^{\beta/4}$.
Decreasing $h_0$ absorbs this term and gives
$f(r)\le Cr^{(\beta-1)/4}\mathsf d_p$.
The heat-square convolution and Young's inequality give
\[
 \|\mathcal V_r[\bar a^\phi]-\mathcal V_r[\bar a^{\widetilde\phi}]\|_p
 \le C\mathsf d_p\int_0^r(r-s)^{-1/2}s^{(\beta-1)/4}\,ds.
\]
The same substitution shows that this integral is a constant times
$r^{(1+\beta)/4}$, proving
\begin{equation}
 \|\mathcal V_r[\bar a^\phi]-\mathcal V_r[\bar a^{\widetilde\phi}]\|_p
 \le Cr^{(1+\beta)/4}\mathsf d_p.
 \label{eq:self-profile-variance-lp}
\end{equation}

\emph{H\"older stability.}
For fixed variance, the pointwise mean difference is at most
$[a]_{C^\beta}|\phi_r-\widetilde\phi_r|^\beta$.
Since the torus has volume one,
$\||\phi_r-\widetilde\phi_r|^\beta\|_{L^2(\mathbb T)}
\le\|\phi_r-\widetilde\phi_r\|_{L^2(\mathbb T)}^\beta$.
Replace the source in \eqref{eq:self-profile-volterra} by
$C\mathsf d_2^\beta$. Its integral operator has unweighted
supremum norm $Ch_0^{\beta/4}<1$, proving
\eqref{eq:self-profile-lp}. Heat contraction also covers damped means.

For heat means $\phi_s=S_sw$ and $\widetilde\phi_s=S_s\widetilde w$,
the same argument in $L^p$ has mean term
$C\|S_r(w-\widetilde w)\|_{p\beta}^\beta$.
The unit volume of $\mathbb T$ and heat contraction therefore give
\begin{align}
 \|\bar a_r^\phi-\bar a_r^{\widetilde\phi}\|_p
 &\le C\|w-\widetilde w\|_p^\beta,
 \label{eq:self-profile-holder-lp}\\
 \|\bar a_r^\phi-\bar a_r^{\widetilde\phi}\|_p
 &\le C\|w-\widetilde w\|_q^\beta,
 \qquad 1\le p,q\le\infty,\quad p\beta\le q.
 \label{eq:self-profile-holder-mixed-lp}
\end{align}
For $\beta\le1/2$ this includes $p=2,4$ with $q=2$.
All constants are uniform in the spatial exponents and in smooth
approximations with the common coefficient bounds.
\end{proof}
\subsection{Centering and coupling the Gaussian references}
\label{sec:compatible-reference-proof}
Fix $0<h\le h_0$, $w,\widetilde w\in H=L^2(\mathbb T)$ and let
$\zeta=w-\widetilde w$. The reference fields $X,\widetilde X$
use the same noises, heat means $S_tw,S_t\widetilde w$, and
profiles $\bar a=\bar a[w]$, $\widetilde{\bar a}=\bar a[\widetilde w]$.
Here $L(\kappa)$ is the exact Haar covariance factor from
Section~3. We prove exact centering and control the difference of
the two Gaussian fields, using the profile stability just established.

\begin{proof}[Proof of Proposition~\ref{prop:self-profile-compatible-gaussian}]
Proposition~\ref{prop:self-profile} constructs the unique reference
profiles. For either reference, write $\phi$ for its prescribed mean
and $\bar a^\phi$ for its profile. At each $(t,x)$ its mean is $\phi_t(x)$
and variance $\mathcal V_t[\bar a^\phi](x)$. Consequently the fixed
point equation says exactly
\[
 \mathbb E[a(X_t^\phi(x))]
 =(\mathcal P_{\mathcal V_t[\bar a^\phi](x)}a)(\phi_t(x))
 =\bar a_t^\phi(x).
\]

For $f=\bar a-\widetilde{\bar a}$, the signed heat estimate
\eqref{eq:haar-heat-signed} at $\rho=0$ gives
$\|S_t(L(\bar a)-L(\widetilde{\bar a}))\|_{\mathrm{HS}}^2
\le Ct^{-1/2}\|f\|_{L^2(\mathbb T)}^2$.
It\^o isometry and \eqref{eq:self-profile-lp} bound the centered difference by
\[
 C\|\zeta\|_{L^2(\mathbb T)}^2\int_0^h(h-z)^{-1/2}z^{(\beta-1)/2}dz
 =CB(1/2,(1+\beta)/2)h^{\beta/2}\|\zeta\|_{L^2(\mathbb T)}^2.
\]
The centered stochastic difference has zero expectation, so its
$L^2$ cross term with the deterministic mean difference vanishes.
Adding the squared mean and centered bounds proves the full estimate.
The damped mean case uses $e^{-\Lambda t}\le1$ and heat contraction.
\end{proof}
\subsection{Spatial gain and comparison of Gaussian bands}
\label{sec:gaussian-band-estimates-proofs}

\subsubsection{The spatial gain for an annular Gaussian observation}
\label{sec:annular-spatial-proof}
We first prove the annular spatial estimate that supplies the
$\lambda^{-2\rho}$ gain in the one-band comparison. The next
statement fixes its additive field, annular tests and spatial envelope.

\begin{lemma}[Annular covariance and weighted norm]
\label{lem:annular-spatial-gain}
Let $W$ be a cylindrical Brownian motion on $L^2(\mathbb T)$,
let $0<h\le h_0$ and $m_0>0$, and define the centered additive
Gaussian heat field by
\[
 \mathcal G_h^{\mathrm{add}}:=m_0\int_0^hS_{h-s}\,dW_s.
\]
For $x,y\in\mathbb T$, let $\operatorname{dist}_{\mathbb T}(x,y)$ denote
their torus distance. For a measurable family of bounded functions $f_x:\mathbb R\to\mathbb R$,
assume, for some $\lambda>0$ and fixed $0<c_{\mathrm{ann}}<C_{\mathrm{ann}}<\infty$,
\begin{equation}
 \operatorname{supp}\bigl(\mathcal F_{z\to\zeta}f_x\bigr)
 \subset\bigl\{\zeta\in\mathbb R:
            c_{\mathrm{ann}}\lambda\le|\zeta|
                            \le C_{\mathrm{ann}}\lambda\bigr\},
 \qquad x\in\mathbb T.
 \label{eq:annular-state-support}
\end{equation}
Here $\mathcal F_{z\to\zeta}$ is the real-line Fourier transform in
the scalar state variable $z$, understood in the sense of tempered
distributions. The annular constants are independent of $x$ and $\lambda$.
With constants allowed to
depend on $c_{\mathrm{ann}},C_{\mathrm{ann}}$, one has
\begin{equation}
 |\mathbb E[f_x(\mathcal G_h^{\mathrm{add}}(x))f_y(\mathcal G_h^{\mathrm{add}}(y))]|
 \le C\|f_x\|_\infty\|f_y\|_\infty
                       e^{-c\lambda^2\min(\operatorname{dist}_{\mathbb T}(x,y),\sqrt h)}.
 \label{eq:sharp-two-point-annular}
\end{equation}
Let $\mathsf A_{\mathrm{band}}\in L^2(\mathbb T)$ be a nonnegative
measurable spatial envelope such that
$\sup_{z\in\mathbb R}|f_x(z)|\le\mathsf A_{\mathrm{band}}(x)$
for almost every $x$. For $0<\rho<1/2$ and $\Lambda\ge1$,
with $\|\mathsf A_{\mathrm{band}}\|_{L^2(\mathbb T)}$ taken over the spatial variable,
\begin{equation}
 \mathbb E\|f_\cdot(\mathcal G_h^{\mathrm{add}}(\cdot))\|_{H^{-\rho}_\Lambda}^2
 \le C[\lambda^{-4\rho}+e^{-c\lambda^2\sqrt h}]\|\mathsf A_{\mathrm{band}} \|_{L^2(\mathbb T)}^2.
 \label{eq:sharp-annular-weighted-norm}
\end{equation}
\end{lemma}

\begin{proof}[Proof of Lemma~\ref{lem:annular-spatial-gain}]
Throughout this proof, write $r:=\operatorname{dist}_{\mathbb T}(x,y)$
for the torus distance between the two spatial points.

\emph{Two-point covariance.}
The smaller eigenvalue of the covariance matrix of
$(\mathcal G_h^{\mathrm{add}}(x),\mathcal G_h^{\mathrm{add}}(y))$ is
\[
 \frac{m_0^2}{2}\int_0^h
 \|p_s(x-\cdot)-p_s(y-\cdot)\|_{L^2(\mathbb T)}^2\,ds.
\]
To bound it from below, the heat semigroup identity gives
\[
 \|p_s(x-\cdot)-p_s(y-\cdot)\|_{L^2(\mathbb T)}^2
 =2\bigl[p_{2s}(0)-p_{2s}(x-y)\bigr]
 \ge c s^{-1/2}
 \quad\text{when }0<s\le r^2/8,\ s\le h_0.
\]
Indeed, on the line the two kernel values are
$(8\pi s)^{-1/2}$ and $(8\pi s)^{-1/2}e^{-r^2/(8s)}$,
so their difference is at least
$(8\pi s)^{-1/2}(1-e^{-1})$ in this range.
The torus heat kernel is the sum of the Gaussian kernel over integer
translates. Since $r\le1/2$, the additional terms are exponentially
small as $s\downarrow0$, uniformly in $x,y$; choosing the fixed
$h_0$ sufficiently small preserves the lower bound.

If $0<r\le\sqrt h$, integrate this bound over
$s\in[r^2/16,r^2/8]\subset[0,h]$ to obtain
$\int_0^h\|p_s(x-\cdot)-p_s(y-\cdot)\|_{L^2(\mathbb T)}^2\,ds
\ge cr$.
If $r>\sqrt h$, use instead
$s\in[h/16,h/8]$, where still $s\le r^2/8$, to obtain the lower
bound $c\sqrt h$.
Thus the smaller covariance eigenvalue is at least
$c\min(r,\sqrt h)$.
Split off an independent
Gaussian pair of covariance $c\min(r,\sqrt h)I_2$, where $I_2$ is
the $2\times2$ identity, and condition on the remaining pair.
The annular multiplier argument from the proof of
Lemma~\ref{lem:bands-gaussian-damping} gives an $L^\infty$ Gaussian
smoothing bound $Ce^{-c\lambda^2v}$ at variance $v$.
Applying it in each coordinate proves \eqref{eq:sharp-two-point-annular}.
Scalar state translations preserve the support and supremum bounds,
so the same estimate holds after conditioning on any fixed translations.

\emph{Spatial norm.}
The positive Bessel kernel is bounded near zero by $Cr^{2\rho-1}$.
Its product with the covariance bound has mass at most
\[
 C\int_0^{\sqrt h}r^{2\rho-1}e^{-c\lambda^2r}dr
 +Ce^{-c\lambda^2\sqrt h}
 \le C\bigl[\lambda^{-4\rho}+e^{-c\lambda^2\sqrt h}\bigr].
\]
Schur's inequality with the envelope $\mathsf A_{\mathrm{band}}$
proves \eqref{eq:sharp-annular-weighted-norm}.
\end{proof}

\subsubsection{Comparing and summing the Gaussian coefficient bands}
\label{sec:gaussian-band-comparison-proof}
Fix a horizon $T$ and a window $[t_0,t_0+h]\subset[0,T]$.
Use the common-noise references $X^U,X^V$ of
Definition~\ref{def:coupled-gaussian-references}.
The input is measured by $\|U-V\|_{\mathscr L_0^2(H)}$, with
$H=L^2(\mathbb T)$ and $D=U-V$ the physical solution difference. We apply the spatial
gain of Section~\ref{sec:annular-spatial-proof} to one block
$a_\lambda$, then sum over dyadic $\lambda$.

\begin{corollary}[Comparison of one Gaussian state-frequency band]
\label{cor:coupled-gaussian-band}
Use the coupled pair and reference fields of
Definition~\ref{def:coupled-gaussian-references}, with the integrated norms of
Definition~\ref{def:window-integrated-norms}. Take
$0<h\le\min(h_0,T)$. For
$\lambda\in\mathscr S_{\mathrm z}$ and $0<\rho<1/2$,
\begin{equation}
 \|a_\lambda(X_h^U)-a_\lambda(X_h^V)\|_{\mathscr L_h^2(H_\Lambda^{-\rho})}
 \le C\bigl[\lambda^{-\beta-2\rho}
       +\lambda^{-\beta}e^{-c\lambda^2\sqrt h}\bigr]
                   \min(1,\lambda \|U-V\|_{\mathscr L_0^2(H)}).
 \label{eq:sharp-main-Gaussian-band}
\end{equation}
The constants may depend on $T$ and the common coefficient bounds.
The same estimate holds for the pair $X_h^U$ and its translate
$X_h^U-e^{-\Lambda h}S_hD_{t_0}$.
\end{corollary}

\begin{proof}[Proof of Corollary~\ref{cor:coupled-gaussian-band}]
\emph{Separate the common additive field from the two translations.}
Keep the window start $t_0$ and its length $h$ explicit. At physical
time $t_0+h$, the references in
\eqref{eq:coupled-gaussian-reference-fields} split as
\[
 \begin{aligned}
 \mathcal G_h^{\mathrm{add},t_0}
   &:=m_0\int_0^h S_{h-s}\,dW_s^{0,t_0},\\
 z_h^{i,t_0}
   &:=\phi_h^{i,t_0}+\int_0^h S_{h-s}L(\bar a_s^{i,t_0})\,dB_s^{t_0},\\
 X_h^{i,t_0}&=\mathcal G_h^{\mathrm{add},t_0}+z_h^{i,t_0},
                         \qquad i\in\{U,V\}.
 \end{aligned}
\]
Here $\phi_h^{U,t_0}=S_hU_{t_0}$ and
$\phi_h^{V,t_0}=S_hU_{t_0}-e^{-\Lambda h}S_hD_{t_0}$.
Conditional on $\mathcal F_{t_0}$, these mean-plus-residual
translations depend only on the future $B$ increments and are
independent of $\mathcal G_h^{\mathrm{add},t_0}$.

Condition further on $B$. The positive amplitude $m_0$ leaves
the independent additive variance needed in
Lemma~\ref{lem:annular-spatial-gain}; this is where nondegeneracy
supplies the spatial gain $\lambda^{-2\rho}$.
Apply that lemma with
\[
 f_x^{t_0,h,\lambda}(y)
 :=a_\lambda(y+z_h^{U,t_0}(x))-a_\lambda(y+z_h^{V,t_0}(x)).
\]
At $y=\mathcal G_h^{\mathrm{add},t_0}(x)$ this is the desired
difference. Translations preserve its annular support, and
\eqref{eq:bands-bernstein-tail} gives
\[
 \|f_x^{t_0,h,\lambda}\|_\infty
 \le C\lambda^{-\beta}
          \min\{1,\lambda|z_h^{U,t_0}(x)-z_h^{V,t_0}(x)|\}.
\]
Use the supremum bound and the mean-value theorem for the two caps.

\emph{Control the difference of the translations.}
Their explicit difference is
\[
 z_h^{U,t_0}-z_h^{V,t_0}
 =e^{-\Lambda h}S_hD_{t_0}+
 \int_0^h S_{h-s}
       [L(\bar a_s^{U,t_0})-L(\bar a_s^{V,t_0})]\,dB_s^{t_0}.
\]
Apply the damped-mean version of
\eqref{eq:self-profile-noise-stability}, conditional on $\mathcal F_{t_0}$.
The stochastic integral has conditional mean zero, and heat contraction
bounds the deterministic term. Consequently
\[
 \mathbb E\bigl[\|z_h^{U,t_0}-z_h^{V,t_0}\|_{L^2(\mathbb T)}^2
                                      \mid\mathcal F_{t_0}\bigr]
 \le C\|D_{t_0}\|_{L^2(\mathbb T)}^2.
\]

\emph{Apply the spatial gain and integrate the window start.}
Use the displayed supremum bound as the spatial envelope in
\eqref{eq:sharp-annular-weighted-norm}, with the translations fixed.
Then average over the residual noise and the initial past. The squared
envelope is controlled by
\[
 \begin{aligned}
 &\int_0^{T-h}\mathbb E\int_{\mathbb T}
       \min\{1,\lambda|z_h^{U,t_0}(x)-z_h^{V,t_0}(x)|\}^2\,dx\,dt_0\\
 &\qquad\le\min\{T,C\lambda^2\|U-V\|_{\mathscr L_0^2(H)}^2\}
       \le C_T\min\{1,\lambda\|U-V\|_{\mathscr L_0^2(H)}\}^2.
 \end{aligned}
\]
Taking square roots gives \eqref{eq:sharp-main-Gaussian-band}.
For the translated comparison, replace $z_h^{V,t_0}$ by
$z_h^{U,t_0}-e^{-\Lambda h}S_hD_{t_0}$. The common additive field is the
same, and the translation difference is exactly
$e^{-\Lambda h}S_hD_{t_0}$, so the same estimates apply.
\end{proof}

\begin{corollary}[Gaussian-band sum for general window and cutoff exponents]
\label{cor:gaussian-band-general-exponents}
Under the hypotheses of Theorem~\ref{thm:coupled-gaussian-band-sum},
the bound \eqref{eq:sharp-main-Gaussian-sum} also holds
for any fixed real time scaling exponent $0<\gamma_{\mathrm{time}}<4$
and cutoff exponent $r>0$
when one instead takes
\begin{equation}
 h_\lambda=\lambda^{-\gamma_{\mathrm{time}}},\qquad
 \lambda_*(\Lambda)=2^{\lceil r\log_2\Lambda\rceil},\qquad
 \Lambda^r\le\lambda_*(\Lambda)<2\Lambda^r.
 \label{eq:gaussian-band-general-cutoff}
\end{equation}
Here $\alpha=\beta+2\rho\in(0,1)$, and $p_{\mathrm{tail}}>0$ is the
desired tail power. The constants and the threshold for $\Lambda$
may depend on these fixed choices.
The theorem uses $r=1/4$ and $\gamma_{\mathrm{time}}=3\beta$; the general form lets the
later parameter calculations display the dependence on both exponents.

\end{corollary}

\begin{proof}[Proof of Theorem~\ref{thm:coupled-gaussian-band-sum} and Corollary~\ref{cor:gaussian-band-general-exponents}]
We prove the estimate for \eqref{eq:gaussian-band-general-cutoff}.
For each fixed $\lambda$, change variables $t_0=t-h_\lambda$ and
apply \eqref{eq:sharp-main-Gaussian-band} with $h=h_\lambda$.
We sum its exponential and algebraic terms separately.
Since $\gamma_{\mathrm{time}}<4$, the exponent $2-\gamma_{\mathrm{time}}/2$ is positive. With
$h_\lambda=\lambda^{-\gamma_{\mathrm{time}}}$, the exponential in
\eqref{eq:sharp-main-Gaussian-band} is $e^{-c\lambda^{2-\gamma_{\mathrm{time}}/2}}$. Using $\min(1,\lambda\|U-V\|_{\mathscr L_0^2(H)})\le1$,
its contribution is bounded by
\[
 C\sum_{\substack{\lambda\in\mathscr S_{\mathrm z}\\\lambda\ge\lambda_*}}
       \lambda^{-\beta}e^{-c\lambda^{2-\gamma_{\mathrm{time}}/2}}
 \le C e^{-c'\lambda_*^{2-\gamma_{\mathrm{time}}/2}}
 \le C e^{-c'\Lambda^{r(2-\gamma_{\mathrm{time}}/2)}}.
\]
For fixed $r,p_{\mathrm{tail}}>0$ this is at most
$\Lambda^{-p_{\mathrm{tail}}}$ for large $\Lambda$.
Also $h_{\lambda_*}\to0$, so the window restrictions hold.
For the algebraic part, the two geometric sums below and above
$\lambda=\|U-V\|_{\mathscr L_0^2(H)}^{-1}$ give
\[
 \sum_{\substack{\lambda\in\mathscr S_{\mathrm z}\\\lambda\ge\lambda_*}}
 \lambda^{-\alpha}\min(1,\lambda\|U-V\|_{\mathscr L_0^2(H)})
 \le C_\alpha\|U-V\|_{\mathscr L_0^2(H)}^\alpha,
 \qquad 0<\alpha<1.
\]
The lower sum uses $1-\alpha>0$ and the upper sum uses $\alpha>0$;
if the norm is zero every summand vanishes.
The translated-pair estimate follows by the same summation.
\end{proof}

\begin{remark}[Gaussian-band sums at and above the Lipschitz threshold]
\label{rem:gaussian-band-other-powers}
The same argument also applies when $\alpha=\beta+2\rho\ge1$.
For $\alpha>1$, the same two geometric sums give
$C\lambda_*^{1-\alpha}\|U-V\|_{\mathscr L_0^2(H)}$.
For $\alpha=1$, the first sum counts dyadic bands and gives
$C\|U-V\|_{\mathscr L_0^2(H)}[1+\log_+(1/(\lambda_*\|U-V\|_{\mathscr L_0^2(H)}))]$, where
$\log_+z:=\max(0,\log z)$ and the expression is zero when
$\|U-V\|_{\mathscr L_0^2(H)}=0$. These replace $C\|U-V\|_{\mathscr L_0^2(H)}^\alpha$ in
\eqref{eq:sharp-main-Gaussian-sum}; the Gaussian tail is unchanged.
The final parameter choice uses $0<\alpha<1$.

\end{remark}

\subsection{Deriving the exact finite expansion by It\^o's formula}
\label{sec:finite-expansion-proof}
Fix a coefficient block $a_\lambda$ and a window $[t_0,t_0+h]$.
As in Definition~\ref{def:coefficient-band-data}, relabel the
physical solution $U_{t_0+t}$ as $U_t$ and condition on the initial past.
For $Z=((T_i,x_i))_{i=1}^d$ with $0<T_i\le h$ and $x_i\in\mathbb T$,
the backward observation is
$\mathbf Y_t^Z=((S_{T_i-t}U_t)(x_i))_{i=1}^d$.
The Gaussian convolution $\mathcal P^Z_{t,s}$ uses covariance
\eqref{eq:observation-increment-covariance};
$\mathscr I^Z_{s,\xi}$ is the covariance-insertion operator in
\eqref{eq:recursive-covariance-insertion}. We prove one It\^o step
and iterate it to obtain the words of
Definition~\ref{def:finite-coefficient-band-words}.

\begin{proof}[Proof of Lemma~\ref{lem:finite-expansion}]

\emph{One insertion is an application of It\^o's formula.}
Fix a deterministic smooth bounded test $F:\mathbb R^d\to\mathbb R$
with bounded derivatives and a terminal time $0<s<\min_iT_i$.
We apply the time-dependent, finite-dimensional It\^o formula
\cite[Chapter~IV]{RevuzYor1999} to
$t\mapsto(\mathcal P^Z_{t,s}F)(\mathbf Y_t^Z)$ on $[0,s]$.
Here the observation vector of
Definition~\ref{def:gaussian-observation-operators} is a continuous
martingale for the uncontrolled solution $U$, with bracket
\[
 \frac{d}{dt}\langle\mathbf Y^Z\rangle_t
 =\int_{\mathbb T}[m_0^2+a(U_t(\xi))]
                \mathbf q_t^Z(\xi)\mathbf q_t^Z(\xi)^*\,d\xi.
\]
The Gaussian convolution solves the backward equation
\begin{equation}
 \partial_t\mathcal P^Z_{t,s}F
 =-\frac12\int_{\mathbb T}\bar{\mathfrak a}_t(\xi)
       (\mathbf q_t^Z(\xi)\cdot\nabla_z)^2
                         \mathcal P^Z_{t,s}F\,d\xi,
 \qquad \mathcal P^Z_{s,s}F=F.
 \label{eq:reference-gaussian-backward-equation}
\end{equation}
The difference between the actual and reference covariance is
\[
 [m_0^2+a(U_t(\xi))]-\bar{\mathfrak a}_t(\xi)
                 =a(U_t(\xi))-\bar a_t(\xi),
\]
so integration of It\^o's formula gives
\begin{align}
 F(\mathbf Y_s^Z)={}&(\mathcal P ^Z_{0,s}F)(\mathbf Y_0^Z)
 +\int_0^s\nabla \mathcal P ^Z_{t,s}F(\mathbf Y_t^Z)\cdot d\mathbf Y_t^Z\nonumber\\
 &+\int_0^s\int_{\mathbb T}(\mathscr I^Z_{t,\xi}\mathcal P ^Z_{t,s}F)
                         (\mathbf Y_t^Z,U_t(\xi))\,d\xi\,dt.
 \label{eq:recursive-one-step}
\end{align}
A drift $\mathbf b_t^Z\,dt$ adds
\begin{equation}
 \int_0^s\nabla\mathcal P^Z_{t,s}F(\mathbf Y_t^Z)
                         \cdot\mathbf b_t^Z\,dt,
 \label{eq:controlled-observation-drift-term}
\end{equation}
and the stochastic integral uses the martingale part of the observation.
Lemma~\ref{lem:exact-bracket-cutoff-limits} extends the identity
from smooth tests and positive gaps to the rough coefficient words.

\emph{Iteration.} Apply \eqref{eq:recursive-one-step} to each history
in the definitions of $R_n,\mathsf P_{n+1},N_{n+1}$. Spatial and
stochastic Fubini give \eqref{eq:recursive-master-identity}; its finite
telescoping sum gives \eqref{eq:recursive-finite-expansion}.
\end{proof}

\subsection{Gaussian and time-kernel estimates used below}
\label{sec:gaussian-input-estimates}
The following statements collect the Gaussian derivative bounds and time kernels
used in Sections~5--6. Their definitions and usable estimates are
kept together here; the longer Gaussian verifications are in
Appendix~\ref{app:gaussian-geometry}, and the time-integration
proofs are in Appendix~\ref{app:time-kernels}.

\subsubsection{Covariance and mean changes: the two spatial input bounds}
\label{sec:multipoint-score-package}
Two estimates from this subsection will be used below:
Lemma~\ref{lem:gaussian-one-input-scores}(iii) keeps a covariance-profile
change in $L^2(\mathbb T)$, while
Lemma~\ref{lem:physical-heat-score}(ii) does the same for an initial-state
change in the Gaussian mean. These are the covariance and mean inputs
in the comparison of Gaussian average terms of Proposition~\ref{prop:all-order-predictable-pair}.
The other parts identify the Gaussian derivative weights, their
spatial localization and the common Gaussian shift used by the
finite expansion. Section~\ref{sec:conditional-gaussian-inputs} then
treats heat-smoothed observations and bounded profile errors;
the supporting heat-span geometry and general leverage identity
are collected in Appendix~\ref{app:gaussian-geometry}.

\begin{definition}[The reference field for the Gaussian input estimates]
\label{def:gaussian-input-reference}
Fix $[t_0,t_0+h]$, $0<h\le h_0$, and write
$\mathbb E_{t_0}=\mathbb E[\,\cdot\mid\mathcal F_{t_0}]$.
Choose $i\in\{U,V\}$ in Definition~\ref{def:coupled-gaussian-references}
and abbreviate $\phi=\phi^{i,t_0}$, $\bar a=\bar a^{i,t_0}$.
Thus $\phi_t=S_tU_{t_0}$ for $i=U$ and
$\phi_t=S_tU_{t_0}-e^{-\Lambda t}S_tD_{t_0}$ for $i=V$.
Using the same shifted noises as in
\eqref{eq:coupled-gaussian-reference-fields}, the chosen field is
\begin{equation}
 X_t:=X_t^{i,t_0}=\phi_t+m_0\int_0^t S_{t-r}\,dW_r^{0,t_0}
                 +\int_0^t S_{t-r}L(\bar a_r)\,dB_r^{t_0},
 \qquad 0\le t\le h.
 \label{eq:score-reference-field}
\end{equation}
The profile is deterministic after conditioning, and
$L(\bar a_r)L(\bar a_r)^*=\mathsf M_{\bar a_r}$ by
\eqref{eq:covariance-factor-identity}. The following estimates, including
Section~\ref{sec:conditional-gaussian-inputs}, hold for both choices
with the same constants; each uses its own mean and covariance.
All times are elapsed from $t_0$.
\end{definition}

\begin{definition}[Observation vector, heat responses and covariance]
\label{def:gaussian-input-observations}
Let $Z=((T_i,x_i))_{i=1}^q$, where $q\ge1$,
$0<T_i\le h$ and $x_i\in\mathbb T$. The scalar field $X_t(x)$
produces the finite observation vector and its conditional mean
\[
 \mathbf X^Z=(X_{T_i}(x_i))_{i=1}^q,\qquad
 \mu^Z=\mathbb E_{t_0}\mathbf X^Z=(\phi_{T_i}(x_i))_{i=1}^q.
\]
Observation $i$ is made at physical time $t_0+T_i$.
All probability norms below are conditional on the initial past.

\emph{Covariance determined by the history.}
Use the heat rows of Definition~\ref{def:gaussian-observation-operators},
with a coordinate set to zero after its observation time:
\begin{equation}
 Q_u^Z(\xi)_i=
 \begin{cases}
 p_{T_i-u}(x_i-\xi),&0\le u<T_i,\\
 0,&u\ge T_i.
 \end{cases}
 \label{eq:history-heat-response}
\end{equation}
This vector records how noise at time $u$ and spatial point $\xi$
reaches the observations in $Z$. It\^o isometry applied to
\eqref{eq:score-reference-field} gives their full covariance:
\begin{equation}
 \begin{aligned}
 \Sigma^Z&:=\mathbb E_{t_0}[(\mathbf X^Z-\mu^Z)(\mathbf X^Z-\mu^Z)^*]\\
 &=\int_0^{h}\int_{\mathbb T}
       (m_0^2+\bar a_u(\xi))Q_u^Z(\xi)Q_u^Z(\xi)^*\,d\xi\,du.
 \end{aligned}
 \label{eq:history-full-observation-covariance}
\end{equation}
In entry $(i,k)$, the time integral stops at $\min(T_i,T_k)$.
The earlier covariance $\Sigma^Z_{t,s}$ in
\eqref{eq:observation-increment-covariance} records only the noise
on the common interval $(t,s)$; \eqref{eq:history-full-observation-covariance}
includes all noise received by each observation up to its own target time.

The unit-variance covariance for the same history is
\begin{equation}
 \mathsf G^Z:=\int_0^{h}\int_{\mathbb T}
                        Q_u^Z(\xi)Q_u^Z(\xi)^*\,d\xi\,du,
 \qquad c\mathsf G^Z\le\Sigma^Z\le C\mathsf G^Z.
 \label{eq:history-unit-covariance}
\end{equation}
The symbols $(\Sigma^Z)^+$ and
$(\mathsf G^Z)^+$ denote the inverses on their common range,
extended by zero on the null space. The superscript $Z$ is displayed
throughout the estimates, including when the centers $x_i$ are integrated.

\end{definition}
The covariance comparison in \eqref{eq:history-unit-covariance}
follows immediately from \eqref{eq:history-full-observation-covariance}:
the scalar variance $m_0^2+\bar a_u(\xi)$ lies between the common
positive lower and upper bounds.

\paragraph{\textbf{Covariance changes: preserve one spatial input.}}
Gaussian integration by parts differentiates an average by
multiplying its integrand by an explicit Gaussian weight.
For a standard normal variable $G$, the first and second weights
are $G$ and $G^2-1$. The formulas below give their counterparts
for the observation vector and its covariance. The first part of the
lemma identifies the weights for changing the mean or changing
the variance at one time--space point.

\begin{lemma}[Gaussian derivative bounds with one spatial $L^2$ input]
\label{lem:gaussian-one-input-scores}
The estimates below hold simultaneously for the two Gaussian references
$X^{U,t_0}$ and $X^{V,t_0}$, with the same constants; in each case use
its own observation mean $\mu^Z$ and covariance $\Sigma^Z$ from
\eqref{eq:history-full-observation-covariance}.
For the history $Z=((T_i,x_i))_{i=1}^q$ above, set
$\displaystyle s_Z:=\min_{1\le i\le q}T_i>0$.
This is the elapsed time from the window start $t_0$ to the earliest
observation, which occurs at physical time $t_0+s_Z$.
All observations receive noise throughout the common interval $(0,s_Z)$.
\begin{enumerate}[label=(\roman*),leftmargin=*]
\item \emph{Size of the Gaussian derivative weights.}
For a mean direction $\mathbf d\in\operatorname{Ran}\Sigma^Z$,
and a covariance direction $Q_u^Z(\xi)Q_u^Z(\xi)^*$, the respective Gaussian weights are
\begin{equation}
\begin{aligned}
 \mathsf S_{\mathrm{mean}}^Z(\mathbf d)
   &=\mathbf d^*(\Sigma^Z)^+(\mathbf X^Z-\mu^Z),\\
 \mathsf S_{\mathrm{cov}}^Z(u,\xi)
   &=\tfrac12\bigl[(Q_u^Z(\xi)^*(\Sigma^Z)^+(\mathbf X^Z-\mu^Z))^2
                              -Q_u^Z(\xi)^*(\Sigma^Z)^+Q_u^Z(\xi)\bigr].
\end{aligned}
\label{eq:gaussian-score-definitions}
\end{equation}
Define $\mathsf{Lev}^Z(u,\xi)$ as the covariance-normalized size of the
heat row. For $1\le p<\infty$,
\begin{equation}
 \begin{aligned}
 \mathsf{Lev}^Z(u,\xi)&:=Q_u^Z(\xi)^*(\mathsf G^Z)^+Q_u^Z(\xi),\\
 \|\mathsf S_{\mathrm{mean}}^Z(\mathbf d)\|_{L^p}
    &\le C_p(\mathbf d^*(\Sigma^Z)^+\mathbf d)^{1/2},\\
 \|\mathsf S_{\mathrm{cov}}^Z(u,\xi)\|_{L^p}
    &\le C_p\mathsf{Lev}^Z(u,\xi).
 \end{aligned}
 \label{eq:multipoint-score-products}
\end{equation}
At any fixed finite derivative order, higher and mixed derivative weights obey
the product of the corresponding direction bounds in
\eqref{eq:multipoint-score-products}, with a constant depending on
that order and $p$. Signed covariance inputs use their absolute
values in this positive bound.

\item \emph{Integrated covariance-derivative bound.}
The normalized row size has total mass at most $q$. It also admits
a singular time weight: for $-1<\gamma_{\mathrm{wt}}<0$,
\begin{equation}
 \begin{aligned}
 \int_0^{h}\int_{\mathbb T}
                       \mathsf{Lev}^Z(u,\xi)\,d\xi\,du&\le q,\\
 \int_0^{h}\int_{\mathbb T}
           u^{\gamma_{\mathrm{wt}}}\mathsf{Lev}^Z(u,\xi)\,d\xi\,du
     &\le C_{q,\gamma_{\mathrm{wt}}}s_Z^{\gamma_{\mathrm{wt}}}.
 \end{aligned}
 \label{eq:multipoint-weighted-leverage}
\end{equation}

\begin{samepage}
\item \emph{A spatial coefficient difference stays in $L^2$.}
Keep $x_1=x$ as the output point. For each $i=2,\ldots,q$, choose
an earlier-center index $1\le\mathrm{par}(i)<i$ and a heat time
$r_i>0$. The history inside the integral is
$Z=((T_1,x),(T_2,x_2),\ldots,(T_q,x_q))$.
For a measurable input $f:(0,h)\times\mathbb T\to\mathbb R$,
write $f_u=f(u,\cdot)$ and define
\begin{equation}
 \begin{aligned}
 I_f(x):=\int_{\mathbb T^{q-1}}
 &\left[\int_0^h\int_{\mathbb T}|f_u(\xi)|\,
                  \mathsf{Lev}^Z(u,\xi)\,d\xi\,du\right]\\
 &\qquad\times\prod_{i=2}^q p_{r_i}(x_{\mathrm{par}(i)}-x_i)
                                       \,dx_2\cdots dx_q,
 \end{aligned}
 \label{eq:heat-averaged-covariance-input}
\end{equation}
Each heat kernel has mass one; when $q=1$, the product and the
integrals in $x_2,\ldots,x_q$ are absent.
If $\|f_u\|_{L^2(\mathbb T)}\le\mathsf A_{\mathrm{in}}u^{\gamma_{\mathrm{wt}}}$,
where $\mathsf A_{\mathrm{in}}\ge0$ and $-1<\gamma_{\mathrm{wt}}\le0$, then,
uniformly over the heat times $r_i$ and indices $\mathrm{par}(i)$,
\begin{equation}
 \|I_f\|_{L^2(\mathbb T)}
 \le C_{q,\gamma_{\mathrm{wt}}}s_Z^{\gamma_{\mathrm{wt}}}
                                      \mathsf A_{\mathrm{in}}.
 \label{eq:multipoint-weighted-input-L2}
\end{equation}
\end{samepage}
\end{enumerate}
\end{lemma}
The proof is given in Appendix~\ref{app:score-moments}.
Only the common variance bounds enter the estimates, so they also
apply to the reserved Gaussian fractions and interpolated covariances
used in the comparisons below.

\paragraph{\textbf{Mean changes: preserve an initial spatial difference.}}
Part~(ii) of the next lemma is the $L^2$ mean-input estimate.
Part~(i) keeps a point input localized, and part~(iii) supplies a
common Gaussian shift for simultaneous state smoothing.
In the notation of
Lemma~\ref{lem:gaussian-one-input-scores}, the size of a linear Gaussian weight
in direction $\mathbf d$ is $(\mathbf d^*(\Sigma^Z)^+\mathbf d)^{1/2}$.

\begin{lemma}[Gaussian derivative bounds for physical heat directions]
\label{lem:physical-heat-score}
Use the history $Z=((T_i,x_i))_{i=1}^q$, its covariance
$\Sigma^Z$ from \eqref{eq:history-full-observation-covariance},
and $s_Z=\min_iT_i>0$, the earliest elapsed observation time.
The constants below depend only on $q$ and the ellipticity bounds.
\begin{enumerate}[label=(\roman*),leftmargin=*]
\item \emph{A point input has a heat-kernel envelope.}
For $\xi\in\mathbb T$, a unit point input at $\xi$ changes the
observed means in direction
$\mathbf d^Z(\xi)=(p_{T_i}(x_i-\xi))_{i=1}^q\in\mathbb R^q$.
Its Gaussian mean-derivative weight satisfies
\begin{equation}
 (\mathbf d^Z(\xi)^*(\Sigma^Z)^+\mathbf d^Z(\xi))^{1/2}
 \le C_q\sum_{i=1}^qT_i^{-1/4}p_{C_qT_i}(x_i-\xi).
 \label{eq:multipoint-first-score-envelope}
\end{equation}
Thus the mean-derivative bound remains localized around the
observation points, on their respective heat scales.

\item \emph{An initial $L^2$ difference stays in $L^2$.}
For $d\in L^2(\mathbb T)$ and $\Lambda\ge0$, its damped heat
evolution is $d_t=e^{-\Lambda t}S_td$. The corresponding vector
of observed mean changes is
$\mathbf d_{\mathrm{obs}}^Z=(d_{T_i}(x_i))_{i=1}^q\in\mathbb R^q$.
Average its squared mean-derivative bound over the observation centers using
the heat-kernel products of Lemma~\ref{lem:gaussian-one-input-scores}(iii):
for each $i\ge2$, choose an earlier index $\mathrm{par}(i)<i$ and a
heat time $r_i>0$. With
$\mathbf x=(x_1,\ldots,x_q)$ and $d\mathbf x=dx_1\cdots dx_q$,
\begin{equation}
 \int_{\mathbb T^q}(\mathbf d_{\mathrm{obs}}^Z)^*(\Sigma^Z)^+\mathbf d_{\mathrm{obs}}^Z
       \prod_{i=2}^q p_{r_i}(x_{\mathrm{par}(i)}-x_i)\,d\mathbf x
                  \le C_qs_Z^{-1/2}\|d\|_{L^2(\mathbb T)}^2.
 \label{eq:multipoint-mean-energy}
\end{equation}
The estimate is uniform in $\Lambda$ and in the averaging times $r_i$.

\item \emph{All observations have a common Gaussian shift.}
Let $\mathbf1=(1,\ldots,1)\in\mathbb R^q$. Then
\begin{equation}
 \Sigma^Z\ge c_q\sqrt{s_Z}\,\mathbf1\mathbf1^*.
 \label{eq:multipoint-common-translation}
\end{equation}
In particular, the centered observation vector can be decomposed
in law into a common scalar Gaussian shift of variance $c_q\sqrt{s_Z}$
and an independent Gaussian residual. This is the variance used
for simultaneous smoothing in the scalar state variable.
\end{enumerate}
\end{lemma}
The proof is given in Appendix~\ref{app:physical-directions}.
For a positive lower time bound $\ell$, when all $T_i\ge\ell$, the first envelope is also bounded by
$C_q\ell^{-1/2}\sum_{i=1}^qT_i^{1/4}p_{C_qT_i}$; positive old heat weights
therefore preserve a physical output heat time. Observations with an
additional heat-smoothing step are treated in the next subsection.

\subsubsection{Conditioning, heat smoothing and bounded profile errors}
\label{sec:conditional-gaussian-inputs}
\label{sec:finite-heat-span-localization}
This subsection supplies two extensions of the preceding input bounds.
First, a point-input derivative weight still has a spatial heat-kernel envelope
when one conditions at an intermediate time and smooths each observation.
Second, a bounded covariance-profile error is controlled by its supremum
norm, uniformly in the observation times and centers.

\paragraph{\textbf{Point inputs after conditioning and heat smoothing.}}
For each observation count $q$, let $h_q^{\mathrm{span}}>0$ be the
time threshold in Lemma~\ref{lem:finite-line-heat-span} of
Appendix~\ref{app:heat-spans}. Below this threshold, a span of $q$
torus heat kernels has a uniform spatial envelope, including when
centers or times collide. The following corollary is the derivative weight
estimate obtained from that supporting geometric fact.

\begin{corollary}[Point-input derivative weights for heat-smoothed observations]
\label{cor:heat-smoothed-observation-score}
Use the same reference $X$ as above, and let $t_{\mathrm{in}}\in[0,h)$ be
the elapsed conditioning time, measured from the window start $t_0$. Record the observation data as
$Z=((\sigma_i,T_i,x_i))_{i=1}^q$, where
$t_{\mathrm{in}}<\sigma_i\le T_i$, $\sigma_i\le h$ and $x_i\in\mathbb T$.
The observation
$\mathcal G_i^{Z,\mathrm{obs}}=(S_{T_i-\sigma_i}X_{\sigma_i})(x_i)$
records the field at time $\sigma_i$ and then smooths it for time
$T_i-\sigma_i$. Conditional on $\mathcal F_{t_0+t_{\mathrm{in}}}$, its covariance
matrix $\Sigma_{t_{\mathrm{in}}}^Z$ has entries
\begin{equation}
 (\Sigma_{t_{\mathrm{in}}}^Z)_{ik}=\int_{t_{\mathrm{in}}}^{\min(\sigma_i,\sigma_k)}\int_{\mathbb T}
 (m_0^2+\bar a_u(\xi))p_{T_i-u}(x_i-\xi)p_{T_k-u}(x_k-\xi)
 \,d\xi\,du.
 \label{eq:history-smoothed-observation-covariance}
\end{equation}
Write $\tau_{\min}^{Z,t_{\mathrm{in}}}:=\min_i(\sigma_i-t_{\mathrm{in}})>0$ for the shortest
fresh-noise interval and $t_i:=T_i-t_{\mathrm{in}}$ for the total heat time.
Assume $t_i\le h_q^{\mathrm{span}}$ for every $i$.

A unit point input at $\xi\in\mathbb T$ at time $t_{\mathrm{in}}$ changes the
observed means in direction
$Q_{t_{\mathrm{in}}}^Z(\xi)=(p_{t_i}(x_i-\xi))_{i=1}^q\in\mathbb R^q$.
Its linear Gaussian weight satisfies
\begin{equation}
 (Q_{t_{\mathrm{in}}}^Z(\xi)^*(\Sigma_{t_{\mathrm{in}}}^Z)^+Q_{t_{\mathrm{in}}}^Z(\xi))^{1/2}
 \le C_q(\tau_{\min}^{Z,t_{\mathrm{in}}})^{-1/2}
             \sum_{i=1}^q t_i^{1/4}p_{C_qt_i}(x_i-\xi).
 \label{eq:contracted-first-score-envelope}
\end{equation}
Thus the available Gaussian evolution time controls the derivative weight
factor $(\tau_{\min}^{Z,t_{\mathrm{in}}})^{-1/2}$, while the total heat times $t_i$
determine the spatial kernels.
\end{corollary}

\begin{proof}[Proof of Corollary~\ref{cor:heat-smoothed-observation-score}]
\emph{The common interval supplies a covariance lower bound.}
Define the heat-row Gram matrix
$\mathsf G_{t_{\mathrm{in}}}^Z:=\int_{\mathbb T}Q_{t_{\mathrm{in}}}^Z(\xi)Q_{t_{\mathrm{in}}}^Z(\xi)^*\,d\xi$.
Every observation receives noise throughout
$[t_{\mathrm{in}},t_{\mathrm{in}}+\tau_{\min}^{Z,t_{\mathrm{in}}}]$. For $\mathbf a\in\mathbb R^q$, ellipticity gives
\[
 \mathbf a^*\Sigma_{t_{\mathrm{in}}}^Z\mathbf a
 \ge c\int_0^{\tau_{\min}^{Z,t_{\mathrm{in}}}}
 \left\|\sum_{i=1}^q a_i p_{t_i-u}(x_i-\cdot)
                                      \right\|_{L^2(\mathbb T)}^2\,du.
\]
Applying $S_u$ to the function inside the norm gives
$\sum_i a_i p_{t_i}(x_i-\cdot)$. By $L^2$ contraction, the integrand
is therefore at least
$\|\sum_i a_i p_{t_i}(x_i-\cdot)\|_{L^2(\mathbb T)}^2
=\mathbf a^*\mathsf G_{t_{\mathrm{in}}}^Z\mathbf a$,
independently of $u$. Integrating this constant lower bound over the
common interval of length $\tau_{\min}^{Z,t_{\mathrm{in}}}$ supplies
the time factor and gives
\[
 \Sigma_{t_{\mathrm{in}}}^Z\ge c\tau_{\min}^{Z,t_{\mathrm{in}}}\mathsf G_{t_{\mathrm{in}}}^Z.
\]
\emph{The heat span supplies the spatial envelope.}
For a direction in the range of $\mathsf G_{t_{\mathrm{in}}}^Z$, the inverse variational
formula gives the following quadratic-form bound:
\[
 Q_{t_{\mathrm{in}}}^Z(\xi)^*(\Sigma_{t_{\mathrm{in}}}^Z)^+Q_{t_{\mathrm{in}}}^Z(\xi)
 \le(c\tau_{\min}^{Z,t_{\mathrm{in}}})^{-1}Q_{t_{\mathrm{in}}}^Z(\xi)^*(\mathsf G_{t_{\mathrm{in}}}^Z)^+Q_{t_{\mathrm{in}}}^Z(\xi).
\]
The evaluation identity for the finite heat span is
\[
 Q_{t_{\mathrm{in}}}^Z(\xi)^*(\mathsf G_{t_{\mathrm{in}}}^Z)^+Q_{t_{\mathrm{in}}}^Z(\xi)
 =\sup_{\substack{0\ne f\in\operatorname{span}
                     \{p_{t_i}(x_i-\cdot)\}_{i=1}^q}}
          \frac{|f(\xi)|^2}{\|f\|_{L^2(\mathbb T)}^2}.
\]
Indeed, writing $f=\sum_i a_i p_{t_i}(x_i-\cdot)$ gives
$f(\xi)=\mathbf a^*Q_{t_{\mathrm{in}}}^Z(\xi)$ and
$\|f\|_{L^2(\mathbb T)}^2=\mathbf a^*\mathsf G_{t_{\mathrm{in}}}^Z\mathbf a$.
Cauchy--Schwarz in this Gram form bounds the ratio by the left-hand
side, with equality for
$\mathbf a=(\mathsf G_{t_{\mathrm{in}}}^Z)^+Q_{t_{\mathrm{in}}}^Z(\xi)$.
This remains valid for dependent heat rows: a vector in the null
space of the Gram matrix gives a continuous heat combination that
vanishes identically, so $Q_{t_{\mathrm{in}}}^Z(\xi)$ lies in its range.
Apply \eqref{eq:finite-torus-span-envelope} and take square roots.
The resulting Gaussian envelopes are bounded by heat kernels:
\[
 t_i^{-1/4}e^{-c_q\operatorname{dist}_{\mathbb T}(x_i,\xi)^2/(2t_i)}
 \le C_qt_i^{1/4}p_{C_qt_i}(x_i-\xi).
\]
This yields \eqref{eq:contracted-first-score-envelope}.
\end{proof}

\paragraph{\textbf{Bounded covariance-profile errors.}}
Use the same history $Z=((\sigma_i,T_i,x_i))_{i=1}^q$ and conditioning
time $t_{\mathrm{in}}$ as above. For $t_{\mathrm{in}}\le u\le h$, define
the active noise-response vector by
\[
 Q_u^{Z,\mathrm{obs}}(\xi)_i=
 \begin{cases}
 p_{T_i-u}(x_i-\xi),&u<\sigma_i,\\
 0,&u\ge\sigma_i.
 \end{cases}
\]
Its $i$th coordinate is zero once observation $i$ has stopped receiving
noise. The covariance is $\Sigma_{t_{\mathrm{in}}}^Z$ from
\eqref{eq:history-smoothed-observation-covariance}.
The reference variance satisfies
$\bar{\mathfrak a}_u(\xi)=m_0^2+\bar a_u(\xi)\ge m_0^2$,
where $m_0>0$ is the common additive-noise amplitude.

\begin{corollary}[Bounded covariance-profile errors]
\label{cor:bounded-profile-score}
For every bounded measurable profile error
$f:[t_{\mathrm{in}},h]\times\mathbb T\to\mathbb R$,
\begin{equation}
 \begin{aligned}
 &\int_{t_{\mathrm{in}}}^h\int_{\mathbb T}|f(u,\xi)|\,
    Q_u^{Z,\mathrm{obs}}(\xi)^*(\Sigma_{t_{\mathrm{in}}}^Z)^+
                        Q_u^{Z,\mathrm{obs}}(\xi)\,d\xi\,du\\
 &\hspace{20mm}\le\frac q{m_0^2}
    \|f\|_{L^\infty([t_{\mathrm{in}},h]\times\mathbb T)}.
 \end{aligned}
 \label{eq:bounded-input-leverage}
\end{equation}
The constant is independent of the observation times and centers.
\end{corollary}

\begin{proof}[Proof of Corollary~\ref{cor:bounded-profile-score}]
The covariance formula gives
\[
 \begin{aligned}
 &\int_{t_{\mathrm{in}}}^h\int_{\mathbb T}\bar{\mathfrak a}_u(\xi)
 Q_u^{Z,\mathrm{obs}}(\xi)^*(\Sigma_{t_{\mathrm{in}}}^Z)^+
                         Q_u^{Z,\mathrm{obs}}(\xi)\,d\xi\,du\\
 &\qquad=\operatorname{tr}\bigl((\Sigma_{t_{\mathrm{in}}}^Z)^+
                                      \Sigma_{t_{\mathrm{in}}}^Z\bigr)\le q.
 \end{aligned}
\]
Use $\bar{\mathfrak a}\ge m_0^2$ and bound $|f|$ by its supremum.
For a reserved Gaussian fraction with variance bounded below by
$\kappa_->0$, the same calculation gives the factor $q/\kappa_-$.
\end{proof}

By the covariance-derivative bound in
Lemma~\ref{lem:gaussian-one-input-scores}(i),
\eqref{eq:bounded-input-leverage} controls the Gaussian derivative
caused by a uniformly bounded profile error. It supplies the
total-mass bound in the higher-martingale and paired-remainder
comparisons, Propositions~\ref{prop:higher-martingale-comparison}
and~\ref{prop:paired-actual-remainder}, and in the coefficient
replacement of Lemma~\ref{lem:exact-bracket-cutoff-limits}.
When the profile input is measured in $L^2$ instead, the relevant
bound is Lemma~\ref{lem:gaussian-one-input-scores}(iii).
The general response-row version of the total-mass calculation is
Lemma~\ref{lem:total-leverage} in Appendix~\ref{app:leverage}.

\begin{samepage}
\subsubsection{Integration of positive time kernels}

\begin{definition}[Clipped root time weight and logarithmic factors]
\label{def:clipped-root-time-weights}
Only one clipped root weight is needed. For $\lambda\ge1$ define
\begin{equation}
\begin{gathered}
 k_\lambda(t)=t^{-1}\min(1,\lambda t^{1/4}),\qquad t>0,\\
 L_\lambda(w)=\int_0^w k_\lambda(t)\,dt
  =\begin{cases}
    4\lambda w^{1/4},&0<w\le\lambda^{-4},\\
    4+\log(\lambda^4w),&w\ge\lambda^{-4}.
   \end{cases}
\end{gathered}
\label{eq:clipped-root-time-weights}
\end{equation}
Write $\mathcal L(\lambda,h)=1+\log(2+\lambda)+|\log h|$.
\end{definition}

A coefficient factor receiving centered differences supplies a power of
its smallest assigned gap. Bounding
that minimum by the geometric mean distributes the same total gain
among those gaps; the resulting powers can be integrated directly.

\begin{lemma}[Integration of insertion gaps]
\label{lem:insertion-gap-integration}
For an integer $m\ge1$, $\eta>0$ and positive gaps $\delta_1,\ldots,\delta_m$,
\begin{equation}
 (\min_i\delta_i)^\eta\prod_{i=1}^m\delta_i^{-1}
 \le\prod_{i=1}^m\delta_i^{\eta/m-1},\qquad
 \int_{(0,h)^m}(\min_i\delta_i)^\eta
                 \prod_{i=1}^m\frac{d\delta_i}{\delta_i}
 \le(m/\eta)^m h^\eta.
 \label{eq:minimum-gap-cube-bound}
\end{equation}
More generally, fix an integer $n\ge0$, a root index set
$I_0\subset\{1,\ldots,n\}$, and positive exponents $\alpha_j$
for $j\notin I_0$. Put $r=|I_0|$ and
$\gamma=\sum_{j\notin I_0}\alpha_j$.
For $\sigma>-1$, $h>0$, $s=h-\sum_{j=1}^n\delta_j$ and
$d\boldsymbol\delta=\prod_{j=1}^n d\delta_j$,
\begin{equation}
 \int_{\substack{\delta_j>0\\\sum_j\delta_j<h}}
 \prod_{j\in I_0}k_\lambda(\delta_j)
 \prod_{j\notin I_0}\delta_j^{\alpha_j-1}
                       s^\sigma\,d\boldsymbol\delta
 \le C_{n,\alpha,\sigma}
                      h^{\gamma+\sigma}L_\lambda(h)^r.
 \label{eq:insertion-gap-simplex-bound}
\end{equation}
Empty products are one; for $n=0$ the left side is $h^\sigma$.
The constant may depend on the positive exponents $\alpha_j$,
but is independent of $\lambda$ and $h$.
\end{lemma}
\end{samepage}
The proof is given in Appendix~\ref{app:positive-gaps}.

\begin{lemma}[Squared time integrals]
\label{lem:outer-square-integrals}
For $\lambda\ge1$, an integer $r\ge0$, $\vartheta>0$ and $h>0$,
\begin{equation}
 \int_0^h w^{\vartheta-1}L_\lambda(w)^{2r}\,dw
 \le\vartheta^{-1}h^\vartheta L_\lambda(h)^{2r}.
 \label{eq:positive-square-integral}
\end{equation}
For $r\ge1$ the critical integral satisfies
\begin{equation}
 \int_0^h w^{-1}L_\lambda(w)^{2r}\,dw
 \le C_r[1+L_\lambda(h)]^{2r+1}.
 \label{eq:root-critical-square-integral}
\end{equation}
For the current-state comparison, let $b,r\ge0$ be integers with
$b+r\ge1$, $0<\beta<1$, $0<\rho<1/2$ and
$0<\chi<\min(1,\beta+2\rho)$. Writing $z_+=\max(z,0)$, for every
$\varepsilon>0$ and $\lambda^{-4}\le h\le1$,
\begin{equation}
 \left(\int_0^h w^{b\beta/2+\rho-\chi/2-1}
                         L_\lambda(w)^{2r}\,dw\right)^{1/2}
 \le C_{b,r,\beta,\rho,\chi,\varepsilon}
 \lambda^{(\chi-2\rho)_++\varepsilon}
 h^{b\beta/4+(2\rho-\chi)_+/4}.
 \label{eq:current-mean-square-integral}
\end{equation}
All constants are independent of $\lambda$ and $h$.
\end{lemma}
The proof is given in Appendix~\ref{app:clipped-convolutions}.

\paragraph{Differentiating the Gaussian references.}
Changing the reference mean, covariance and subtracted profiles
produces three corresponding contributions. The product
identity, including its hypotheses and proof, is
Lemma~\ref{lem:three-leg-product} in
Appendix~\ref{sec:full-product-identity}.

\clearpage
\flushbottom
\section{Comparison estimates for the finite expansion}
We compare coefficient observations from the coupled solutions $U$
and $V$ by separating an average determined by the initial past,
a stochastic fluctuation and a remainder. We first describe these
objects through a common Gaussian interpolation, then estimate their
changes between the two solutions.
\begin{enumerate}[label=(\roman*),leftmargin=*]
\item \textbf{The initial-past averages.} Proposition~\ref{prop:all-order-predictable-pair} controls their change when the initial mean, reference covariance and subtracted profile values vary together.
\item \textbf{The leading stochastic fluctuation.} Proposition~\ref{prop:complete-martingale-comparison} compares the first martingale after subtracting its Gaussian-reference fluctuation.
\item \textbf{The higher stochastic fluctuations.} Proposition~\ref{prop:higher-martingale-comparison} directly compares the martingale words with at least one covariance insertion.
\item \textbf{The unexpanded remainder.} Proposition~\ref{prop:paired-actual-remainder} shows that increasing the fixed expansion depth improves its comparison bound.
\end{enumerate}
\paragraph{Reading route.}
Read the definitions and four statements below, then the common
centering and spatial estimates in Section~\ref{sec:centered-word-inputs}.
The Gaussian-average proof is in Section~\ref{sec:predictable-comparison-proof};
the stochastic comparisons are in Section~\ref{sec:martingale-comparisons}.
The remainder, initial-window and frequency estimates follow in
Sections~\ref{sec:remainder-comparisons}--\ref{sec:frequency-summation-energy}.
Section~\ref{sec:approximation-and-transfer} justifies the approximation
limits; Section~6 checks the change of coordinates. Supporting calculations are in
Appendices~\ref{app:finite-word-details}--\ref{app:approximation-details}.

\subsection{Main results}
\label{sec:comparison-main-results}

\subsubsection{Data and comparison norms}
Let $U$ be the given weak solution and $V$ its smooth-coefficient
companion with extra drift $\Lambda(U-V)$, where $\Lambda\ge1$
is the coupling strength. Fix a time horizon $T>0$ and work on a
window $[t_0,t_0+h]\subset[0,T]$, where $t_0$ is the physical starting
time and $h$ the duration. At elapsed time $s\in[0,h]$, the physical
state is $U_{t_0+s}$. In the formulas below, $t,q,r$ also denote local
elapsed times when so specified; ordered histories use
$h=s_0>s_1>\cdots>s_n>0$.
Condition on the initial past $\mathcal F_{t_0}$ and write
$\mathbb E_{t_0}=\mathbb E[\,\cdot\mid\mathcal F_{t_0}]$.
Here $g$ is the diffusion coefficient, $m_0>0$ is the common
additive-noise amplitude, and $a(z)=g(z)^2-m_0^2$ is the residual
variance. The heat semigroup on $\mathbb T$ and its kernel are
$S_t$ and $p_t$.
Choose a dyadic state-frequency label
$\lambda\in\mathscr S_{\mathrm z}=\{1,2,4,\ldots\}$:
$a_\lambda$ is the Fourier component of $a$ with oscillation scale
$\lambda^{-1}$ in its scalar argument $z$.
Assume $a$ is bounded, uniformly positive and $\beta$-H\"older,
$0<\beta<1$. In the main comparison statements we use the same
band duration as Theorem~\ref{thm:coupled-gaussian-band-sum}:
$h=h_\lambda:=\lambda^{-3\beta}$. We write $h$ in the formulas to
show the time gain from each insertion. Here $3\beta$ is the chosen
value of $\gamma_{\mathrm{time}}$ in the relation
$h_\lambda=\lambda^{-\gamma_{\mathrm{time}}}$: it determines how
the duration shrinks with frequency, while $\beta$ is the
coefficient's H\"older exponent. All sums over high bands
start at the fixed cutoff \eqref{eq:gaussian-band-dyadic-cutoff},
$\lambda_*\asymp\Lambda^{1/4}$.
The proofs establish the estimates for general admissible durations;
the initial-time extensions use each window's actual shorter duration.

\paragraph{\textbf{Comparison norms and time thresholds.}}
Use a spatial regularity exponent $0<\rho<1/2$ and
$H=L^2(\mathbb T)$. The coupling strength $\Lambda$ also fixes the
scale in the negative Sobolev norm.
The spatial norm is
$\|f\|_{H_\Lambda^{-\rho}}=\|(\Lambda-\Delta)^{-\rho/2}f\|_{L^2(\mathbb T)}$.
For a window output $Q(t_0)$, an unqualified output norm means
\[
 \|Q\|=\|Q\|_{\mathscr L_h^2(H_\Lambda^{-\rho})}
 =\left(\int_0^{T-h}\mathbb E
             \|Q(t_0)\|_{H_\Lambda^{-\rho}}^2\,dt_0\right)^{1/2}.
\]
The case $h=0$ gives the full-horizon norm.
Let $a_{n_{\mathrm{sm}}}$ be the state-mollified coefficient used by
$V$; its index $n_{\mathrm{sm}}$ records coefficient smoothing,
independently of the expansion depth $n$.
The solution and coefficient differences are measured by
$\|U-V\|_{\mathscr L_0^2(H)}$ and
$\mathsf d_{\mathrm{coef}}(U,V)^2=\|a(U)-a_{n_{\mathrm{sm}}}(V)\|_{\mathscr L_0^2(H_\Lambda^{-\rho})}^2
+\varepsilon_{\mathrm{app}}$, where $\varepsilon_{\mathrm{app}}\ge0$
is a deterministic approximation source.
The common short-time bound $h_0\in(0,1]$ comes from
Propositions~\ref{prop:self-profile-compatible-gaussian}
and~\ref{prop:self-profile}: on $0<h\le h_0$, the Gaussian profiles
exist uniquely and satisfy the stated stability estimates.
It depends only on the coefficient bounds, $\beta$ and $m_0$,
uniformly over the initial states, block labels and coupling strengths.
The order-dependent heat-span threshold is specified immediately
before Proposition~\ref{prop:all-order-predictable-pair}.

\subsubsection{Gaussian interpolation and average terms}
\label{sec:all-order-covariance-martingale}
\label{sec:all-order-covariance-strong-paired}
\begin{definition}[Gaussian average terms for the coupled Gaussian references]
\label{def:predictable-mean-functional}
For $0\le\theta\le1$, let $F_0^\theta:\mathbb R\to\mathbb R$
denote the initial scalar test function; throughout this section
we take $F_0^\theta=a_\lambda$ for the fixed block $\lambda$,
independently of $\theta$. Write $d=U_{t_0}-V_{t_0}\in H$ and define
$\phi_t^\theta=S_tU_{t_0}-(1-\theta)e^{-\Lambda t}S_td$,
$0\le t\le h$.
This function is the spatial mean of the Gaussian reference at
elapsed time $t$. Its endpoints are explicitly
$\phi_t^0=S_tU_{t_0}-e^{-\Lambda t}S_t(U_{t_0}-V_{t_0})$ and
$\phi_t^1=S_tU_{t_0}$; they are the $V$-side and $U$-side means,
respectively. Each path lies in $C([0,h];H)$ and is smooth in $x$
for $t>0$. The superscripts label the references, while $t$ is elapsed time.
Let $\bar a^\theta$ be the self-consistent profile for this mean,
whose existence is given by Proposition~\ref{prop:self-profile}.
For a list of $d_{\mathrm{obs}}\ge1$ targets
$Z=((T_i,x_i))_{i=1}^{d_{\mathrm{obs}}}$, with $0<T_i\le h$ and
$x_i\in\mathbb T$, define $\mathcal P^{Z,\theta}_{q,s}$ as centered
Gaussian convolution on $\mathbb R^{d_{\mathrm{obs}}}$ with covariance
\begin{equation}
 (\Sigma^{Z,\theta}_{q,s})_{ik}
 =\int_q^s\int_{\mathbb T}(m_0^2+\bar a_r^\theta(\eta))
 p_{T_i-r}(x_i-\eta)p_{T_k-r}(x_k-\eta)\,d\eta\,dr,
 \qquad 0\le q<s\le\min_iT_i.
 \label{eq:predictable-mean-covariance}
\end{equation}
In particular, the self-consistency condition says
$\bar a_t^\theta(x)=
(\mathcal P^{((t,x)),\theta}_{0,t}a)(\phi_t^\theta(x))$.
This displays both dependencies: the profile determines the
Gaussian covariance, and the prescribed mean is the argument at
which the Gaussian average is evaluated.

For a history $h=s_0>s_1>\cdots>s_n>0$ and points
$\xi_0=x,\xi_1,\ldots,\xi_n\in\mathbb T$, write
$Z_j=((s_0,\xi_0),\ldots,(s_j,\xi_j))$.
Starting from this $F_0^\theta$, for $1\le j\le n$ define the
inserted tests by
\begin{equation}
 \begin{aligned}
 F_j^\theta(z_0,\ldots,z_j)
 :={}&\frac12[a(z_j)-\bar a_{s_j}^\theta(\xi_j)]\\
 &\times\left[
 \left(\sum_{i=0}^{j-1}p_{s_i-s_j}(\xi_i-\xi_j)\partial_{z_i}\right)^2
 \mathcal P^{Z_{j-1},\theta}_{s_j,s_{j-1}}F_{j-1}^\theta
 \right](z_0,\ldots,z_{j-1}).
 \end{aligned}
 \label{eq:predictable-mean-test-recursion}
\end{equation}
The dependence of $F_j^\theta$ on the chosen time--space history
is understood. Each step first Gaussian-smooths the previous test,
then differentiates twice in the heat-row direction, and finally
multiplies by $\tfrac12[a(z_j)-\bar a_{s_j}^\theta(\xi_j)]$.
The new variable $z_j$ is included in the Gaussian smoothing at
the next step. The factor $1/2$ comes from the covariance term in
It\^o's formula.
The initial block $a_\lambda$ is the \emph{root}; the inserted
discrepancies are \emph{rough coefficient factors}. The value
subtracted in the $j$th factor is $\bar a_{s_j}^\theta(\xi_j)$.
Each factor uses the full coefficient $a$.
The general Gaussian average term is explicitly
\begin{equation}
 \begin{aligned}
 \mathsf P_{n+1,\lambda}^\theta(h,x)
 :=\int_{0<s_n<\cdots<s_1<h}\int_{\mathbb T^n}
 &\bigl(\mathcal P^{Z_n,\theta}_{0,s_n}F_n^\theta\bigr)
 \bigl(\phi_h^\theta(x),\phi_{s_1}^\theta(\xi_1),\ldots,
                         \phi_{s_n}^\theta(\xi_n)\bigr)\\
 &\hspace{3em}\,d\xi_1\cdots d\xi_n\,ds_n\cdots ds_1.
 \end{aligned}
 \label{eq:predictable-mean-word-formula}
\end{equation}
For $n=0$, the history integral is absent and $s_0=h$.
The tests $F_j^\theta$ depend on $\lambda$ through the root
$a_\lambda$; the means, profiles and Gaussian covariances use the
full coefficient $a$ and are common to all blocks at a fixed window.
We call $\mathsf P_{n+1,\lambda}^\theta$ a \emph{Gaussian average term}:
all future Gaussian increments have been averaged out, so its value
is determined by the initial past $\mathcal F_{t_0}$.
\end{definition}

\paragraph{\textbf{Endpoint cases and a short example.}}
For the $U$ expansion ($\theta=1$), write $\bar a=\bar a^1$,
$\mathcal P^Z=\mathcal P^{Z,1}$ and
$\mathsf P_{j,\lambda}=\mathsf P_{j,\lambda}^1$; $\theta=0$ gives
the $V$-side reference.
For a target list $Z$, the physical backward observation vector is
\[
 \mathbf Y_q^Z=
 \bigl((S_{T_i-q}U_{t_0+q})(x_i)\bigr)_{i=1}^{d_{\mathrm{obs}}},
 \qquad 0\le q\le\min_iT_i.
\]
Applying It\^o's formula to these physical observations with the
$\theta=1$ reference gives the martingale terms $N_{j,\lambda}$ and
remainders $R_{n,\lambda}$.

For example, with $Z_0=((h,x))$ and
$Z_1=((h,x),(s_1,\xi))$, the first two Gaussian average terms are
\begin{equation}
 \begin{aligned}
 \mathsf P_{1,\lambda}^\theta(h,x)
 &=\bigl(\mathcal P^{Z_0,\theta}_{0,h}a_\lambda\bigr)
                                      (\phi_h^\theta(x)),\\
 \mathsf P_{2,\lambda}^\theta(h,x)
 &=\int_0^h\int_{\mathbb T}
 \bigl(\mathcal P^{Z_1,\theta}_{0,s_1}F_1^\theta\bigr)
       (\phi_h^\theta(x),\phi_{s_1}^\theta(\xi))\,d\xi\,ds_1.
 \end{aligned}
 \label{eq:predictable-first-two-interpolated}
\end{equation}
Here $F_1^\theta$ is the one-insertion test in
\eqref{eq:predictable-mean-test-recursion}. Thus the first term is
one Gaussian average of $a_\lambda$; the second averages the test
with one inserted covariance discrepancy, then integrates its
insertion point. The detailed formulas for $N_{1,\lambda}$,
$N_{2,\lambda}$ and the first two remainders appear in
Appendix~\ref{sec:first-two-expansion-orders}.

\paragraph{\textbf{Expansion order and continuous time variables.}}
The full expansion at depth $n$ is
\[
 a_\lambda(U_{t_0+h}(x))
 =\sum_{j=1}^n(\mathsf P_{j,\lambda}+N_{j,\lambda})+R_{n,\lambda}.
\]
Here $\mathsf P_{j,\lambda}$ is the Gaussian average fixed by the initial
past, $N_{j,\lambda}$ is the conditionally centered stochastic integral,
and $R_{n,\lambda}$ is the term left after $n$ insertions. The words
$\mathsf P_{j,\lambda},N_{j,\lambda}$ each contain $j-1$ insertions, so their subscripts
record expansion order. At depth $n$, the times
$h>s_1>\cdots>s_n>0$ range over the full ordered simplex and each
position $\xi_i$ ranges over $\mathbb T$. These are continuous
integration variables. Higher orders repeat the same It\^o
decomposition with the inserted tests of
Definition~\ref{def:predictable-mean-functional} at $\theta=1$.
The terms $N_{j,\lambda}$ are stochastic integrals.
Each comparison below specifies which means and profiles are used.

The block indices make explicit the dependence suppressed in
$\mathsf P_j,N_j,R_n$ of Section~\ref{sec:recursive-gaussian-expansion}.
For the companion, the terminal test is still $a_\lambda$, while
the physical states and Gaussian reference are those of $V$ and
each local insertion uses $a_{n_{\mathrm{sm}}}(V)-\bar a^V$.
Its control drift $\Lambda(U-V)$ also contributes to It\^o's formula.
The joint coordinates of Section~\ref{sec:single-species-joint-flow-recursion}
cancel that drift and record the additional covariance terms
separately; \eqref{eq:joint-V-finite-expansion} displays the resulting
$V$ expansion before subtraction from the $U$ expansion.

\paragraph{\textbf{Continuous interpolation between the endpoints.}}
The auxiliary parameter $\theta$ lets us compare endpoints by
integration of a derivative, with
$\partial_\theta\phi_t^\theta=e^{-\Lambda t}S_td$.
Recomputing the profile and its subtracted values along this path preserves
their centering. This gives the linear dependence on
$\|d\|_{L^2(\mathbb T)}$ in the next result.

\paragraph{\textbf{Admissible window lengths.}}
For $q$ observations, the heat-span estimate has a positive
short-time bound $h_q^{\mathrm{span}}$: whenever all target ages are
at most this bound, the normalized observation rows have the uniform
Gaussian spatial envelopes used to estimate the derivative weights.
These are the bounds of Lemma~\ref{lem:finite-line-heat-span}.
A word with $n$ insertions has at most $n+1$ observations, so set
\begin{equation}
 h_n:=\min\{h_0,h_q^{\mathrm{span}}:1\le q\le n+1\}.
 \label{eq:finite-order-window-threshold}
\end{equation}
Thus $0<h_{n+1}\le h_n\le h_0$. For each fixed order $n$, $h_n$
is a constant independent of $\lambda$ and $\Lambda$.
Once the maximal depth $N_\beta$ is fixed, the single threshold
$h_{N_\beta}$ works for all its orders. This depth and its positive
time threshold are fixed before the frequency and coupling
parameters tend to infinity.

\begin{samepage}
\begin{proposition}[Comparison of Gaussian average terms for one coefficient block]
\label{prop:all-order-predictable-pair}
Fix a dyadic block $a_\lambda$ and $n\ge1$. Take
$h=h_\lambda=\lambda^{-3\beta}\le h_n$, with $h_n$ defined by
\eqref{eq:finite-order-window-threshold}. Put
$\mathcal L(\lambda,h):=1+\log(2+\lambda)+|\log h|$.
For the words of Definition~\ref{def:predictable-mean-functional},
both with initial test $F_0^\theta=a_\lambda$, and
$d=U_{t_0}-V_{t_0}$, the spatial comparison bounds are
\begin{align}
 \|\mathsf P_{n+1,\lambda}^1-\mathsf P_{n+1,\lambda}^0\|_{L^2(\mathbb T)}
 &\le C_n\mathcal L(\lambda,h)^{C_n}
       \lambda^{-\beta}h^{n\beta/4-1/4}\|d\|_{L^2(\mathbb T)},
 \label{eq:all-order-predictable-clipped-linear-pair}\\
 \|\mathsf P_{n+1,\lambda}^1-\mathsf P_{n+1,\lambda}^0\|_{L^2(\mathbb T)}
 &\le C_n\mathcal L(\lambda,h)^{C_n}\lambda^{-\beta}h^{n\beta/4}
                      \min\{1,h^{-1/4}\|d\|_{L^2(\mathbb T)}\}.
 \label{eq:all-order-predictable-clipped-capped-pair}
\end{align}
The constants depend on $n$ and the common coefficient bounds,
and are independent of the block label $\lambda$, window length
$h$ and coupling strength $\Lambda$.

\begin{samepage}
For $0<\rho<1/2$, the form used in the later comparison is
\begin{equation}
 \begin{aligned}
 &\|\mathsf P_{n+1,\lambda}^1-\mathsf P_{n+1,\lambda}^0
                           \|_{\mathscr L_h^2(H_\Lambda^{-\rho})}\\
 &\qquad\le C_n\Lambda^{-\rho/2}\mathcal L(\lambda,h)^{C_n}
             \lambda^{-\beta}h^{n\beta/4-1/4}\|U-V\|_{\mathscr L_0^2(H)}.
 \end{aligned}
 \label{eq:all-order-predictable-integrated-sobolev}
\end{equation}
\end{samepage}
\end{proposition}
\end{samepage}

\paragraph{\textbf{Reading guide for Proposition~\ref{prop:all-order-predictable-pair}.}}
The right-hand side of the comparison estimate
\eqref{eq:all-order-predictable-clipped-linear-pair} is the endpoint
bound \eqref{eq:centered-Gaussian-average-size} multiplied by
$h^{-1/4}\|d\|_{L^2(\mathbb T)}$.
The Gaussian law and the subtracted profiles change together,
so all three derivative contributions must satisfy this bound. In Section~\ref{sec:predictable-comparison-proof},
the proof first keeps the old centering increments fixed, then
bounds the mean, covariance and subtraction terms in $L^2(\mathbb T)$. Integrating the interpolation gives the comparison;
the separate endpoint size gives its cap. The later frequency form
is \eqref{eq:ordinary-predictable-comparison-row}.

\subsubsection{The first and higher martingales}
We compare the leading actual stochastic term after subtracting
its Gaussian-reference fluctuation.
\label{sec:direct-actual-gamma-negative-input}
\label{sec:capped-actual-gamma-independent}
\label{sec:actual-gamma-current-state-cap}

\begin{definition}[Backward fields and complete martingales for the coupled pair]
\label{def:coupled-backward-martingales}
Use the window, block $a_\lambda$, conditional expectation
$\mathbb E_{t_0}$ and shifted noises fixed above. All times $s,r\in[0,h]$
are elapsed from $t_0$.
\begin{enumerate}[label=(\roman*),leftmargin=*]
\item \emph{The actual backward fields and their noises.}
Put $d=U_{t_0}-V_{t_0}$. Let $X^U,X^V$ be the common-noise
Gaussian references of Definition~\ref{def:coupled-gaussian-references},
with means $\phi^U=\phi^1$, $\phi^V=\phi^0$ and profiles
$\bar a^U=\bar a^1$, $\bar a^V=\bar a^0$ from
Definition~\ref{def:predictable-mean-functional}.
The actual backward fields are
\[
 \widehat U_s^U=S_{h-s}U_{t_0+s},\qquad
 \widehat U_s^V=S_{h-s}V_{t_0+s}
                 +(e^{-\Lambda s}-e^{-\Lambda h})S_hd.
\]
Both end at the corresponding physical state at time $t_0+h$.
Their shared shifted Brownian motions are
$W_s^{0,t_0}=W_{t_0+s}^0-W_{t_0}^0$ and
$B_s^{t_0}=B_{t_0+s}-B_{t_0}$.
Here $m_0>0$ is the common additive-noise amplitude,
$a=g^2-m_0^2$ is the residual variance, and
$a_{n_{\mathrm{sm}}}$ is its state smoothing used by $V$.
The actual residual-noise operators are
$L_U(s)=L(a(U_{t_0+s}))$ and
$L_V(s)=L(a_{n_{\mathrm{sm}}}(V_{t_0+s}))$.
The factor satisfies $L(\kappa)L(\kappa)^*=\mathsf M_\kappa$, and
$\mathsf M_\kappa$ is multiplication by $\kappa$.

\item \emph{The first actual stochastic integral.}
For $i\in\{U,V\}$ and $x\in\mathbb T$, the remaining reference
variance at time $s$ is
\begin{equation}
 \upsilon_s^i(x):=\int_s^h\int_{\mathbb T}p_{h-r}(x-y)^2
                  (m_0^2+\bar a_r^i(y))\,dy\,dr.
 \label{eq:complete-martingale-remaining-variance}
\end{equation}
Set $F_s^i(x,z):=(\mathcal P_{\upsilon_s^i(x)}a_\lambda)(z)$,
where $(\mathcal P_vf)(z)=\mathbb E[f(z+\sqrt v\,\zeta_{\mathrm G})]$
for $\zeta_{\mathrm G}\sim N(0,1)$.
Thus $F_s^i$ is the fixed block $a_\lambda$ smoothed by the
Gaussian noise still to come, and $\partial_zF_s^i$ denotes its
state derivative. Define the stochastic-integral process by
\begin{equation}
 \begin{aligned}
 N_{1,\lambda}^i(s,x):={}&m_0\int_0^s
 \partial_zF_r^i(x,\widehat U_r^i(x))
                         (S_{h-r}\,dW_r^{0,t_0})(x)\\
 &+\int_0^s\partial_zF_r^i(x,\widehat U_r^i(x))
                         (S_{h-r}L_i(r)\,dB_r^{t_0})(x).
 \end{aligned}
 \label{eq:explicit-first-martingale}
\end{equation}
These are the noise terms in the backward fields; the companion's
remaining control drift is treated separately in the expansion.
Write $N_{1,\lambda}^i=N_{1,\lambda}^i(h,\cdot)$ for the terminal
spatial field. The subscript $1$ means that no covariance has yet
been inserted, and $\lambda$ records the terminal test block.

\item \emph{Subtract the Gaussian-reference fluctuation.}
For the same terminal test, define the centered Doob martingale
\begin{equation}
 N_{1,\lambda}^{i,\mathrm{ref}}(s,x)
 :=\mathbb E[a_\lambda(X_h^i(x))\mid\mathcal F_{t_0+s}]
                  -\mathbb E_{t_0}[a_\lambda(X_h^i(x))].
 \label{eq:reference-first-doob-martingale}
\end{equation}
The complete correction process is
$\Gamma_{i,\lambda}(s,x):=
N_{1,\lambda}^i(s,x)-N_{1,\lambda}^{i,\mathrm{ref}}(s,x)$.
Its terminal spatial field, used in the comparison below, is
\[
 \Gamma_{i,\lambda}:=\Gamma_{i,\lambda}(h,\cdot)
 =N_{1,\lambda}^i-
       \bigl(a_\lambda(X_h^i)-\mathbb E_{t_0}[a_\lambda(X_h^i)]\bigr).
\]
Thus the quantity being compared is the actual first stochastic
fluctuation after its Gaussian-reference fluctuation has been removed.
\end{enumerate}
\end{definition}

\begin{lemma}[The complete correction is a martingale]
\label{lem:complete-correction-martingale}
For fixed $h,x$ and $i\in\{U,V\}$, the processes
$N_{1,\lambda}^i(s,x)$, $N_{1,\lambda}^{i,\mathrm{ref}}(s,x)$ and
$\Gamma_{i,\lambda}(s,x)$ of
Definition~\ref{def:coupled-backward-martingales} are square-integrable
martingales for $(\mathcal F_{t_0+s})_{0\le s\le h}$, starting at zero.
\end{lemma}

\begin{proof}[Proof of Lemma~\ref{lem:complete-correction-martingale}]
The integrands in \eqref{eq:explicit-first-martingale} are predictable.
By It\^o isometry, $|\partial_zF_r^i|\le\|a_\lambda'\|_\infty$
and bounded noise covariance give
$\mathbb E|N_{1,\lambda}^i(h,x)|^2
\le C\|a_\lambda'\|_\infty^2\int_0^h
\|p_{h-r}\|_{L^2(\mathbb T)}^2\,dr<\infty$.
The reference term is the Doob martingale of the bounded variable
$a_\lambda(X_h^i(x))$. Their difference is therefore a
square-integrable martingale, and both terms start at zero.
\end{proof}

\begin{definition}[Bandwise coefficient-replacement error]
\label{def:bandwise-coefficient-error}
Fix an expansion depth $N_\beta\ge1$. The approximation error
measures the effect of replacing $a$ by $a_{n_{\mathrm{sm}}}$
in the finite expansion, while keeping the physical processes
$(U,V^{n_{\mathrm{sm}},\Lambda})$ and their noise operators fixed.
It includes three kinds of correction:
\begin{itemize}[leftmargin=*]
\item replacing the terminal test $a_\lambda$ by
      $(a_{n_{\mathrm{sm}}})_\lambda$;
\item replacing one inserted coefficient $a$ by $a_{n_{\mathrm{sm}}}$;
\item replacing the Gaussian reference constructed from $a$ by
      the one constructed from $a_{n_{\mathrm{sm}}}$, at the same
      prescribed mean path $\phi^i$. The profile subtracted in each
      centered insertion and the covariance used in each Gaussian
      average are replaced by those of this recomputed reference.
\end{itemize}
Write $\operatorname{Err}_{n_{\mathrm{sm}},\lambda}$ for the sum
of the $\mathscr L_h^2(H_\Lambda^{-\rho})$ norms of these
single-replacement corrections. Its exact finite-sum definition is
\eqref{eq:bandwise-approximation-error} in
Appendix~\ref{app:coefficient-replacement-words}.
\end{definition}
At fixed $\lambda,\Lambda,h,N_\beta$ and cutoffs,
Lemma~\ref{lem:exact-bracket-cutoff-limits} gives
\[
 \operatorname{Err}_{n_{\mathrm{sm}},\lambda}\longrightarrow0
 \qquad\text{as}\qquad
 \epsilon_{n_{\mathrm{sm}}}
 :=\|a_{n_{\mathrm{sm}}}-a\|_\infty\longrightarrow0.
\]
The lemma also justifies removing the cutoff tails. The separate
amplitude-truncation tail $\lambda^{-p_{\mathrm{tail}}}$ is not
included in this error.

\begin{proposition}[Comparison of the first martingale corrections]
\label{prop:complete-martingale-comparison}
Let $U,V$ be the coupled solutions of
Section~\ref{sec:comparison-main-results}. Fix a dyadic
coefficient block $a_\lambda$ and its duration
$h=h_\lambda=\lambda^{-3\beta}\le h_0$.
For each window $[t_0,t_0+h]$ and $i\in\{U,V\}$, compare the
actual first stochastic fluctuation with its own Gaussian-reference
fluctuation by writing
$\Gamma_{i,\lambda}=N_{1,\lambda}^i-N_{1,\lambda}^{i,\mathrm{ref}}$,
as in Definition~\ref{def:coupled-backward-martingales}.
Here $N_{1,\lambda}^i$ is the terminal stochastic integral in
\eqref{eq:explicit-first-martingale}, and
$N_{1,\lambda}^{i,\mathrm{ref}}
=a_\lambda(X_h^i)-\mathbb E_{t_0}[a_\lambda(X_h^i)]$ is the
terminal Gaussian fluctuation from
\eqref{eq:reference-first-doob-martingale}.
Both use the same test block $a_\lambda$.

Set $\alpha=\beta+2\rho$, with $0<\rho<1/2$.
Measure the coefficient and solution differences by
$\mathsf d_{\mathrm{coef}}(U,V)$ and
$\|U-V\|_{\mathscr L_0^2(H)}$ from
\eqref{eq:coefficient-discrepancy-norm}.
Recall that $\|\cdot\|=\|\cdot\|_{\mathscr L_h^2(H_\Lambda^{-\rho})}$
is the $L^2$ norm over probability and window starts
$t_0\in[0,T-h]$, with spatial $H_\Lambda^{-\rho}$ norm
(Definition~\ref{def:window-integrated-norms}).
Let $\operatorname{Err}_{n_{\mathrm{sm}},\lambda}$ be the bandwise
coefficient-replacement error of
Definition~\ref{def:bandwise-coefficient-error}, explicitly
\eqref{eq:bandwise-approximation-error}, and write
$\mathcal L(\lambda,h)=1+\log(2+\lambda)+|\log h|$.
Then the difference of the two complete corrections satisfies
\begin{equation}
 \begin{aligned}
 \|\Gamma_{U,\lambda}-\Gamma_{V,\lambda}\|
 \le C\mathcal L(\lambda,h)^C\bigl[&
 \underbrace{\lambda^{-\beta}\mathsf d_{\mathrm{coef}}(U,V)}
              _{\mathrm{actual\ noise\ amplitudes}}\\
 &+\underbrace{\lambda^{-\alpha}\min(1,\lambda \|U-V\|_{\mathscr L_0^2(H)})}
              _{\mathrm{evaluation\ points}}\\
 &+\underbrace{\lambda^{-\alpha}h^{(\beta-1)/4}\|U-V\|_{\mathscr L_0^2(H)}}
              _{\mathrm{Gaussian\ reference\ profiles}}\bigr]
 +\operatorname{Err}_{n_{\mathrm{sm}},\lambda}.
 \end{aligned}
 \label{eq:current-cap-complete-Gamma}
\end{equation}
The constant is uniform in $\lambda$, $\Lambda$ and
$n_{\mathrm{sm}}$, with the other parameters fixed.
\end{proposition}

\paragraph{\textbf{Mechanism and reading route.}}
This is a first-order stability estimate with three sources of
variation. The representation
\eqref{eq:complete-martingale-representation} identifies the
actual and reference stochastic integrals.
Although the terminal test $a_\lambda$ is common, its backward
smoothings $F_s^U,F_s^V$ use different remaining reference variances.
The telescoping identity
\eqref{eq:complete-integrand-telescoping} changes their noise
amplitude, the argument of the smoothed test, and the variance
inside that test one at a time.
It\^o isometry then reduces each contribution to
a spatial trace and a time integral.

The detailed estimates preserve two features needed later:
the signed coefficient difference in $H_\Lambda^{-\rho}$,
and the additional frequency factor $\lambda^{-2\rho}$ in
the evaluation-point and reference terms.
The reference-profile variation also has an integrable
$s^{(\beta-1)/4}$ singularity at the window start, while the
remaining Gaussian variance vanishes at the terminal time.
Checking these spatial and temporal bounds accounts for most
of the calculation in
Section~\ref{sec:complete-martingale-bound-proof}.
For a first reading, follow the exact split
\eqref{eq:complete-integrand-telescoping} and the three resulting
bounds \eqref{eq:complete-coefficient-contribution},
\eqref{eq:current-cap-root-row} and
\eqref{eq:complete-reference-contribution}; their sum is the
asserted comparison. The intervening calculations verify these
three bounds using the estimates from Sections~3 and~4.

\begin{definition}[Higher stochastic terms for the coupled pair]
\label{def:higher-coupled-stochastic-terms}
The first martingale contains no covariance insertion and is
compared through its complete correction $\Gamma_{i,\lambda}$.
For $j\ge2$, put $n=j-1$: $N_{j,\lambda}^i$ contains $n$
insertions. We spell out its formula for both $i=U,V$.
Fix $x\in\mathbb T$ and an ordered history
$h=s_0>s_1>\cdots>s_n>0$, with $\xi_0=x$ and
$Z_\ell=((s_k,\xi_k))_{k=0}^\ell$.
The integration domain is
\begin{equation}
 \mathcal H_n(h)=\{(\mathbf s,\boldsymbol\xi)\in(0,h)^n\times\mathbb T^n:
                       h>s_1>\cdots>s_n>0\},\qquad
 d\mathbf s\,d\boldsymbol\xi=\prod_{\ell=1}^n ds_\ell\,d\xi_\ell.
 \label{eq:actual-history-domain}
\end{equation}
All these times are elapsed from $t_0$ and range continuously.
Write $a_U=a$, $a_V=a_{n_{\mathrm{sm}}}$, and use
$\mathcal P^{Z,U}=\mathcal P^{Z,1}$,
$\mathcal P^{Z,V}=\mathcal P^{Z,0}$ from
\eqref{eq:predictable-mean-covariance}.
The reference profiles remain $\bar a^U=\bar a^1$ and
$\bar a^V=\bar a^0$ from Definition~\ref{def:predictable-mean-functional}.
For this history, start with $F_0^{i,\lambda}=a_\lambda$ and, for $1\le\ell\le n$, set
\begin{equation}
 \begin{aligned}
 F_\ell^{i,\lambda}(z_0,\ldots,z_\ell)
 :={}&\frac12[a_i(z_\ell)-\bar a_{s_\ell}^i(\xi_\ell)]\\
 &\times\left[
 \left(\sum_{k=0}^{\ell-1}p_{s_k-s_\ell}(\xi_k-\xi_\ell)\partial_{z_k}\right)^2
 \mathcal P^{Z_{\ell-1},i}_{s_\ell,s_{\ell-1}}F_{\ell-1}^{i,\lambda}
 \right](z_0,\ldots,z_{\ell-1}).
 \end{aligned}
 \label{eq:actual-coupled-inserted-tests}
\end{equation}
Thus the root is $a_\lambda$ for both processes; the inserted
coefficient is $a$ for $U$ and $a_{n_{\mathrm{sm}}}$ for $V$.
Only the root carries the block label $\lambda$.

For $0\le r\le s_n$, the actual observation vector has coordinates
\begin{equation}
 \begin{aligned}
 (\mathbf Y_r^{U,Z_n})_k
   &=(S_{s_k-r}U_{t_0+r})(\xi_k),\\
 (\mathbf Y_r^{V,Z_n})_k
   &=(S_{s_k-r}V_{t_0+r})(\xi_k)
        +(e^{-\Lambda r}-e^{-\Lambda s_k})(S_{s_k}d)(\xi_k),
       \qquad 0\le k\le n,
 \end{aligned}
 \label{eq:actual-coupled-history-observations}
\end{equation}
where $d=U_{t_0}-V_{t_0}$, as in
Definition~\ref{def:coupled-backward-martingales}.
Let $\mathbf M^{i,Z_n}$ be its martingale part, starting at zero.
Its noise coordinates are explicitly
\begin{equation}
 d(\mathbf M_r^{i,Z_n})_k
 =m_0(S_{s_k-r}\,dW_r^{0,t_0})(\xi_k)
   +(S_{s_k-r}L_i(r)\,dB_r^{t_0})(\xi_k),\qquad 0\le k\le n,
 \label{eq:actual-coupled-history-noise}
\end{equation}
with $L_U(r)=L(a(U_{t_0+r}))$ and
$L_V(r)=L(a_{n_{\mathrm{sm}}}(V_{t_0+r}))$.
The higher stochastic word is
\begin{equation}
 \begin{aligned}
 N_{j,\lambda}^i(t_0;h,x)
 =\int_{\mathcal H_n(h)}\left[
 \int_0^{s_n}\nabla\bigl(\mathcal P^{Z_n,i}_{r,s_n}F_n^{i,\lambda}\bigr)
       (\mathbf Y_r^{i,Z_n})\cdot d\mathbf M_r^{i,Z_n}
 \right]d\mathbf s\,d\boldsymbol\xi,
 \qquad n=j-1.
 \end{aligned}
 \label{eq:actual-higher-martingale-formula}
\end{equation}
The gradient acts on all $n+1$ state coordinates.
For $n=0$, the outer history integral is absent and the formula
reduces to \eqref{eq:explicit-first-martingale}.
The control drift is separate from this stochastic integral.
As before, $N_{j,\lambda}^i$ denotes the resulting terminal spatial
field, with $t_0,h,x$ suppressed. Its additional insertions supply
the time gains used in the next proposition. Coefficient and
reference-profile replacements are accounted for by
$\operatorname{Err}_{n_{\mathrm{sm}},\lambda}$.

\end{definition}

For the chosen duration $h=h_\lambda=\lambda^{-3\beta}$,
$j\ge2$ and $0<\chi<\min(1,\beta+2\rho)$, define the frequency
decay exponents
\begin{equation}
 \begin{aligned}
 w_j&=\beta+\frac34(j-1)\beta^2,\\
 q_{j,\chi}&=w_j+\min\left(2\rho-\chi,\frac{3\beta}{4}(2\rho-\chi)\right),\\
 q_{j,\mathrm{prof}}&=w_j-\frac{3\beta}{4}(1-2\rho).
 \end{aligned}
 \label{eq:higher-martingale-degrees}
\end{equation}
They measure the noise, current-state and profile contributions,
respectively. Each additional insertion increases all three
exponents by $3\beta^2/4$.

\begin{proposition}[Higher-order martingale comparison]
\label{prop:higher-martingale-comparison}
Fix one dyadic block $a_\lambda$, an order $j\ge2$ and
$0<\chi<\min(1,\beta+2\rho)$, with the exponents in
\eqref{eq:higher-martingale-degrees}. For
$h=h_\lambda=\lambda^{-3\beta}\le h_{j-1}$, a sufficiently small
fixed loss $\varepsilon>0$ and any prescribed tail exponent
$p_{\mathrm{tail}}>0$, the stochastic words above satisfy
\begin{equation}
 \begin{aligned}
 \|N_{j,\lambda}^U-N_{j,\lambda}^V\|
 \le C_{j,\varepsilon,p_{\mathrm{tail}}}\big[&
     \lambda^{-w_j+\varepsilon}\mathsf d_{\mathrm{coef}}(U,V)
    +\lambda^{-q_{j,\chi}+\varepsilon}
                         \|U-V\|_{\mathscr L_0^2(H)}^{\chi-\varepsilon}\\
    &+\lambda^{-q_{j,\mathrm{prof}}+\varepsilon}
                         \|U-V\|_{\mathscr L_0^2(H)}^{1-\varepsilon}
    +\lambda^{-p_{\mathrm{tail}}}\big]
    +\operatorname{Err}_{n_{\mathrm{sm}},\lambda}.
 \end{aligned}
 \label{eq:higher-martingale-comparison}
\end{equation}
The constant may depend on the fixed order, coefficient bounds
and fixed exponents, but is independent of $\lambda$ and $\Lambda$.
\end{proposition}
The final power is a separately controlled truncation tail;
fixed logarithmic factors are included in the small loss $\varepsilon$.
Section~\ref{sec:higher-martingale-main-proof} checks the time
integral for each of the three contributions.

\subsubsection{The paired actual remainder}
\label{sec:paired-actual-covariance-remainder}
\begin{proposition}[Paired actual remainder]
\label{prop:paired-actual-remainder}
Fix a dyadic block $a_\lambda$, a depth $N_\beta\ge1$ and
$0<\chi<\beta<1$. The paired estimate below concerns the
common-coefficient remainders, both with root $a_\lambda$ and
inserted coefficient $a$. Explicitly, for $n=N_\beta$,
\begin{equation}
 \begin{aligned}
 R_{n,\lambda}^U(t_0;h,x)
 &=\int_{\mathcal H_n(h)}F_n^1(\mathbf Y_{s_n}^{U,Z_n})
                                      \,d\mathbf s\,d\boldsymbol\xi,\\
 R_{n,\lambda}^V(t_0;h,x)
 &=\int_{\mathcal H_n(h)}F_n^0(\mathbf Y_{s_n}^{V,Z_n})
                                      \,d\mathbf s\,d\boldsymbol\xi,
 \end{aligned}
 \label{eq:common-coefficient-remainders}
\end{equation}
where $F_n^1,F_n^0$ are the endpoint tests of
\eqref{eq:predictable-mean-test-recursion}, and the actual
observations are \eqref{eq:actual-coupled-history-observations}.
Replacing the $V$ factors by $a_{n_{\mathrm{sm}}}$ gives its actual
expansion remainder; the resulting corrections are included in
$\operatorname{Err}_{n_{\mathrm{sm}},\lambda}$.
Take $h=h_\lambda=\lambda^{-3\beta}\le h_{N_\beta}$,
where $h_{N_\beta}$ is the fixed time threshold in
\eqref{eq:finite-order-window-threshold}.
Their difference satisfies, for a sufficiently small fixed exponent
loss $\varepsilon>0$,
\begin{equation}
 \|R_{N_\beta,\lambda}^U-R_{N_\beta,\lambda}^V\|_{\mathscr L_h^2(H_\Lambda^{-\rho})}
 \le C_{N_\beta,\chi,\varepsilon}\Lambda^{-\rho/2}
 \lambda^{-\beta+\chi+\varepsilon}
 h^{\beta N_\beta/4-\varepsilon}\|U-V\|_{\mathscr L_0^2(H)}^{\chi-\varepsilon}.
 \label{eq:paired-actual-remainder}
\end{equation}
The main gain is $h^{\beta N_\beta/4}$, which improves with the
expansion depth. The exponent $\chi$ trades some frequency decay
for the factor $\|U-V\|_{\mathscr L_0^2(H)}^{\chi}$, measuring the difference
between the solutions; later we take $\chi=\beta/2$.
For each fixed $\beta$, we choose $N_\beta$ large but fixed,
independently of $\lambda$ and $\Lambda$, to make the remainder
negligible in the closure argument. The parameter $\varepsilon>0$
records arbitrarily small losses and is chosen last.

Its individual cap from \eqref{eq:capped-root-actual-strong-remainder}
remains available.

\end{proposition}

\begin{remark}[The standard windows and their frequency degrees]
\label{rem:standard-band-windows}
The duration $h_\lambda=\lambda^{-3\beta}$ shrinks as the band
frequency grows. Since $0<\beta<1$, it satisfies
$\lambda^{-4}\le h_\lambda$ for $\lambda\ge1$.
The reference Gaussian damping is
$e^{-c\lambda^2\sqrt{h_\lambda}}
=e^{-c\lambda^{2-3\beta/2}}$, which decays faster than every
inverse power of $\lambda$.
For each fixed order $n$,
\begin{equation}
 \frac{h_\lambda}{h_n}=\frac{1}{h_n\lambda^{3\beta}}\longrightarrow0,
 \qquad h_\lambda\le h_n\ \text{whenever}\
 \lambda\ge h_n^{-1/(3\beta)}.
 \label{eq:window-threshold-high-band}
\end{equation}
Once the depth $N_\beta$ is fixed, $h_{N_\beta}$ is a fixed positive
threshold. Thus all bands $\lambda\ge\lambda_*\asymp\Lambda^{1/4}$
have admissible durations for sufficiently large $\Lambda$.
The initial strip $0\le t<h_\lambda$ uses duration $h=t$;
Section~\ref{sec:initial-window-comparison-proof} proves the corresponding
true-age versions of the martingale and remainder bounds.

For $j=n+1$, the predictable factor
$\lambda^{-\beta}h_\lambda^{(j-1)\beta/4-1/4}$ is
$\lambda^{-p_j}$, with
\begin{equation}
 p_j=\beta+\frac34(j-1)\beta^2-\frac34\beta
     =w_j-\frac34\beta.
 \label{eq:all-order-predictable-clipped-degree}
\end{equation}
In particular, the lowest predictable degree is
$p_2=\beta-\gamma_{\mathrm{time}}(1-\beta)/4
=\beta/4+3\beta^2/4>0$ for $\gamma_{\mathrm{time}}=3\beta$.
This positivity is one requirement behind the choice of duration;
Lemma~\ref{lem:parameter-energy-closure} checks the full set of
parameter inequalities.
The general window exponents are stated with the proofs, where
the separate time factors are estimated before this substitution.
\end{remark}

\paragraph{\textbf{Roadmap for the comparison proofs.}}
The four main estimates share one requirement: a term with $n$
insertions must preserve its factor
$\lambda^{-\beta}h^{n\beta/4}$, up to the stated logarithms and
small losses. Gaussian differentiation produces singular time
factors; centering turns them into coefficient differences whose
H\"older bounds make the time integrals finite.

The comparison separates changes of the actual noise amplitude,
the current observations and the Gaussian reference profile.
The noise change is measured by the signed coefficient norm
$\mathsf d_{\mathrm{coef}}(U,V)$; the other changes are measured by
$\|U-V\|_{\mathscr L_0^2(H)}$, both defined in
\eqref{eq:coefficient-discrepancy-norm}.
Keeping the coefficient norm linear allows the comparison to close
with the damped energy estimate in Section~\ref{sec:all-beta-joint-complete-closure}.

\subsection{Working estimates for the comparison proofs}
\label{sec:centered-word-inputs}

The comparison proofs use the expansion terms of
Section~\ref{sec:comparison-main-results}, with one fixed block $a_\lambda$.
We first establish the time gains from centering, then the spatial
bounds for one changed input. Calculations begin with finite noise
series, smooth tests and strictly positive gaps. The resulting uniform
majorants permit the limits of Lemma~\ref{lem:exact-bracket-cutoff-limits},
which preserves the physical brackets. The first two orders are written
out in Appendix~\ref{sec:first-two-expansion-orders}.

Throughout the following comparison proofs, an unqualified output
norm is $\mathscr L_h^2(H_\Lambda^{-\rho})$ from Section~5.1.
The inputs $\|U-V\|_{\mathscr L_0^2(H)}$ and $\mathsf d_{\mathrm{coef}}(U,V)$ are the global
solution and coefficient distances in \eqref{eq:coefficient-discrepancy-norm};
$\mathcal L(\lambda,h)=1+\log(2+\lambda)+|\log h|$ records logarithms.
For $1\le p<\infty$, the conditional moment norm is
\[
 \|F\|_{L^p(\mathbb P_{t_0};L^2(\mathbb T))}
       =(\mathbb E_{t_0}\|F\|_{L^2(\mathbb T)}^p)^{1/p}.
\]
Here $\mathbb E_{t_0}=\mathbb E[\,\cdot\mid\mathcal F_{t_0}]$;
$\mathbb P_{t_0}$ in the norm records this conditioning.
The parameters $\varepsilon>0$ and $p_{\mathrm{tail}}>0$ denote a
small fixed loss and a prescribed tail exponent. For the higher
native current-mean estimate, $0<\chi<\min(1,\beta+2\rho)$.

\subsubsection{Turning centered differences into an integrable time kernel}
\label{sec:centered-word-calculus}
The comparison proofs below use one common calculation:
Gaussian integration by parts produces a mean-zero quadratic weight;
subtracting a value independent of its increment then gives a clipped
root factor or a positive H\"older gain. We record the assignments and
estimates here; Appendix~\ref{app:centered-word-kernel} verifies them.

Fix a block $a_\lambda$, an insertion count $n\ge1$, a window
$[t_0,t_0+h]$ with $0<h\le h_n$, and $0\le\theta\le1$.
Use the mean $\phi^\theta$ and profile $\bar a^\theta$ of
Definition~\ref{def:predictable-mean-functional}.
Fix a history $h=s_0>s_1>\cdots>s_n>0$,
with centers $\xi_0=x,\xi_1,\ldots,\xi_n$, and put
$\delta_j=s_{j-1}-s_j$ and $Z_j=((s_i,\xi_i))_{i=0}^j$.
The scalar state variable $w_\ell$ below is the coordinate at
$(s_\ell,\xi_\ell)$ in the insertion rule
\eqref{eq:predictable-mean-test-recursion}.

The starting identity is ordinary Gaussian integration by parts:
for $G\sim N(0,1)$ and a smooth bounded function $f$ with bounded
derivatives,
\[
 \mathbb E f''(G)=\mathbb E[(G^2-1)f(G)],\qquad
 \mathbb E(G^2-1)=0.
\]
Thus a second derivative of a Gaussian average can be expressed by
multiplying its integrand by a mean-zero quadratic weight. The
following is the precise weight for interval $(s_j,s_{j-1})$.

Choose $\kappa_*>0$ below half the common lower bound on
$m_0^2+\bar a^\theta$, and split this variance profile into the
constant part $\kappa_*$ and the positive remainder
$m_0^2+\bar a^\theta-\kappa_*$. Realize their contributions on
each interval as independent Gaussian vectors. For the fixed
history $Z_{j-1}=((s_i,\xi_i))_{0\le i<j}$, let
$\boldsymbol\eta_j$ be the centered observation increment from
the constant part on $(s_j,s_{j-1})$. Its covariance $C_j$, with
rows and columns indexed by $0,\ldots,j-1$, is explicitly
\[
 (C_j)_{ik}=\kappa_*\int_{s_j}^{s_{j-1}}\int_{\mathbb T}
 p_{s_i-r}(\xi_i-\zeta)p_{s_k-r}(\xi_k-\zeta)\,d\zeta\,dr.
\]
Both $\boldsymbol\eta_j$ and $C_j$ depend on this fixed history.
Define the insertion direction
$b_j=(p_{s_i-s_j}(\xi_i-\xi_j))_{0\le i<j}\in\mathbb R^j$.
Every vector in $\ker C_j$ annihilates the heat rows throughout
this interval and, by continuity, at $s_j$; hence
$b_j\in\operatorname{Ran}C_j$.
\begin{lemma}[Centering one interval increment]
\label{lem:interval-increment-centering}
Use $\boldsymbol\eta_j,C_j,b_j$ above, and let $C_j^+$ denote the
inverse on the covariance range, extended by zero on its null space.
The weight representing the half second derivative in direction $b_j$ is
\begin{equation}
 \mathsf S_j
 =\frac12\left[
   (b_j^*C_j^+\boldsymbol\eta_j)^2-b_j^*C_j^+b_j\right].
 \label{eq:interval-gaussian-derivative-weight}
\end{equation}
Let $\mathbb E_j$ average only $\boldsymbol\eta_j$, with all other
Gaussian increments fixed. For the coefficient product $H$ with root
$a_\lambda$ and factors $a(w_\ell)-\bar a_{s_\ell}^\theta(\xi_\ell)$,
let $H^{(j,0)}$ be its value with this interval increment suppressed
in every state argument and all subtracted values fixed. Then
\begin{equation}
 \mathbb E_j \mathsf S_j=0,
 \qquad
 \mathbb E_j[\mathsf S_jH]
        =\mathbb E_j[\mathsf S_j(H-H^{(j,0)})].
 \label{eq:procedure-score-centering}
\end{equation}
\end{lemma}
\begin{proof}[Proof of Lemma~\ref{lem:interval-increment-centering}]
Apply \eqref{eq:gaussian-score-definitions} to the constant Gaussian
part. The identity
\[
 \mathbb E_j(b_j^*C_j^+\boldsymbol\eta_j)^2=b_j^*C_j^+b_j
\]
gives $\mathbb E_j\mathsf S_j=0$.
Since $H^{(j,0)}$ is independent of $\boldsymbol\eta_j$,
$\mathbb E_j[\mathsf S_jH^{(j,0)}]=0$, proving the centering identity.
\end{proof}
The constant Gaussian part has a law fixed as $\theta$ varies;
the independent remainder is available when differentiating the
reference mean and covariance in the comparison proofs.

The product difference assigns the increment to individual factors.
Explicitly, if $H_\ell$ and $\widetilde H_\ell$ are their values
before and after suppressing the increment, then
\begin{equation}
 \prod_{\ell=0}^n H_\ell-\prod_{\ell=0}^n\widetilde H_\ell
 =\sum_{\ell=0}^n(H_\ell-\widetilde H_\ell)
       \prod_{k<\ell}H_k\prod_{k>\ell}\widetilde H_k.
 \label{eq:procedure-product-difference}
\end{equation}

\begin{definition}[Centered assignments and their trees]
\label{def:centered-word-assignments}
Fix a block $a_\lambda$, one coupled reference calculation and a
history $h=s_0>\cdots>s_n>0$, $n\ge1$, with centers
$\xi_0=x,\xi_1,\ldots,\xi_n\in\mathbb T$.
Put $(T_i,x_i)=(s_i,\xi_i)$ for $0\le i\le n$,
$\delta_j=s_{j-1}-s_j$ and
$\boldsymbol\delta=(\delta_1,\ldots,\delta_n)$.
\begin{enumerate}[label=(\roman*),leftmargin=*,itemsep=5pt]
\item \textbf{Label the coefficient factors.}
Writing $w_i$ for the scalar state at observation $(s_i,\xi_i)$,
the factors in \eqref{eq:predictable-mean-test-recursion} have labels
\[
 \begin{array}{ccl}
 0&:&a_\lambda(w_0)\quad\text{(the root)},\\[2pt]
 i=1,\ldots,n&:&a(w_i)-\bar a_{s_i}(\xi_i).
 \end{array}
\]
\item \begin{minipage}[t]{\linewidth}
\textbf{Assign each interval difference to one factor.}
A \emph{centered assignment} is a map
\[
 \boxed{\displaystyle
 \pi:\{1,\ldots,n\}\longrightarrow\{0,\ldots,n-1\},
 \qquad 0\le\pi(j)<j.}
\]
Its input is an interval label; its output is a coefficient-factor label:
\[
 \underbrace{j}_{\text{interval }(s_j,s_{j-1})}
 \ \xmapsto{\ \pi\ }\ 
 \underbrace{i=\pi(j)}_{\text{factor receiving the difference}}.
\]
In particular, $\pi(j)=0$ assigns the difference to the root
$a_\lambda$; $\pi(j)=i>0$ assigns it to inserted factor $i$.
\end{minipage}
\item \textbf{Draw the tree.}
The tree has vertices $0,\ldots,n$, root $0$, and the edge rule
\[
 \boxed{\pi(j)=i\quad\Longleftrightarrow\quad i\longrightarrow j.}
\]
Parent labels decrease when followed toward the root, so every
vertex is connected to $0$. The children of vertex $i$ are
$I_i:=\{j\in\{1,\ldots,n\}:\pi(j)=i\}$, $0\le i\le n$.
\end{enumerate}
\end{definition}

\noindent\textbf{One assignment, written as a map and drawn as a tree.}
\par\nobreak\smallskip
\noindent
\begin{minipage}[c]{0.43\linewidth}
\centering
\includegraphics[width=184pt]{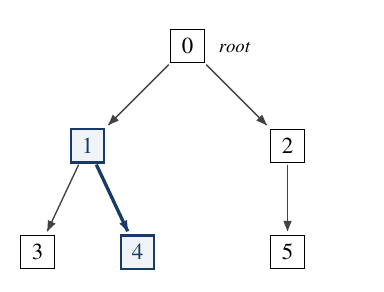}
\end{minipage}\hfill
\begin{minipage}[c]{0.54\linewidth}
For $n=5$, choose
\[
 \begin{array}{c|ccccc}
 j&1&2&3&4&5\\\hline
 \pi(j)&0&0&1&1&2.
 \end{array}
\]
The thicker edge illustrates
\[
 \pi(4)=1\quad\Longleftrightarrow\quad1\to4:
\]
the difference from interval $4$ acts on coefficient factor $1$.
\end{minipage}
\par\smallskip

\paragraph{Read off the index sets and counts.}
The nonempty child sets of nonroot vertices are
$\mathscr B_{\mathrm{rough}}:=\{I_i:1\le i<n,\ I_i\ne\varnothing\}$.
Together with $I_0$, they partition $\{1,\ldots,n\}$; child labels
need not be consecutive. The quantities used in the estimates are
\[
 \begin{aligned}
 n_{\mathrm r}&:=|I_0|
       &&\text{(root children)},\\
 n_{\mathrm b}&:=|\mathscr B_{\mathrm{rough}}|
       &&\text{(nonroot vertices with children)},\\
 u&:=n-1-n_{\mathrm b}
       &&\text{(older nonroot vertices without children)},\\
 \gamma_{\mathrm b}&:=n_{\mathrm b}\beta/4
       &&\text{(total time gain from nonroot child sets)}.
 \end{aligned}
\]
Vertex $n$ has no children and is counted separately from $u$.
In the example, $I_0=\{1,2\}$, $I_1=\{3,4\}$ and $I_2=\{5\}$,
so $n_{\mathrm r}=n_{\mathrm b}=u=2$.

\paragraph{Why centering produces these assignments.}
For interval $j$, Lemma~\ref{lem:interval-increment-centering} replaces
the coefficient product by its difference with that increment suppressed.
The product identity \eqref{eq:procedure-product-difference} assigns this
difference to one affected factor $i<j$. Recording the choice as
$\pi(j)=i$ and repeating over the intervals gives exactly the trees above.

\paragraph{How the tree determines the time weights.}
The root supplies $\lambda^{-\beta}$. The second derivative on
interval $j$ contributes the time factor $\delta_j^{-1}$ to its
edge $i\to j$. Centering supplies the
following additional factors, estimated in
Lemma~\ref{lem:centered-word-time-kernel}:
\begin{center}
\begin{tabular}{@{}ll@{}}
\toprule
Tree feature & Centering or moment factor\\
\midrule
Root edge $0\to j$ & $\min(1,\lambda\delta_j^{1/4})$\\[2pt]
Nonroot vertex with children $I_i\ne\varnothing$
 & $\displaystyle(\min_{j\in I_i}\delta_j)^{\beta/4}$\\[4pt]
Older nonroot vertex without children & $h^{\beta/4}$\\[2pt]
Newest vertex $n$ & $s_n^{\beta/4}$\\
\bottomrule
\end{tabular}
\end{center}
Thus each root edge gives
$k_\lambda(\delta_j)=\delta_j^{-1}\min(1,\lambda\delta_j^{1/4})$;
its integral is $L_\lambda(w)=\int_0^w k_\lambda(t)\,dt$, as in
\eqref{eq:clipped-root-time-weights}. Each nonroot vertex with
children shares one total gain $\beta/4$ among its child gaps.
This tree records where centered differences act. The heat-kernel
parent choices in Appendix~\ref{app:centered-word-kernel} are separate.
All empty products equal one.

\paragraph{Differences of the root factor.}
Write $\delta_z f(y)=f(y+z)-f(y)$ for a state translation difference. For every deterministic increment vector $\mathbf z=(z_1,\ldots,z_{n_{\mathrm r}})$, Bernstein gives
\begin{equation}
 |\delta_{z_1}\cdots\delta_{z_{n_{\mathrm r}} }a_\lambda(y)|
 \le C_{n_{\mathrm r}} \lambda^{-\beta}\prod_{i=1}^{n_{\mathrm r}}\min(1,\lambda|z_i|).
 \label{eq:capped-root-mixed-difference}
\end{equation}

Gaussian differentiation and integration of the spatial heat kernels
leave the factor $\delta_j^{-1}$ on each centered interval.
The detailed spatial bounds are proved in
Corollary~\ref{cor:interval-derivative-spatial-envelopes} in
Appendix~\ref{app:centered-word-kernel},
\eqref{eq:all-order-second-score-envelope}--\eqref{eq:all-order-first-score-envelope}.

\begin{lemma}[Common time kernel after centering]
\label{lem:centered-word-time-kernel}
Consider $\mathsf P_{n+1,\lambda}^\theta$ from
Definition~\ref{def:predictable-mean-functional}, with $n\ge1$
and $0<h\le h_n$ as in \eqref{eq:finite-order-window-threshold}.
Fix an assignment $\pi$ from
Definition~\ref{def:centered-word-assignments}, and put
\[
 s=s_n=h-\sum_{j=1}^n\delta_j,\qquad
 \alpha_j=\frac{\beta}{4|I_{\pi(j)}|}\quad(j\notin I_0).
\]
\begin{enumerate}[label=(\roman*),leftmargin=*,itemsep=4pt]
\item \emph{Before time integration.}
In \eqref{eq:predictable-mean-word-formula}, fix the gaps
$\delta_j>0$, $\sum_j\delta_j<h$. After centering, select the
term indexed by $\pi$ and integrate the insertion centers
$\xi_1,\ldots,\xi_n$. The resulting function of $x$ has
$L^2(\mathbb T)$ norm
at most
\begin{equation}
 C_n\lambda^{-\beta}h^{u\beta/4}
 \prod_{j\in I_0}k_\lambda(\delta_j)
 \prod_{j\notin I_0}\delta_j^{\alpha_j-1}s^{\beta/4}.
 \label{eq:all-order-predictable-common-ledger}
\end{equation}
\item \emph{After time integration.}
Integrating over the gaps and summing the assignments gives
\begin{equation}
 \|\mathsf P_{n+1,\lambda}^\theta(h,\cdot)\|_{L^2(\mathbb T)}
 \le C_n\mathcal L(\lambda,h)^{C_n}
                   \lambda^{-\beta}h^{n\beta/4}.
 \label{eq:centered-Gaussian-average-size}
\end{equation}
\end{enumerate}
\end{lemma}
The proof is given in Appendix~\ref{app:centered-word-kernel}.

\paragraph{One insertion: where the time gain appears.}
For $n=1$, the old interval can only be assigned to the root:
$n_{\mathrm r}=1$, $n_{\mathrm b}=u=0$. The fixed-gap bound
\eqref{eq:all-order-predictable-common-ledger} is then
$C\lambda^{-\beta}k_\lambda(h-s)s^{\beta/4}$.
The clipped factor comes from centering the old derivative weight; the positive
power comes from the newest coefficient discrepancy. In particular,
\[
 \|\mathsf P_{2,\lambda}^\theta\|_{L^2(\mathbb T)}
 \le C\lambda^{-\beta}\int_0^h k_\lambda(h-s)s^{\beta/4}\,ds
 \le C\lambda^{-\beta}h^{\beta/4}L_\lambda(h).
\]
At higher orders, each nonempty rough block contributes another
$\beta/4$ and root-assigned intervals contribute logarithms.
The lemma records this accounting once for all the comparison proofs.

\subsubsection{Spatial bounds for coefficient and solution differences}
\label{sec:capped-root-convolution}
\label{sec:capped-root-native-mean}
\label{sec:capped-root-native-noise}
\label{sec:direct-paired-native-martingale}
\label{sec:hybrid-native-actual-mean}
The next two lemmas keep one changed coefficient or solution state
in its spatial norm, and preserve the time gains of the other factors.
Work on $[t_0,t_0+h]$ with one block $a_\lambda$, order $j\ge2$
and $0<h\le h_{j-1}$. Write $D_t=U_t-V_t$ and $d=D_{t_0}$.
The global input and output norms are those recalled at the start
of Section~\ref{sec:centered-word-inputs}.

\begin{lemma}[Spatial estimates for one changed input]
\label{lem:comparison-spatial-inputs}
Under the preceding setup, the following estimates hold.
\begin{enumerate}[label=(\roman*),leftmargin=*,itemsep=5pt]
\item \emph{A coefficient difference stays in $L^2$.}
The residual-noise factors $L$ of
Definition~\ref{def:coupled-backward-martingales} satisfy, for $t>0$,
\begin{equation}
 \|S_t(L(a(U))-L(a(V)))\|_{\rm HS}
       \le Ct^{-1/4}\|a(U)-a(V)\|_{L^2(\mathbb T)}.
 \label{eq:native-integrator-heat-l2}
\end{equation}
\item \emph{A solution difference stays in $L^2$.}
For $s,t>0$, $s+t\le1$, the spatial trace is
\begin{equation}
 \|Q_\Lambda^{1/2}S_t\mathsf M_dS_s\|_{\rm HS}
       \le C_\rho(t+s)^{\rho/2-1/4}\|d\|_{L^2(\mathbb T)},
 \qquad Q_\Lambda=(\Lambda-\Delta)^{-\rho}.
 \label{eq:one-weight-heat-tree}
\end{equation}
For an insertion history $h=s_0>\cdots>s_{j-1}>q\ge0$, put
$\delta=s_{j-1}-q$. The same trace holds with $d$ replaced by
$D_{t_0+q}$ or a positive $L^2$ envelope. After the heat averaging in
\eqref{eq:multipoint-mean-energy}, the mean-derivative bound is at most
$C_j\delta^{-1/4}$ times the spatial $L^2$ norm of its input.
This also applies to the compensated initial-state direction
$(e^{-\Lambda T}-e^{-\Lambda q})S_Td$, for target times
$s_{j-1}\le T\le h$, uniformly in
$\Lambda$.
\end{enumerate}
\end{lemma}
The proof is given in Appendix~\ref{app:physical-input-details}.

\begin{lemma}[Interpolation preserves the coefficient time gain]
\label{lem:consistent-current-interpolation}
Use $\phi^\theta,\bar a^\theta$ from
Definition~\ref{def:predictable-mean-functional}, and set
$U_q^\theta=(1-\theta)V_{t_0+q}+\theta U_{t_0+q}$,
$0\le\theta\le1$, $0\le q\le h$.
\begin{enumerate}[label=(\roman*),leftmargin=*,itemsep=5pt]
\item \emph{The future mean.}
At any target time $q<T\le h$, the interpolated Gaussian continuation has mean
\begin{equation}
 S_{T-q}U_q^\theta+\phi_T^\theta-S_{T-q}\phi_q^\theta
       =\phi_T^\theta+S_{T-q}(U_q^\theta-\phi_q^\theta).
 \label{eq:hybrid-controlled-future-mean}
\end{equation}
\item \emph{The coefficient moment.}
For every $1\le p<\infty$ and $0<q\le h$,
\begin{equation}
 \sup_{\substack{0\le\theta\le1\\\xi\in\mathbb T}}
 \left(\mathbb E_{t_0}
   |a(U_q^\theta(\xi))-\bar a_q^\theta(\xi)|^p\right)^{1/p}
       \le C_pq^{\beta/4}.
 \label{eq:consistent-current-coefficient-moment}
\end{equation}
Adding the bounded future Gaussian increments up to $h$ gives
the same $h^{\beta/4}$ moment for the inserted coefficient factors.
The constants are uniform in $\Lambda$.
\end{enumerate}
\end{lemma}
\begin{proof}[Proof of Lemma~\ref{lem:consistent-current-interpolation}]
The mean identity is algebraic. Put
$\eta_q=D_{t_0+q}-e^{-\Lambda q}S_qd$. The damped mild equation gives
\[
 \eta_q=\int_0^q e^{-\Lambda(q-r)}S_{q-r}
                       (L_U(r)-L_V(r))\,dB_r^{t_0},
\]
with $L_U,L_V$ from Definition~\ref{def:coupled-backward-martingales}.
Bounded noise covariances, heat-kernel square integrability and
Burkholder--Davis--Gundy give conditional pointwise moment bounds
$C_pq^{1/4}$ for both $U_{t_0+q}-\phi_q^1$ and $\eta_q$.
Since $V_{t_0+q}-\phi_q^0=(U_{t_0+q}-\phi_q^1)-\eta_q$,
the bound holds uniformly for $U_q^\theta-\phi_q^\theta$.
H\"older continuity of $a$ and
\eqref{eq:self-profile-displacement} prove
\eqref{eq:consistent-current-coefficient-moment}.
Heat contraction and the future Gaussian variance bound give the
asserted $h^{\beta/4}$ moment at later targets.
\end{proof}

\subsection{Stability of the Gaussian average terms}
\label{sec:comparison-proof-details}
\label{sec:predictable-comparison-proof}
Fix one block $a_\lambda$, an insertion depth $n\ge1$ and a window
$[t_0,t_0+h]$ with $\lambda^{-4}\le h\le h_n$; the proposition takes
$h=h_\lambda=\lambda^{-3\beta}$.
Write $d=U_{t_0}-V_{t_0}$. For $0\le\theta\le1$ and $0\le t\le h$,
use $\phi_t^\theta=S_tU_{t_0}-(1-\theta)e^{-\Lambda t}S_td$ and
its self-consistent profile $\bar a^\theta$ from
Definition~\ref{def:predictable-mean-functional}.

\begin{proof}[Proof of Proposition~\ref{prop:all-order-predictable-pair}]
Condition on $\mathcal F_{t_0}$ and differentiate the interpolation of
Definition~\ref{def:predictable-mean-functional}, with $d\in L^2(\mathbb T)$ fixed.
We prove
\begin{equation}
 \boxed{\displaystyle
 \|\partial_\theta\mathsf P_{n+1,\lambda}^\theta\|_{L^2(\mathbb T)}
 \le C_n\mathcal L(\lambda,h)^{C_n}
       \lambda^{-\beta}h^{n\beta/4-1/4}\|d\|_{L^2(\mathbb T)}.}
 \label{eq:predictable-path-derivative-bound}
\end{equation}
Start with smooth coefficients and data; all bounds below are
uniform under the common coefficient assumptions.

\medskip
\noindent\emph{1. Preserve the centering identities along the interpolation.}\par
Use the constant Gaussian fraction $\kappa_*$ of
Lemma~\ref{lem:interval-increment-centering} for the old weights and
differences. Its law is independent of $\theta$, preserving
\eqref{eq:procedure-score-centering}. Conditional on this fraction,
\eqref{eq:three-leg-product} separates the mean, covariance and
subtracted-value derivatives, denoted by $\mathrm{mean}$,
$\mathrm{cov}$ and $\mathrm{sub}$. Put each new derivative on a
separate elliptic fraction. The heat-span bounds apply since
$h\le h_n$ and there are at most $n+1$ observations.

For an assignment of Definition~\ref{def:centered-word-assignments},
put $s=s_n$. Its $n-1$ older factors comprise $n_{\mathrm b}$
with assigned differences and $u$ without, so $u+n_{\mathrm b}=n-1$.
In \eqref{eq:all-order-predictable-common-ledger}, the unassigned
older factors are bounded by a constant times $h^{u\beta/4}$,
the sum of the positive gap exponents is $n_{\mathrm b}\beta/4$,
and the bound for the newest factor contains $s^{\beta/4}$.
Lemma~\ref{lem:insertion-gap-integration} bounds the corresponding
integral by $C_n\mathcal L(\lambda,h)^{C_n}$ times $h$ to the power
\[
 \frac{u\beta}{4}+\frac{n_{\mathrm b}\beta}{4}+\frac\beta4
 =\frac{n\beta}{4}.
\]
We now estimate the three differentiated contributions in the spatial $L^2$ norm.

\medskip
\noindent\emph{2. Bound the mean, covariance and subtraction terms.}\par
The profile estimates \eqref{eq:self-profile-displacement} and
\eqref{eq:self-profile-lp} give
\begin{equation}
 |\bar a_t^\theta-a(\phi_t^\theta)|\le Ct^{\beta/4},\qquad
 \|\partial_\theta \bar a_t^\theta\|_{L^2(\mathbb T)}
       \le Ct^{(\beta-1)/4}\|d\|_{L^2(\mathbb T)}.
 \label{eq:all-order-profile-path-bounds}
\end{equation}
H\"older continuity and Gaussian moments then bound each
$a(w)-\bar a_t^\theta(\xi)$ by $C_pt^{\beta/4}$ in every fixed
finite $L^p$ norm.

\emph{Mean contribution.}
The direction of the mean derivative at observation $i$ is
$e^{-\Lambda s_i}(S_{s_i}d)(\xi_i)$.
By \eqref{eq:multipoint-mean-energy} and
\eqref{eq:multipoint-score-products}, its Gaussian derivative
weight has heat-averaged norm at most
$C_ns^{-1/4}\|d\|_{L^2(\mathbb T)}$.
The newest coefficient factor remains bounded in $L^p$ by
$C_ps^{\beta/4}$, so the final time exponent is $(\beta-1)/4$.
Lemma~\ref{lem:insertion-gap-integration}, with
$\sigma=(\beta-1)/4>-1$, therefore gives the explicit spatial bound
\begin{equation}
 \begin{aligned}
 \big\|(\partial_\theta\mathsf P_{n+1,\lambda}^\theta)_{\mathrm{mean}}
                                                   \big\|_{L^2(\mathbb T)}
 &\le C_n\mathcal L(\lambda,h)^{C_n}\lambda^{-\beta}
       h^{(u+n_{\mathrm b})\beta/4+(\beta-1)/4}
                                      \|d\|_{L^2(\mathbb T)}\\
 &=C_n\mathcal L(\lambda,h)^{C_n}\lambda^{-\beta}
       h^{n\beta/4-1/4}\|d\|_{L^2(\mathbb T)}.
 \end{aligned}
 \label{eq:predictable-mean-contribution-bound}
\end{equation}

\emph{Derivative of the subtracted profile values.}
For a subtracted constant $b$ independent of the state variable,
$\delta_z(a-b)=\delta_z a$, where
$\delta_z f(y)=f(y+z)-f(y)$. Thus the derivative of a subtracted value is zero
on an assigned factor. At the newest factor,
\eqref{eq:all-order-profile-path-bounds} gives
$\|\partial_\theta\bar a_s^\theta\|_{L^2(\mathbb T)}
\le Cs^{(\beta-1)/4}\|d\|_{L^2(\mathbb T)}$;
the gap integral is bounded exactly as in
\eqref{eq:predictable-mean-contribution-bound}.
At an unassigned older factor, $u\ge1$ and $s_i\ge s$, so
\[
 \|\partial_\theta\bar a_{s_i}^\theta\|_{L^2(\mathbb T)}
 \le Cs_i^{(\beta-1)/4}\|d\|_{L^2(\mathbb T)}
 \le Cs^{(\beta-1)/4}\|d\|_{L^2(\mathbb T)}.
\]
There are now $u-1$ undifferentiated unassigned older factors.
Their bound contains $h^{(u-1)\beta/4}$; the newest factor
still has $L^p$ norm at most $C_ps^{\beta/4}$.
Apply the gap-integral estimate with $\sigma=(2\beta-1)/4>-1$.
The resulting power is
\[
 \frac{(u-1+n_{\mathrm b})\beta}{4}+\frac{2\beta-1}{4}
 =\frac{n\beta}{4}-\frac14.
\]
Thus \eqref{eq:predictable-mean-contribution-bound} also bounds the
subtracted-value contribution.

\emph{Covariance contribution.}
Apply \eqref{eq:multipoint-weighted-input-L2} with spatial input
$f_t=\partial_\theta\bar a_t^\theta$ and time exponent
$(\beta-1)/4$. The hypothesis on $f_t$ follows from
\eqref{eq:all-order-profile-path-bounds}.
Combining this estimate with
\eqref{eq:all-order-predictable-common-ledger} preserves all $n$
coefficient-factor bounds. The final time exponent is
$(2\beta-1)/4>-1$; the gap integration gives
\begin{equation}
 \begin{aligned}
 \big\|(\partial_\theta\mathsf P_{n+1,\lambda}^\theta)_{\mathrm{cov}}
                                                   \big\|_{L^2(\mathbb T)}
 &\le C_n\mathcal L(\lambda,h)^{C_n}\lambda^{-\beta}
       h^{(u+n_{\mathrm b})\beta/4+(2\beta-1)/4}
                                      \|d\|_{L^2(\mathbb T)}\\
 &=C_n\mathcal L(\lambda,h)^{C_n}\lambda^{-\beta}
       h^{n\beta/4-1/4+\beta/4}\|d\|_{L^2(\mathbb T)}.
 \end{aligned}
 \label{eq:predictable-covariance-contribution-bound}
\end{equation}
Since $h^{\beta/4}\le1$, this right-hand side is also bounded
by that in \eqref{eq:predictable-mean-contribution-bound}.
Summing the three contributions, with the finitely many assignments
and differentiated factors at order $n$, proves
\eqref{eq:predictable-path-derivative-bound} uniformly in $\theta$.

\medskip
\noindent\emph{3. Recover the endpoint difference, its cap and the window norm.}\par
Integrate \eqref{eq:predictable-path-derivative-bound} using
\begin{equation}
 \mathsf P_{n+1,\lambda}^1-\mathsf P_{n+1,\lambda}^0
 =\int_0^1\partial_\theta\mathsf P_{n+1,\lambda}^\theta\,d\theta.
\label{eq:procedure-predictable-interpolation}
\end{equation}
This proves \eqref{eq:all-order-predictable-clipped-linear-pair}.
Taking its minimum with the sum of the endpoint bounds
\eqref{eq:centered-Gaussian-average-size} gives
\eqref{eq:all-order-predictable-clipped-capped-pair}.

Finally, the spectral inequality
\[
 \|f\|_{H_\Lambda^{-\rho}}
 \le\Lambda^{-\rho/2}\|f\|_{L^2(\mathbb T)}
\]
and \eqref{eq:all-order-predictable-clipped-linear-pair} give the
window estimate \eqref{eq:all-order-predictable-integrated-sobolev}
after squaring, taking expectation and integrating over $t_0$:
\[
 \int_0^{T-h}\mathbb E\|U_{t_0}-V_{t_0}\|_{L^2(\mathbb T)}^2\,dt_0
 \le\|U-V\|_{\mathscr L_0^2(H)}^2.
\]
Lemma~\ref{lem:exact-bracket-cutoff-limits} justifies smooth
approximation and removal of the positive-gap restrictions.
The three asserted bounds follow for the original coefficient and data.
\end{proof}

\begin{samepage}
\subsection{Martingale comparisons}
\label{sec:martingale-comparisons}

\subsubsection{A stochastic-integral formula for the first correction}
\label{sec:complete-martingale-representation-proof}
Fix one block $a_\lambda$ and a terminal window length $h$.
For $i\in\{U,V\}$, $\Gamma_{i,\lambda}$ is the complete first
martingale correction of Definition~\ref{def:coupled-backward-martingales}.
This subsection writes it as a difference of stochastic integrals,
with every backward field and remaining variance specified below.
Its quantitative comparison is proved in
Section~\ref{sec:complete-martingale-bound-proof}.
There are two successive subtractions: first subtract each
solution's Gaussian-reference fluctuation to obtain
$\Gamma_{i,\lambda}$; then compare
$\Gamma_{U,\lambda}-\Gamma_{V,\lambda}$.
The first subtraction separates a change of noise amplitude from
a change of evaluation point, as displayed in
\eqref{eq:complete-martingale-representation}.

\end{samepage}
\begin{proof}[Proof of Proposition~\ref{prop:complete-martingale-comparison} (representation)]
\emph{1. Form the actual and Gaussian backward fields.}
Use the means, profiles, shifted noises and actual amplitudes
$L_U(s)=L(a(U_{t_0+s}))$, $L_V(s)=L(a_{n_{\mathrm{sm}}}(V_{t_0+s}))$
of Definition~\ref{def:coupled-backward-martingales}.
Write $U_s^U=U_{t_0+s}$ and $U_s^V=V_{t_0+s}$.
The identity \eqref{eq:covariance-factor-identity} fixes their residual
covariances; $m_0dW^0$ is their common additive channel.

Write $D_t=U_t-V_t$ and set $d=D_{t_0}$,
$\tau=h-s$, and $\eta_s=D_{t_0+s}-e^{-\Lambda s}S_sd$.
The reference drifts are $q_s^U=0$ and
$q_s^V=\Lambda e^{-\Lambda s}S_sd$. The actual and Gaussian backward fields are
\begin{equation}
 \widehat U_s^i=S_{h-s}U_s^i+\int_s^hS_{h-r}q_r^i\,dr,\qquad
 \mathcal G_s^i=S_{h-s}X_s^i+\int_s^hS_{h-r}q_r^i\,dr.
\label{eq:complete-backward-fields}
\end{equation}
The formula for $\widehat U^V$ agrees with
Definition~\ref{def:coupled-backward-martingales}, since
$\int_s^h S_{h-r}q_r^V\,dr=(e^{-\Lambda s}-e^{-\Lambda h})S_hd$.
The physical equations give
\begin{equation}
 d\widehat U_s^i
 =\mathbf1_{\{i=V\}}\Lambda S_{h-s}\eta_s\,ds
       +m_0S_{h-s}\,dW_s^{0,t_0}
       +S_{h-s}L_i(s)\,dB_s^{t_0}.
 \label{eq:complete-actual-backward-differential}
\end{equation}
The heat evolution cancels; for $V$, the correction subtracts
$\Lambda e^{-\Lambda s}S_hd\,ds$ from the control drift.
The reference mild equation similarly gives
\begin{equation}
 \mathcal G_s^i=\phi_h^i
       +m_0\int_0^sS_{h-r}\,dW_r^{0,t_0}
       +\int_0^sS_{h-r}L(\bar a_r^i)\,dB_r^{t_0},
 \qquad 0\le s\le h.
 \label{eq:gaussian-backward-field-martingale}
\end{equation}
Here $\mathcal G_0^i=\phi_h^i$ and $\mathcal G_h^i=X_h^i$.
With mean and profile fixed at $t_0$, this is a Gaussian martingale.
Define its amplitudes and the actual-minus-reference differences by
\begin{equation}
 L_i^0(s):=L(\bar a_s^i),\qquad
 \mathcal B_i(s):=L_i(s)-L_i^0(s),\qquad i\in\{U,V\}.
 \label{eq:complete-noise-factor-differences}
\end{equation}
The superscript $0$ denotes the reference operator, fixed by
the window-start data; time arguments are suppressed below.

\emph{2. Write the reference fluctuation as an It\^o integral.}
Use the remaining variance $\upsilon_s^i(x)$ from
\eqref{eq:complete-martingale-remaining-variance} and the smoothed block
$F_s^i(x,z)$ of Definition~\ref{def:coupled-backward-martingales}.
We suppress $x$ and write $(F_s^i)^{(r)}$ for its $r$th state
derivative, where $r\ge0$ is an integer.
Its reference evaluation is the conditional terminal test:
\[
 F_s^i(x,\mathcal G_s^i(x))
       =\mathbb E[a_\lambda(X_h^i(x))\mid\mathcal F_{t_0+s}].
\]
The scalar Gaussian backward equation and It\^o's formula hence
identify the Doob martingale in
\eqref{eq:reference-first-doob-martingale} with
$\int_0^s(F_r^i)'(\mathcal G_r^i)\,d\mathcal G_r^i$.

\emph{3. Separate the amplitude and evaluation-point differences.}
Subtract this reference integral from the actual stochastic term
\eqref{eq:explicit-first-martingale}. Add and subtract
$\int_0^h(F_s^i)'(\widehat U_s^i)\,d\mathcal G_s^i$.
The actual and reference noises share the additive channel;
their remaining difference is
$S_\tau\mathcal B_i(s)\,dB_s^{t_0}$.
This gives the complete correction as
\begin{equation}
 \Gamma_{i,\lambda}=\int_0^h(F_s^i)'(\widehat U_s^i)S_\tau \mathcal B_i\,dB_s^{t_0}
 +\int_0^h[(F_s^i)'(\widehat U_s^i)-(F_s^i)'(\mathcal G_s^i)]\,d\mathcal G_s^i.
 \label{eq:complete-martingale-representation}
\end{equation}
The integrals change the noise amplitude and evaluation point,
respectively. The $V$ control in
\eqref{eq:complete-actual-backward-differential} is a separate
finite-variation term in the expansion.
\end{proof}

\subsubsection{Estimating the first martingale correction}
\label{sec:complete-martingale-bound-proof}
The proof establishes \eqref{eq:current-cap-complete-Gamma} for any
$\lambda^{-4}\le h\le h_0$. Substitution of
$h=h_\lambda=\lambda^{-3\beta}$ gives the main statement.
Fix $a_\lambda$ and such a duration $h$.
Write $D_t=U_t-V_t$, $d=D_{t_0}$, $\tau=h-s$ and $\alpha=\beta+2\rho$.
Use the actual and Gaussian backward fields
$\widehat U_s^i,\mathcal G_s^i$ in
\eqref{eq:complete-backward-fields}--\eqref{eq:gaussian-backward-field-martingale},
and the reference-noise operators $L_i^0=L(\bar a_s^i)$ in
\eqref{eq:complete-noise-factor-differences}.
Use the smoothed block $F_s^i$ of
Definition~\ref{def:coupled-backward-martingales}, with remaining
variance \eqref{eq:complete-martingale-remaining-variance}; a prime
still denotes a state derivative.
Below $\delta$ means the $U$-minus-$V$ difference,
$Q_\Lambda=(\Lambda-\Delta)^{-\rho}$, and $\mathsf M_b$ is
multiplication by a spatial function $b$.
Output norms below are $\mathscr L_h^2(H_\Lambda^{-\rho})$,
over probability and window starts.

\begin{proof}[Proof of Proposition~\ref{prop:complete-martingale-comparison} (estimate)]
Use \eqref{eq:complete-martingale-representation} on the shared
noises $(W^{0,t_0},B^{t_0})$. Subtract the actual integrands by
\begin{equation}
 \begin{aligned}
 &(F_s^U)'(\widehat U_s^U)S_\tau L_U
                -(F_s^V)'(\widehat U_s^V)S_\tau L_V\\
 ={}&(F_s^U)'(\widehat U_s^U)S_\tau(L_U-L_V)
                       &&\text{(noise amplitude)}\\
 &+\big[(F_s^U)'(\widehat U_s^U)
                   -(F_s^U)'(\widehat U_s^V)\big]S_\tau L_V
                       &&\text{(evaluation point)}\\
 &+\big[(F_s^U)'(\widehat U_s^V)
                   -(F_s^V)'(\widehat U_s^V)\big]S_\tau L_V
                       &&\text{(remaining variance)}.
 \end{aligned}
 \label{eq:complete-integrand-telescoping}
\end{equation}
The scalar fields in this identity act by pointwise multiplication.
Apply the same identity to the reference integrands, replacing
$\widehat U_s^i$ by $\mathcal G_s^i$ and $L_i$ by $L_i^0$.
For the shared additive channel, only the evaluation-point and
remaining-variance rows occur.
It\^o isometry reduces the output norm to squared spatial traces
integrated in $s,t_0$; higher moments use BDG
\cite[Chapter~4]{DaPratoZabczyk2014}.

\medskip
\noindent\emph{1. The actual-noise change gives the coefficient-discrepancy term.}\par
Here $L_U-L_V=L(a(U_{t_0+s}))
-L(a_{n_{\mathrm{sm}}}(V_{t_0+s}))$.
Apply \eqref{eq:haar-negative-negative-trace} to the first derivative
of $L$ and integrate along the straight segment between these two
coefficient profiles. For any bounded multiplier $b$, this gives
\[
 \|Q_\Lambda^{1/2}\mathsf M_b S_\tau(L_U-L_V)\|_{\mathrm{HS}}^2
 \le C\|b\|_\infty^2\tau^{-1/2}
       \|a(U_{t_0+s})-a_{n_{\mathrm{sm}}}(V_{t_0+s})
                                      \|_{H_\Lambda^{-\rho}}^2.
\]
Use $b=(F_s^U)'(\widehat U_s^U)$.
The remaining variance is bounded below by $c\sqrt\tau$, so
annular Gaussian damping gives, for every integer $r\ge0$,
\begin{equation}
 \|(F_s^i)^{(r)}\|_\infty
 \le C_r\lambda^{r-\beta}e^{-c\lambda^2\sqrt\tau},
 \qquad i\in\{U,V\}.
 \label{eq:complete-smoothed-block-derivatives}
\end{equation}
The time integral needed here, and in the next two parts, is
\begin{equation}
 \begin{aligned}
 \int_0^h\tau^{\gamma-1}e^{-c\lambda^2\sqrt\tau}\,d\tau
 &=2\lambda^{-4\gamma}
       \int_0^{\lambda^2\sqrt h}u^{2\gamma-1}e^{-cu}\,du\\
 &\le C_\gamma\lambda^{-4\gamma},\qquad\gamma>0,
 \end{aligned}
 \label{eq:complete-annular-time-integral}
\end{equation}
by $u=\lambda^2\sqrt\tau$.
For the first derivative, its squared time factor is therefore
\[
 \lambda^{2-2\beta}\int_0^h\tau^{-1/2}
 e^{-c\lambda^2\sqrt\tau}\,d\tau\le C\lambda^{-2\beta}.
\]
After integrating over $t_0$, the coefficient input is bounded by
$\mathsf d_{\mathrm{coef}}(U,V)$ from
\eqref{eq:coefficient-discrepancy-norm}.
Consequently the first actual row of
\eqref{eq:complete-integrand-telescoping} has stochastic-integral
norm at most
\begin{equation}
 C\lambda^{-\beta}\mathsf d_{\mathrm{coef}}(U,V).
 \label{eq:complete-coefficient-contribution}
\end{equation}

\medskip
\noindent\emph{2. The evaluation-point changes give the clipped solution-difference term.}\par
We keep the smoothed test $F_s^U$ fixed and change its actual
and reference arguments. The two differences of arguments are
controlled by the global solution norm.
For the actual fields, the exact formula is
\begin{align}
 \delta \widehat U_s&=S_\tau D_{t_0+s}
       -e^{-\Lambda s}(1-e^{-\Lambda\tau})S_hD_{t_0},
 \label{eq:current-cap-exact-Z}\\
 \|\delta \widehat U_s\|_{\mathscr L_h^2(H)}&\le2\|U-V\|_{\mathscr L_0^2(H)}.
 \label{eq:current-cap-global-Z}
\end{align}
The second line follows from heat contraction and time translation.
For the references, define the difference of their centered
increments by
$G_s^0=(X_s^U-\phi_s^U)-(X_s^V-\phi_s^V)$. Then
\begin{align}
 \delta \mathcal G_s&=e^{-\Lambda h}S_hd+S_\tau G_s^0,
 \label{eq:actual-Gamma-delta-M}\\
 \|\delta \mathcal G_s\|_{L^p(\mathbb P_{t_0};L^2(\mathbb T))}&\le C_p\|d\|_{L^2(\mathbb T)}.
 \label{eq:actual-Gamma-Gaussian-moments}
\end{align}
Use \eqref{eq:self-profile-noise-stability}, Gaussian Hilbert
moments and heat contraction for the last line. The variance
integral is
$\int_0^s(s-r)^{-1/2}r^{(\beta-1)/2}\,dr=Cs^{\beta/2}$,
which is finite for $\beta>0$.
Integrating the conditional bound over $t_0$ gives
$\|\delta\mathcal G_s\|_{\mathscr L_h^2(H)}
\le C\|U-V\|_{\mathscr L_0^2(H)}$.

For $z,u\in\mathbb R$, the first two derivative bounds give
\[
 |(F_s^U)'(z+u)-(F_s^U)'(z)|
 \le C\lambda^{1-\beta}e^{-c\lambda^2\sqrt\tau}
                                      \min(1,\lambda|u|).
\]
For any bounded noise operator with $LL^*\le CI$, the output trace is
\begin{equation}
 \|Q_\Lambda^{1/2}\mathsf M_bS_\tau L\|_{\mathrm{HS}}^2
       \le C\tau^{\rho-1/2}\|b\|_{L^2(\mathbb T)}^2.
 \label{eq:complete-output-trace}
\end{equation}
This is the pathwise Bessel estimate
\eqref{eq:bessel-output-trace}, using
$2|b(x)b(y)|\le|b(x)|^2+|b(y)|^2$.
It permits dependent $b,L$ and includes the additive channel.

Apply \eqref{eq:complete-output-trace} first with
$b=(F_s^U)'(\widehat U_s^U)-(F_s^U)'(\widehat U_s^V)$,
and then with
$b=(F_s^U)'(\mathcal G_s^U)-(F_s^U)'(\mathcal G_s^V)$.
For each choice, the clipped estimate and the argument bounds
above give an integrated squared input at most
\[
 C\lambda^{2-2\beta}e^{-c\lambda^2\sqrt\tau}
       \min\{1,\lambda^2\|U-V\|_{\mathscr L_0^2(H)}^2\}.
\]
The time integration now uses
$\lambda^{2-2\beta}\int_0^h\tau^{\rho-1/2}
e^{-c\lambda^2\sqrt\tau}\,d\tau
\le C\lambda^{-2\beta-4\rho}$, by
\eqref{eq:complete-annular-time-integral}.
Taking square roots shows that the evaluation-point contributions
together have norm at most
\begin{equation}
 C\lambda^{-\alpha}\min(1,\lambda \|U-V\|_{\mathscr L_0^2(H)}),\qquad\alpha=\beta+2\rho.
 \label{eq:current-cap-root-row}
\end{equation}

\medskip
\noindent\emph{3. The reference changes give the linear reference-profile term.}\par
The remaining changes are $L_U^0-L_V^0$ and
$\delta\upsilon_s:=\upsilon_s^U-\upsilon_s^V$. Both depend on
$\delta\bar a_s:=\bar a_s^U-\bar a_s^V$.
Integrating \eqref{eq:all-order-profile-path-bounds} over the
interpolation parameter $\theta$ gives
\begin{equation}
 \|\delta\bar a_s\|_{L^2(\mathbb T)}
       \le Cs^{(\beta-1)/4}\|d\|_{L^2(\mathbb T)}.
 \label{eq:complete-profile-input}
\end{equation}

For the reference-noise amplitude, integrate
\eqref{eq:haar-one-l2-trace} along the two profiles and use
\eqref{eq:complete-smoothed-block-derivatives} with $r=1$.
The resulting squared trace is bounded by
\[
 C\lambda^{2-2\beta}e^{-c\lambda^2\sqrt\tau}
        \tau^{\rho-1/2}s^{(\beta-1)/2}\|d\|_{L^2(\mathbb T)}^2.
\]

For the two remaining-variance rows, hold the evaluation point
$z$ fixed. Differentiating Gaussian convolution in its variance gives
\begin{equation}
 (F_s^U)'(z)-(F_s^V)'(z)
 =\frac{\delta\upsilon_s}{2}\int_0^1
   \partial_z^3\mathcal P_{(1-\theta)\upsilon_s^V+\theta\upsilon_s^U}
                       a_\lambda(z)\,d\theta.
 \label{eq:complete-variance-derivative}
\end{equation}
At the same spatial point $x$,
\eqref{eq:complete-martingale-remaining-variance} and
\eqref{eq:complete-profile-input} give, for $0<s<h$,
\begin{equation}
 \begin{aligned}
 \|\delta\upsilon_s\|_{L^2(\mathbb T)}
 &\le C\|d\|_{L^2(\mathbb T)}
       \int_s^h(h-r)^{-1/2}r^{(\beta-1)/4}\,dr\\
 &\le C\sqrt\tau\,s^{(\beta-1)/4}\|d\|_{L^2(\mathbb T)}.
 \end{aligned}
 \label{eq:complete-remaining-variance-input}
\end{equation}
Each interpolated variance in
\eqref{eq:complete-variance-derivative} is at least $c\sqrt\tau$.
The third-derivative bound and
\eqref{eq:complete-output-trace} therefore bound each variance
row's squared trace by
\[
 C\lambda^{6-2\beta}e^{-c\lambda^2\sqrt\tau}
        \tau^{\rho+1/2}s^{(\beta-1)/2}\|d\|_{L^2(\mathbb T)}^2.
\]

The two time kernels satisfy
\begin{equation}
 \begin{aligned}
 \lambda^{2-2\beta}\int_0^h
   (h-s)^{\rho-1/2}e^{-c\lambda^2\sqrt{h-s}}
                          s^{(\beta-1)/2}\,ds
 &\le C\lambda^{-2\alpha}h^{(\beta-1)/2},\\
 \lambda^{6-2\beta}\int_0^h
   (h-s)^{\rho+1/2}e^{-c\lambda^2\sqrt{h-s}}
                          s^{(\beta-1)/2}\,ds
 &\le C\lambda^{-2\alpha}h^{(\beta-1)/2}.
 \end{aligned}
 \label{eq:complete-profile-time-integrals}
\end{equation}
To verify this, on $s\ge h/2$ bound
$s^{(\beta-1)/2}$ by $Ch^{(\beta-1)/2}$ and apply
\eqref{eq:complete-annular-time-integral}, with exponents
$\gamma=\rho+1/2$ and $\gamma=\rho+3/2$, respectively.
On $s<h/2$, use $h-s\asymp h$ and
$\int_0^{h/2}s^{(\beta-1)/2}\,ds\le Ch^{(\beta+1)/2}$.
After division by the asserted right-hand side, the two bounds
reduce to
$(\lambda^4h)^{\rho+1/2}e^{-c\sqrt{\lambda^4h}}\le C$
and
$(\lambda^4h)^{\rho+3/2}e^{-c\sqrt{\lambda^4h}}\le C$.
These hold uniformly for $\lambda^4h\ge1$.
Taking square roots and integrating the initial difference
$\|d\|_{L^2(\mathbb T)}^2$ over window starts gives a combined
reference contribution of norm at most
\begin{equation}
 C\lambda^{-\alpha}h^{(\beta-1)/4}
                              \|U-V\|_{\mathscr L_0^2(H)}.
 \label{eq:complete-reference-contribution}
\end{equation}

Adding \eqref{eq:complete-coefficient-contribution},
\eqref{eq:current-cap-root-row} and
\eqref{eq:complete-reference-contribution} proves
\eqref{eq:current-cap-complete-Gamma}.
Lemma~\ref{lem:exact-bracket-cutoff-limits} handles the additive
coefficient error and approximation limits.
\end{proof}

\subsubsection{Comparing the higher-order stochastic words}
\label{sec:higher-martingale-main-proof}
\label{sec:general-word-parameters}
Fix one coefficient block $a_\lambda$, an order $j\ge2$ and
$n=j-1$ insertions. Take $\lambda^{-4}\le h\le h_{j-1}$ and
$0<\chi<\min(1,\beta+2\rho)$. The target is
\eqref{eq:higher-martingale-comparison}; we first keep its bounds
in terms of the duration $h$. Section~\ref{sec:word-frequency-decay}
then substitutes $h=\lambda^{-\gamma_{\mathrm{time}}}$, including the main choice
$\gamma_{\mathrm{time}}=3\beta$, in all three bounds at once.
The words are the stochastic integrals
\eqref{eq:actual-higher-martingale-formula}, and the output norm is
$\mathscr L_h^2(H_\Lambda^{-\rho})$.

\begin{proof}[Proof of Proposition~\ref{prop:higher-martingale-comparison}]
\emph{1. Subtract the noise operators before taking norms.}
For $0\le q<h$ and $f\in H=L^2(\mathbb T)$, stochastic Fubini
in \eqref{eq:actual-higher-martingale-formula}, first at finite
cutoffs, gives the noise operator
\begin{equation}
 \begin{aligned}
 (\mathscr K_q^if)(x)=\int_{\mathcal H_n(h)\cap\{s_n>q\}}
 \sum_{k=0}^n
 &\partial_{z_k}\bigl(\mathcal P^{Z_n,i}_{q,s_n}F_n^{i,\lambda}\bigr)
                   (\mathbf Y_q^{i,Z_n})\\
 &\times(S_{s_k-q}f)(\xi_k)\,d\mathbf s\,d\boldsymbol\xi.
 \end{aligned}
 \label{eq:procedure-noise-kernel}
\end{equation}
The tests and actual observation vectors here are
\eqref{eq:actual-coupled-inserted-tests} and
\eqref{eq:actual-coupled-history-observations}.
Thus $N_{j,\lambda}^i=m_0\int_0^h\mathscr K_q^i\,dW_q^{0,t_0}
+\int_0^h\mathscr K_q^iL_i(q)\,dB_q^{t_0}$.
Use $L_U(q),L_V(q)$ and the shifted noises from
\eqref{eq:actual-coupled-history-noise}. Future Gaussian increments
are averaged out, so the kernel is predictable.

Replace only the inserted $V$ coefficient by $a$, and denote
the resulting integral and operator by $N_{j,\lambda}^{V;a}$ and
$\mathscr K^{V;a}$. The observations, $L_V$ and profile stay fixed;
the error is bounded by \eqref{eq:bandwise-approximation-error}.
Suppressing the time argument of $L_U,L_V$, subtract by
\begin{equation}
 \mathscr K_q^UL_U-\mathscr K_q^{V;a}L_V
   =\mathscr K_q^U(L_U-L_V)
                  +(\mathscr K_q^U-\mathscr K_q^{V;a})L_V.
 \label{eq:procedure-stochastic-telescoping}
\end{equation}
The additive channel has only the operator difference.
Stochastic Fubini gives
\begin{equation}
 \begin{aligned}
 N_{j,\lambda}^U-N_{j,\lambda}^{V;a}
 ={}&\underbrace{\int_0^h\mathscr K_q^U(L_U(q)-L_V(q))\,dB_q^{t_0}}
                         _{N_{j,\lambda}^{\mathrm{noise}}}\\
 &+m_0\int_0^h(\mathscr K_q^U-\mathscr K_q^{V;a})\,dW_q^{0,t_0}\\
 &+\int_0^h(\mathscr K_q^U-\mathscr K_q^{V;a})L_V(q)\,dB_q^{t_0}.
 \end{aligned}
 \label{eq:higher-martingale-integrand-split}
\end{equation}
In the last two integrals, separate observation and reference changes
by the interpolation $U_q^\theta=(1-\theta)V_{t_0+q}+\theta U_{t_0+q}$
and the consistent reference $\phi^\theta,\bar a^\theta$ of
Definition~\ref{def:predictable-mean-functional}. For a fixed history
$Z_n$, its interpolated current observation is explicitly
\begin{equation}
 (\mathbf Y_q^{\theta,Z_n})_k
 =(S_{s_k-q}U_q^\theta)(\xi_k)
   +\phi_{s_k}^\theta(\xi_k)
   -(S_{s_k-q}\phi_q^\theta)(\xi_k),\qquad 0\le k\le n.
 \label{eq:higher-current-observation-path}
\end{equation}
Its endpoints are \eqref{eq:actual-coupled-history-observations}.
With $D_t=U_t-V_t$ and $d=D_{t_0}$, its derivative is
\[
 (\partial_\theta\mathbf Y_q^{\theta,Z_n})_k
 =(S_{s_k-q}D_{t_0+q})(\xi_k)
       +(e^{-\Lambda s_k}-e^{-\Lambda q})(S_{s_k}d)(\xi_k).
\]
Lemma~\ref{lem:comparison-spatial-inputs}(ii) controls both
physical directions. Formula~\eqref{eq:three-leg-product} supplies
the covariance and subtraction changes, while
Lemma~\ref{lem:consistent-current-interpolation} preserves their
rough moments. Denote the current and profile contributions by
$\delta_{\mathrm{current}}N_{j,\lambda}$ and
$\delta_{\mathrm{profile}}N_{j,\lambda}$; step~3 separates their caps.

For the time integrations, let $q$ be the stochastic integration
time, $w=h-q$ the remaining time to the terminal observation, and
$\delta=s_n-q$ the gap to the earliest insertion. In a centered
assignment, $n_{\mathrm r}$ counts root-assigned intervals and
$n_{\mathrm b}$ counts rough coefficient factors receiving a nonempty block of
intervals, as in Definition~\ref{def:centered-word-assignments}.
There are $n-n_{\mathrm b}$ other rough factors, including the
newest one. Their moments supply $h^{(n-n_{\mathrm b})\beta/4}$.
Only finitely many assignments occur at this fixed order.

\emph{2. The changed noise: combine the Gaussian derivative and the signed trace.}
At the current physical time $t_0+q$, put
$f=a(U_{t_0+q})-a_{n_{\mathrm{sm}}}(V_{t_0+q})$ and
$B_\Delta=L_U(q)-L_V(q)$.
Here $L$ is the covariance factor in
\eqref{eq:covariance-factor-identity}. Its signed heat estimate
\eqref{eq:haar-heat-signed} gives
\begin{equation}
 \|S_tB_\Delta\|_{\mathrm{HS}}
 \le C\bigl[\Lambda^{\rho/2}t^{-1/4}+t^{-1/4-\rho/2}\bigr]
                                  \|f\|_{H_\Lambda^{-\rho}}.
 \label{eq:procedure-signed-input}
\end{equation}
The integrated input is $\mathsf d_{\mathrm{coef}}(U,V)$ from
\eqref{eq:coefficient-discrepancy-norm}.

The changed-noise term in \eqref{eq:higher-martingale-integrand-split}
contains $\mathscr K_q^U(L_U(q)-L_V(q))$. To see which smoothing
reaches this noise difference, first fix a history
$Z_n=((s_k,\xi_k))_{k=0}^n$, with $s_0=h$, $\xi_0=x$ and $q<s_n$.
Its contribution to \eqref{eq:procedure-noise-kernel} has the
explicit spatial kernel
\begin{equation}
 \begin{aligned}
 \mathscr K_q^{U,Z_n}(x,\zeta)
 :=\sum_{k=0}^n
 &\partial_{z_k}\bigl(\mathcal P^{Z_n,U}_{q,s_n}F_n^{U,\lambda}\bigr)
                    (\mathbf Y_q^{U,Z_n})\\
 &\times p_{s_k-q}(\xi_k-\zeta).
 \end{aligned}
 \label{eq:higher-fixed-history-noise-kernel}
\end{equation}
Here $\zeta$ is the noise-input location; history integration gives
$\mathscr K_q^U$.

Every heat age $s_k-q$ in this formula is at least
$\delta=s_n-q$. The semigroup identity
\[
 S_{s_k-q}=S_{s_k-q-\delta/2}S_{\delta/2},\qquad 0\le k\le n,
\]
factors $S_{\delta/2}$ from every row. Replacing
$p_{s_k-q}$ by $p_{s_k-q-\delta/2}$ defines
$\widetilde{\mathscr K}_q^{U,Z_n}$ and gives
$\mathscr K_q^{U,Z_n}=\widetilde{\mathscr K}_q^{U,Z_n}S_{\delta/2}$
before history integration.

Keep the coefficient and old-gap amplitudes outside the spatial
bound. Gaussian integration by parts on a reserved last-half-interval
fraction, followed by
Corollary~\ref{cor:interval-derivative-spatial-envelopes} and
positive heat integration, gives the following bound for a
contribution normalized to have random amplitude at most one:
\begin{equation}
 \mathscr K=\widetilde{\mathscr K} S_{\delta/2},\qquad
 |\widetilde{\mathscr K}|\le C_j\delta^{-1/2}w^{1/4}p_{C_jw},
 \quad w=h-q.
 \label{eq:native-integrator-exact-factor}
\end{equation}
The heat kernel is evaluated at $x-\zeta$ and its mass bounds
the operator norm. Lemma~\ref{lem:normalized-kernel-truncation}
handles unbounded amplitudes. Apply
\eqref{eq:haar-scale-separated-negative-input} with the right
$S_{\delta/2}$ next to $B_\Delta$, preserving the signed
$H_\Lambda^{-\rho}$ input. The time factors multiply to
\begin{equation}
 \underbrace{\delta^{-1/2}w^{1/4}}_{\text{last Gaussian derivative}}
 \underbrace{w^{\rho/2-1/4}\delta^{-\rho/2}
               \sqrt{1+\log(w/\delta)}}_{\text{signed spatial estimate}}
 =\delta^{-(1+\rho)/2}w^{\rho/2}
               \sqrt{1+\log(w/\delta)}.
 \label{eq:higher-noise-time-factor}
\end{equation}

At fixed $q$, integrate the old gaps by
\eqref{eq:insertion-gap-simplex-bound}, with final variable $\delta$,
exponent $-(1+\rho)/2>-1$ and rough degree $n_{\mathrm b}\beta/4$.
Including $w^{\rho/2}$ but omitting the logarithm gives
$Cw^{n_{\mathrm b}\beta/4-1/2}L_\lambda(w)^{n_{\mathrm r}}$.
For the logarithm choose $0<\eta<(1-\rho)/2$ and use
$\sqrt{1+\log(w/\delta)}\le C_\eta(w/\delta)^\eta$.
The last-gap exponent remains greater than $-1$; the extra
$w^\eta$ cancels the loss $w^{-\eta}$ in the same simplex integral.
Thus the same time power holds with the logarithm included.

For the stochastic integral use
$Q_\Lambda=(\Lambda-\Delta)^{-\rho}$ and Hilbert-space It\^o isometry:
\begin{equation}
 \mathbb E\left\|\int_0^h\mathscr K_qL_q\,dB_q\right\|_{H_\Lambda^{-\rho}}^2
 =\mathbb E\int_0^h
           \|Q_\Lambda^{1/2}\mathscr K_qL_q\|_{\mathrm{HS}}^2\,dq.
 \label{eq:procedure-final-isometry}
\end{equation}
Take $B_q=B_q^{t_0}$, $\mathscr K_q=\mathscr K_q^U$ and
$L_q=B_\Delta(q)$, and integrate over window starts.
If $n_{\mathrm b}>0$, \eqref{eq:positive-square-integral} integrates
$w^{n_{\mathrm b}\beta/2-1}L_\lambda(w)^{2n_{\mathrm r}}$.
If $n_{\mathrm b}=0$, all $n\ge1$ old intervals are assigned to the
root, and \eqref{eq:root-critical-square-integral} integrates
$w^{-1}L_\lambda(w)^{2n}$. Both cases give
\begin{equation}
 \|N_{j,\lambda}^{\rm noise}\|
       \le C_j\mathcal L(\lambda,h)^{C_j}\lambda^{-\beta}h^{n\beta/4}\mathsf d_{\mathrm{coef}}(U,V).
 \label{eq:capped-native-actual-noise-ideal}
\end{equation}
for bounded normalized amplitudes. The moment truncation in
Appendix~\ref{app:physical-input-details} removes this restriction:
for every prescribed small loss $\varepsilon>0$ and tail power
$p_{\mathrm{tail}}>0$,
\begin{equation}
 \|N_{j,\lambda}^{\rm noise}\|
 \le C_{j,\varepsilon,p_{\mathrm{tail}}}
       \lambda^{-\beta+\varepsilon}h^{n\beta/4}
             \mathsf d_{\mathrm{coef}}(U,V)
       +C_{j,\varepsilon,p_{\mathrm{tail}}}\lambda^{-p_{\mathrm{tail}}}.
 \label{eq:higher-noise-duration-bound}
\end{equation}
The right heat factor appears in both pieces of that truncation,
and the constants are uniform in $\Lambda$.

\emph{3. The current observations: preserve the root cap before squaring.}
The two physical directions in
\eqref{eq:higher-current-observation-path} have global spatial
$L^2$ size at most $C_j\|U-V\|_{\mathscr L_0^2(H)}$ and the mean-derivative
bounds of Lemma~\ref{lem:comparison-spatial-inputs}(ii).
The fractional mean comparison uses $\chi$ of this input and
has last-gap factor
$\delta^{-1/2}\delta^{-\chi/4}$: the first factor is the final
Gaussian derivative, and the second bounds the fractional mean derivative.
Its exponent $-1/2-\chi/4$ is greater than $-1$.

Here the cap is applied to the full consistent endpoint difference.
To keep the profile correction separate, let $A_{\mathrm{cur}}$
be the bound for the current-observation change before capping,
$A_{\mathrm{prof}}$ its separate profile bound, and
$C_{\mathrm{cap}}$ the full endpoint amplitude bound. The elementary
inequality
\[
 \min(A_{\mathrm{cur}}+A_{\mathrm{prof}},C_{\mathrm{cap}})
 \le\min(2A_{\mathrm{cur}},C_{\mathrm{cap}})+A_{\mathrm{prof}}
\]
preserves each endpoint's rough factor together with its own subtracted value.
This is the current contribution denoted by
$\delta_{\mathrm{current}}N_{j,\lambda}$; the remaining profile
bound is estimated in step~4. The positive spatial-input estimate
and moment interpolation are those of
Appendix~\ref{app:physical-input-details}.

Integrating the old gaps with the final exponent
$-1/2-\chi/4$, and using the output heat gain, gives the time kernel
\begin{equation}
 C\lambda^{-\beta}h^{(n-n_{\mathrm b})\beta/4}
 w^{n_{\mathrm b}\beta/4-1/2-\chi/4+\rho/2}L_\lambda(w)^{n_{\mathrm r}}.
 \label{eq:capped-native-mean-integrand}
\end{equation}
After squaring, its $w$-power is
$n_{\mathrm b}\beta/2+\rho-\chi/2-1$.
Apply \eqref{eq:current-mean-square-integral} with
$b=n_{\mathrm b}$ and $r=n_{\mathrm r}$. Their sum is positive
because there is at least one old interval. Including the
unassigned-factor power in
\eqref{eq:capped-native-mean-integrand} proves
\begin{equation}
 \|\delta_{\mathrm{current}}N_{j,\lambda}\|
 \le C\lambda^{-\beta+(\chi-2\rho)_++\varepsilon}
       h^{n\beta/4+(2\rho-\chi)_+/4}
                         \|U-V\|_{\mathscr L_0^2(H)}^{\chi-\varepsilon},
 \qquad x_+:=\max(x,0).
 \label{eq:higher-current-duration-bound}
\end{equation}

\emph{4. The reference profile: differentiate the subtracted values and covariance.}
An interval difference cancels the deterministic value subtracted in its factor.
Consequently a derivative of a subtracted value occurs only at a factor with no
assigned interval. It replaces one $h^{\beta/4}$ by
$q^{(\beta-1)/4}\|U-V\|_{\mathscr L_0^2(H)}$, using
\eqref{eq:all-order-profile-path-bounds} and the fact that the
insertion time is at least $q$.
The remaining old-gap integration and output trace have squared
power $(h-q)^{n_{\mathrm b}\beta/2+\rho-1}$.
The required time integral is explicitly
\[
 \int_0^h(h-q)^{n_{\mathrm b}\beta/2+\rho-1}
                       q^{(\beta-1)/2}\,dq
 =h^{n_{\mathrm b}\beta/2+\rho+(\beta-1)/2}
       B(n_{\mathrm b}\beta/2+\rho,(\beta+1)/2),
\]
with both beta-integral parameters positive. Root-assigned intervals add only
$L_\lambda(h)^{n_{\mathrm r}}$, by monotonicity, also when
$n_{\mathrm b}=0$. Taking the square root and restoring the
coefficient factors without children gives the duration factor
\[
 h^{(n-n_{\mathrm b}-1)\beta/4}
 h^{n_{\mathrm b}\beta/4+\rho/2+(\beta-1)/4}
       =h^{n\beta/4+\rho/2-1/4}.
\]

For a covariance derivative use
\eqref{eq:multipoint-weighted-input-L2} and the bounded total leverage
of Corollary~\ref{cor:bounded-profile-score}.
Its profile input at time $z\ge q$ obeys
$z^{(\beta-1)/4}\le q^{(\beta-1)/4}$ and preserves every rough factor.
It therefore has an additional $h^{\beta/4}$ compared with the
derivative of a subtracted value. The spatial input is kept in $L^2$ throughout.
The moment bounds in Appendix~\ref{app:physical-input-details} yield
\begin{equation}
 \|\delta_{\mathrm{profile}}N_{j,\lambda}\|
 \le C\lambda^{-\beta+\varepsilon}
       h^{n\beta/4+\rho/2-1/4}
                         \|U-V\|_{\mathscr L_0^2(H)}^{1-\varepsilon}.
 \label{eq:higher-profile-duration-bound}
\end{equation}
Conditioning companion moments at $t_0$ gives the same linear
profile input with the integrable $q$-weight above.

\emph{5. Add the contributions and restore the actual coefficient.}
Sum \eqref{eq:higher-noise-duration-bound},
\eqref{eq:higher-current-duration-bound} and
\eqref{eq:higher-profile-duration-bound} in
\eqref{eq:higher-martingale-integrand-split}, adding the separate
coefficient error. Distribute the prescribed loss among the
finitely many moment, power and truncation estimates. Since
$\lambda^{-4}\le h\le1$, the fixed powers of $\mathcal L(\lambda,h)$
are absorbed into this loss. Substituting $h=\lambda^{-3\beta}$
as in Section~\ref{sec:word-frequency-decay} gives exactly
\eqref{eq:higher-martingale-comparison}.
Lemma~\ref{lem:exact-bracket-cutoff-limits} justifies removal of the
finite cutoffs and supplies the stated coefficient-replacement error.
\end{proof}

\subsection{Actual remainders and initial time windows}
\label{sec:remainder-comparisons}

\subsubsection{Comparing the unexpanded actual remainder}
\label{sec:paired-remainder-proof}
The remainder is the part of the finite expansion still evaluated
at an actual observation. For a history with $n$ insertions,
integrate the one-step It\^o identity \eqref{eq:recursive-one-step}
over the old history. For the ordinary $U$ expansion this gives
\begin{equation}
 R_{n,\lambda}^U
   =\mathsf P_{n+1,\lambda}^U+N_{n+1,\lambda}^U+R_{n+1,\lambda}^U,
 \qquad R_{0,\lambda}^U=a_\lambda(U_{t_0+h}(x)).
 \label{eq:procedure-recursion}
\end{equation}
The three terms are the next Gaussian average, stochastic integral
and unexpanded remainder. For a controlled observation, the
first-derivative drift term
\eqref{eq:controlled-observation-drift-term} is also kept.

We establish \eqref{eq:paired-actual-remainder} for every
$\lambda^{-4}\le h\le h_{N_\beta}$; its main-statement form follows
by taking $h=\lambda^{-3\beta}$.
Fix a depth $N_\beta\ge1$, one block $a_\lambda$,
$0<\chi<\beta$ and $\lambda^{-4}\le h\le h_{N_\beta}$.
The terms $R_{N_\beta,\lambda}^U,R_{N_\beta,\lambda}^V$ are the
remainders with $N_\beta$ insertions, given by the history integral
in Proposition~\ref{prop:paired-actual-remainder}; their norm is
$\mathscr L_h^2(H_\Lambda^{-\rho})$. Write $D_t=U_t-V_t$.
Put $n=N_\beta$ in the following comparison.
For the common-coefficient remainders explicitly defined in
\eqref{eq:common-coefficient-remainders}, form
$R_{n,\lambda}^{V\mid U}$ from $R_{n,\lambda}^V$ by replacing
every Gaussian convolution $\mathcal P^{Z,0}$ in
\eqref{eq:predictable-mean-test-recursion} by $\mathcal P^{Z,1}$.
Keep the inserted coefficient $a$, the subtracted $V$ profile
$\bar a^0$ and the actual vector $\mathbf Y^{V,Z}$ fixed.
Thus only the future covariance changes. The exact split is
\begin{equation}
 R_{n,\lambda}^U-R_{n,\lambda}^V
  =\bigl(R_{n,\lambda}^U-R_{n,\lambda}^{V\mid U}\bigr)
    +\bigl(R_{n,\lambda}^{V\mid U}-R_{n,\lambda}^V\bigr).
 \label{eq:procedure-remainder-split}
\end{equation}
The first difference pairs current values and their subtracted values at
fixed covariance; the second changes only that covariance at the
fixed $V$ observation vector and subtracted value.
We compare one factor at a time while preserving the time gains
of every unchanged factor.

\begin{proof}[Proof of Proposition~\ref{prop:paired-actual-remainder}]
Use the split \eqref{eq:procedure-remainder-split}, telescoping
one coefficient factor at a time. This keeps exactly one changed
factor in each summand and preserves the centering of the others.

\emph{Compare the current values before the covariance.}
Fix a history $h=s_0>\cdots>s_{N_\beta}=s$, write
$\delta_i=s_{i-1}-s_i$ for its gaps, and use the centered assignments
of Definition~\ref{def:centered-word-assignments}. Condition on the
sigma field at $t_0+s$. The $N_\beta$ centered interval derivative weights and $N_\beta-1$
old rough factors are inside the auxiliary future expectation; the
newest $a(U_{t_0+s})-\bar a_s^U$ is outside and fixed there. Its
$s^{\beta/4}$ moment is conditional at the old past $t_0$, applied
after the future expectation. The corrected
controlled-mean identity \eqref{eq:hybrid-controlled-future-mean}
gives the same moment for $V$.

During the covariance change, bounded elliptic future increments
remain of size $h^{1/4}$ about the fixed $V$ mean, so its own subtracted value
preserves the required rough moments.
Physical current directions are finite sums of positive heat averages
of $D_{t_0+s}$ and damped $D_{t_0}$, including the compensated remaining
bridge. Their global $L^2$ norms are at most $C_N\|U-V\|_{\mathscr L_0^2(H)}$, uniformly in
times and $\Lambda$, by time translation and heat contraction.
The profile direction has global bound $C\|U-V\|_{\mathscr L_0^2(H)}^\beta$: use $L^2$ heat
contraction, finite-measure norm comparison and the profile Volterra
inequality.

\emph{Pair one centered factor.} Reserve a constant elliptic fraction for every old interval derivative weight, below
all profiles on the comparison path, and keep its inverse covariance matrices
fixed. At common future covariance, expand the same finite differences
at both endpoints and telescope, with exactly one paired slot. Write
$\delta_zf(y)=f(y+z)-f(y)$ for a state translation difference.
For $r_{\mathrm{inc}}\ge0$ old translation increments
$\mathbf z=(z_1,\ldots,z_{r_{\mathrm{inc}}})$ and one paired increment
$e$, a root slot obeys
\begin{equation}
 |\delta_e\delta_{z_1}\cdots\delta_{z_{r_{\mathrm{inc}}}}a_\lambda(y)|
 \le C_{r_{\mathrm{inc}}}\lambda^{-\beta}\min(1,\lambda|e|)^\chi
                   \prod_{i=1}^{r_{\mathrm{inc}}}\min(1,\lambda|z_i|).
 \label{eq:paired-remainder-root-slot}
\end{equation}
It follows from \eqref{eq:capped-root-mixed-difference} by weakening
the $e$ minimum to exponent $\chi$. The positive one-input estimate and
moment interpolation multiply the ordinary bound by at most
$C\lambda^\chi \|U-V\|_{\mathscr L_0^2(H)}^{\chi-\varepsilon}$
and preserve every old clip. A paired rough block $I$ instead uses
\begin{equation}
 \min\{\min_{i\in I}\delta_i^{\beta/4},|e|^\beta\}
 \le(\min_{i\in I}\delta_i)^{(\beta-\chi)/4}|e|^\chi.
 \label{eq:paired-remainder-rough-block}
\end{equation}
Only that block loses $\chi/4$ degree and keeps positive degree
$(\beta-\chi)/4$. An unassigned paired factor has
$h^{(\beta-\chi)/4}\|U-V\|_{\mathscr L_0^2(H)}^{\chi-\varepsilon}$; the newest actual paired
factor has $s^{(\beta-\chi)/4}\|U-V\|_{\mathscr L_0^2(H)}^{\chi-\varepsilon}$. Apply these
caps to the whole centered endpoint factor, including its own subtracted value.
All other factors keep endpoint rough moments. H\"older is applied
only after the positive input estimate, with losses smaller than the
fixed margins $\chi,\beta-\chi$.

\emph{Change the reference covariance.} At this step keep the
current $V$ observation vector and subtracted value fixed, and interpolate only the future
covariance through the explicit profile
$(1-\theta)\bar a^V+\theta\bar a^U$, $0\le\theta\le1$.
Keep old derivative weight fractions fixed and put the covariance derivative
on a separate reserved Gaussian fraction.
Only the total covariance variation matters, so these local
constructions patch along the uniformly elliptic path. The weighted
total-mass bound of Corollary~\ref{cor:bounded-profile-score}, together
with the spatial input estimate
\eqref{eq:multipoint-weighted-input-L2}, gives one $L^2$ mark
$\bar a^U-\bar a^V$ of bounded total mass and temporal weight zero. Old centerings and clips
survive; hence this leg is bounded by
$C_{N_\beta}\lambda^{-\beta}h^{\beta N_\beta/4}\mathcal L(\lambda,h)^C\|U-V\|_{\mathscr L_0^2(H)}^{\beta-\varepsilon}$,
which is no larger than the asserted $\|U-V\|_{\mathscr L_0^2(H)}^{\chi-\varepsilon}$ row on
the bounded range of $\|U-V\|_{\mathscr L_0^2(H)}$ given by preliminary damping.

\emph{Integrate the time gains.} A root slot multiplies the ordinary
cap by $\lambda^\chi$. A paired rough block, unassigned factor, or
newest factor loses only $h^{\chi/4}$, and
$h^{-\chi/4}\le\lambda^\chi$. For a paired rough block $I$,
distribute its remaining positive degree as
$(\beta-\chi)/(4|I|)$ per assigned gap; all other gap exponents
stay unchanged. Lemma~\ref{lem:insertion-gap-integration} applies
also when all old intervals go to the root, since the newest
factor has exponent at least $(\beta-\chi)/4>-1$.
It yields the ordinary $h^{N_\beta\beta/4}$ cap with at most
the stated factor $\lambda^\chi$. The positive majorants remove
the time cutoffs and prove the paired remainder bound.
\end{proof}

\subsubsection{Extending the comparisons to initial time windows}
\label{sec:initial-window-comparison-proof}
For a block with nominal duration $h_\lambda$, the window ending
at physical time $t$ is $[\max(0,t-h_\lambda),t]$.
When $0<t<h_\lambda$, this is $[0,t]$: its duration is $t$.
The preceding proofs assumed a fixed duration at least
$\lambda^{-4}$, whereas $t$ can be arbitrarily small.
The next result handles this endpoint by integrating the time
kernels over $0<s<t<h_\lambda$. Equal initial data also simplify
the comparison: the two Gaussian references constructed from $a$
coincide.

\begin{corollary}[Comparison bounds on initial windows]
\label{cor:initial-window-comparisons}
Use the coupled pair of Section~\ref{sec:comparison-main-results},
with $U_0=V_0=u_0$, $0<\beta<1$ and $0<\rho<1/2$.
Fix one coefficient block $a_\lambda$ and a depth $N_\beta\ge1$.
Take $h_\lambda=\lambda^{-\gamma_{\mathrm{time}}}\le\min(T,h_{N_\beta})$, $0<\gamma_{\mathrm{time}}<4$,
where $h_{N_\beta}$ is defined in
\eqref{eq:finite-order-window-threshold}.
For $0<t<h_\lambda$, let
$\Gamma_{i,\lambda}^{\mathrm{init}}(t)$,
$N_{j,\lambda}^{i,\mathrm{init}}(t)$ and
$R_{N_\beta,\lambda}^{i,\mathrm{init}}(t)$, $i\in\{U,V\}$,
denote the terminal spatial fields from
Definition~\ref{def:coupled-backward-martingales},
\eqref{eq:actual-higher-martingale-formula} and
\eqref{eq:common-coefficient-remainders} on $[0,t]$, respectively.
Thus $t_0=0$ and $h=t$. Their terminal-time norm is
\begin{equation}
 \|A\|_{\mathrm{init},\lambda}
 :=\left(\int_0^{h_\lambda}\mathbb E
              \|A(t)\|_{H_\Lambda^{-\rho}}^2\,dt\right)^{1/2},
 \qquad
 \|f\|_{H_\Lambda^{-\rho}}
 =\|(\Lambda-\Delta)^{-\rho/2}f\|_{L^2(\mathbb T)}.
 \label{eq:initial-window-output-norm}
\end{equation}
The input norms are the full-horizon quantities
$\|U-V\|_{\mathscr L_0^2(H)}
=(\int_0^T\mathbb E\|U_s-V_s\|_{L^2(\mathbb T)}^2ds)^{1/2}$
and $\mathsf d_{\mathrm{coef}}(U,V)$ from
\eqref{eq:coefficient-discrepancy-norm}.
Let $\operatorname{Err}_{n_{\mathrm{sm}},\lambda}^{\mathrm{init}}$
be the finite replacement-error sum
\eqref{eq:bandwise-approximation-error}, with each correction
evaluated on $[0,t]$ and each output norm replaced by
\eqref{eq:initial-window-output-norm}.

The following bounds hold. Recall
$\mathcal L(\lambda,h)=1+\log(2+\lambda)+|\log h|$.
\begin{enumerate}[label=(\roman*),leftmargin=*,itemsep=5pt]
\item \emph{The first complete correction satisfies}
\begin{equation}
 \begin{aligned}
 \|\Gamma_{U,\lambda}^{\mathrm{init}}
       -\Gamma_{V,\lambda}^{\mathrm{init}}\|_{\mathrm{init},\lambda}
 \le C\mathcal L(\lambda,h_\lambda)^C\bigl[&
        \lambda^{-\beta}\mathsf d_{\mathrm{coef}}(U,V)\\
 &+\lambda^{-\beta-2\rho}
        \min(1,\lambda\|U-V\|_{\mathscr L_0^2(H)})\bigr].
 \end{aligned}
 \label{eq:initial-first-martingale-comparison}
\end{equation}
\item \emph{Higher stochastic terms.} For $2\le j\le N_\beta$,
put $n=j-1$ and choose $0<\chi<\min(1,\beta+2\rho)$.
For sufficiently small fixed $\varepsilon>0$ and any
$p_{\mathrm{tail}}>0$, with $x_+=\max(x,0)$,
\begin{equation}
 \begin{aligned}
 \|N_{j,\lambda}^{U,\mathrm{init}}
       -N_{j,\lambda}^{V,\mathrm{init}}\|_{\mathrm{init},\lambda}
 \le C\bigl[&\lambda^{-\beta+\varepsilon}h_\lambda^{n\beta/4}
                         \mathsf d_{\mathrm{coef}}(U,V)\\
 &+\lambda^{-\beta+(\chi-2\rho)_++\varepsilon}
       h_\lambda^{n\beta/4+(2\rho-\chi)_+/4}
                         \|U-V\|_{\mathscr L_0^2(H)}^{\chi-\varepsilon}\\
 &+\lambda^{-p_{\mathrm{tail}}}\bigr]
       +\operatorname{Err}_{n_{\mathrm{sm}},\lambda}^{\mathrm{init}}.
 \end{aligned}
 \label{eq:initial-higher-martingale-comparison}
\end{equation}
\item \emph{Common-coefficient remainders.} Both remainders have
inserted coefficient $a$. For $0<\chi_{\mathrm{rem}}<\beta$ and
sufficiently small fixed $\varepsilon>0$,
\begin{equation}
 \begin{aligned}
 \|R_{N_\beta,\lambda}^{U,\mathrm{init}}
       -R_{N_\beta,\lambda}^{V,\mathrm{init}}\|_{\mathrm{init},\lambda}
 \le{}&C\Lambda^{-\rho/2}
       \lambda^{-\beta+\chi_{\mathrm{rem}}+\varepsilon}
       h_\lambda^{N_\beta\beta/4-\varepsilon}\\
 &\times\|U-V\|_{\mathscr L_0^2(H)}^{\chi_{\mathrm{rem}}-\varepsilon}.
 \end{aligned}
 \label{eq:initial-paired-remainder-comparison}
\end{equation}
Using the inserted coefficient $a_{n_{\mathrm{sm}}}$ in the $V$
remainder adds at most
$\operatorname{Err}_{n_{\mathrm{sm}},\lambda}^{\mathrm{init}}$.
\end{enumerate}
The constants are uniform in $\lambda$, $\Lambda$ and $n_{\mathrm{sm}}$.
The amplitude-truncation tail in (ii) is separate from coefficient
replacement. Taking $\gamma_{\mathrm{time}}=3\beta$ gives the main statements' frequency
powers, with the reference-profile terms zero.
\end{corollary}

\begin{proof}[Proof of Corollary~\ref{cor:initial-window-comparisons}]
\emph{1. Identify what remains when the initial data agree.}
Write $D_s=U_s-V_s$. Since $D_0=0$, the interpolated means in
Definition~\ref{def:predictable-mean-functional} are
$\phi_s^\theta=S_su_0$ for every $\theta$.
Uniqueness of the profile gives $\bar a^U=\bar a^V$;
the references and their Doob martingales therefore coincide.
The current-observation direction in
\eqref{eq:higher-current-observation-path} reduces to
$(S_{s_k-q}D_q)(\xi_k)$, and the compensated difference is
$\eta_q=D_q$. Only the actual-noise and current-observation changes remain;
in \eqref{eq:procedure-remainder-split}, $R^{V\mid U}=R^V$.

\emph{2. Integrate over the terminal time.}
We use the positive kernel bounds in the preceding proofs before
their final duration integration. For a nonnegative integrable
kernel $\omega$ and a nonnegative input $f$, Tonelli gives
\begin{equation}
 \begin{aligned}
 \int_0^{h_\lambda}\int_0^t\omega(t-s)f(s)\,ds\,dt
 &=\int_0^{h_\lambda}f(s)
             \int_0^{h_\lambda-s}\omega(r)\,dr\,ds\\
 &\le\|\omega\|_{L^1(0,h_\lambda)}
                  \int_0^{h_\lambda}f(s)\,ds.
 \end{aligned}
 \label{eq:initial-triangle-kernel-bound}
\end{equation}
After It\^o isometry, $\omega$ is the squared stochastic kernel and
$f(s)$ the mean squared spatial input. For ordinary integrals,
weighted Cauchy--Schwarz and \eqref{eq:initial-triangle-kernel-bound} give
\begin{equation}
 \int_0^{h_\lambda}
       \left(\int_0^t\omega(t-s)f(s)\,ds\right)^2dt
 \le\|\omega\|_{L^1(0,h_\lambda)}^2
                       \int_0^{h_\lambda}f(s)^2\,ds.
 \label{eq:initial-triangle-convolution-bound}
\end{equation}
Heat contraction and
Remark~\ref{rem:actual-input-moment-interpolation} bound the input
by the full-horizon norm. The kernel mass is taken up to
$h_\lambda\ge\lambda^{-4}$, even when the actual duration
$t<\lambda^{-4}$.

\emph{3. Estimate the first correction.}
In \eqref{eq:explicit-first-martingale} the two smoothed tests
$F_s^U,F_s^V$ coincide. Subtract the actual integrands, separating
the noise difference $L_U(s)-L_V(s)$ from the evaluation-point
difference $S_{t-s}D_s$.
The signed-input and output trace estimates in
Section~\ref{sec:complete-martingale-bound-proof} give the two
squared time-kernel masses, respectively,
\begin{equation}
 \begin{aligned}
 \lambda^{2-2\beta}\int_0^{h_\lambda}
             r^{-1/2}e^{-c\lambda^2\sqrt r}\,dr
       &\le C\lambda^{-2\beta},\\
 \lambda^{2-2\beta}\int_0^{h_\lambda}
             r^{\rho-1/2}e^{-c\lambda^2\sqrt r}\,dr
       &\le C\lambda^{-2\beta-4\rho}.
 \end{aligned}
 \label{eq:initial-first-kernel-masses}
\end{equation}
Both integrals converge at zero, by the same substitution
$u=\lambda^2\sqrt r$ used there.
For the evaluation-point term, heat contraction and the
fixed-variance clipped difference bound give an integrated
squared input at most
$C\min(h_\lambda,\lambda^2\|U-V\|_{\mathscr L_0^2(H)}^2)$.
Apply \eqref{eq:initial-triangle-kernel-bound}, take square roots
and use $h_\lambda\le1$ to obtain
\eqref{eq:initial-first-martingale-comparison}.

\emph{4. Estimate the higher stochastic terms.}
Use \eqref{eq:higher-martingale-integrand-split}, whose profile
contribution is zero. For an assignment with counts
$n_{\mathrm b},n_{\mathrm r}$ from
Definition~\ref{def:centered-word-assignments}, the changed-noise
kernel after the old-gap integration is, with $w=t-q$,
\[
 C\lambda^{-\beta}h_\lambda^{(n-n_{\mathrm b})\beta/4}
       w^{n_{\mathrm b}\beta/4-1/2}L_\lambda(w)^{n_{\mathrm r}}.
\]
Here $n=j-1$, and only the outside duration factors have been
bounded by $h_\lambda$. Its squared mass is controlled by
\eqref{eq:positive-square-integral} for $n_{\mathrm b}>0$ and
\eqref{eq:root-critical-square-integral} otherwise; in the latter
case $n_{\mathrm r}=n\ge1$ ensures integrability at zero.
For the current observations use the kernel
\eqref{eq:capped-native-mean-integrand} and its squared integral
\eqref{eq:current-mean-square-integral}. These estimates include
$w<\lambda^{-4}$. Thus It\^o isometry and
\eqref{eq:initial-triangle-kernel-bound} give
\eqref{eq:higher-noise-duration-bound} and
\eqref{eq:higher-current-duration-bound} with $h=h_\lambda$.
Lemma~\ref{lem:normalized-kernel-truncation} supplies the separate tail.
Absorb logarithms into the prescribed loss and restore the $V$
coefficient to obtain \eqref{eq:initial-higher-martingale-comparison}.

\emph{5. Estimate the remainder and coefficient replacements.}
Apply the current-value comparison of
Section~\ref{sec:paired-remainder-proof}; the covariance change is zero.
Lemma~\ref{lem:insertion-gap-integration} integrates each assigned
rough block with degree $\beta/4$, and each root interval contributes
$L_\lambda(h_\lambda)$. The unassigned and newest factors complete
the total degree to $N_\beta\beta/4$.
Pairing the root adds $\lambda^{\chi_{\mathrm{rem}}}$ by
\eqref{eq:paired-remainder-root-slot}. Pairing any rough factor
replaces just its degree $\beta/4$ by
$(\beta-\chi_{\mathrm{rem}})/4>0$, using
\eqref{eq:paired-remainder-rough-block}, and hence adds at most
$h_\lambda^{-\chi_{\mathrm{rem}}/4}\le\lambda^{\chi_{\mathrm{rem}}}$.
All gap kernels remain integrable. Apply the ordinary-integral
bound \eqref{eq:initial-triangle-convolution-bound} to $D_s$, with
Lemma~\ref{lem:comparison-spatial-inputs}(ii) and
Remark~\ref{rem:actual-input-moment-interpolation} for its spatial
and moment norms. Finally the embedding
$\|f\|_{H_\Lambda^{-\rho}}\le\Lambda^{-\rho/2}\|f\|_{L^2(\mathbb T)}$
and the prescribed small losses give
\eqref{eq:initial-paired-remainder-comparison}.

Apply the same integrations to \eqref{eq:bandwise-approximation-error}.
The limits follow from Corollary~\ref{cor:integrated-coefficient-error}
and Lemma~\ref{lem:exact-bracket-cutoff-limits}.
\end{proof}

\begin{samepage}
\subsection{Frequency summation and damped energy}
\label{sec:frequency-summation-energy}

\subsubsection{Converting word bounds into frequency decay}
\label{par:predictable-endpoint-comparison}
\label{sec:word-frequency-decay}
Fix a block $a_\lambda$ and set $h=\lambda^{-\gamma_{\mathrm{time}}}$, $0<\gamma_{\mathrm{time}}<4$.
The main statements take $\gamma_{\mathrm{time}}=3\beta$; the following result keeps
the dependence on $\gamma_{\mathrm{time}}$ explicit.
\end{samepage}

\begin{corollary}[Frequency decay of the higher-martingale terms]
\label{cor:higher-martingale-frequency-bounds}
For the coupled pair and stochastic integrals of
Definition~\ref{def:higher-coupled-stochastic-terms}, fix one
state-frequency block $a_\lambda$, an order $j\ge2$,
$0<\rho<1/2$ and $0<\chi<\min(1,\beta+2\rho)$.
Take $h=\lambda^{-\gamma_{\mathrm{time}}}$ with $0<\gamma_{\mathrm{time}}<4$ and
$\lambda^{-4}\le h\le h_{j-1}$. Define
\begin{equation}
 \begin{aligned}
 w_j&=\beta+\gamma_{\mathrm{time}}(j-1)\beta/4,\\
 q_{j,\chi}&=w_j+\min(2\rho-\chi,\gamma_{\mathrm{time}}(2\rho-\chi)/4),\\
 q_{j,\mathrm{prof}}&=w_j+\gamma_{\mathrm{time}}\rho/2-\gamma_{\mathrm{time}}/4.
 \end{aligned}
 \label{eq:higher-martingale-degrees-general}
\end{equation}
For a sufficiently small fixed loss $\varepsilon>0$ and any
prescribed tail exponent $p_{\mathrm{tail}}>0$,
\begin{equation}
 \begin{aligned}
 \|N_{j,\lambda}^U-N_{j,\lambda}^V\|_{\mathscr L_h^2(H_\Lambda^{-\rho})}
 \le C_{j,\varepsilon,p_{\mathrm{tail}}}\big[&
     \lambda^{-w_j+\varepsilon}\mathsf d_{\mathrm{coef}}(U,V)
    +\lambda^{-q_{j,\chi}+\varepsilon}
                         \|U-V\|_{\mathscr L_0^2(H)}^{\chi-\varepsilon}\\
    &+\lambda^{-q_{j,\mathrm{prof}}+\varepsilon}
                         \|U-V\|_{\mathscr L_0^2(H)}^{1-\varepsilon}
    +\lambda^{-p_{\mathrm{tail}}}\big]
    +\operatorname{Err}_{n_{\mathrm{sm}},\lambda}.
 \end{aligned}
 \label{eq:higher-martingale-frequency-combined}
\end{equation}
The first three terms on the right correspond, in order, to:
\begin{enumerate}[label=(\roman*),leftmargin=*,itemsep=3pt]
\item \emph{The actual noise amplitude:} the difference
$L(a(U))-L(a_{n_{\mathrm{sm}}}(V))$ in
\eqref{eq:higher-martingale-integrand-split}.
\item \emph{The current observations:} the change from the $V$-side
to the $U$-side observation vector in
\eqref{eq:actual-coupled-history-observations}, including its control correction.
\item \emph{The Gaussian reference profile:} the change from $\bar a^V$
to $\bar a^U$, which changes both the Gaussian covariances in
\eqref{eq:predictable-mean-covariance} and the subtracted centering
values in \eqref{eq:actual-coupled-inserted-tests}.
\end{enumerate}
The global input quantities are defined in
\eqref{eq:coefficient-discrepancy-norm}. The term
$\lambda^{-p_{\mathrm{tail}}}$ is the truncation tail, and
$\operatorname{Err}_{n_{\mathrm{sm}},\lambda}$ is the
coefficient-replacement error defined in
\eqref{eq:bandwise-approximation-error}.
The constant is independent of $\lambda$, $\Lambda$ and
$n_{\mathrm{sm}}$; fixed logarithmic factors are included in
the loss $\varepsilon$.
\end{corollary}
\begin{proof}[Proof of Corollary~\ref{cor:higher-martingale-frequency-bounds}]
Put $n=j-1$ and $x_+=\max(x,0)$. Adding the three duration estimates
\eqref{eq:higher-noise-duration-bound},
\eqref{eq:higher-current-duration-bound} and
\eqref{eq:higher-profile-duration-bound} gives
\begin{equation}
 \begin{aligned}
 &\|N_{j,\lambda}^U-N_{j,\lambda}^V\|_{\mathscr L_h^2(H_\Lambda^{-\rho})}\\
 &\quad\le C_{j,\varepsilon,p_{\mathrm{tail}}}\big[
       \lambda^{-\beta+\varepsilon}h^{n\beta/4}
                                 \mathsf d_{\mathrm{coef}}(U,V)\\
 &\qquad\quad
       +\lambda^{-\beta+(\chi-2\rho)_++\varepsilon}
        h^{n\beta/4+(2\rho-\chi)_+/4}
                           \|U-V\|_{\mathscr L_0^2(H)}^{\chi-\varepsilon}\\
 &\qquad\quad
       +\lambda^{-\beta+\varepsilon}
        h^{n\beta/4+\rho/2-1/4}
                           \|U-V\|_{\mathscr L_0^2(H)}^{1-\varepsilon}
       +\lambda^{-p_{\mathrm{tail}}}\big]
       +\operatorname{Err}_{n_{\mathrm{sm}},\lambda}.
 \end{aligned}
 \label{eq:higher-martingale-combined-duration}
\end{equation}
For $h=\lambda^{-\gamma_{\mathrm{time}}}$, the three factors become
\[
 \begin{aligned}
 \lambda^{-\beta}h^{n\beta/4}&=\lambda^{-w_j},\\
 \lambda^{-\beta+(\chi-2\rho)_+}
 h^{n\beta/4+(2\rho-\chi)_+/4}
 &=\lambda^{-w_j+(\chi-2\rho)_+
                   -\gamma_{\mathrm{time}}(2\rho-\chi)_+/4}
   =\lambda^{-q_{j,\chi}},\\
 \lambda^{-\beta}h^{n\beta/4+\rho/2-1/4}
 &=\lambda^{-w_j-\gamma_{\mathrm{time}}\rho/2
                   +\gamma_{\mathrm{time}}/4}
   =\lambda^{-q_{j,\mathrm{prof}}}.
 \end{aligned}
\]
The middle equality uses $0<\gamma_{\mathrm{time}}/4<1$:
its added decay is $\gamma_{\mathrm{time}}(2\rho-\chi)/4$ when
$\chi\le2\rho$, and $2\rho-\chi$ otherwise. These are exactly the
minimum in \eqref{eq:higher-martingale-degrees-general}.
Also $\mathcal L(\lambda,\lambda^{-\gamma_{\mathrm{time}}})
\le C_{\gamma_{\mathrm{time}}}(1+\log\lambda)$, so fixed logarithmic
powers are included in the prescribed small loss. The separate tail
and coefficient error in \eqref{eq:higher-martingale-combined-duration}
remain additive. This proves \eqref{eq:higher-martingale-frequency-combined};
setting $\gamma_{\mathrm{time}}=3\beta$ gives the main comparison.
\end{proof}

\begin{remark}[Dyadic sums for the first martingale correction]
\label{rem:first-martingale-dyadic-sums}
For a general duration $h=\lambda^{-\gamma_{\mathrm{time}}}$, $0<\gamma_{\mathrm{time}}<4$, the high-state
sum of \eqref{eq:current-cap-complete-Gamma} needs $\beta>0$ and
$q_\Gamma=\alpha-\gamma_{\mathrm{time}}(1-\beta)/4>0$; the main statement uses $\gamma_{\mathrm{time}}=3\beta$.
Here $\alpha=\beta+2\rho$, $q_\Gamma$ is the decay exponent of the reference-profile term,
and $\lambda_*\asymp\Lambda^{1/4}$ is the dyadic lower cutoff from
\eqref{eq:gaussian-band-dyadic-cutoff}. The root-minimum sum is
$C\|U-V\|_{\mathscr L_0^2(H)}^\alpha$ for $\alpha<1$, logarithmic at equality, and
$C\lambda_*^{1-\alpha}\|U-V\|_{\mathscr L_0^2(H)}$ over $\lambda\ge \lambda_*$ for $\alpha>1$.
The actual-noise sum is $C\lambda_*^{-\beta}\mathsf d_{\mathrm{coef}}(U,V)$.
\end{remark}

\begin{corollary}[Frequency bounds for Gaussian averages and actual remainders]
\label{cor:average-remainder-frequency-bounds}
Use the common-coefficient terms of
Definition~\ref{def:predictable-mean-functional} and
\eqref{eq:common-coefficient-remainders}. Let $n\ge1$ denote the
number of insertions in both parts, and take
$h=\lambda^{-\gamma_{\mathrm{time}}}\le h_n$,
$0<\gamma_{\mathrm{time}}<4$ and $1\le p<\infty$.
\begin{enumerate}[label=(\roman*),leftmargin=*,itemsep=4pt]
\item \emph{Gaussian average terms.}
The spatial $L^2$ norm of either endpoint $\mathsf P_{n+1,\lambda}$
is bounded, up to a fixed-order constant, by
\begin{equation}
 \lambda^{-\beta}h^{n\beta/4}\mathcal L(\lambda,h)^n,
 \label{eq:capped-root-predictable-size}
\end{equation}
Its frequency degree is
$w_{n+1}=\beta+\gamma_{\mathrm{time}}n\beta/4$.
Set $p_{n+1}=w_{n+1}-\gamma_{\mathrm{time}}/4$. For every
$\varepsilon>0$, the endpoint difference satisfies
\begin{equation}
 \|\mathsf P_{n+1,\lambda}^U-\mathsf P_{n+1,\lambda}^V\|_{\mathscr L_h^2(H_\Lambda^{-\rho})}
 \le C_{n+1,\varepsilon}\Lambda^{-\rho/2}
       \lambda^{-p_{n+1}+\varepsilon}\|U-V\|_{\mathscr L_0^2(H)}.
 \label{eq:ordinary-predictable-comparison-row}
\end{equation}
\item \emph{Actual remainders.}
For $i=U,V$, let $R_{n,\lambda}^i$ be the remainder in
\eqref{eq:common-coefficient-remainders}, using its matching
reference. Then
\begin{equation}
 \|R_{n,\lambda}^i\|_{L^p(\mathbb P_{t_0};L^2(\mathbb T))}
 \le C_{n ,p}\mathcal L(\lambda,h)^{C_{n} }\lambda^{-\beta}h^{n \beta/4}.
 \label{eq:capped-root-actual-strong-remainder}
\end{equation}
The remainder degree is
$\beta+\gamma_{\mathrm{time}}n\beta/4=w_{n+1}$; embedding in
$H_\Lambda^{-\rho}$ adds the factor $\Lambda^{-\rho/2}$.
\end{enumerate}
\end{corollary}
\begin{proof}[Proof of Corollary~\ref{cor:average-remainder-frequency-bounds}]

\medskip
\noindent\emph{1. The size of each Gaussian average: \eqref{eq:capped-root-predictable-size}.}\par
\nopagebreak[4]
For an assignment of Definition~\ref{def:centered-word-assignments},
write $s=s_n$, $\delta_j=s_{j-1}-s_j$, $I_0$ for its root block,
and $\mathscr B_{\mathrm{rough}}$ for its nonempty rough blocks.
Set $n_{\mathrm r}=|I_0|$,
$n_{\mathrm b}=|\mathscr B_{\mathrm{rough}}|$ and
$u=n-1-n_{\mathrm b}$. The integrated assignment bound
\eqref{eq:procedure-final-predictable-integral} is
$C_n\lambda^{-\beta}h^{n\beta/4}L_\lambda(h)^{n_{\mathrm r}}$.
Since $n_{\mathrm r}\le n$ and
$L_\lambda(h)\le C\mathcal L(\lambda,h)$, summing assignments
proves \eqref{eq:capped-root-predictable-size}.
Its frequency degree follows from the identity
\[
 \lambda^{-\beta}h^{n\beta/4}
 =\lambda^{-\beta-\gamma_{\mathrm{time}}n\beta/4}
 =\lambda^{-w_{n+1}}.
\]

\medskip
\noindent\emph{2. The endpoint comparison: \eqref{eq:ordinary-predictable-comparison-row}.}\par
\nopagebreak[4]
In \eqref{eq:all-order-predictable-integrated-sobolev}, substitute
\[
 \lambda^{-\beta}h^{n\beta/4-1/4}
 =\lambda^{-\beta-\gamma_{\mathrm{time}}(n\beta-1)/4}
 =\lambda^{-p_{n+1}}.
\]
For the fixed time scaling exponent, the logarithmic factor obeys
\[
 \mathcal L(\lambda,\lambda^{-\gamma_{\mathrm{time}}})
 \le C_{\gamma_{\mathrm{time}}}(1+\log\lambda),\qquad
 \mathcal L(\lambda,h)^{C_n}
 \le C_{n,\varepsilon}\lambda^\varepsilon.
\]
This gives \eqref{eq:ordinary-predictable-comparison-row}.

\medskip
\noindent\emph{3. The actual remainder: \eqref{eq:capped-root-actual-strong-remainder}.}\par
\nopagebreak[4]
Use the actual history integrals
\eqref{eq:common-coefficient-remainders}, and again fix their
earliest insertion time $s=s_n$. The newest coefficient factor
is now evaluated at the actual solution. At the two endpoints
these factors are explicitly
\[
 a(U_{t_0+s}(\xi_n))-\bar a_s^U(\xi_n),\qquad
 a(V_{t_0+s}(\xi_n))-\bar a_s^V(\xi_n).
\]
These factors are fixed under the future Gaussian expectation and
random under $\mathbb E_{t_0}$.
Lemma~\ref{lem:consistent-current-interpolation}, at $\theta=1,0$,
gives, for every finite moment order $p\ge1$,
\begin{equation}
 \begin{aligned}
 &\sup_{\xi\in\mathbb T}
  \bigl(\mathbb E_{t_0}|a(U_{t_0+s}(\xi))-\bar a_s^U(\xi)|^p\bigr)^{1/p}\\
 &\quad+\sup_{\xi\in\mathbb T}
  \bigl(\mathbb E_{t_0}|a(V_{t_0+s}(\xi))-\bar a_s^V(\xi)|^p\bigr)^{1/p}
 \le C_ps^{\beta/4}.
 \end{aligned}
 \label{eq:frequency-remainder-newest-moment}
\end{equation}
The same lemma and \eqref{eq:hybrid-controlled-future-mean}
give $h^{\beta/4}$ moments for the unassigned older factors,
uniformly in $\Lambda$.

For the later Gaussian intervals, apply the centering identity
\eqref{eq:procedure-score-centering} and the spatial bounds used
in \eqref{eq:all-order-predictable-common-ledger}.
An interval assigned to $a_\lambda$ gives
$k_\lambda(\delta_j)$; a nonempty block $I$ assigned to one older
coefficient gives
$(\min_{j\in I}\delta_j)^{\beta/4}\prod_{j\in I}\delta_j^{-1}$.
For $p\ge2$, conditional Minkowski and H\"older apply with
the factor bounds above and \eqref{eq:multipoint-score-products}
at fixed higher moments depending on $n,p$. For $p<2$, use
the order-two bound and norm monotonicity.

Consequently, the contribution of each assignment to the
conditional $L^p(\mathbb P_{t_0};L^2(\mathbb T))$ norm of the
remainder is bounded by
\begin{equation}
 \begin{aligned}
 C_{n,p}\lambda^{-\beta}h^{u\beta/4}
 &\int_{\substack{\delta_j>0\\\sum_j\delta_j<h}}
 s^{\beta/4}\prod_{j\in I_0}k_\lambda(\delta_j)\\
 &\quad\times\prod_{I\in\mathscr B_{\mathrm{rough}}}
 \left[(\min_{j\in I}\delta_j)^{\beta/4}
                  \prod_{j\in I}\delta_j^{-1}\right]
                         d\boldsymbol\delta .
 \end{aligned}
 \label{eq:frequency-remainder-gap-integral}
\end{equation}
Here $s=h-\sum_j\delta_j$.
Apply \eqref{eq:minimum-gap-cube-bound} to each rough block,
then \eqref{eq:insertion-gap-simplex-bound} with total block
degree $n_{\mathrm b}\beta/4$ and newest-factor exponent
$\sigma=\beta/4$. The result is
\[
 C_{n,p}\lambda^{-\beta}
 h^{(u+n_{\mathrm b}+1)\beta/4}L_\lambda(h)^{n_{\mathrm r}}
 =C_{n,p}\lambda^{-\beta}h^{n\beta/4}L_\lambda(h)^{n_{\mathrm r}}.
\]
Summing the assignments and using $n_{\mathrm r}\le n$ proves
\eqref{eq:capped-root-actual-strong-remainder} for each $i=U,V$.
The same identity $\lambda^{-\beta}h^{n\beta/4}=\lambda^{-w_{n+1}}$
gives its asserted frequency degree. Finally, the spatial
embedding acts inside the conditional moment norm:
\begin{equation}
 \begin{aligned}
 \|R_{n,\lambda}^i\|_{L^p(\mathbb P_{t_0};H_\Lambda^{-\rho})}
 &\le\Lambda^{-\rho/2}
       \|R_{n,\lambda}^i\|_{L^p(\mathbb P_{t_0};L^2(\mathbb T))}\\
 &\le C_{n,p}\Lambda^{-\rho/2}\mathcal L(\lambda,h)^{C_n}
                       \lambda^{-w_{n+1}}.
 \end{aligned}
 \label{eq:frequency-remainder-negative-norm}
\end{equation}
This proves the last assertion of part~(ii).
\end{proof}

\subsubsection{Damped energy and the frequency sums used in closure}
\label{sec:haar-decoupled-mass}
\label{sec:finite-low-covariance-comparison}
\begin{lemma}[Damped energy and low-frequency absorption]
\label{lem:damped-energy-low-frequency}
Let $U,V$ be the controlled pair satisfying
\eqref{eq:companion-original-rotated} and \eqref{eq:companion-controlled},
driven by the same Brownian noises $W^0,B$ and using the same fixed
covariance-factor map $L$. Assume $U_0=V_0$, write $D=U-V$, and take
$\Lambda\ge1$, $0<\rho<1/2$. Thus $D$ satisfies the damped difference
equation \eqref{eq:companion-difference}. Use the global norms in
\eqref{eq:coefficient-discrepancy-norm}.
\begin{enumerate}[label=(\roman*),leftmargin=*,itemsep=4pt]
\item \emph{Energy from the coefficient difference.}
\begin{equation}
 \|U-V\|_{\mathscr L_0^2(H)}^2\le C\Lambda^{\rho-1/2}\mathsf d_{\mathrm{coef}}(U,V)^2.
 \label{eq:decoupled-MJ}
\end{equation}
\item \emph{Low state frequencies.}
For $\lambda_*\ge1$, set
$a_{\le\lambda_*}=a_{\mathrm{low}}+
\sum_{\lambda\in\mathscr S_{\mathrm z},\,\lambda\le\lambda_*}a_\lambda$.
Then
\begin{equation}
 \|a_{\le \lambda_* }(U)-a_{\le \lambda_* }(V)\|_{\mathscr L_0^2(H_\Lambda^{-\rho})}
       \le C\Lambda^{-\rho/2}\lambda_*^{1-\beta}\|U-V\|_{\mathscr L_0^2(H)}.
 \label{eq:current-cap-low-raw}
\end{equation}
For $\lambda_*\asymp\Lambda^r$, $r>0$, the low-frequency bound
is at most
$C\Lambda^{-1/4+r(1-\beta)}\mathsf d_{\mathrm{coef}}(U,V)$.
This coefficient tends to zero when $r(1-\beta)<1/4$.
\end{enumerate}
\end{lemma}
\begin{proof}[Proof of Lemma~\ref{lem:damped-energy-low-frequency}]
\emph{Energy.}
Apply \eqref{eq:companion-signed-energy} with $h=T$ and include
the nonnegative allowance in \eqref{eq:coefficient-discrepancy-norm}.
This proves \eqref{eq:decoupled-MJ}. Its hypotheses are exactly the
shared-noise coupling and damped difference equation assumed here.
If $a$ is constant, then $a_{n_{\mathrm{sm}}}=a$ and the noise-factor
difference in \eqref{eq:companion-difference} is zero; hence $D=0$.

\emph{Low frequencies.}
The low-block convolution kernel has integrable derivative, and
\eqref{eq:bands-bernstein-tail} bounds the dyadic derivatives. Since
$0<\beta<1$,
\[
 \|\partial_z a_{\le\lambda_*}\|_\infty
 \le C+C\sum_{\substack{\lambda\in\mathscr S_{\mathrm z}\\\lambda\le\lambda_*}}
                    \lambda^{1-\beta}
 \le C_\beta\lambda_*^{1-\beta}.
\]
The mean-value theorem and
$\|f\|_{H_\Lambda^{-\rho}}\le\Lambda^{-\rho/2}\|f\|_{L^2(\mathbb T)}$
give \eqref{eq:current-cap-low-raw} after probability--time integration.
Combining it with the energy bound gives
\begin{equation}
 \begin{aligned}
 &\|a_{\le\lambda_*}(U)-a_{\le\lambda_*}(V)\|_{\mathscr L_0^2(H_\Lambda^{-\rho})}\\
 &\qquad\le C\Lambda^{-\rho/2}\lambda_*^{1-\beta}
                  \Lambda^{\rho/2-1/4}\mathsf d_{\mathrm{coef}}(U,V)\\
 &\qquad=C\Lambda^{-1/4}\lambda_*^{1-\beta}
                                      \mathsf d_{\mathrm{coef}}(U,V).
 \end{aligned}
 \label{eq:low-frequency-energy-substitution}
\end{equation}
For $\lambda_*\asymp\Lambda^r$, the coefficient is
$C\Lambda^{-1/4+r(1-\beta)}$, which vanishes when $r(1-\beta)<1/4$.
For $r=1/4$ it is $C\Lambda^{-\beta/4}$, giving the stated absorption.
\end{proof}

\begin{remark}[The remainder contribution to feedback]
\label{rem:remainder-feedback-degree}
For the exponent $0<\chi<\beta$ in
Proposition~\ref{prop:paired-actual-remainder}, put
$W_{N_\beta}=\beta+\gamma_{\mathrm{time}}\beta N_\beta/4>\chi$.
Its dyadic high-band contribution is
\[
 C\Lambda^{-\rho/2}\lambda_*^{-(W_{N_\beta}-\chi)+o(1)}
                         \|U-V\|_{\mathscr L_0^2(H)}^{\chi-o(1)},
\]
where $o(1)$ denotes exponent losses made arbitrarily small at fixed
order. Combining this particular term with
$\|U-V\|_{\mathscr L_0^2(H)}^2\le C\Lambda^{\rho-1/2}\mathsf d_{\mathrm{coef}}(U,V)^2$ and Young's inequality
gives the following exponent for this contribution. Consider the
corresponding isolated distance inequality at zero loss:
\[
 \|U-V\|_{\mathscr L_0^2(H)}
 \le C\Lambda^{-1/4-r(W_{N_\beta}-\chi)}
                           \|U-V\|_{\mathscr L_0^2(H)}^{\chi}.
\]
For a positive distance, divide by $\|U-V\|_{\mathscr L_0^2(H)}^{\chi}$ and
raise to the power $2/(1-\chi)$. This gives the squared-distance
decay exponent
\begin{equation}
 \frac{1/2+2r(W_{N_\beta}-\chi)}{1-\chi}.
 \label{eq:decoupled-feedback}
\end{equation}
Section~\ref{sec:all-beta-joint-complete-closure} performs the same
calculation for the other two nonlinear terms and chooses parameters
so that all three exponents have the required positive margin.
\end{remark}

\subsection{Passing to the coefficient and noise limits}
\label{sec:approximation-and-transfer}

The estimates proved at finite cutoffs pass to the actual coefficients
and noise series by the following lemma. The coefficient-error calculation
at colliding times, including the integrable majorants, is in
Appendix~\ref{sec:uniform-coefficient-diagonal}. Section~6 verifies the
centering and spatial bounds directly for its joint coordinates.

\subsubsection{Passing to one common coefficient approximation}
\label{sec:common-approximation-proof}
The next lemma controls the errors defined in
\eqref{eq:bandwise-approximation-error} for the fixed physical pair
$(U,V^{n_{\mathrm{sm}},\Lambda})$, finite depth $N_\beta$,
root-frequency cutoff $\lambda_{\mathrm{cut}}$ and positive-gap
cutoff $\eta>0$, which restricts every included time gap to be at
least $\eta$. It gives a single approximation for all included expansion terms.

\begin{lemma}[Exact brackets and approximation limits]
\label{lem:exact-bracket-cutoff-limits}
Assume the boundedness, uniform ellipticity and common $C^\beta$
coefficient bounds of Sections 2--3. Fix a finite expansion order and
a coupled pair $(U,V)$ with its exact Haar innovation. The finite-word
identities and comparisons hold for this pair under the parameter
ranges stated for each estimate. Noise-series approximation preserves the physical bracket in the limit. Coefficient-approximation
errors are removed at finite root-frequency and positive-gap cutoffs,
with uniform tails as those restrictions are removed.
There is a modulus
$\omega_{\Lambda,\lambda_{\mathrm{cut}},\eta}:[0,\infty)\to[0,\infty)$,
vanishing at zero, such that
\[
 \sum_{\substack{\lambda\in\mathscr S_{\mathrm z}\\\lambda\le\lambda_{\mathrm{cut}}}}
 \operatorname{Err}_{n_{\mathrm{sm}},\lambda}
 \le\omega_{\Lambda,\lambda_{\mathrm{cut}} ,\eta}(\|a_{n_{\mathrm{sm}}}-a\|_\infty),
 \qquad \omega_{\Lambda,\lambda_{\mathrm{cut}} ,\eta}(r)\longrightarrow0\quad(r\downarrow0).
\]
\end{lemma}

The proof is given in Appendix~\ref{app:approximation-limits}.

\clearpage
\section{A stopped finite-order comparison in joint coordinates}
\label{sec:single-species-joint-flow-recursion}
The controlled companion $V$ has an extra drift $\Lambda(U-V)$.
Backward coordinates built from the coupled pair cancel this drift,
so the finite expansion can be compared through its covariance terms.
We define the recursive operations and construct these coordinates.
Lemma~\ref{lem:joint-coordinate-transfer} verifies the centering,
spatial input and profile bounds needed to repeat the Section~5
comparisons in these coordinates.
Three results complete the comparison:
\begin{enumerate}[label=(\roman*),leftmargin=*]
\item \textbf{The finite identity.} Theorem~\ref{thm:joint-finite-comparison} compares the initial-past averages, stochastic integrals and remainder, recording the extra covariance corrections separately.
\item \textbf{The extra covariance terms.} Proposition~\ref{prop:all-order-first-hit-joint-boundary} bounds their total contribution by a small multiple of the coefficient-comparison norm.
\item \textbf{The leading stochastic term.} Proposition~\ref{prop:joint-first-pair} combines the first average and stochastic integral into a Gaussian-reference difference and the complete martingale estimated in Section~5.
\end{enumerate}
The proofs and supporting estimates are collected in the last subsection.

\subsection{Main results}

\subsubsection{The coupled pair and the expansion terms}
Use the coupled solutions $U,V$ on $[0,T]$, with the same initial
state, and write $D=U-V$. As in Section~5, $a=g^2-m_0^2$ is the
residual variance, $a_{n_{\mathrm{sm}}}$ its state mollification,
and $a_\lambda$ a dyadic state-frequency band.
The heat semigroup and kernel are $S_t$ and $p_t$.
In this section's main statements use
$h=h_\lambda=\lambda^{-3\beta}$ and the cutoff
\eqref{eq:gaussian-band-dyadic-cutoff},
$\lambda_*\asymp\Lambda^{1/4}$, as in Section~5.1.
The residual noise factors are $L_U=L(a(U))$ and
$L_V=L(a_{n_{\mathrm{sm}}}(V))$, driven by the shared independent
Brownian motions $W^0,B$:
\[
 dU=\Delta U\,dt+m_0dW^0+L_UdB,\qquad
 dV=\Delta V\,dt+\Lambda D\,dt+m_0dW^0+L_VdB.
\]
Write $B_\Delta=L_U-L_V$ and $f=a(U)-a_{n_{\mathrm{sm}}}(V)$.
On $[t_0,t_0+h]\subset[0,T]$, $t_0$ is the physical starting time and
$h$ the duration. Local times are measured from $t_0$, so an insertion
at $s_j$ occurs at physical time $t_0+s_j$.
The notation $\mathbb E_{t_0}$ means conditioning on
$\mathcal F_{t_0}$, and the shifted noises are
$W_s^{0,t_0}=W_{t_0+s}^0-W_{t_0}^0$, $B_s^{t_0}=B_{t_0+s}-B_{t_0}$.

Recall that $\mathsf P_{j,\lambda}$ is the initial-past Gaussian average,
$N_{j,\lambda}$ the centered stochastic integral, and $R_{n,\lambda}$ the remainder
after $n$ covariance insertions. An inserted factor
$a(w)-\bar a_t(\xi)$ is a rough coefficient factor, with subtracted value
$\bar a_t(\xi)$; the original factor $a_\lambda$ is the root.
Here $\mathsf P_{j,\lambda},N_{j,\lambda}$ contain $j-1$ insertions: these integer
subscripts record expansion order. The insertion times are continuous
integration variables, as in Definition~\ref{def:finite-coefficient-band-words}.
The reference profile is fixed by the initial past and chosen so
that $\mathbb E_{t_0}[a(X_t(x))]=\bar a_t(x)$ for its Gaussian reference $X$.

\subsubsection{The drift-cancelling backward coordinates}
\begin{definition}[Scalar joint backward coordinate]
\label{def:scalar-joint-backward-coordinate}
For a terminal point $x\in\mathbb T$ on $[t_0,t_0+h]$, set
\begin{equation}
 Y_s^{V,h}(x):=(S_{h-s}U_{t_0+s})(x)
       -e^{-\Lambda(h-s)}(S_{h-s}D_{t_0+s})(x),\qquad 0\le s\le h.
 \label{eq:joint-scalar-backward-coordinate}
\end{equation}
\end{definition}
\begin{lemma}[Drift cancellation in the backward coordinate]
\label{lem:joint-backward-drift-cancellation}
The coordinate in Definition~\ref{def:scalar-joint-backward-coordinate}
is a martingale, ends at $V_{t_0+h}(x)$, and satisfies
\begin{equation}
 \begin{aligned}
 dY_s^{V,h}(x)={}&m_0(S_{h-s}\,dW_s^{0,t_0})(x)\\
 &+\bigl(S_{h-s}[L_V+(1-e^{-\Lambda(h-s)})B_\Delta]\,dB_s^{t_0}\bigr)(x),
 \qquad s<h,
 \end{aligned}
 \label{eq:joint-scalar-drift-cancellation}
\end{equation}
where the noise factors are evaluated at time $t_0+s$.
\end{lemma}
\begin{proof}[Proof of Lemma~\ref{lem:joint-backward-drift-cancellation}]
Fix the window start $t_0$, terminal point $x$ and duration $h$.
Write $\tau=h-s$ and $\vartheta_s=e^{-\Lambda(h-s)}$; all noise
operators below are evaluated at physical time $t_0+s$.
The coupled equations recalled above, or equivalently
\eqref{eq:companion-difference}, give
\[
 dD_{t_0+s}=(\Delta-\Lambda)D_{t_0+s}\,ds
                         +B_\Delta\,dB_s^{t_0}.
\]
Backward heat transport cancels the Laplacian in both equations:
\[
 \begin{aligned}
 d(S_\tau U_{t_0+s})
   &=m_0S_\tau\,dW_s^{0,t_0}+S_\tau L_U\,dB_s^{t_0},\\
 d(S_\tau D_{t_0+s})
   &=-\Lambda S_\tau D_{t_0+s}\,ds
                            +S_\tau B_\Delta\,dB_s^{t_0}.
 \end{aligned}
\]
Since $d\vartheta_s=\Lambda\vartheta_s\,ds$, the product rule gives
\begin{equation}
 \begin{aligned}
 d(\vartheta_sS_\tau D_{t_0+s})
 &=\bigl[\Lambda\vartheta_sS_\tau D_{t_0+s}
       -\Lambda\vartheta_sS_\tau D_{t_0+s}\bigr]ds
       +\vartheta_sS_\tau B_\Delta\,dB_s^{t_0}\\
 &=\vartheta_sS_\tau B_\Delta\,dB_s^{t_0}.
 \end{aligned}
 \label{eq:joint-backward-product-cancellation}
\end{equation}
Subtract this differential from that of $S_\tau U_{t_0+s}$ and use
$L_U-\vartheta_sB_\Delta=L_V+(1-\vartheta_s)B_\Delta$.
Evaluation at $x$ proves \eqref{eq:joint-scalar-drift-cancellation}.
At $s=h$, $\vartheta_h=1$ and $S_0=I$, so
$Y_h^{V,h}(x)=U_{t_0+h}(x)-D_{t_0+h}(x)=V_{t_0+h}(x)$.
The computations can first be made for $s\le h-\epsilon$ in the
mild equations. The uniform covariance bound and
$\int_0^h\|p_{h-s}\|_{L^2(\mathbb T)}^2ds<\infty$ give square-integrable
terminal limits as $\epsilon\downarrow0$. Thus the local stochastic
integrals extend to a martingale on the whole closed interval.
\end{proof}
Its residual noise factor is
$L_V+(1-e^{-\Lambda(h-s)})B_\Delta$ rather than $L_V$;
the resulting difference of covariance brackets is the
\emph{bracket defect} tracked below. Here the bracket is the
quadratic covariation of the observation coordinates.
The zero drift removes the control integral from the backward
It\^o expansion. Its replacement is this covariance correction,
whose factor $1-e^{-\Lambda(h-s)}\le\min(1,\Lambda(h-s))$
vanishes at the terminal time. This factor is used explicitly in
\eqref{eq:first-defect-integrated-pair-kernel}.
A branch is a sequence of local coefficient insertions. It is
expanded until this extra covariance term first appears, and that
term is recorded as a \emph{first-defect mark}.

\begin{definition}[Joint backward observation fields]
\label{def:joint-backward-observations}
For an ordered list
$Z=((\sigma_i,T_i,x_i))_{i=1}^{d_{\mathrm{obs}}}$, let $\sigma_i$
be the observation time, $T_i$ its heat-target time, and $x_i$ its
spatial center. For each fixed history these are deterministic elapsed times.
Assume $0<\sigma_i\le T_i\le h$,
$x_i\in\mathbb T$ and a fixed finite observation count
$d_{\mathrm{obs}}\ge1$. Define
\[
 Y_t^{V,Z,i}:=(S_{T_i-t}U_{t_0+t})(x_i)
       -e^{-\Lambda(\sigma_i-t)}(S_{T_i-t}D_{t_0+t})(x_i),
       \qquad 0\le t\le\sigma_i.
\]
The vector of joint observations is explicitly
\[
 \mathbf Y_t^{V,Z}:=(Y_t^{V,Z,1},\ldots,Y_t^{V,Z,d_{\mathrm{obs}}}),
 \qquad 0\le t\le\min_i\sigma_i.
\]
This extends \eqref{eq:joint-scalar-backward-coordinate} to several
possibly heat-smoothed observations. Each coordinate has the same
drift cancellation as \eqref{eq:joint-scalar-drift-cancellation}.
At its observation time $\sigma_i$, it satisfies
$Y_{\sigma_i}^{V,Z,i}=(S_{T_i-\sigma_i}V_{t_0+\sigma_i})(x_i)$.
The compatible Gaussian references $X^U,X^V$ have prescribed means
$\phi_s^U=S_sU_{t_0}$ and
$\phi_s^V=S_sU_{t_0}-e^{-\Lambda s}S_sD_{t_0}$.
Their profiles $\bar a^i$, $i\in\{U,V\}$, satisfy the centering
property just recalled; the joint $V$ calculation uses $\bar a^V$.
These are the references of
Definition~\ref{def:coupled-gaussian-references}.
Work on $0\le q\le s\le\min_i\sigma_i$.
A \emph{heat row} is the spatial kernel of an observation.
For this list, define
\[
 q_t^{Z,i}(\xi):=p_{T_i-t}(x_i-\xi),\qquad
 \mathbf q_t^Z(\xi):=(q_t^{Z,i}(\xi))_{i=1}^{d_{\mathrm{obs}}}.
\]
For $i\in\{U,V\}$, define the reference covariance and its Gaussian
averaging operator by
\begin{equation}
 \begin{aligned}
 \Sigma^{Z,i}_{q,s}
 &:=\int_q^s\int_{\mathbb T}(m_0^2+\bar a_t^i(\xi))
            \mathbf q_t^Z(\xi)\mathbf q_t^Z(\xi)^*\,d\xi\,dt,\\
 (\mathcal P^{Z,i}_{q,s}F)(y)
 &:=\int_{\mathbb R^{d_{\mathrm{obs}}}}F(y+z)\,
                          \mathcal N(0,\Sigma^{Z,i}_{q,s})(dz),
 \end{aligned}
 \label{eq:joint-gaussian-averaging-operator}
\end{equation}
where $F:\mathbb R^{d_{\mathrm{obs}}}\to\mathbb R$ is bounded and Borel.
When the endpoint is suppressed, $\Sigma^Z_{q,s}=\Sigma^{Z,V}_{q,s}$ and
$\mathcal P^Z_{q,s}=\mathcal P^{Z,V}_{q,s}$.
The mean $y$ is supplied by the observation coordinate at time $q$;
this operator averages only the future centered Gaussian increment.
The covariance and the operator depend on the full list $Z$
and on the selected window-start profile $\bar a^i$.
\end{definition}

\begin{definition}[Covariance defect of the joint fields]
\label{def:joint-covariance-defect}
For $0\le t<\min_i\sigma_i$, abbreviate
$q_i=q_t^{Z,i}(\cdot)\in H=L^2(\mathbb T)$ and set
$r_i=r_i(t):=1-e^{-\Lambda(\sigma_i-t)}$.
In the following formulas the noise operators are evaluated at the
physical time $t_0+t$: $L_V=L(a_{n_{\mathrm{sm}}}(V_{t_0+t}))$ and
$B_\Delta=L(a(U_{t_0+t}))-L_V$.
The $i$th joint coordinate has residual-noise row
$q_i[L_V+r_iB_\Delta]$, whereas the local $V$ row is $q_iL_V$.
Subtracting their covariance products defines the signed correction
\begin{equation}
 \begin{aligned}
 E_{ij}(t):={}&r_i\langle B_\Delta^*q_i,L_V^*q_j\rangle_H\\
 &+r_j\langle L_V^*q_i,B_\Delta^*q_j\rangle_H\\
 &+r_ir_j\langle B_\Delta^*q_i,B_\Delta^*q_j\rangle_H.
 \end{aligned}
 \label{eq:single-species-direct-defect}
\end{equation}
The first two terms are the cross covariances of the local $V$
noise and the added noise; the last is the covariance of the added
noise itself. The common additive-noise channel cancels in the
subtraction. Equivalently, the quadratic covariation satisfies
\begin{equation}
 \begin{aligned}
 \frac{d}{dt}\langle Y^{V,Z,i},Y^{V,Z,j}\rangle_t
 ={}&\int_{\mathbb T}
       [m_0^2+a_{n_{\mathrm{sm}}}(V_{t_0+t}(\xi))]
           q_i(\xi)q_j(\xi)\,d\xi+E_{ij}(t).
 \end{aligned}
 \label{eq:joint-bracket-defect-definition}
\end{equation}
Thus $E_{ij}(t)$ corrects the actual local $V$ covariance. The
separate difference $a_{n_{\mathrm{sm}}}(V)-\bar a^V$ from the
Gaussian reference is the local coefficient insertion.
For a smooth test $F:\mathbb R^{d_{\mathrm{obs}}}\to\mathbb R$,
define the defect operator by its action:
\begin{equation}
 (\mathcal E_tF)(y)
 :=\frac12\sum_{i,j=1}^{d_{\mathrm{obs}}}E_{ij}(t)
                     \frac{\partial^2F}{\partial y_i\partial y_j}(y).
 \label{eq:joint-covariance-defect-operator}
\end{equation}
Thus $\partial_{ij}$ below means differentiation in the $i$th
and $j$th state coordinates. Formula~\eqref{eq:single-species-direct-defect} also
shows that the defect vanishes when $B_\Delta=0$.

\end{definition}

\begin{definition}[Recursive tests, expansion terms and covariance corrections]
\label{def:joint-first-defect-marks}
Fix $t_0,h,x$ and the terminal test $a_\lambda$.
The following operations define the two expansions explicitly.

\smallskip
\noindent\emph{(i) Observation lists and Gaussian averaging.}
For $m\ge0$, let $h=s_0>s_1>\cdots>s_m>0$, set $x_0=x$, and
write $Z_m=((s_k,s_k,x_k))_{k=0}^m$.
The integration domain and its measure are
\[
 \begin{aligned}
 \mathcal H_m(h)&:=\{(s_1,\ldots,s_m):0<s_m<\cdots<s_1<h\}
                                      \times\mathbb T^m,\\
 dZ_m&:=ds_1\cdots ds_m\,dx_1\cdots dx_m.
 \end{aligned}
\]
For $m=0$ there is no history integration and $s_0=h$.
For $i=U,V$ use the observation vectors
\[
 \begin{aligned}
 \mathbf Y_q^{U,Z_m}
   &=\bigl((S_{s_k-q}U_{t_0+q})(x_k)\bigr)_{k=0}^m,\\
 \mathbf Y_q^{V,Z_m}
   &=\bigl((S_{s_k-q}U_{t_0+q})(x_k)
       -e^{-\Lambda(s_k-q)}(S_{s_k-q}D_{t_0+q})(x_k)\bigr)_{k=0}^m,
       \qquad 0\le q\le s_m.
 \end{aligned}
\]
The second is Definition~\ref{def:joint-backward-observations}
for this point list.
Use the endpoint Gaussian operators
\eqref{eq:joint-gaussian-averaging-operator} for this list. In particular,
\begin{equation}
 \mathbf q_u^{Z_m}(\xi)_k=p_{s_k-u}(x_k-\xi),\quad 0\le k\le m,
 \qquad \mathcal P^{Z_m}_{q,t}=\mathcal P^{Z_m,V}_{q,t},
 \quad 0\le q\le t\le s_m.
 \label{eq:joint-endpoint-gaussian-operators}
\end{equation}

\smallskip
\noindent\emph{(ii) Inserting one local coefficient discrepancy.}
Start with $F_0^U=F_0^V=a_\lambda$. When the next point
$(s_{m+1},x_{m+1})$ is inserted, put
\[
 \partial_{\mathrm{heat},m+1}
 :=\sum_{k=0}^m p_{s_k-s_{m+1}}(x_k-x_{m+1})\partial_{z_k}.
\]
For $\mathbf z=(z_0,\ldots,z_m)$ and the new state variable $w$,
define
\begin{equation}
 \begin{aligned}
 F_{m+1}^U(\mathbf z,w)
  :={}&\tfrac12[a(w)-\bar a_{s_{m+1}}^U(x_{m+1})]\\
     &\quad\times
       [\partial_{\mathrm{heat},m+1}^{\,2}
          \mathcal P^{Z_m,U}_{s_{m+1},s_m}F_m^U](\mathbf z),\\
 F_{m+1}^V(\mathbf z,w)
  :={}&\tfrac12[a_{n_{\mathrm{sm}}}(w)-\bar a_{s_{m+1}}^V(x_{m+1})]\\
     &\quad\times
       [\partial_{\mathrm{heat},m+1}^{\,2}
          \mathcal P^{Z_m,V}_{s_{m+1},s_m}F_m^V](\mathbf z).
 \end{aligned}
 \label{eq:joint-local-test-recursion}
\end{equation}
Each step therefore Gaussian-averages the preceding test,
differentiates twice in the displayed heat direction, and
multiplies by half the local coefficient discrepancy. Both
expansions preserve the original terminal test $a_\lambda$.

\smallskip
\noindent\emph{(iii) Initial averages, stochastic integrals and remainders.}
For $j=m+1\ge1$, define the two initial-past averages and the
two stochastic integrals explicitly by
\begin{equation}
 \begin{aligned}
 \mathsf P_{j,\lambda}^{U}
  &:=\int_{\mathcal H_m(h)}
       (\mathcal P^{Z_m,U}_{0,s_m}F_m^U)(\mathbf Y_0^{U,Z_m})\,dZ_m,\\
 \mathsf P_{j,\lambda}^{V,\mathrm{joint}}
  &:=\int_{\mathcal H_m(h)}
       (\mathcal P^{Z_m,V}_{0,s_m}F_m^V)(\mathbf Y_0^{V,Z_m})\,dZ_m,\\
 N_{j,\lambda}^{U}
  &:=\int_{\mathcal H_m(h)}\int_0^{s_m}
       \nabla(\mathcal P^{Z_m,U}_{q,s_m}F_m^U)(\mathbf Y_q^{U,Z_m})
                                  \cdot d\mathbf Y_q^{U,Z_m}\,dZ_m,\\
 N_{j,\lambda}^{V,\mathrm{joint}}
  &:=\int_{\mathcal H_m(h)}\int_0^{s_m}
       \nabla(\mathcal P^{Z_m,V}_{q,s_m}F_m^V)(\mathbf Y_q^{V,Z_m})
                                  \cdot d\mathbf Y_q^{V,Z_m}\,dZ_m.
 \end{aligned}
 \label{eq:joint-ordinary-output-formulas}
\end{equation}
For $m\ge1$, define the two remainders by
\begin{equation}
 \begin{aligned}
 R_{m,\lambda}^{U}
  &:=\int_{\mathcal H_m(h)}F_m^U(\mathbf Y_{s_m}^{U,Z_m})\,dZ_m,\\
 R_{m,\lambda}^{V,\mathrm{joint}}
  &:=\int_{\mathcal H_m(h)}F_m^V(\mathbf Y_{s_m}^{V,Z_m})\,dZ_m.
 \end{aligned}
 \label{eq:joint-remainder-definition}
\end{equation}
The gradient acts on the $m+1$ state variables.
All six quantities depend on $(t_0,h,x)$, suppressed here.
At $q=0$, $\mathbf Y_0^{i,Z_m}=(\phi_{s_k}^i(x_k))_{k=0}^m$;
the initial averages therefore use the same means as
\eqref{eq:predictable-mean-word-formula}. At the last insertion,
$(\mathbf Y_{s_m}^{V,Z_m})_m=V_{t_0+s_m}(x_m)$, so the newest
factor in the $V$ remainder is evaluated at the actual solution.
The superscript $V,\mathrm{joint}$ records evaluation in the
joint $V$ coordinates; the $U$ quantities use the ordinary
backward heat observations. The index $j=m+1$ on
$\mathsf P,N$ means exactly $m$ local insertions, whereas $R_m$
contains $m$ insertions without taking another It\^o step.

\smallskip
\noindent\emph{(iv) The extra covariance correction.}
With $\mathcal E_q$ from
Definition~\ref{def:joint-covariance-defect}, set
\[
 \mathfrak J_{m,\lambda,h}(t_0,x)
 :=\int_{\mathcal H_m(h)}\int_0^{s_m}
   (\mathcal E_q\mathcal P^{Z_m,V}_{q,s_m}F_m^V)
                         (\mathbf Y_q^{V,Z_m})\,dq\,dZ_m.
\]
This is the correction after $m$ local insertions. It is recorded
without further expansion, so each such term contains precisely
one covariance defect.
\par\smallskip
\noindent\emph{(v) Comparing the two constructions.}
For $j\ge1$, define the two comparison differences by
\begin{equation}
 \begin{aligned}
 \delta_{UV}\mathsf P_{j,\lambda}
   &:=\mathsf P_{j,\lambda}^U
                         -\mathsf P_{j,\lambda}^{V,\mathrm{joint}},\\
 \delta_{UV}N_{j,\lambda}
   &:=N_{j,\lambda}^U-N_{j,\lambda}^{V,\mathrm{joint}}.
 \end{aligned}
 \label{eq:joint-comparison-differences}
\end{equation}
The first compares the initial-past Gaussian averages, and the
second compares the stochastic integrals, at the same
$(t_0,h,x,\lambda)$. \end{definition}
Starting from the root $a_\lambda$ and expanding to depth $n\ge1$
gives the following separate $V$ identity, proved with
Theorem~\ref{thm:joint-finite-comparison}:
\begin{equation}
 \begin{aligned}
 a_\lambda(V_{t_0+h})
 ={}&\sum_{j=1}^{n}\bigl(\mathsf P_{j,\lambda}^{V,\mathrm{joint}}
                         +N_{j,\lambda}^{V,\mathrm{joint}}\bigr)
       +R_{n,\lambda}^{V,\mathrm{joint}}\\
 &+\sum_{m=0}^{n-1}\mathfrak J_{m,\lambda,h}.
 \end{aligned}
 \label{eq:joint-V-finite-expansion}
\end{equation}
The last sum records the extra joint covariance. It appears with
a minus sign when this identity is subtracted from the $U$ expansion.

\subsubsection{The finite comparison theorem}
Use the spatial norm
$\|f\|_{H_\Lambda^{-\rho}}=\|(\Lambda-\Delta)^{-\rho/2}f\|_{L^2(\mathbb T)}$,
$0<\rho<1/2$. The comparison sizes are
$\|U-V\|_{\mathscr L_0^2(H)}^2=\mathbb E\int_0^T\|D_t\|_{L^2(\mathbb T)}^2dt$ and
$\mathsf d_{\mathrm{coef}}(U,V)^2=\mathbb E\int_0^T\|f_t\|_{H_\Lambda^{-\rho}}^2dt
+\varepsilon_{\mathrm{app}}$, with deterministic
$\varepsilon_{\mathrm{app}}\ge0$. A window output has squared norm
\[
 \|Q\|^2=\int_0^{T-h}\mathbb E\|Q(t_0)\|_{H_\Lambda^{-\rho}}^2dt_0,
\]
as in Definition~\ref{def:window-integrated-norms}.
All terms below use the same terminal test $a_\lambda$ at both
endpoints; their $\lambda$ indices are displayed.
The estimates below
state the common-coefficient part; replacing the companion coefficient
by $a_{n_{\mathrm{sm}}}$ adds the separate error
$\operatorname{Err}_{n_{\mathrm{sm}},\lambda}$ from
\eqref{eq:bandwise-approximation-error}, the sum of the norms
of coefficient-test, inserted-factor and reference-profile corrections.
Its convergence is given by Lemma~\ref{lem:exact-bracket-cutoff-limits}.
For these joint comparisons, include the same coefficient-slot and
profile differences for the expansion terms in
\eqref{eq:joint-V-finite-expansion}, their complete first correction,
and the marks $\mathfrak J_{m,\lambda,h}$, $0\le m<N_\beta$.
Their actual integrators and the bracket defect
\eqref{eq:single-species-direct-defect} remain fixed.
The full index set is \eqref{eq:joint-approximation-word-index-set}.
This only enlarges the nonnegative error sum.

\begin{theorem}[Finite comparison with first-defect stopping]
\label{thm:joint-finite-comparison}
Assume the coefficient hypotheses of
Theorem~\ref{thm:main-weak-uniqueness}. More explicitly, let
$g:\mathbb R\to\mathbb R$ satisfy, at the effective exponent
$0<\beta<1$,
\[
 c_0\le |g(z)|\le C_0,\qquad
 |g(z)-g(z')|\le L_g|z-z'|^\beta,\qquad z,z'\in\mathbb R,
\]
for $c_0>0$ and finite $C_0,L_g$.
Use the residual variance $a=g^2-m_0^2$, the exact Haar coupling
and the self-consistent profiles above. Assume $D_0=0$ and $\Lambda\ge1$. Fix $N_\beta\ge2$, $0<\rho<1/2$,
$h=h_\lambda=\lambda^{-3\beta}$, and choose
$0<\chi<\min(1,\beta+2\rho)$ and $0<\chi_{\mathrm{rem}}<\beta$.
Let $p_{\mathrm{tail}}>0$ be a prescribed tail exponent and
$\mathcal L(\lambda,h)=1+\log(2+\lambda)+|\log h|$. Take dyadic $\lambda\ge2$ large enough that
$h=h_\lambda\le h_{N_\beta}\le h_0$, with the threshold
from \eqref{eq:finite-order-window-threshold}. Use the point observations and recursive operations of
Definition~\ref{def:joint-first-defect-marks}. For $1\le j\le N_\beta$, use the explicit differences
\eqref{eq:joint-comparison-differences} in that definition. Expand only the local coefficient discrepancies,
recording each extra covariance correction without expanding it.
The remainder $R_{N_\beta,\lambda}^{V,\mathrm{joint}}$ is the second
formula in \eqref{eq:joint-remainder-definition} with $m=N_\beta$.
Then the resulting \emph{finite identity} is
\begin{equation}
 a_\lambda(U_{t_0+h})-a_\lambda(V_{t_0+h})
  =\sum_{j=1}^{N_\beta}(\delta_{UV}\mathsf P_{j,\lambda}+\delta_{UV}N_{j,\lambda})
     +(R_{N_\beta,\lambda}^U-R_{N_\beta,\lambda}^{V,\mathrm{joint}})
     -\sum_{n=0}^{N_\beta-1}\mathfrak J_{n,\lambda,h}.
 \label{eq:joint-finite-stopped-identity}
\end{equation}
For $j\ge2$, use $w_j,q_{j,\chi},q_{j,\mathrm{prof}}$ from
\eqref{eq:higher-martingale-degrees} and $p_j$ from
\eqref{eq:all-order-predictable-clipped-degree}.
For the remainder, set $W_{N_\beta}:=\beta+3\beta^2N_\beta/4$.
These are precisely the Section~5 frequency degrees.

Apart from fixed logarithmic powers and arbitrarily small fixed losses,
the following bounds hold, with $2\le j\le N_\beta$ in both native and
predictable rows (a state power may lose $\varepsilon>0$ and
$\|U-V\|_{\mathscr L_0^2(H)}^\theta$ may become $\|U-V\|_{\mathscr L_0^2(H)}^{\theta-\varepsilon}$):
\begin{align*}
 \|\delta_{UV}N_{j,\lambda}\|&\le C_{N_\beta}\bigl\{
       \lambda^{-w_j}\mathsf d_{\mathrm{coef}}(U,V)
       +\lambda^{-q_{j,\chi}}\|U-V\|_{\mathscr L_0^2(H)}^\chi\\
 &\qquad+\lambda^{-q_{j,\mathrm{prof}}}\|U-V\|_{\mathscr L_0^2(H)}\bigr\}
       +C_{N_\beta}\lambda^{-p_{\mathrm{tail}} },\\
 \|\delta_{UV}\mathsf P_{j,\lambda}\|&\le
       C_{N_\beta}\Lambda^{-\rho/2}\lambda^{-p_j}\|U-V\|_{\mathscr L_0^2(H)},
       \qquad 2\le j\le N_\beta,\\
 \|R_{N_\beta,\lambda}^U-R_{N_\beta,\lambda}^{V,\mathrm{joint}}\|&\le
       C_{N_\beta}\Lambda^{-\rho/2}\lambda^{-(W_{N_\beta}-\chi_{\mathrm{rem}})}\|U-V\|_{\mathscr L_0^2(H)}^{\chi_{\mathrm{rem}}},\\
 \Big\|\sum_{n=0}^{N_\beta-1}\mathfrak J_{n,\lambda,h}\Big\|
       &\le C_{N_\beta}\lambda^{-\beta}\mathcal L(\lambda,h)^{C_{N_\beta}}
                         [1+\log(2+\Lambda h)]\mathsf d_{\mathrm{coef}}(U,V).
\end{align*}
The two pairing parameters have distinct purposes: $\chi$ belongs to
the native current-mean square integral, whereas $\chi_{\mathrm{rem}}$ pairs one
rough factor in the unmarked actual remainder.

The first Gaussian average and stochastic integral are treated separately.
They are regrouped into a terminal compatible Gaussian band and a
complete actual martingale as in \eqref{eq:joint-first-pair-regrouping};
the latter uses \eqref{eq:current-cap-complete-Gamma}. All estimates
hold at the fixed finite depth $N_\beta$.
\end{theorem}

\subsubsection{The first-defect boundary}
The mark $\mathfrak J_{n,\lambda,h}$ defined above contains $n$ local
insertions followed by the first covariance defect. The input norm
$\mathsf d_{\mathrm{coef}}(U,V)$ is the global signed coefficient norm from the common setup.
\begin{proposition}[Finite first-hit boundary]
\label{prop:all-order-first-hit-joint-boundary}
Fix $N_\beta$, $0<\beta<1$, $0<\rho<1/2$ and $0\le n\le N_\beta$.
At each dyadic block $a_\lambda$, use
$h=h_\lambda=\lambda^{-3\beta}\le h_{N_\beta}$.
The marks following $n$ ordinary local insertions satisfy
\begin{equation}
 \|\mathfrak J_{n,\lambda,h}\|_{\mathscr L_h^2(H_\Lambda^{-\rho})}
 \le C_{N_\beta}\lambda^{-\beta}[1+\log(2+\lambda^4h)]^{C_{N_\beta}}
 [1+\log(2+\Lambda h)]\mathsf d_{\mathrm{coef}}(U,V).
 \label{eq:all-order-first-hit-joint-boundary}
\end{equation}
This includes the true initial triangle. With the cutoff
\eqref{eq:gaussian-band-dyadic-cutoff}, summing over
$\lambda\ge\lambda_*\asymp\Lambda^{1/4}$ gives at most
$C_{N_\beta}(1+\log\Lambda)^{C_{N_\beta}}
\Lambda^{-\beta/4}\mathsf d_{\mathrm{coef}}(U,V)=o\mathsf d_{\mathrm{coef}}(U,V)$.
\end{proposition}
\subsubsection{The first predictable--martingale pair}
Take $m=0$ in \eqref{eq:joint-ordinary-output-formulas}.
The compatible references of
Definition~\ref{def:coupled-gaussian-references} give
\[
 \mathsf P_{1,\lambda}^U=\mathbb E_{t_0}[a_\lambda(X_h^U)],\qquad
 \mathsf P_{1,\lambda}^{V,\mathrm{joint}}
                     =\mathbb E_{t_0}[a_\lambda(X_h^V)].
\]
Define the complete first martingale corrections by
\begin{equation}
 \begin{aligned}
 \Gamma_{U,\lambda}^{\mathrm{joint}}
  &:=N_{1,\lambda}^U-\bigl[a_\lambda(X_h^U)-\mathsf P_{1,\lambda}^U\bigr],\\
 \Gamma_{V,\lambda}^{\mathrm{joint}}
  &:=N_{1,\lambda}^{V,\mathrm{joint}}
      -\bigl[a_\lambda(X_h^V)-\mathsf P_{1,\lambda}^{V,\mathrm{joint}}\bigr].
 \end{aligned}
 \label{eq:joint-first-complete-corrections}
\end{equation}
Thus $\Gamma_{i,\lambda}^{\mathrm{joint}}$ is the actual stochastic term
minus its centered Gaussian-reference value. As in
\eqref{eq:reference-first-doob-martingale}, its time-dependent version
subtracts the reference Doob martingale from the truncated actual
stochastic integral; hence the displayed quantity is a terminal
martingale value.

\begin{proposition}[The first predictable--martingale pair]
\label{prop:joint-first-pair}
At the fixed block $a_\lambda$, take
$h=h_\lambda=\lambda^{-3\beta}\le h_0$. For the compatible references
and block-dependent complete martingales just defined,
\begin{equation}
 \delta_{UV}\mathsf P_{1,\lambda}+\delta_{UV}N_{1,\lambda}
 =a_\lambda(X_h^U)-a_\lambda(X_h^V)
       +\Gamma_{U,\lambda}^{\mathrm{joint}}-\Gamma_{V,\lambda}^{\mathrm{joint}}.
 \label{eq:joint-first-pair-regrouping}
\end{equation}
Set $\alpha=\beta+2\rho$ and $\mathcal L(\lambda,h)=1+\log(2+\lambda)+|\log h|$.
The complete joint martingale satisfies
\[
 \begin{aligned}
 \|\Gamma_{U,\lambda}^{\mathrm{joint}}-\Gamma_{V,\lambda}^{\mathrm{joint}}\|
 &\le C\mathcal L(\lambda,h)^C\bigl[\lambda^{-\beta}\mathsf d_{\mathrm{coef}}(U,V)\\
 &\quad+\lambda^{-\alpha}\min(1,\lambda \|U-V\|_{\mathscr L_0^2(H)})
 +\lambda^{-\alpha+3\beta(1-\beta)/4}\|U-V\|_{\mathscr L_0^2(H)}\bigr]
 +\operatorname{Err}_{n_{\mathrm{sm}},\lambda}.
 \end{aligned}
\]
The terminal Gaussian difference is controlled by
\eqref{eq:sharp-main-Gaussian-band} and its dyadic sum
\eqref{eq:sharp-main-Gaussian-sum}, with the same fixed durations and
cutoff as here.
\end{proposition}

\subsection{Proofs and supporting stopped-history estimates}
\label{sec:joint-proof-details}
\begin{remark}[The general parameter form used in the proofs]
\label{rem:joint-general-parameters}
The main statements use $h_\lambda=\lambda^{-3\beta}$ and
$\lambda_*\asymp\Lambda^{1/4}$. For the proofs and the parameter
calculation in Section~\ref{sec:all-beta-joint-complete-closure}, we
also record the general form. Take $h_\lambda=\lambda^{-\gamma_{\mathrm{time}}}$,
$0<\gamma_{\mathrm{time}}<4$, and the cutoff $\lambda_*\asymp\Lambda^r$, $r>0$, of
\eqref{eq:gaussian-band-general-cutoff}. With $\alpha=\beta+2\rho$,
for $j\ge2$ use the following degree symbols in these general calculations:
\begin{equation}
 \begin{aligned}
 w_j&=\beta+\gamma_{\mathrm{time}}(j-1)\beta/4,&p_j&=w_j-\gamma_{\mathrm{time}}/4,\\
 q_{j,\chi}&=w_j+\min(2\rho-\chi,\gamma_{\mathrm{time}}(2\rho-\chi)/4),&
 q_{j,\mathrm{prof}}&=w_j+\gamma_{\mathrm{time}}\rho/2-\gamma_{\mathrm{time}}/4,\\
 q_\Gamma&=\alpha-\gamma_{\mathrm{time}}(1-\beta)/4,&
 W_{N_\beta}&=\beta+\gamma_{\mathrm{time}}\beta N_\beta/4.
 \end{aligned}
 \label{eq:joint-general-degrees}
\end{equation}
The estimates in Theorem~\ref{thm:joint-finite-comparison} preserve the
same form with these degrees. The last linear term in
Proposition~\ref{prop:joint-first-pair} becomes
$\lambda^{-q_\Gamma}\|U-V\|_{\mathscr L_0^2(H)}$, and the summed boundary term
in Proposition~\ref{prop:all-order-first-hit-joint-boundary} has
factor $\Lambda^{-r\beta}$.
Setting $\gamma_{\mathrm{time}}=3\beta$ and $r=1/4$ gives the main statements.

\end{remark}

The three main proofs come first. Their inputs are then established
in the supporting statements: the exact It\^o identity, the stopped
heat-row bounds, the first-defect time integral and the transfer
of the Section~5 estimates. The It\^o calculation, time integration,
and transfer verification are proved here; the longer stopped-row
verifications are in Appendix~\ref{app:stopped-heat-details}. In the boundary argument,
the new defect time is integrated before the old centered gaps.

\subsubsection{The finite comparison}
\nopagebreak[4]
\begin{proof}[Proof of Theorem~\ref{thm:joint-finite-comparison}]
Fix $t_0,h,x,a_\lambda$. We prove the finite identity, then bound
its terms.

\medskip
\noindent\emph{1. One It\^o step and its four outputs.}\par
Use the tests and integrals of
Definition~\ref{def:joint-first-defect-marks}. For each fixed
history $Z_m$, apply the backward It\^o formula
\eqref{eq:single-species-joint-one-step} to $F_m^V$ and integrate
over $\mathcal H_m(h)$. Let
$\mathcal V_0=a_\lambda(V_{t_0+h}(x))$ and
$\mathcal V_m=R_{m,\lambda}^{V,\mathrm{joint}}$ for $m\ge1$.
The initial value and stochastic integral are exactly
\eqref{eq:joint-ordinary-output-formulas}. The local covariance
term is the integral of the $V$ test in
\eqref{eq:joint-local-test-recursion}, hence is $\mathcal V_{m+1}$.
The extra covariance term is
$\mathfrak J_{m,\lambda,h}$. Consequently
\begin{equation}
 \mathcal V_m=\mathsf P_{m+1,\lambda}^{V,\mathrm{joint}}
       +N_{m+1,\lambda}^{V,\mathrm{joint}}
       +\mathcal V_{m+1}+\mathfrak J_{m,\lambda,h}.
 \label{eq:joint-four-output-recursion}
\end{equation}
Iterate for $m=0,\ldots,N_\beta-1$, expanding only $\mathcal V_{m+1}$
and recording each $\mathfrak J_m$ once. This gives
\eqref{eq:joint-V-finite-expansion}. The ordinary $U$ identity is
Lemma~\ref{lem:finite-expansion}, with local coefficient $a$ and
no extra covariance term. Subtraction gives
\eqref{eq:joint-finite-stopped-identity}, including the minus sign
in front of the marks.

\medskip
\noindent\emph{2. The ordinary terms preserve their comparison bounds.}\par
Lemma~\ref{lem:joint-coordinate-transfer} verifies that the local
factors, centering identities, spatial inputs and profile bounds
give the same kernels as in the Section~5 proofs. Its concrete
verification yields the following three estimates:
\begin{itemize}[leftmargin=*]
\item The higher stochastic terms use
\eqref{eq:higher-martingale-frequency-combined}, with the degrees
$w_j,q_{j,\chi},q_{j,\mathrm{prof}}$ of
\eqref{eq:joint-general-degrees}.
\item The initial averages use
\eqref{eq:all-order-predictable-integrated-sobolev}. For $j-1$
insertions its frequency factor is
$\lambda^{-\beta}h^{(j-1)\beta/4-1/4}=\lambda^{-p_j}$.
\item The unexpanded remainder uses
\eqref{eq:paired-actual-remainder}. Its frequency factor is
$\lambda^{-\beta+\chi_{\mathrm{rem}}}h^{N_\beta\beta/4}
=\lambda^{-(W_{N_\beta}-\chi_{\mathrm{rem}})}$.
\end{itemize}
Here $h=\lambda^{-\gamma_{\mathrm{time}}}$ until the final
substitution. The logarithms, small losses and truncation tail are
exactly those of the cited comparisons. The first stochastic term
is treated together with its average by
Proposition~\ref{prop:joint-first-pair}.

\medskip
\noindent\emph{3. The stopped terms and the limits.}\par
Proposition~\ref{prop:all-order-first-hit-joint-boundary} bounds each
of the finitely many marks; the triangle inequality gives the last
row of the theorem. The recursion is first justified with smooth
coefficient tests, finite noise series and strictly positive gaps.
Lemma~\ref{lem:exact-bracket-cutoff-limits} removes these restrictions
with the physical bracket fixed, as specified in that lemma.
The coefficient replacements contribute the additive
$\operatorname{Err}_{n_{\mathrm{sm}},\lambda}$ defined before the
theorem. Setting $\gamma_{\mathrm{time}}=3\beta$ in the degree
formulas now gives all the displayed exponents.
\end{proof}

\subsubsection{The first-defect boundary}
\nopagebreak[4]
\begin{proof}[Proof of Proposition~\ref{prop:all-order-first-hit-joint-boundary}]
Fix $0\le n\le N_\beta$, $a_\lambda$ and
$\lambda^{-4}\le h\le h_{N_\beta}$.
The input is $f_t=a(U_t)-a_{n_{\mathrm{sm}}}(V_t)$, with global
negative norm at most $\mathsf d_{\mathrm{coef}}(U,V)$.

\medskip
\noindent\emph{1. Express all endpoint ages in one new time variable.}\par
Fix the old insertion gaps
$d_i=s_{i-1}-s_i>0$, $i=1,\ldots,n$, where $s_0=h$.
At the defect time $q<s_n$, put
$\tau=s_n-q$ and $b_i=s_i-s_n$ for $i=0,\ldots,n$.
Thus the age of endpoint $i$ is $s_i-q=\tau+b_i$ and $b_n=0$.
For $n=0$ this means $s_n=h$, $b_0=0$, with no old gaps.

Apply Corollary~\ref{cor:first-defect-endpoint-pairs} after the
spatial heat integrations. Apart from the root and old-factor
bounds, each ordered endpoint pair $(i,j)$ contributes the positive
time kernel
\begin{equation}
 \begin{aligned}
 K_{ij}(\tau;\mathbf b):={}&
 r(\tau+b_i)(\tau+b_i)^{-1/2}(\tau+b_j)^{-1/2}\\
 &\times[1+\Lambda^{-\rho/2}(\tau+b_i)^{-\rho/2}],
 \qquad r(u)=\min(1,\Lambda u).
 \end{aligned}
 \label{eq:first-defect-positive-time-kernel}
\end{equation}
Lemma~\ref{lem:first-defect-time-integration}, first with singularity
exponent $0$ and then with $\rho/2$, gives uniformly in the offsets
\begin{equation}
 \int_0^{s_n}K_{ij}(\tau;\mathbf b)\,d\tau
 \le C_\rho[1+\log(2+\Lambda h)].
 \label{eq:first-defect-integrated-pair-kernel}
\end{equation}
The second part uses $\Lambda^{-\rho/2}\Lambda^{\rho/2}=1$;
the defect factor $1-e^{-\Lambda(s_i-q)}$ ensures integrability.

\medskip
\noindent\emph{2. Integrate the old centered differences.}\par
The interval-centering identity of
Lemma~\ref{lem:interval-increment-centering} and the assignments of
Definition~\ref{def:centered-word-assignments} apply to the $n$ old
derivatives. Write $I_0$ for the intervals assigned to the root,
and $\mathscr B_{\mathrm{rough}}$ for the nonempty blocks assigned
to rough coefficient factors. Set
$c_\lambda(I_0)=\min(1,\lambda\min_{i\in I_0}d_i^{1/4})$
when $I_0\ne\varnothing$, and $c_\lambda(\varnothing)=1$.
For each fixed assignment their positive bound is
\begin{equation}
 C_n\lambda^{-\beta}c_\lambda(I_0)
    \prod_{I\in\mathscr B_{\mathrm{rough}}}
                 \min_{i\in I}d_i^{\beta/4}
    \prod_{i=1}^n\frac1{d_i}.
 \label{eq:first-defect-old-gap-bound}
\end{equation}
Unassigned coefficient factors are bounded. The independent
Gaussian fraction used for the new covariance derivative leaves
these old centered differences unchanged; fixed Gaussian moments
are controlled by \eqref{eq:multipoint-score-products}.

Each rough block integrates to $C_{k,\beta}h^{\beta/4}$ by
\eqref{eq:minimum-gap-cube-bound}. For a root block of size $k$,
bound its clip by
$\prod_{i\in I_0}\min(1,\lambda d_i^{1/4})^{1/k}$ and use
\[
 \int_0^h\min(1,\lambda u^{1/4})^{1/k}\frac{du}{u}
 \le4k+\log\max(1,\lambda^4h).
\]
Enlarge the old-gap simplex to these cubes. Since $h\le1$,
\eqref{eq:first-defect-old-gap-bound} integrates to
$C_n\lambda^{-\beta}[1+\log(2+\lambda^4h)]^{C_n}$.

\medskip
\noindent\emph{3. Take the global norm and sum the blocks.}\par
For the global norm, put $t=t_0+h$. The physical input time is
$t_0+q=t-\tau-\sum_{i=1}^n d_i$.
Extend $f$ by zero outside $[0,T]$. Its translated restriction has
probability--time $L^2(H_\Lambda^{-\rho})$ norm no greater than
that of $f$. Minkowski's inequality thus bounds the output by
the positive kernel mass times
$\mathsf d_{\mathrm{coef}}(U,V)$. Tonelli permits the new-time
integration \eqref{eq:first-defect-integrated-pair-kernel} before
the old-gap integration. Summing the finite endpoint pairs and
assignments proves \eqref{eq:all-order-first-hit-joint-boundary}.

On an initial window $[0,t]$, $t<h_\lambda$, use its actual duration
in the same gap variables and enlarge the positive integration
domain to the one with duration $h_\lambda$.
The zero extension and the terminal-time integration just used
give the same bound; no lower bound $t\ge\lambda^{-4}$ is needed
for these positive kernel estimates. This is the initial-triangle
norm of Definition~\ref{def:window-integrated-norms}.

Finally take $h=\lambda^{-\gamma_{\mathrm{time}}}$ and
$\lambda\ge\lambda_*\asymp\Lambda^r$. The dyadic series of
$\lambda^{-\beta}$ times fixed logarithmic powers is bounded by
$C(1+\log\Lambda)^C\Lambda^{-r\beta}$: write
$\lambda=2^k\lambda_*$ and sum the convergent series
$\sum_{k\ge0}2^{-k\beta}(1+k)^C$.
For $r=1/4$ and $\gamma_{\mathrm{time}}=3\beta$ this is the
summed bound in the proposition, and tends to zero as
$\Lambda\to\infty$.
\end{proof}

\subsubsection{The first predictable--martingale pair}
\nopagebreak[4]
\begin{proof}[Proof of Proposition~\ref{prop:joint-first-pair}]
Fix $a_\lambda$ and $\lambda^{-4}\le h\le h_0$. Put
$\tau=h-s$, $\vartheta_s=e^{-\Lambda\tau}$ and
$Y_s^{U,h}:=S_\tau U_{t_0+s}$; use $Y_s^{V,h}$ from
\eqref{eq:joint-scalar-backward-coordinate}. Their residual-noise
operators, evaluated at $t_0+s$, satisfy
\begin{equation}
 \begin{aligned}
 J_s^U&:=L_U(s),\qquad
 J_s^V:=(1-\vartheta_s)L_U(s)+\vartheta_sL_V(s),\\
 J_s^U-J_s^V&=\vartheta_sB_\Delta(s),\qquad
 Y_s^{U,h}-Y_s^{V,h}=\vartheta_sS_\tau D_{t_0+s}.
 \end{aligned}
 \label{eq:joint-first-comparison-inputs}
\end{equation}
The noise formula follows from
\eqref{eq:joint-scalar-drift-cancellation}; convexity gives
$J_s^i(J_s^i)^*\le CI$.

Use the Gaussian backward fields $\mathcal G_s^i$ of
\eqref{eq:gaussian-backward-field-martingale}, the remaining variances
\eqref{eq:complete-martingale-remaining-variance}, and the smoothed
tests $F_s^i$ from Definition~\ref{def:coupled-backward-martingales}.
The stochastic-integral identity
\eqref{eq:complete-martingale-representation} holds with
$\widehat U_s^i$ replaced by $Y_s^{i,h}$ and
$\mathcal B_i(s)$ by $J_s^i-L(\bar a_s^i)$:
subtract the same reference Doob integral from the actual integral
\eqref{eq:joint-ordinary-output-formulas} at $m=0$.
Adding the initial averages proves
\eqref{eq:joint-first-pair-regrouping}.

For the estimate, apply the exact integrand split
\eqref{eq:complete-integrand-telescoping} with actual fields $Y_s^{i,h}$
and noise amplitudes $J_s^i$.
Its three bounds have the following inputs.

\emph{Noise amplitudes.} The first difference in
\eqref{eq:joint-first-comparison-inputs} is the signed input
$B_\Delta$ multiplied by $0\le\vartheta_s\le1$.
The signed trace \eqref{eq:haar-negative-negative-trace} and
annular integration \eqref{eq:complete-annular-time-integral}
therefore give \eqref{eq:complete-coefficient-contribution}.

\emph{Evaluation points.} The second difference in
\eqref{eq:joint-first-comparison-inputs} has global spatial $L^2$
norm at most $\|U-V\|_{\mathscr L_0^2(H)}$, by heat contraction
and time translation. The reference difference still satisfies
\eqref{eq:actual-Gamma-Gaussian-moments}. Since the actual companion
$J^V$ has bounded covariance, the output trace
\eqref{eq:complete-output-trace} gives exactly the clipped bound
\eqref{eq:current-cap-root-row}.

\emph{Reference profiles.} The reference fields and variances are
unchanged, so their inputs
\eqref{eq:complete-profile-input}--\eqref{eq:complete-remaining-variance-input}
and time integrals \eqref{eq:complete-profile-time-integrals}
give \eqref{eq:complete-reference-contribution}.

Sum these bounds, add the coefficient error from
Lemma~\ref{lem:exact-bracket-cutoff-limits}, and set $h=\lambda^{-3\beta}$.
For a general duration the profile degree is $q_\Gamma$ in
\eqref{eq:joint-general-degrees}. The terminal Gaussian difference
uses \eqref{eq:sharp-main-Gaussian-band} and
\eqref{eq:sharp-main-Gaussian-sum}.
\end{proof}

\subsubsection{The covariance defect and the one-step identity}
The covariance calculation uses the profile $\bar a^V$.
The terminal references $X^U,X^V$ use their respective profiles,
with their centered difference controlled by
\eqref{eq:self-profile-noise-stability}.
Outside the common window every row carries its actual indicator
$\mathbf1_{t<\sigma_i}$; the stopped-history estimates below use
these masks. Heat targets $T_i$ and observation times $\sigma_i$ are
kept separate throughout.

\begin{lemma}[Covariance defect and backward It\^o identity]
\label{lem:joint-backward-ito}
Use the joint coordinates and defect matrix of
\eqref{eq:joint-bracket-defect-definition}. Fix a terminal time
$0<s<\min_i\sigma_i$, so all coordinates are active on $[0,s]$.
At time $t\le s$, abbreviate $q_i=q_t^{Z,i}(\cdot)\in H$,
$L_U=L_U(t_0+t)$, $L_V=L_V(t_0+t)$, and $B_\Delta=L_U-L_V$.
Set $\vartheta_i(t)=e^{-\Lambda(\sigma_i-t)}$ and
$f=a(U_{t_0+t})-a_{n_{\mathrm{sm}}}(V_{t_0+t})$.
With $\mathcal A=L_UB_\Delta^*-B_\Delta L_U^*$, the defect also has
the following signed-coefficient decomposition:
\begin{align}
 E_{ij}={}&\left(1-\tfrac{\vartheta_i(t)+\vartheta_j(t)}{2}\right)\langle q_i,\mathsf M_fq_j\rangle\nonumber\\
 &-\left(\tfrac{\vartheta_i(t)+\vartheta_j(t)}{2}-\vartheta_i(t)\vartheta_j(t)\right)
                        \langle q_i,B_\Delta B_\Delta^*q_j\rangle\nonumber\\
 &+\tfrac{\vartheta_i(t)-\vartheta_j(t)}{2}\langle q_i,\mathcal A q_j\rangle.
 \label{eq:single-species-three-defects}
\end{align}
Here $E_{ij}$ denotes $E_{ij}(t)$.
For a smooth bounded test $F:\mathbb R^{d_{\mathrm{obs}}}\to\mathbb R$
with bounded derivatives, put $\mathbf Y_t=\mathbf Y_t^{V,Z}$.
Let $\mathcal P^Z$ use the $V$ reference profile, and use the defect
operator $\mathcal E_t$ defined after
\eqref{eq:joint-bracket-defect-definition}. Then
\begin{align}
 F(\mathbf Y_s)={}&\mathcal P ^Z_{0,s}F(\mathbf Y_0)
       +\int_0^s\nabla \mathcal P ^Z_{t,s}F(\mathbf Y_t)\cdot d\mathbf Y_t\nonumber\\
 &+\tfrac12\int_0^s\int_{\mathbb T}[a_{n_{\mathrm{sm}}}(V_{t_0+t}(\xi))-\bar a_t^V(\xi)]
        (\mathbf q_t^Z(\xi)\cdot\nabla)^2\mathcal P ^Z_{t,s}F(\mathbf Y_t)\,d\xi\,dt
       +\int_0^s\mathcal E_t\mathcal P ^Z_{t,s}F(\mathbf Y_t)\,dt.
 \label{eq:single-species-joint-one-step}
\end{align}
\end{lemma}

\begin{proof}[Proof of Lemma~\ref{lem:joint-backward-ito}]
\emph{1. Compute the covariance correction.}
At the fixed physical time $t_0+t$, let $A=L_U$, $B=B_\Delta$,
so $L_V=A-B$. The actual residual row of observation $i$ is
$q_i(A-\vartheta_iB)$, with
$\vartheta_i=e^{-\Lambda(\sigma_i-t)}$.
Before pairing with $q_i,q_j$, the difference from the local
$V$ covariance is exactly
\[
 \begin{aligned}
 &(A-\vartheta_iB)(A-\vartheta_jB)^*-(A-B)(A-B)^*\\
 &\qquad=(1-\vartheta_j)AB^*+(1-\vartheta_i)BA^*
                              +(\vartheta_i\vartheta_j-1)BB^*.
 \end{aligned}
\]
The factor identity \eqref{eq:covariance-factor-identity} says
$AA^*-(A-B)(A-B)^*=\mathsf M_f$, where
$\mathsf M_f u=fu$ is multiplication by
$f=a(U_{t_0+t})-a_{n_{\mathrm{sm}}}(V_{t_0+t})$.
It follows that $AB^*+BA^*=\mathsf M_f+BB^*$, while
$AB^*-BA^*=\mathcal A$.
Substituting these symmetric and skew parts in the last display
and pairing with $q_i,q_j$ gives
\eqref{eq:single-species-three-defects}. Equivalently, expansion
with $L_V$ gives the three terms in
\eqref{eq:single-species-direct-defect}. The additive covariance
is the same in the two rows and cancels.

\emph{2. Apply the time-dependent It\^o formula.}
Let $u(t,y)=(\mathcal P^Z_{t,s}F)(y)$ and define the reference
covariance density by
\[
 C_{ij}^{Z,V}(t)=\int_{\mathbb T}
       [m_0^2+\bar a_t^V(\xi)]q_t^{Z,i}(\xi)q_t^{Z,j}(\xi)\,d\xi.
\]
Differentiating the Gaussian convolution in
Definition~\ref{def:joint-backward-observations} gives its backward
equation
$\partial_tu+\frac12\sum_{i,j}C_{ij}^{Z,V}(t)\partial_{ij}u=0$,
with $u(s,\cdot)=F$. This is
\eqref{eq:reference-gaussian-backward-equation} for the present list.
Each coordinate of $\mathbf Y$ has zero drift by
Lemma~\ref{lem:joint-backward-drift-cancellation}, so It\^o's formula
takes the explicit form
\[
 du(t,\mathbf Y_t)=\nabla u(t,\mathbf Y_t)\cdot d\mathbf Y_t
 +\frac12\sum_{i,j}
 \left[\frac{d\langle Y^i,Y^j\rangle_t}{dt}
                   -C_{ij}^{Z,V}(t)\right]
                    \partial_{ij}u(t,\mathbf Y_t)\,dt.
\]
By \eqref{eq:joint-bracket-defect-definition}, the matrix in
brackets is
\[
 \int_{\mathbb T}
 [a_{n_{\mathrm{sm}}}(V_{t_0+t}(\xi))-\bar a_t^V(\xi)]
                 q_t^{Z,i}(\xi)q_t^{Z,j}(\xi)\,d\xi+E_{ij}(t).
\]
The first term produces
$\frac12[a_{n_{\mathrm{sm}}}(V)-\bar a^V]
(\mathbf q_t^Z(\xi)\cdot\nabla)^2u$; the second is
$\mathcal E_tu$. Integration from $0$ to $s$ proves
\eqref{eq:single-species-joint-one-step}.
For the recursion one also uses $s=\min_i\sigma_i$.
This follows by approaching that time from below, first for the
smooth, positive-gap tests, using heat-square integrability and
then the limits of Lemma~\ref{lem:exact-bracket-cutoff-limits}.
\end{proof}
The four terms in the identity are the initial predictable value,
the martingale, the next local coefficient insertion, and the
bracket defect. Appending $(t,t,\xi)$ expands the third term.
The fourth is recorded as a first-defect mark, ending that branch.
For a stopped history, each row is extended by zero after its
own $\sigma_i$. Lemma~\ref{lem:exact-bracket-cutoff-limits} justifies
smoothing and stochastic Fubini with these physical brackets.

\subsubsection{Stopped heat rows and Gaussian derivative weights}
\label{sec:stopped-heat-contraction-joint-defect}
\label{sec:comparable-stopped-heat-span-score}
\label{sec:all-age-stopped-point-score-frame}

Three estimates supply the inputs to the main proofs: damping
controls the covariance range, heat localization controls a changing
noise or mean, and age grouping combines observations at different
scales. Inverses always mean inverses on the
actual covariance range; constants may depend on the finite number
$n$ of observations. We use $\Sigma^{-1/2}$ as shorthand for
$(\Sigma^+)^{1/2}$, the inverse square root on that range.

We first compare complete stopped heat histories. The contraction
below then gives a covariance bound for their current rows.

For a list $(\sigma_i,T_i,x_i)_{i=1}^n$ with
$q<\sigma_i\le T_i$ and scalars $z_1,\ldots,z_n$, define, for $u\ge q$,
\[
 F(u)=\sum_{i=1}^n z_i\mathbf1_{u<\sigma_i}p_{T_i-u}(x_i-\cdot),\quad
 G(u)=\sum_{i=1}^n z_i\mathbf1_{u<\sigma_i}
       e^{-\Lambda(\sigma_i-u)}p_{T_i-u}(x_i-\cdot).
\]

Define the stopped heat row
$R(u,x)_i=\mathbf1_{u<\sigma_i}p_{T_i-u}(x_i-x)$ and its Gram matrix
$\mathsf G=\int_q^\infty\int_{\mathbb T}R(u,x)R(u,x)^*\,dx\,du$.
Here $q$ is the current local time, and the list
$Z=((\sigma_i,T_i,x_i))_{i=1}^n$ is fixed throughout this subsection.
For the $V$ reference, the future observation covariance is explicitly
\[
 \Sigma=\int_q^h\int_{\mathbb T}
       [m_0^2+\bar a_u^V(\xi)]R(u,\xi)R(u,\xi)^*\,d\xi\,du.
\]
For a reserved Gaussian fraction use the corresponding fixed
fraction of this covariance. The same arguments apply to either
profile and to uniformly elliptic interpolations of the profiles.
Thus $c\mathsf G\le\Sigma\le C\mathsf G$. The symbols
$\Sigma,\mathsf G,Q$ below all depend on this list and on $q$;
their history superscripts are suppressed only in this subsection.
Regard $Q:H\to\mathbb R^n$ as the operator with rows
$Q_i=p_{T_i-q}(x_i-\cdot)$, and put
$Q_{\mathrm{def},i}=(1-e^{-\Lambda(\sigma_i-q)})Q_i$.
Define the diagonal damping matrix
$\mathsf{Damp}_\Lambda(q)=\operatorname{diag}(\exp(-\Lambda(\sigma_i-q)))_{i=1}^n$.
Thus $\mathsf{Damp}_\Lambda(q)Q$ multiplies the $i$th heat row by its damping weight.
For any age $\tau\ge0$, set $r(\tau)=\min(1,\Lambda\tau)$.

\begin{lemma}[Damping stopped heat histories]
\label{lem:stopped-history-contraction}
Damping the stopped heat histories is an $L^2$ contraction:
\begin{equation}
 \|G\|_{L^2([q,\infty)\times\mathbb T)}\le\|F\|_{L^2([q,\infty)\times\mathbb T)}.
 \label{eq:stopped-heat-damping-contraction}
\end{equation}
This controls the current heat rows by the future observation covariance.
If all observation and target ages are comparable to $w$, then
\begin{equation}
 QQ^*\le C_nw^{-1}\Sigma,\qquad
 Q_{\mathrm{def}}Q_{\mathrm{def}} ^*\le C_nr_w^2w^{-1}\Sigma,\quad r_w=\min(1,\Lambda w).
 \label{eq:comparable-stopped-current-ranges}
\end{equation}
For arbitrary ages, with $\delta=\min_i(\sigma_i-q)$, the damped
rows satisfy
\[
 (\mathsf{Damp}_\Lambda(q)Q)(\mathsf{Damp}_\Lambda(q)Q)^*
 \le C\delta^{-1}\Sigma.
\]
\end{lemma}
The proof is given in Appendix~\ref{app:stopped-damping}.

Two forms of heat localization are useful below. A translated
projection preserves cancellation in a signed noise input; a positive
physical envelope carries a mean input through spatial integration.

\begin{lemma}[Heat localization for signed noise and mean inputs]
\label{lem:comparable-heat-span-Fourier}
Suppose the observation ages $\sigma_i-q$ and target ages $T_i-q$
are all comparable to $w>0$, with fixed comparison constants.
Let $E=\operatorname{span}\{Q_i:i=1,\ldots,n\}$ be the span of the
current heat rows. Let $\mathsf{Proj}_{E,x}$ be the orthogonal projection onto the
space obtained by translating every center in $E$ by $x\in\mathbb T$.
For sufficiently small $w$ and a bounded operator $B:H\to H$,
the translated projection satisfies
\begin{equation}
 \int_{\mathbb T}\|\mathsf{Proj}_{E,x}B\|_{\mathrm{HS}}^2\,dx
       \le C_n\sqrt w\,\|S_{c_nw}B\|_{\mathrm{HS}}^2.
 \label{eq:comparable-span-input-heat}
\end{equation}
For arbitrary ages, set
$\delta=\min_i(\sigma_i-q)>0$; the damped linear Gaussian weight obeys
\begin{equation}
 |\Sigma^{-1/2}\mathsf{Damp}_\Lambda(q)Q(\xi)|
 \le C_n\delta^{-1/2}\sum_{i=1}^n(T_i-q)^{1/4}
                      p_{C_n(T_i-q)}(x_i-\xi),
 \label{eq:damped-first-score-physical-envelope}
\end{equation}
and preserves a genuine right factor $S_{\delta/2}$.
\end{lemma}
The proof is given in Appendix~\ref{app:stopped-localization}.

Point observations can have very different ages. We group comparable
ages so that each new derivative weight acts on one controlled covariance block.
The lemma gives the two outputs used later: a signed covariance bound
and an integrated mean-energy bound.

Fix an observation time $q$ and $n$ genuine point observations,
so their observation and target times coincide. Write
$t_i=\sigma_i-q=T_i-q$ for their ages, with $0<t_i\le h_n$,
$i=1,\ldots,n$, using \eqref{eq:finite-order-window-threshold}. Order the ages
and join adjacent ones whose ratio is at most a fixed
$\mathfrak r_{\mathrm{sep}}>16n^2$. The resulting groups are disjoint
nonempty subsets $g\subset\{1,\ldots,n\}$ partitioning the observation
indices. Set $\tau_g=\min_{i\in g}t_i$; the rows $Q_g$, histories
$F_g$ and Grams $\mathsf G_g$ are restricted to the indices in $g$.
Here $Q_g$ consists of the rows $p_{t_i}(x_i-\cdot)$, $i\in g$;
$F_g(u)=\sum_{i\in g}z_i\mathbf1_{u<t_i}p_{t_i-u}(x_i-\cdot)$
for arbitrary coefficients $z_i$. Here $u$ is elapsed time after $q$;
thus $\mathsf G_g$ is the matrix satisfying
$z^*\mathsf G_gz=\int_0^\infty\|F_g(u)\|_{L^2(\mathbb T)}^2du$.
Write $\Sigma_g$ for the corresponding covariance block,
$\Sigma_g\ge c\mathsf G_g$, and
$E_g=\operatorname{span}\{p_{t_i}(x_i-\cdot):i\in g\}$.
The diagonal matrix $\mathsf{Damp}_{\Lambda,g}(q)$ is the restriction of
$\mathsf{Damp}_\Lambda(q)$ to the indices in $g$.
A sum over $i\in g$ ranges over the observations in that group.

The derivative weights in the next lemma are the explicit
Gaussian integration-by-parts factors from
\eqref{eq:gaussian-score-definitions}. For a centered Gaussian
vector $\zeta$ of covariance $\Sigma_g$ and deterministic directions
$v,w$ in its range, the first and second weights are, respectively,
\[
 v^*\Sigma_g^+\zeta,\qquad
 (v^*\Sigma_g^+\zeta)(w^*\Sigma_g^+\zeta)-v^*\Sigma_g^+w.
\]
Their $L^p$ bounds are $C_p|\Sigma_g^{-1/2}v|$ and
$C_p|\Sigma_g^{-1/2}v||\Sigma_g^{-1/2}w|$.
These formulas explain what is being bounded: a changed mean
supplies one direction, and a changed covariance supplies two
noise-response directions. Lemma~\ref{lem:all-age-score-frame}
preserves one such physical input in its spatial norm.

\begin{lemma}[Gaussian derivative bounds for observations of different ages]
\label{lem:all-age-score-frame}
\emph{A signed covariance input.} For $\mathcal B =L(a(U))-L(a_{n_{\mathrm{sm}}}(V))$ and a bounded companion $L$, a covariance
derivative weight with a nonnegative amplitude bound $\mathsf A_0$ for the
remaining factors, signed row $Q_{\mathrm{def},g}\mathcal B $ and row $Q_{g'}L$
has global output bound
\begin{equation}
 C_n\mathsf A_0 r(\tau _g)\tau _g^{-1/2}\tau _{g'}^{-1/2}
               [1+\Lambda^{-\rho/2}\tau _g^{-\rho/2}]\mathsf d_{\mathrm{coef}}(U,V).
 \label{eq:all-age-one-input-score}
\end{equation}
For equal group ages $\tau_g=\tau_{g'}=w$, this reads
\begin{equation}
 C_n\mathsf A_0r_w[w^{-1}+\Lambda^{-\rho/2}w^{-1-\rho/2}]\mathsf d_{\mathrm{coef}}(U,V).
 \label{eq:comparable-one-input-covariance-score}
\end{equation}
A second bounded $\mathcal B$ may be used as the companion $L$.

\emph{A mean input.} For one comparable group of scale $\tau$, let $x\in\mathbb T$
translate all observation centers together, with their relative centers
fixed. A mean input $d\in H$ satisfies
\begin{equation}
 \int_{\mathbb T}|\Sigma_g^{-1/2}\mathsf{Damp}_{\Lambda,g}(q)Q_gd|^2\,dx\le C_n\tau ^{-1/2}\|S_{c_n\tau }d\|_{L^2(\mathbb T)}^2,
 \label{eq:joint-native-mean-one-input}
\end{equation}

The corresponding positive pointwise envelope is
\begin{equation}
 |\Sigma_g^{-1/2}\mathsf{Damp}_{\Lambda,g}(q)Q_gd|\le C_n\tau ^{-1/4}
                   \sum_{i\in g}(p_{C_nt_i}*|d|)(x_i).
 \label{eq:joint-native-mean-positive-row}
\end{equation}
Thus the spatial norm of a linear Gaussian weight over ages at least $s$
is bounded by $C_ns^{-1/4}\|d\|_{L^2(\mathbb T)}$. All estimates permit joint positive integration
of the relative centers, and use only the eventual global second
moment of the distinguished input.
\end{lemma}
The proof is given in Appendix~\ref{app:stopped-age-groups}.

\subsubsection{Integrating the first defect}
\label{sec:all-order-first-hit-joint-mark}

Ordinary local insertions starting from $a_\lambda(V_{t_0+h}(x))$
produce the full point list $h=s_0>\cdots>s_n>q$.
Every local insertion appends its physical point and time. The
full list is kept when a centered difference becomes constant
in one coordinate. Old centerings use their reserved Gaussian
fractions and subtracted values at the old past; a separate future fraction
still observes the full point list. The new derivative weight is normalized
on the covariance range of that list.

\begin{corollary}[The first-defect endpoint-pair bound]
\label{cor:first-defect-endpoint-pairs}
For the point list $h=s_0>\cdots>s_n>q$ above, put
$f_t=a(U_t)-a_{n_{\mathrm{sm}}}(V_t)$ and $\tau_i=s_i-q$.
Let $\mathsf A_0$ bound the remaining root and old coefficient
factors in the all-age frame. Then the whole covariance defect has
the positive endpoint-pair majorant
\begin{equation}
 C_n\mathsf A_0\sum_{i,j=0}^n r(\tau_i)\tau_i^{-1/2}\tau_j^{-1/2}
 [1+\Lambda^{-\rho/2}\tau_i^{-\rho/2}]\|f_{t_0+q}\|_{H_\Lambda^{-\rho}},
 \quad \tau_i=s_i-q.
 \label{eq:first-hit-pair-age-kernel}
\end{equation}
Here $r(\tau)=\min(1,\Lambda\tau)$, and the spatial norm is taken
after the one-input integration with the relative centers in their
original joint heat-tree measure.
\end{corollary}
\begin{proof}[Proof of Corollary~\ref{cor:first-defect-endpoint-pairs}]
Condition on $\mathcal F_{t_0+q}$, so $L_V$ and $B_\Delta$ at
time $t_0+q$ are fixed for the future Gaussian average.
Formula~\eqref{eq:single-species-direct-defect} consists of two
cross covariances and one covariance quadratic in $B_\Delta$.
Apply \eqref{eq:all-age-one-input-score} to their second Gaussian
derivative weights, using the separate future covariance fraction
specified above. For the first cross covariance, the signed row
is $(1-e^{-\Lambda\tau_i})Q_iB_\Delta$ and the companion is
$Q_jL_V$. Its time factor is
\[
 r(\tau_i)\tau_i^{-1/2}\tau_j^{-1/2}
                       [1+\Lambda^{-\rho/2}\tau_i^{-\rho/2}].
\]
For the second cross covariance interchange $i,j$; the sum over
all ordered pairs in \eqref{eq:first-hit-pair-age-kernel} already
contains this term. In the quadratic covariance use its $i$th
$B_\Delta$ as the signed input and the other as the companion.
The latter is bounded since $\|B_\Delta\|\le\|L_U\|+\|L_V\|\le C$,
and its additional damping factor is at most one.
Thus this term has the same majorant. The age groups of
Lemma~\ref{lem:all-age-score-frame} have internally comparable
ages, so their group bounds give the displayed endpoint bounds
with constants depending only on the observation count.
Restoring the common amplitude $\mathsf A_0$ and the signed
spatial input $\|f_{t_0+q}\|_{H_\Lambda^{-\rho}}$ proves
\eqref{eq:first-hit-pair-age-kernel}. The spatial averaging here
is exactly the positive heat-kernel integration in that lemma;
the global probability--time norm is taken later in the proof of
Proposition~\ref{prop:all-order-first-hit-joint-boundary}.
\end{proof}
The next lemma integrates its new time before any of the old gaps.

\begin{lemma}[New-time integration before old centered gaps]
\label{lem:first-defect-time-integration}
Fix $\Lambda\ge1$, $0<h\le1$ and $r(u)=\min(1,\Lambda u)$.
Let $H_{\mathrm{int}}\in[0,h]$ be an integration length,
$b_i,b_j\in[0,h]$ endpoint offsets, and $s\in[0,1/2)$ a singularity exponent.
Then
\begin{equation}
 \int_0^{H_{\mathrm{int}}} r(\tau+b_i)(\tau+b_i)^{-1/2-s}
              (\tau+b_j)^{-1/2}d\tau
 \le\begin{cases}C[1+\log(2+\Lambda h)],&s=0,\\
 C_s\Lambda^s,&s>0.
 \end{cases}
 \label{eq:first-hit-new-time-integral}
\end{equation}
\end{lemma}
\begin{proof}[Proof of Lemma~\ref{lem:first-defect-time-integration}]
\emph{The equal-age bound.}
Write $r(u)=\min(1,\Lambda u)$. Splitting at $u=\Lambda^{-1}$ gives
\[
 \int_0^{2h}r(u)u^{-1-s}\,du
 \le\begin{cases}
 C[1+\log(2+\Lambda h)],&s=0,\\
 C_s\Lambda^s,&0<s<1/2.
 \end{cases}
\]
For $s>0$, the two pieces can be bounded explicitly by extending
the upper limit to infinity:
\[
 \Lambda\int_0^{\Lambda^{-1}}u^{-s}\,du
       +\int_{\Lambda^{-1}}^\infty u^{-1-s}\,du
 =\Lambda^s\left(\frac1{1-s}+\frac1s\right).
\]
At $s=0$, the first piece is at most one and the second is
logarithmic up to the actual upper limit $2h$.

\emph{Unequal endpoint ages.}
If $b_i\le b_j$, the integrand in
\eqref{eq:first-hit-new-time-integral} is at most
$r(\tau+b_i)(\tau+b_i)^{-1-s}$, so the preceding bound applies.
If $b_i>b_j$, put $u=\tau+b_j$ and $b=b_i-b_j>0$.
On $u\ge b$, $r(u+b)\le2r(u)$ and the same bound applies.
On $u<b$, use $r(u+b)\le2r(b)$ and integrate $u^{-1/2}$:
more explicitly, this contribution is at most
\[
 2r(b)b^{-1/2-s}\int_0^b u^{-1/2}\,du
       =4r(b)b^{-s}.
\]
For $s=0$ this is at most $4$. For $s>0$, if $b\le\Lambda^{-1}$,
then $r(b)b^{-s}=\Lambda b^{1-s}\le\Lambda^s$; if
$b\ge\Lambda^{-1}$, then $r(b)b^{-s}=b^{-s}\le\Lambda^s$.
This proves the estimate with the damping factor attached to its
original signed endpoint.
\end{proof}

\subsubsection{Transferring the higher-order comparison estimates}
\label{sec:joint-native-and-remainder-transfer}
\label{sec:predictable-damped-mean-history-frame}
The following verification is shared by the higher-order stochastic integrals,
the Gaussian average terms and the paired remainder. It is independent
of the separate first-order regrouping in Proposition~\ref{prop:joint-first-pair}.

\paragraph{\textbf{The comparison inputs after the coordinate change.}}
The Section~5 proofs apply with the same centering identities and time
integrals when the changed physical observation rows satisfy the
same spatial input estimates. At a current window time $q$,
observation $i$ is recorded at time $\sigma_i$ and heat-smoothed
to its target time $T_i\ge\sigma_i>q$, at center $x_i$.
Its noise-response row is
$Q_if=(S_{T_i-q}f)(x_i)$ for $f\in H$.
For one terminal observation the drift cancellation is
\eqref{eq:joint-scalar-drift-cancellation}.
With $L_U=L(a(U_{t_0+q}))$, $L_V=L(a_{n_{\mathrm{sm}}}(V_{t_0+q}))$
and $B_\Delta=L_U-L_V$, the changed noise row is exactly
\[
 Q_iL_U-Q_i\bigl[L_V+(1-e^{-\Lambda(\sigma_i-q)})B_\Delta\bigr]
                    =e^{-\Lambda(\sigma_i-q)}Q_iB_\Delta.
\]
The lemma below checks these noise and mean inputs, the coefficient
moments and the profile changes directly against the kernels in the
Section~5 proofs. The coordinate change also produces a
quadratic-covariation defect, given by
\eqref{eq:single-species-direct-defect}. That term is recorded at
its first occurrence and estimated separately in
Proposition~\ref{prop:all-order-first-hit-joint-boundary}.

\begin{lemma}[Transfer to the joint backward coordinates]
\label{lem:joint-coordinate-transfer}
For the coupled solutions of Proposition~\ref{prop:construct-companion}
and the profiles of Definition~\ref{def:coupled-gaussian-references},
the joint backward coordinates
obey the higher-martingale, Gaussian-average and paired-remainder
estimates \eqref{eq:higher-martingale-frequency-combined},
\eqref{eq:ordinary-predictable-comparison-row} and
\eqref{eq:paired-actual-remainder}, respectively. Use the joint
outputs of \eqref{eq:joint-ordinary-output-formulas} and
\eqref{eq:joint-remainder-definition}, signed input
$\mathsf d_{\mathrm{coef}}(U,V)$ and solution input
$\|U-V\|_{\mathscr L_0^2(H)}$, with the parameter ranges of the cited
estimates and $h\le h_{N_\beta}$. The frequency degrees, small losses
and separate approximation errors are unchanged.
\end{lemma}
\begin{proof}[Proof of Lemma~\ref{lem:joint-coordinate-transfer}]
We check the centering, spatial inputs and profile changes used in
the three comparison proofs. Fix $a_\lambda$
and $Z=((s_i,s_i,x_i))_{i=0}^n$, with
$h=s_0>\cdots>s_n>q$, $x_0=x$.
Set $\delta=s_n-q$, $w=h-q$ and $\tau_i=s_i-q\in[\delta,w]$.
Constants are uniform in $t_0,\Lambda$.

\medskip
\noindent\emph{1. The local factors, centerings and their moments.}\par
The tests \eqref{eq:joint-local-test-recursion} have the root,
factors and insertion rule of \eqref{eq:predictable-mean-test-recursion};
the separately stopped defect leaves their assignments and centerings
unchanged. Work with common coefficient $a$; restoring $a_{n_{\mathrm{sm}}}$
in the $V$ slots adds the specified replacement error.

Put $d=D_{t_0}$ and
$\eta_q=D_{t_0+q}-e^{-\Lambda q}S_qd$. The proof of
Lemma~\ref{lem:consistent-current-interpolation} gives conditional
pointwise $L^p$ bounds $C_pq^{1/4}$ for both
$U_{t_0+q}-S_qU_{t_0}$ and $\eta_q$.
Use the same interpolation convention as
Definition~\ref{def:predictable-mean-functional}: $\theta=0$ at $V$
and $\theta=1$ at $U$. Its old-past mean is $\phi_t^\theta$ from that
definition. The mean continued from $q$ to a target $q<t\le h$ is
$\phi_t^{\theta,q}:=S_{t-q}U_{t_0+q}
-(1-\theta)e^{-\Lambda(t-q)}S_{t-q}D_{t_0+q}$.
Their exact difference is
\begin{equation}
 \phi_t^{\theta,q}-\phi_t^\theta
 =S_{t-q}(U_{t_0+q}-S_qU_{t_0})
       -(1-\theta)e^{-\Lambda(t-q)}S_{t-q}\eta_q.
 \label{eq:joint-native-consistent-rough-moments}
\end{equation}
Heat contraction, the two moment bounds and the future Gaussian
variance bound $C\sqrt{t-q}$ give an $L^p$ displacement of size
$C_pt^{1/4}$ about $\phi_t^\theta$.
H\"older continuity of $a$ and
$|\bar a_t^\theta-a(\phi_t^\theta)|\le Ct^{\beta/4}$ from
\eqref{eq:self-profile-displacement} therefore give the required
$C_pt^{\beta/4}$ moment for each factor with its own subtracted
value. This includes $q^{\beta/4}$ for the newest actual factor.

Use fixed elliptic fractions for the old interval weights and
a separate one for the new derivative. Their centering and
moments \eqref{eq:multipoint-score-products} are unchanged,
giving the centering and moment bounds used in
Lemma~\ref{lem:centered-word-time-kernel}.

\medskip
\noindent\emph{2. The changed noise and mean preserve their spatial norms.}\par
Set $\vartheta_i=e^{-\Lambda\tau_i}$ and let
$Q_i f=(S_{\tau_i}f)(x_i)$, so $Q:H\to\mathbb R^{n+1}$ is
the vector of current heat rows. Direct subtraction gives
\begin{align}
 Y_q^{U,Z,i}-Y_q^{V,Z,i}
   &=\vartheta_i(S_{\tau_i}D_{t_0+q})(x_i),\nonumber\\
 Q_iL_U-Q_i[L_V+(1-\vartheta_i)B_\Delta]
   &=\vartheta_i Q_iB_\Delta.
 \label{eq:joint-native-exact-differences}
\end{align}
Write $\mathsf{Damp}_\Lambda(q)=\operatorname{diag}(\vartheta_i)$.
The unchanged residual and additive companions are finite sums of
\begin{equation}
 QL_V+(Q-\mathsf{Damp}_\Lambda(q)Q)B_\Delta,\qquad m_0Q.
 \label{eq:joint-native-companion-split}
\end{equation}
Their right-hand operators are uniformly bounded. The stopped
history and localization bounds
\eqref{eq:stopped-heat-damping-contraction}--\eqref{eq:damped-first-score-physical-envelope}
apply to each of $Q$, $\mathsf{Damp}_\Lambda(q)Q$ and their
difference before these companions are summed.

The required right heat factor is explicit: for every $i$,
\[
 \vartheta_iS_{\tau_i}
       =(\vartheta_iS_{\tau_i-\delta/2})S_{\delta/2},\qquad
 (1-\vartheta_i)S_{\tau_i}
       =((1-\vartheta_i)S_{\tau_i-\delta/2})S_{\delta/2}.
\]
Normalize the left rows using the reserved future covariance
after $q+\delta/2$, as in
Lemma~\ref{lem:comparable-heat-span-Fourier}. Their observation ages
are still at least $\delta/2$ and their target ages are comparable
to the original ages. Its positive envelope and the signed
spatial transport of Proposition~\ref{prop:unified-dyadic-column-transport}
therefore give the same fixed-history kernel
\eqref{eq:native-integrator-exact-factor}. The signed trace therefore gives exactly
\eqref{eq:higher-noise-time-factor}, with input $f_{t_0+q}$ in
$H_\Lambda^{-\rho}$. Keeping that input signed until transport,
use the same old-gap integral and final It\^o isometry as in the
higher-martingale proof.

For the changed mean, \eqref{eq:joint-native-mean-one-input} and
\eqref{eq:joint-native-mean-positive-row} preserve $D_{t_0+q}$ in
spatial $L^2$. Jensen, unit heat mass and time translation give
the global norm $\|U-V\|_{\mathscr L_0^2(H)}$.
For an old-past direction, use $d=D_{t_0}$ with its damping factor.
These give the noise and mean inputs to
\eqref{eq:higher-noise-duration-bound} and
\eqref{eq:higher-current-duration-bound}, with the same distinct inner
gap $\delta$ and outer age $w$.

\medskip
\noindent\emph{3. Profiles, the last actual factor and the limits.}\par
The mean direction along the old-past interpolation is
$e^{-\Lambda t}S_td$. Thus \eqref{eq:all-order-profile-path-bounds} gives the profile
source $Ct^{(\beta-1)/4}\|d\|_{L^2(\mathbb T)}$.
Use it directly for subtracted values and under
\eqref{eq:multipoint-weighted-input-L2} for covariance changes,
placing the new derivative on the separate fraction from step~1.
For the paired remainder, the alternative bound in
\eqref{eq:self-profile-lp} gives
$\|\bar a_t^U-\bar a_t^V\|_{L^2(\mathbb T)}
\le C\|d\|_{L^2(\mathbb T)}^\beta$.
After integration over probability and window starts, these give
the linear profile input in \eqref{eq:higher-profile-duration-bound}
and the $\beta$-power profile input of the paired-remainder proof.

At the last insertion $q$, the factor
$a_{n_{\mathrm{sm}}}(V_{t_0+q})-\bar a_q^V$ is measurable at
$\mathcal F_{t_0+q}$ and remains outside the auxiliary future
Gaussian expectation. Its conditional $q^{\beta/4}$ moment was
checked in step~1, up to the separate coefficient error.
Compare current means and their own subtracted values first,
then the future covariance at the fixed $V$ mean. The changed
mean is the first line of
\eqref{eq:joint-native-exact-differences}; the two profile bounds
above control the remaining changes. Thus the root-slot and
rough-factor estimates
\eqref{eq:paired-remainder-root-slot}--\eqref{eq:paired-remainder-rough-block}
and their time integrations prove
\eqref{eq:paired-actual-remainder} with
$\chi=\chi_{\mathrm{rem}}$. Branches already stopped at a defect
are absent from this remainder.

On an initial window, $t_0=0$ and $d=D_0=0$, so the old-past
mean and common-profile differences vanish. Keep the actual ages
and use the initial-triangle integrations of
Corollary~\ref{cor:initial-window-comparisons}.
Finally Lemma~\ref{lem:exact-bracket-cutoff-limits} removes the
finite series and positive-gap restrictions, preserving the actual
brackets. Its coefficient-slot and profile corrections are the
additive $\operatorname{Err}_{n_{\mathrm{sm}},\lambda}$ and are
kept separate from both comparison distances.
The Gaussian-average proof uses the same old-past means and tests,
so \eqref{eq:all-order-predictable-integrated-sobolev} applies directly.
For the stochastic terms the verified rows give the same kernels
\eqref{eq:higher-noise-time-factor} and
\eqref{eq:capped-native-mean-integrand}, and the same profile integral
\eqref{eq:higher-profile-duration-bound}. Their time integrations and
the already checked remainder comparison prove all three assertions.
\end{proof}

\section{Closing the comparison and proving weak uniqueness}
\label{sec:uniform-small-beta-direct-ledger}
\label{sec:all-beta-joint-complete-closure}

The uniqueness argument has one quantitative requirement: for
each prescribed solution, its controlled smooth companion must
approach it with vanishing control energy. Once this is known,
Proposition~\ref{prop:law-comparison} compares both prescribed
solutions with the same uncontrolled smooth law.
We give that complete uniqueness proof in
Section~\ref{sec:closure-local-target}.

The remaining subsections establish its quantitative input.
Section~\ref{sec:closure-coefficient-inequality} sums the
Section~6 comparisons; Section~\ref{sec:closure-parameters-energy}
obtains the energy rate $\Lambda^{-5/2}$; and
Section~\ref{sec:closure-approximation-laws} chooses the common
approximation sequence. Section~\ref{sec:conditional-global} then applies standard conditional
restart to obtain uniqueness on every finite horizon.

\subsection{The uniqueness argument}
\label{sec:closure-local-target}
Assume the coefficient hypotheses of
Theorem~\ref{thm:main-weak-uniqueness}: the scalar diffusion $g$ is
bounded, uniformly nonzero and globally $\beta_0$-H\"older for
some $\beta_0>0$. Fix
$0<\beta\le\min(\beta_0,4/25)$.
The same coefficient is globally $\beta$-H\"older: use its original
H\"older bound at state distances at most one and boundedness at
larger distances. The smaller exponent is used only to choose
parameters in the comparison proof.

\begin{theorem}[Local weak uniqueness]
\label{thm:all-positive-beta-joint}
Under the hypotheses of Theorem~\ref{thm:main-weak-uniqueness},
any two given weak solutions from the same deterministic
$H=L^2(\mathbb T)$ state have the same state-path law on every
horizon $0<T\le1$. A sufficient expansion depth is
$N_\beta=\lceil6/\beta^2\rceil$ for any effective
$0<\beta\le\min(\beta_0,4/25)$.
\end{theorem}

Recall the coupling notation from \eqref{eq:framework-residual}:
$0<m_0<c_0$, where $c_0$ is the lower bound for $|g|$,
$a=g^2-m_0^2$ is the residual variance, and
$a_{n_{\mathrm{sm}}}$ is its state mollification.
The proof below uses only the conclusion of
Lemma~\ref{lem:common-approximation-control}: the companions can
be chosen with both their path error and their control energy
tending to zero. Their construction and decay estimates are
verified in the three subsections that follow the proof.

\begin{proof}[Proof of Theorem~\ref{thm:all-positive-beta-joint}]
Fix $0<T\le1$ and two weak solutions $U^{(1)},U^{(2)}$ from
the same deterministic initial state. For any coupled pair use
$\|U-V\|_{\mathscr L_0^2(H)}^2
=\mathbb E\int_0^T\|U_t-V_t\|_{L^2(\mathbb T)}^2dt$.
Lemma~\ref{lem:common-approximation-control} supplies one
deterministic smoothing index $n_{\mathrm{sm}}(\Lambda)$ and
controlled companions $V^{(1)},V^{(2)}$ such that, for $i=1,2$,
\[
 \|U^{(i)}-V^{(i)}\|_{\mathscr L_0^2(H)}\longrightarrow0,
 \qquad
 \mathcal C_\Lambda^{(i)}:=m_0^{-2}\Lambda^2
          \|U^{(i)}-V^{(i)}\|_{\mathscr L_0^2(H)}^2\longrightarrow0.
\]
The first quantity measures the error between the paths; the
second is the energy needed to remove the companion's control
drift. The common smoothing index makes the uncontrolled smooth
law the same for both comparisons.

Let $\nu_{n_{\mathrm{sm}}(\Lambda)}$ be the canonical law on
$\mathcal X_T=C([0,T];H)$ of the uncontrolled smooth equation
with diffusion $\sqrt{m_0^2+a_{n_{\mathrm{sm}}(\Lambda)}}$ and
the common initial state. This is the law in
Proposition~\ref{prop:law-comparison}; it is independent of the
space used to construct either companion.
Recall the inclusion
$\iota:\mathcal X_T\to\mathcal Y_T=L^2([0,T];H)$.
For any test $F$ on $\mathcal Y_T$ bounded by one and
one-Lipschitz for the truncated norm metric,
\eqref{eq:companion-law-bound} gives, for $i=1,2$,
\[
 \left|\mathbb EF(U^{(i)})
       -\nu_{n_{\mathrm{sm}}(\Lambda)}(F\circ\iota)\right|
 \le\|U^{(i)}-V^{(i)}\|_{\mathscr L_0^2(H)}
       +2\sqrt{\mathcal C_\Lambda^{(i)}}
       +2\mathcal C_\Lambda^{(i)}\longrightarrow0.
\]
Subtract the two comparisons. Their middle reference law is
identical, so the triangle inequality gives
$\mathbb EF(U^{(1)})=\mathbb EF(U^{(2)})$.
Bounded Lipschitz tests determine probability laws on the Polish
space $\mathcal Y_T$; hence the integrated path laws agree.
Lemma~\ref{lem:integrated-laws} then identifies their laws on
$\mathcal X_T$, proving local weak uniqueness.
\end{proof}

The local theorem holds for every deterministic initial state.
Section~\ref{sec:conditional-global} justifies its use after
conditioning on a state history and extends the conclusion to
arbitrary finite horizons.

\subsection{Summing the expansion into one coefficient inequality}
\label{sec:closure-coefficient-inequality}

We now establish the quantitative estimate used in the uniqueness
proof. There are two inputs: the signed-energy bound from
Section~5 and the coefficient comparison obtained by summing
Section~6. This subsection collects those inequalities;
Section~\ref{sec:closure-parameters-energy} combines them.

Fix a horizon $0<T\le1$. For each prescribed weak solution $U$,
take $\Lambda\ge1$, $n_{\mathrm{sm}}\ge1$ and let
$V=V^{n_{\mathrm{sm}},\Lambda}$ be its controlled companion
from Proposition~\ref{prop:construct-companion}, on the same noise
pair and with $V_0=U_0$. Thus the difference $D=U-V$ satisfies
\eqref{eq:companion-difference}.
For $0<\rho<1/2$, use the spatial norm
$\|f\|_{H_\Lambda^{-\rho}}
=\|(\Lambda-\Delta)^{-\rho/2}f\|_{L^2(\mathbb T)}$.
Write $a=g^2-m_0^2$ and let $a_{n_{\mathrm{sm}}}$ be its state
mollification. Let $\varepsilon_{\mathrm{app}}\ge0$ be the
deterministic allowance for regularization tails in
\eqref{eq:coefficient-discrepancy-norm}.
The two comparison norms are
\begin{equation}
 \begin{aligned}
 \|U-V\|_{\mathscr L_0^2(H)}^2
   &:=\mathbb E\int_0^T\|U_t-V_t\|_{L^2(\mathbb T)}^2\,dt,\\
 \mathsf d_{\mathrm{coef}}(U,V)^2
   &:=\mathbb E\int_0^T
       \|a(U_t)-a_{n_{\mathrm{sm}}}(V_t)\|_{H_\Lambda^{-\rho}}^2\,dt
                  +\varepsilon_{\mathrm{app}}.
 \end{aligned}
 \label{eq:closure-comparison-norms}
\end{equation}
Lemma~\ref{lem:damped-energy-low-frequency}(i)
and the bounded-noise isometry give
\begin{equation}
 \begin{aligned}
 \|U-V\|_{\mathscr L_0^2(H)}
   &\le C\Lambda^{\rho/2-1/4}\mathsf d_{\mathrm{coef}}(U,V),\\
 \|U-V\|_{\mathscr L_0^2(H)}^2&\le C\Lambda^{-1/2}.
 \end{aligned}
 \label{eq:joint-closure-energy}
\end{equation}
The first estimate converts coefficient bounds into solution
energy bounds. The second is the preliminary estimate, which
will be improved to $C\Lambda^{-5/2}$ after the coefficient
limit. This suffices because
$\Lambda^2\Lambda^{-5/2}=\Lambda^{-1/2}\to0$.

We use the same cutoff and durations as the main statements of
Sections~4--6, for dyadic state frequencies
$\lambda\in\mathscr S_{\mathrm z}=\{2^k:k=0,1,\ldots\}$:
\[
 \lambda_*=2^{\lceil\frac14\log_2\Lambda\rceil}
      \asymp\Lambda^{1/4},\qquad h_\lambda=\lambda^{-3\beta}.
\]
The negative spatial exponent $\rho$ and the two comparison
powers $\chi,\chi_{\mathrm{rem}}$ will be chosen in
Lemma~\ref{lem:parameter-energy-closure}. For now require
$\alpha:=\beta+2\rho<1$,
$0<\chi<\alpha$ and $0<\chi_{\mathrm{rem}}<\beta$.

\paragraph{\textbf{The frequency degrees used in the sum.}}
Only the lowest higher-order degrees are needed. In the notation
of Theorem~\ref{thm:joint-finite-comparison} and
\eqref{eq:joint-general-degrees}, with
$\gamma_{\mathrm{time}}=3\beta$, they are
\begin{align*}
 w_2&=\beta+\tfrac34\beta^2,&
 p_2&=w_2-\tfrac34\beta,\\
 q_{2,\chi}&=w_2+\min(2\rho-\chi,\tfrac34\beta(2\rho-\chi)),&
 q_{2,\mathrm{prof}}&=w_2-\tfrac34\beta(1-2\rho),\\
 q_\Gamma&=\beta+2\rho-\tfrac34\beta(1-\beta),&
 W_{N_\beta}&=\beta+\tfrac34\beta^2N_\beta.
\end{align*}
At each further order, $w_j,p_j,q_{j,\chi},q_{j,\mathrm{prof}}$
increase by $3\beta^2/4$; in particular
$p_j=p_2+3(j-2)\beta^2/4$ for $j\ge2$.
The degrees $w_2,q_{2,\chi},q_{2,\mathrm{prof}}$ belong to the
three higher-martingale contributions; $p_j$ belongs to the
Gaussian averages, $q_\Gamma$ to the first complete martingale,
and $W_{N_\beta}$ to the actual remainder.

\paragraph{\textbf{Linear coefficients and approximation errors.}}
Let $\mathfrak L_{N_\beta,\Lambda}\ge1$ be a fixed power of
$1+\log\Lambda$, collecting the logarithms in the first-martingale
noise and first-defect estimates. Define
\begin{align}
 L_{\mathrm{coef}}:={}&\lambda_*^{-w_2}
       +\mathfrak L_{N_\beta,\Lambda}\lambda_*^{-\beta},\nonumber\\
 L_{\mathrm{dist}}:={}&1+\Lambda^{-\rho/2}\lambda_*^{1-\beta}
 +\lambda_*^{-q_\Gamma}+\lambda_*^{-q_{2,\mathrm{prof}}}
 +\Lambda^{-\rho/2}\sum_{j=2}^{N_\beta+1}\lambda_*^{-p_j}.
 \label{eq:all-beta-joint-linear-slots}
\end{align}
These are deterministic coefficients of, respectively,
$\mathsf d_{\mathrm{coef}}(U,V)$ and
$\|U-V\|_{\mathscr L_0^2(H)}$.
The term $\Lambda^{-\rho/2}\lambda_*^{1-\beta}$ is the
low-frequency estimate of Lemma~\ref{lem:damped-energy-low-frequency}(ii).
The remaining decaying terms are the profile and Gaussian-average
bounds. The common upper bound includes the additional positive order
$j=N_\beta+1$ used in the common upper bound; the actual expansion
ends at order $N_\beta$.

At each fixed $\Lambda$, let $\lambda_{\mathrm{cut}}$ be the
upper state-frequency cutoff and $\eta>0$ the lower
gap cutoff, as in Lemma~\ref{lem:exact-bracket-cutoff-limits}.
Define the summed approximation error by
\begin{equation}
 \operatorname{Err}_{n_{\mathrm{sm}},\Lambda}
 :=C_{\mathrm{app}}
      \sum_{\substack{\lambda\in\mathscr S_{\mathrm z}\\
                             \lambda\le\lambda_{\mathrm{cut}}}}
                  \operatorname{Err}_{n_{\mathrm{sm}},\lambda}
           +\sqrt{\varepsilon_{\mathrm{app}}}.
 \label{eq:summed-approximation-error}
\end{equation}
Each summand is the bandwise error defined by
\eqref{eq:bandwise-approximation-error}.
The deterministic constant $C_{\mathrm{app}}$ absorbs the finite
summation coefficients and depends on the fixed parameters and
cutoffs, but is independent of $n_{\mathrm{sm}}$ and of the
prescribed solution pair. Section~\ref{sec:closure-approximation-laws}
chooses these cutoffs and the common smoothing index.

\begin{proposition}[Master coefficient-comparison inequality]
\label{prop:master-coefficient-comparison}
Under the coefficient assumptions of
Theorem~\ref{thm:main-weak-uniqueness}, use the coupled pair and
norms in \eqref{eq:closure-comparison-norms}, the cutoff and
durations above, and the degrees just specified. Take $\Lambda$
sufficiently large that the high-frequency windows satisfy
Theorem~\ref{thm:joint-finite-comparison}. Suppose
$w_2,q_\Gamma,q_{2,\mathrm{prof}},p_2,q_{2,\chi}$ and
$W_{N_\beta}-\chi_{\mathrm{rem}}$ are positive. Then
\begin{equation}
 \begin{aligned}
 \mathsf d_{\mathrm{coef}}(U,V)\le C\big[&L_{\mathrm{coef}}\mathsf d_{\mathrm{coef}}(U,V)
          +L_{\mathrm{dist}}\|U-V\|_{\mathscr L_0^2(H)}\\
       &+\|U-V\|_{\mathscr L_0^2(H)}^\alpha
          +\lambda_*^{-q_{2,\chi}}\|U-V\|_{\mathscr L_0^2(H)}^\chi\\
       &+\Lambda^{-\rho/2}\lambda_*^{-(W_{N_\beta}-\chi_{\mathrm{rem}})}
                    \|U-V\|_{\mathscr L_0^2(H)}^{\chi_{\mathrm{rem}}}\big]\\
       &+\operatorname{Err}_{n_{\mathrm{sm}},\Lambda}
       +\lambda_*^{-p_{\mathrm{tail}}}.
 \end{aligned}
 \label{eq:all-beta-joint-master}
\end{equation}
Here $L_{\mathrm{coef}},L_{\mathrm{dist}}$ and the summed error are
\eqref{eq:all-beta-joint-linear-slots} and
\eqref{eq:summed-approximation-error}.
The prescribed exponent $p_{\mathrm{tail}}>0$ controls the
separate truncation tail. Fixed logarithmic factors and
arbitrarily small power losses from the Section~6 comparisons
are suppressed in this display; they are restored in the proof
of Lemma~\ref{lem:parameter-energy-closure}.
\end{proposition}

\begin{proof}[Proof of Proposition~\ref{prop:master-coefficient-comparison}]
\emph{Split the state-frequency bands and use the complete time norm.}
The Littlewood--Paley decomposition of $a$ and the triangle
inequality split $a(U)-a(V)$ into the bands below $\lambda_*$
and those at or above it. Replacing the actual companion
coefficient $a_{n_{\mathrm{sm}}}$ by $a$ contributes the root
coefficient error already included in
\eqref{eq:summed-approximation-error}.
The low-band contribution is
$C\Lambda^{-\rho/2}\lambda_*^{1-\beta}
\|U-V\|_{\mathscr L_0^2(H)}$ by
Lemma~\ref{lem:damped-energy-low-frequency}(ii).
For each high band, use the finite identity
\eqref{eq:joint-finite-stopped-identity}.
Use $t=t_0+h_\lambda$ for the window-start norm.
Corollary~\ref{cor:initial-window-comparisons} and the initial
part of Proposition~\ref{prop:all-order-first-hit-joint-boundary}
cover $0<t<h_\lambda$, giving the full $[0,T]$ norm.

\emph{Identify the linear and the three nonlinear contributions.}
The higher stochastic-integral estimate
\eqref{eq:higher-martingale-frequency-combined} contributes
$\lambda_*^{-w_2}\mathsf d_{\mathrm{coef}}(U,V)$,
$\lambda_*^{-q_{2,\chi}}\|U-V\|_{\mathscr L_0^2(H)}^\chi$ and
$\lambda_*^{-q_{2,\mathrm{prof}}}\|U-V\|_{\mathscr L_0^2(H)}$.
Higher orders have strictly larger degrees, so order two bounds
their finite sum. The Gaussian-average estimate
\eqref{eq:all-order-predictable-integrated-sobolev} gives the
sum of $\Lambda^{-\rho/2}\lambda_*^{-p_j}$ multiplying the
solution distance. The actual remainder gives the last
nonlinear term of \eqref{eq:all-beta-joint-master}, by
\eqref{eq:paired-actual-remainder}.

For the first average and stochastic integral use
\eqref{eq:joint-first-pair-regrouping}.
The terminal Gaussian sum is bounded by
$C\|U-V\|_{\mathscr L_0^2(H)}^\alpha$ by
\eqref{eq:sharp-main-Gaussian-sum}.
The complete first martingale contributes
$C\mathfrak L_{N_\beta,\Lambda}\lambda_*^{-\beta}
\mathsf d_{\mathrm{coef}}(U,V)$ and
$C\lambda_*^{-q_\Gamma}\|U-V\|_{\mathscr L_0^2(H)}$.
Its clipped mean term is bounded by
$C\|U-V\|_{\mathscr L_0^2(H)}^\alpha$, using the dyadic sum in
Remark~\ref{rem:first-martingale-dyadic-sums} for $0<\alpha<1$.
Finally the first-defect estimate
\eqref{eq:all-order-first-hit-joint-boundary} is linear in
$\mathsf d_{\mathrm{coef}}(U,V)$ and supplies the remaining
logarithmic term in $L_{\mathrm{coef}}$.

Sum the positive degrees geometrically, keeping the logarithms
specified in the statement. The coefficient corrections give
\eqref{eq:summed-approximation-error}; the discarded terms give
$\lambda_*^{-p_{\mathrm{tail}}}$ by
Lemma~\ref{lem:exact-bracket-cutoff-limits}.
Their sum is \eqref{eq:all-beta-joint-master}.
\end{proof}

\subsection{From the coefficient inequality to energy decay}
\label{sec:closure-parameters-energy}
The coefficient inequality has two kinds of terms to handle.
The linear ones must have coefficients tending to zero after
the first energy estimate in \eqref{eq:joint-closure-energy}.
The three powers below one are handled by Young's inequality;
their resulting squared-energy decay rates must exceed two.
We use the rate $5/2$ to allow all fixed logarithmic factors
and small losses.

\begin{lemma}[Parameter choice and quantitative energy decay]
\label{lem:parameter-energy-closure}
For $0<\beta\le\min(\beta_0,4/25)$, choose
\begin{gather}
 \rho=\tfrac12-\tfrac35\beta,\qquad r=\tfrac14,
 \qquad \gamma_{\mathrm{time}}=3\beta,\nonumber\\
 \chi=1-\tfrac{11}{50}\beta,\qquad
 \chi_{\mathrm{rem}}=\beta/2,\qquad
 N_\beta=\lceil6/\beta^2\rceil.
 \label{eq:uniform-direct-family}
\end{gather}
Here $r$ and $\gamma_{\mathrm{time}}$ record the fixed cutoff and
duration exponents used above.
All linear terms in
Proposition~\ref{prop:master-coefficient-comparison} are absorbable
for sufficiently large $\Lambda$. Choose $p_{\mathrm{tail}}>5$
and sufficiently small fixed losses. There is a constant
$C_{\beta,T}$ such that, with the energy-level approximation error
\[
 \operatorname{Err}^{\mathrm{en}}_{n_{\mathrm{sm}},\Lambda}
 :=C_{\beta,T}\Lambda^{-3\beta/5}
                       \operatorname{Err}_{n_{\mathrm{sm}},\Lambda}^{\,2},
\]
the coupled pair satisfies
\begin{equation}
 \|U-V\|_{\mathscr L_0^2(H)}^2
 \le C_{\beta,T}\Lambda^{-5/2}
       +\operatorname{Err}^{\mathrm{en}}_{n_{\mathrm{sm}},\Lambda}.
 \label{eq:all-beta-final-energy}
\end{equation}
\end{lemma}

\begin{proof}[Proof of Lemma~\ref{lem:parameter-energy-closure}]
\emph{1. Check the parameters and absorb the linear terms.}
The chosen parameters satisfy
\[
 \tfrac12-\rho=\tfrac35\beta,\qquad
 1-\alpha=\tfrac15\beta,\qquad
 \alpha-\chi=\tfrac1{50}\beta,\qquad
 1-\chi=\tfrac{11}{50}\beta.
\]
Thus $0<\rho<1/2$, $0<\chi<\alpha<1$,
$0<\chi_{\mathrm{rem}}<\beta$ and $0<3\beta<4$.
The lowest degrees become
\begin{align*}
 w_2&=\beta+3\beta^2/4>0,&
 p_2&=\beta/4+3\beta^2/4>0,\\
 q_\Gamma&=1-19\beta/20+3\beta^2/4>0,&
 q_{2,\mathrm{prof}}&=\beta-3\beta^2/20>0,\\
 q_{2,\chi}&=\beta/50+3\beta^2/4>0,&
 W_{N_\beta}&\ge\beta+9/2.
\end{align*}
Here $2\rho-\chi=-49\beta/50<0$ selects the first argument of the
minimum, and $N_\beta\ge6/\beta^2$ gives the last bound and
$W_{N_\beta}-\chi_{\mathrm{rem}}>0$.

The signed-energy estimate now reads
$\|U-V\|_{\mathscr L_0^2(H)}
\le C\Lambda^{-3\beta/10}\mathsf d_{\mathrm{coef}}(U,V)$.
Substitute this in the linear part of
\eqref{eq:all-beta-joint-master}. Its coefficient in front of
$\mathsf d_{\mathrm{coef}}(U,V)$ is
$C(L_{\mathrm{coef}}+\Lambda^{-3\beta/10}L_{\mathrm{dist}})$.
The positive degrees above imply $L_{\mathrm{coef}}\to0$,
including its logarithmic factor. The only growing frequency
factor in $L_{\mathrm{dist}}$ is the low-band term, and it satisfies
\[
 \Lambda^{-3\beta/10}\Lambda^{-\rho/2}\lambda_*^{1-\beta}
 =\Lambda^{-1/4}\lambda_*^{1-\beta}
 \asymp\Lambda^{-\beta/4}\longrightarrow0.
\]
The remaining terms are negative powers of $\Lambda$ times
bounded or decaying frequency factors, including
$\Lambda^{-1/4}\sum_j\lambda_*^{-p_j}$.
For large $\Lambda$, the whole linear coefficient is at most
$1/2$ and is absorbed.

\emph{2. Apply Young's inequality to the three nonlinear powers.}
First use the zero-loss exponents. Define the total
additive norm error
\[
 \mathcal R_{n_{\mathrm{sm}},\Lambda}
 :=\operatorname{Err}_{n_{\mathrm{sm}},\Lambda}
                      +\lambda_*^{-p_{\mathrm{tail}}}.
\]
After the preceding absorption, multiply the coefficient
inequality by the signed-energy factor $C\Lambda^{-3\beta/10}$.
Since $\lambda_*\asymp\Lambda^{1/4}$, this gives
\begin{equation}
 \begin{aligned}
 \|U-V\|_{\mathscr L_0^2(H)}\le C\big[&
    \Lambda^{-3\beta/10}\|U-V\|_{\mathscr L_0^2(H)}^\alpha\\
 &+\Lambda^{-3\beta/10-q_{2,\chi}/4}
                              \|U-V\|_{\mathscr L_0^2(H)}^\chi\\
 &+\Lambda^{-3\beta/10-\rho/2-(W_{N_\beta}-\chi_{\mathrm{rem}})/4}
                   \|U-V\|_{\mathscr L_0^2(H)}^{\chi_{\mathrm{rem}}}\big]\\
 &+C\Lambda^{-3\beta/10}\mathcal R_{n_{\mathrm{sm}},\Lambda}.
 \end{aligned}
 \label{eq:closure-three-power-feedback}
\end{equation}
For $x\ge0$, $0<z<1$ and $\epsilon>0$, Young's inequality gives
$cx^z\le\epsilon x+C_{z,\epsilon}c^{1/(1-z)}$.
Apply this to $x=\|U-V\|_{\mathscr L_0^2(H)}$, absorb the
three $\epsilon x$ terms and square.
Each $C\Lambda^{-a}x^z$ gives squared-energy rate
$2a/(1-z)$. In \eqref{eq:closure-three-power-feedback} these are
\begin{align*}
 \text{terminal Gaussian and first mean:}\quad&
 \frac{3\beta/5}{1-\alpha}=3,\\
 \text{higher stochastic integrals:}\quad&
 \frac{3\beta/5+q_{2,\chi}/2}{1-\chi}
     =\frac{61}{22}+\frac{75}{44}\beta,\\
 \text{actual remainder:}\quad&
 \frac{3\beta/5+\rho+(W_{N_\beta}-\chi_{\mathrm{rem}})/2}
      {1-\chi_{\mathrm{rem}}}
 \ge\frac{11/4+\beta/4}{1-\beta/2}>\frac{11}{4}.
\end{align*}
All exceed $5/2$; only the remainder rate uses $N_\beta$.
The additive error contributes
$C\Lambda^{-3\beta/5}\mathcal R_{n_{\mathrm{sm}},\Lambda}^{\,2}$.

\emph{3. Restore the logarithms, losses and the truncation tail.}
Fix $\beta,N_\beta$ and all moment orders before taking
$\Lambda\to\infty$. There are finitely many losses at this depth.
For every fixed $C$ and every $\delta>0$,
$(1+\log\Lambda)^C\le C_\delta\Lambda^\delta$.
Choose the losses and $\delta$ small enough that the three
exponents just computed remain above $5/2$, and the coefficients
in step~1 still decay.
If a linear solution-distance term has the small loss
$\|U-V\|_{\mathscr L_0^2(H)}^{1-\varepsilon}$, its contribution after the
signed-energy substitution is
$C\Lambda^{-a}\|U-V\|_{\mathscr L_0^2(H)}^{1-\varepsilon}$ with some fixed $a>0$
from step~1. The same Young inequality gives squared-energy
exponent $2a/\varepsilon$. Taking $\varepsilon$ sufficiently
small keeps this exponent above $5/2$ as well.

For $p_{\mathrm{tail}}>5$, the squared tail has rate
$3\beta/5+p_{\mathrm{tail}}/2>5/2$. Separating it from the coefficient
error by $(u+v)^2\le2u^2+2v^2$ gives \eqref{eq:all-beta-final-energy}.
\end{proof}

\subsection{Choosing a common approximation sequence}
\label{sec:closure-approximation-laws}
Lemma~\ref{lem:parameter-energy-closure} gives an energy estimate
with an explicit residual approximation error. The next lemma
chooses the same smoothing sequence for two prescribed solutions
and removes that error. It is the quantitative input used in
the proof of Theorem~\ref{thm:all-positive-beta-joint}.

\begin{lemma}[A common approximation with vanishing control energy]
\label{lem:common-approximation-control}
Assume the coefficient hypotheses of
Theorem~\ref{thm:main-weak-uniqueness}. For the fixed effective
$0<\beta\le\min(\beta_0,4/25)$,
fix $0<T\le1$ and two prescribed weak solutions $U^{(1)},U^{(2)}$
with the same deterministic initial state. Use the parameter
family \eqref{eq:uniform-direct-family}. There is a deterministic
choice $n_{\mathrm{sm}}(\Lambda)$ for which the controlled companions
$V^{(i)}=V^{(i),n_{\mathrm{sm}}(\Lambda),\Lambda}$ of
Proposition~\ref{prop:construct-companion} satisfy, for $i=1,2$,
\begin{equation}
 \begin{aligned}
 \|U^{(i)}-V^{(i)}\|_{\mathscr L_0^2(H)}^2
    &\le C_{\beta,T}\Lambda^{-5/2},\\
 \mathcal C_\Lambda^{(i)}
    :=m_0^{-2}\Lambda^2
              \|U^{(i)}-V^{(i)}\|_{\mathscr L_0^2(H)}^2
    &\le C_{\beta,T,m_0}\Lambda^{-1/2}\longrightarrow0.
 \end{aligned}
 \label{eq:closure-common-diagonal-control}
\end{equation}
The index is the same for both solutions; each companion is
constructed on its prescribed solution's probability space.
\end{lemma}

\begin{proof}[Proof of Lemma~\ref{lem:common-approximation-control}]
Use \eqref{eq:uniform-direct-family} and construct the two companions
by Proposition~\ref{prop:construct-companion}.
For sufficiently large $\Lambda$, all high-band durations satisfy
\eqref{eq:window-threshold-high-band}; the low bands and initial
windows are already covered by
Proposition~\ref{prop:master-coefficient-comparison}.

At each fixed $\Lambda$, Lemma~\ref{lem:exact-bracket-cutoff-limits}
and \eqref{eq:diagonal-rms-first} allow the following deterministic
choice: first choose common finite frequency and positive-gap cutoffs
so their tails are small, then choose one $n_{\mathrm{sm}}(\Lambda)$
large enough for both pairs that
$\operatorname{Err}_{n_{\mathrm{sm}},\Lambda}\le
\Lambda^{-p_{\mathrm{tail}}}$, including the source
$\sqrt{\varepsilon_{\mathrm{app}}}$ and the fixed summation constants
in \eqref{eq:summed-approximation-error}.
There are only two sets of requirements, and the replacement-error
moduli are uniform in the companions, so a common index suffices.

In \eqref{eq:all-beta-final-energy} the energy-level error is then
at most $C\Lambda^{-3\beta/5-2p_{\mathrm{tail}}}$.
Since $p_{\mathrm{tail}}>5$, it is bounded by
$C\Lambda^{-5/2}$. Multiplication by $m_0^{-2}\Lambda^2$ proves
\eqref{eq:closure-common-diagonal-control} for both pairs.
\end{proof}

\subsection{Conditional restart and global weak uniqueness}
\label{sec:conditional-global}
The passage from local to global uniqueness uses standard conditional
restart on the Polish path space $\mathcal X_h=C([0,h];H)$.

\begin{lemma}[Conditional restart]
\label{lem:conditional-restart}
Assume the coefficient hypotheses of
Theorem~\ref{thm:main-weak-uniqueness}, and let $U$ be a weak mild
solution on $[0,T]$. For deterministic $t_0\ge0$ and $h>0$ with
$t_0+h\le T$, condition on $\sigma(U_s:0\le s\le t_0)$.
Almost surely, the conditional law of
$(U_{t_0+s})_{0\le s\le h}$ is the state-path law of a weak mild
solution from $U_{t_0}$. If $h\le1$, this law is uniquely determined
by $U_{t_0}$, by Theorem~\ref{thm:all-positive-beta-joint}.
\end{lemma}
\begin{proof}[Proof of Lemma~\ref{lem:conditional-restart}]
Use the equivalent cylinder martingale problem for the weak mild
equation~\cite{Kunze2013}. Its martingale identities are preserved by
disintegration at a deterministic time: the conditional future starts
at $U_{t_0}$ and solves the same problem; see
\cite[Chapter~4]{EthierKurtz1986}. The equivalence with weak mild
solutions gives the first assertion, and local uniqueness gives the second.
\end{proof}

\begin{proof}[Proof of Theorem~\ref{thm:main-weak-uniqueness}]
Weak existence is standard~\cite{GatarekGoldys1994,DaPratoZabczyk2014}.
For uniqueness, choose an integer $J\ge T$ and set $t_j=jT/J$.
Suppose two solutions have the same law on $C([0,t_j];H)$.
Under this common law of the past, Lemma~\ref{lem:conditional-restart}
and Theorem~\ref{thm:all-positive-beta-joint} identify their conditional
laws on $[t_j,t_{j+1}]$, since the interval has length at most one
and both restart from the same state. Integrating against the common
law of the past gives equality on $C([0,t_{j+1}];H)$.
Starting from the common deterministic initial state, finite induction
proves equality of the full state-path laws on $[0,T]$.
\end{proof}

\section{Scalar extensions: unbounded diffusion and measurable drift}
\label{sec:scalar-extensions}
We prove Theorem~\ref{thm:main-extensions} in two steps: localize the
drift-free equation to allow unbounded diffusion, then add a measurable
linear-growth drift by Girsanov.

\subsection{Unbounded diffusion from bounded initial data}
\label{sec:unbounded-coefficient-linf-localization}

This subsection proves the scalar unbounded drift-free input to
Theorem~\ref{thm:main-extensions}, using the bounded scalar
Theorem~\ref{thm:main-weak-uniqueness}. The initial state is bounded
and deterministic.

\begin{proposition}[Scalar unbounded diffusion from bounded initial data]
\label{thm:unbounded-coefficient-linf-localization}
On the one-dimensional torus, suppose
$g:\mathbb R\to\mathbb R$ is globally
$\beta$-H\"older for some $0<\beta<1$ and
$|g(z)|\ge c_0>0$ for all $z$. No uniform upper bound on
$g$ is imposed. For every deterministic
$u_0\in L^\infty(\mathbb T)$, a weak mild solution
with continuous $L^2$ paths exists and its law is unique on each finite
time interval. The initial function need not have a continuous representative.
\end{proposition}

The globally Lipschitz case $\beta=1$ follows from the usual mild
uniqueness argument. Here we treat every $0<\beta<1$. The coefficient satisfies the sublinear bound
\begin{equation}
 |g(z)|\le C(1+|z|^\beta).
 \label{eq:linf-localization-sublinear-growth}
\end{equation}

\subsubsection{The stochastic remainder is continuous through time zero}

\begin{lemma}[Stochastic-remainder continuity at the initial time]
\label{lem:unbounded-remainder-continuity}
Under the hypotheses of
Proposition~\ref{thm:unbounded-coefficient-linf-localization}, for every
finite $T$ and every sufficiently large fixed $p$, any weak mild
solution has a stochastic remainder $Z_t=U_t-S_tu_0$, with $Z_0=0$,
and obeys
\begin{equation}
 \mathbb E\sup_{t\le T}\|U_t\|_{L^2(\mathbb T)}^p<\infty,
 \qquad
 \mathbb E\|Z\|_{C([0,T]\times\mathbb T)}^p<\infty,
 \quad Z_t=U_t-S_tu_0,\quad Z_0=0.
 \label{eq:linf-localization-remainder-moments}
\end{equation}
These bounds are uniform over coefficients with the same growth
constant in \eqref{eq:linf-localization-sublinear-growth}. Constants
may depend on $p,T,\beta$ and $\|u_0\|_{L^2(\mathbb T)}$.

\end{lemma}
\begin{proof}[Proof of Lemma~\ref{lem:unbounded-remainder-continuity}]
The growth bound gives
\[
 \|S_r\mathsf M_{g(v)}\|_{\mathrm{HS}}^2
 \le Cr^{-1/2}(1+\|v\|_{L^2(\mathbb T)}^{2\beta}).
\]
Stop at $\tau_n=\inf\{t:\|U_t\|_{L^2(\mathbb T)}\ge n\}$.
Hilbert-space stochastic factorization, applied to the stopped
integrands on $[0,T]$ with exponent in $(1/p,1/4)$, gives
\cite[Chapter~5]{DaPratoZabczyk2014}
\[
 \mathbb E\sup_{t\le T\wedge\tau_n}\|U_t\|_{L^2(\mathbb T)}^p
 \le C(1+\|u_0\|_{L^2(\mathbb T)}^p)
   +C\left(\mathbb E\sup_{t\le T\wedge\tau_n}
                                  \|U_t\|_{L^2(\mathbb T)}^p\right)^\beta.
\]
Young's inequality absorbs the sublinear term uniformly in $n$;
letting $n\to\infty$ proves the first moment bound.

For the uniform spatial norm, H\"older's inequality and heat scaling give
\begin{equation}
 \int_{\mathbb T}p_r(x-y)^2|g(v(y))|^2dy
 \le C r^{-(1+\beta)/2}(1+\|v\|_{L^2(\mathbb T)}^{2\beta}).
 \label{eq:linf-localization-pointwise-kernel}
\end{equation}
Choose $0<\eta<(1-\beta)/4$ and $p>3/(2\eta)$.
Then $2\eta+(1+\beta)/2<1$ makes the inner stochastic convolution
in the factorization integrable in $L^p(\Omega\times[0,T]\times\mathbb T)$.
The heat estimate $\|S_rf\|_\infty\le Cr^{-1/(2p)}\|f\|_p$ and
$\eta>3/(2p)$ make the outer fractional convolution continuous
into $C([0,T]\times\mathbb T)$, with value zero at time zero.
The factorization theorem therefore gives the asserted supremum
moment of $Z$. Increasing $p$ covers every required finite moment.
\end{proof}

In particular $U$ is continuous with values in $C(\mathbb T)$
on every interval $[\varepsilon,T]$. At time zero we use
\begin{equation}
 \|S_tu_0\|_\infty\le\|u_0\|_\infty,
 \qquad \|Z_t\|_\infty\longrightarrow0.
 \label{eq:linf-localization-initial-gap}
\end{equation}
This is the use of bounded initial data in the localization below:
the state starts strictly inside a fixed uniform-norm cutoff.
The Girsanov transfer in Section~\ref{sec:linear-growth-drift}
requires only $L^2$ initial data; the bounded-data restriction in the
unbounded-diffusion case of Theorem~\ref{thm:main-extensions} comes
from this drift-free localization step.

\subsubsection{Cutoffs and a measurable positive stopping time}

Let ${R_{\mathrm{cut}}}>\|u_0\|_\infty+1$ and define the scalar clipping map
$\pi_{R_{\mathrm{cut}}}(z)=\max(-R_{\mathrm{cut}},\min(z,R_{\mathrm{cut}}))$. Put
$g_{R_{\mathrm{cut}}}=g\circ\pi_{R_{\mathrm{cut}}}$. The projection is $1$-Lipschitz and
$|\pi_{R_{\mathrm{cut}}}z|\le|z|$, so $g_{R_{\mathrm{cut}}}$ has the same global H\"older seminorm
bound and a cutoff-independent sublinear growth bound. Moreover
\[
 |g_{R_{\mathrm{cut}}}(z)|=|g(\pi_{R_{\mathrm{cut}}}z)|\ge c_0.
\]
Its uniform upper bound can depend on ${R_{\mathrm{cut}}}$. Thus every fixed $g_{R_{\mathrm{cut}}}$
is covered by Theorem~\ref{thm:main-weak-uniqueness} for all deterministic
$L^2$ initial states.

\begin{lemma}[Positive exit times and uniform nonexplosion]
\label{lem:unbounded-cutoff-exit}
For $R_{\mathrm{cut}}>\|u_0\|_\infty+1$, define
\begin{equation}
 \tau_{R_{\mathrm{cut}}}=\inf\{t>0:\|U_t\|_\infty\ge {R_{\mathrm{cut}}}\}.
 \label{eq:linf-localization-exit}
\end{equation}
This is an almost surely positive stopping time. On a finite hit,
$U_{\tau_{R_{\mathrm{cut}}}}$ is continuous and has uniform norm
$R_{\mathrm{cut}}$. For every fixed $T$ and sufficiently large fixed $p$,
\begin{equation}
 \mathbb P(\tau_{R_{\mathrm{cut}}}\le T)
 \le \frac{C_{p,T}}{(R_{\mathrm{cut}}-\|u_0\|_\infty)^p}
 \longrightarrow0.
 \label{eq:linf-localization-exit-probability}
\end{equation}
The estimate also holds uniformly for solutions with any cutoff
coefficient $g_{R_{\mathrm{cut}}'}$.
\end{lemma}
\begin{proof}[Proof of Lemma~\ref{lem:unbounded-cutoff-exit}]
Positive-time continuity in $C(\mathbb T)$ and the strict initial gap
\eqref{eq:linf-localization-initial-gap} make $\tau_{R_{\mathrm{cut}}}$
a positive stopping time, with norm $R_{\mathrm{cut}}$ at a finite hit.
A hit by $T$ forces
$\sup_{t\le T}\|U_t-S_tu_0\|_\infty\ge R_{\mathrm{cut}}-\|u_0\|_\infty$.
Markov's inequality and Lemma~\ref{lem:unbounded-remainder-continuity}
give \eqref{eq:linf-localization-exit-probability}, uniformly over the
cutoffs because their growth bounds are common.
\end{proof}

\subsubsection{Stopped uniqueness follows by pasting laws}

\begin{proof}[Proof of Proposition~\ref{thm:unbounded-coefficient-linf-localization}]
For each fixed $R_{\mathrm{cut}}$, Theorem~\ref{thm:main-weak-uniqueness}
makes the $g_{R_{\mathrm{cut}}}$ equation well posed from every
deterministic $H$ state. Its laws form a Borel continuation kernel
by the martingale-problem formulation
\cite[Chapter~4]{EthierKurtz1986} and \cite{Kunze2013}.
Continue an original solution at $\tau_{R_{\mathrm{cut}}}$ with this
cutoff law from its exit state. The coefficients agree before the
stop, so the pasted path solves the cutoff martingale problem.
Uniqueness identifies its law and therefore the original stopped law.

For two original solutions, expectations of a bounded Borel path
functional $F$ consequently differ by at most $2\|F\|_\infty$ times
the sum of their exit probabilities. These vanish as
$R_{\mathrm{cut}}\to\infty$ by
\eqref{eq:linf-localization-exit-probability}, proving uniqueness.

Weak existence follows from the continuous-coefficient, linear-growth
martingale-solution theorem
\cite{GatarekGoldys1994,DaPratoZabczyk2014}. The growth bound
\eqref{eq:linf-localization-sublinear-growth} and the heat-smoothed
Hilbert--Schmidt estimate in Lemma~\ref{lem:unbounded-remainder-continuity}
verify its hypotheses; continuity of $g$ gives continuity of the
heat-smoothed coefficient.
\end{proof}

\subsection{Measurable drifts of linear growth via Girsanov}
\label{sec:linear-growth-drift}

This subsection completes Theorem~\ref{thm:main-extensions} by adding the
drift to the scalar drift-free laws of Theorem~\ref{thm:main-weak-uniqueness}
and Proposition~\ref{thm:unbounded-coefficient-linf-localization}.
The argument has three steps: a moment bound uniform under stopped
changes of measure; a proof that the full Girsanov density has expectation
one; and a canonical expression for that density, which identifies the
state-path law.

The use of linear growth to justify Girsanov goes back to
Bene\v{s}~\cite{Benes1971}; see \cite{KlebanerLiptser2014} for extensions
of this method beyond Novikov's condition. For cylindrical noise and
measurable linear-growth drift with additive diffusion, a detailed
existence-and-uniqueness argument is given in
\cite[Appendix~A.1]{DaPratoFlandoliPriolaRockner2013}.
We use the cylindrical Girsanov theorem in
\cite[Section~10.2]{DaPratoZabczyk2014}.
The distinction between state-law uniqueness and uniqueness of the
joint state-noise law is treated in
\cite{Ondrejat2004,CriensRitter2022}.
Here the density is constructed directly from the canonical state
martingales.

\subsubsection{The equation and the transfer statement}
\label{sec:drift-transfer-statement}
Fix $T<\infty$ and use $H=L^2(\mathbb T)$ with its real inner product
$\langle\cdot,\cdot\rangle_{H}$. For the Girsanov step it is enough to
assume that $g:\mathbb R\to\mathbb R$ is continuous and
\begin{equation}
 \mathfrak a(z):=g(z)^2\ge c_0^2,\qquad
 |g(z)|\le L_g(1+|z|),\qquad z\in\mathbb R,
 \label{eq:drift-diffusion-assumptions}
\end{equation}
where $c_0>0$ and $L_g<\infty$ are fixed. The drift
$b:[0,T]\times\mathbb T\times\mathbb R\to\mathbb R$ is Borel and
\begin{equation}
 |b(t,x,z)|\le L_b(1+|z|).
 \label{eq:drift-linear-growth}
\end{equation}
The constant $L_b$ may depend on the fixed horizon $T$.

\begin{definition}[Mild equation with a measurable drift]
\label{def:linear-growth-drift-solution}
For $v\in H$, write $b_t(v)(x)=b(t,x,v(x))$ and
$g(v)(x)=g(v(x))$. A weak mild solution from a deterministic $u_0\in H$
is an adapted continuous $H$-valued process $U$ on a filtered
probability space carrying a cylindrical Brownian motion $W$ on $H$,
such that
\begin{equation}
 U_t=S_tu_0+\int_0^t S_{t-s}b_s(U_s)\,ds
              +\int_0^t S_{t-s}\mathsf M_{g(U_s)}\,dW_s.
 \label{eq:drift-mild-equation}
\end{equation}
The drift integral is a Bochner integral in $H$, and
$\mathsf M_{g(U_s)}$ means pointwise multiplication before heat smoothing. The solution path space is
$\mathcal X_T:=C([0,T];H)$, with its Borel sigma field.
\end{definition}

The growth and nondegeneracy bounds make
$v\mapsto b_t(v),g(v),b_t(v)/g(v),b_t(v)/g(v)^2$ Borel maps into $H$,
jointly Borel in $(t,v)$. Their compositions with a continuous
adapted state process are predictable.

\begin{proposition}[Transfer of weak well-posedness by Girsanov]
\label{prop:linear-growth-girsanov-transfer}
Assume \eqref{eq:drift-diffusion-assumptions}--\eqref{eq:drift-linear-growth}.
For a fixed deterministic $u_0\in H$, suppose that the drift-free
equation has a weak mild solution and a unique state-path law
$\mu_0$ on $\mathcal X_T$.
Then \eqref{eq:drift-mild-equation} has a weak mild solution and a
unique state-path law $\mu_b$ on $\mathcal X_T$.
The measures $\mu_b$ and $\mu_0$ are mutually absolutely continuous.
\end{proposition}

Only the drift-free well-posedness input uses the H\"older hypotheses
of Theorem~\ref{thm:main-extensions}. The transfer proposition applies
to every continuous diffusion satisfying
\eqref{eq:drift-diffusion-assumptions} for which that input is known.

\subsubsection{Uniform moments under the stopped changes of measure}
\label{sec:drift-uniform-moments}
Both directions of Girsanov produce a drift of the form
$\kappa_t b_t(U_t)$, where the predictable scalar process $\kappa_t$
takes values in $[0,1]$. We therefore prove the required moment bound
once, uniformly in this multiplier.

\begin{lemma}[A uniform second-moment bound]
\label{lem:drift-uniform-second-moment}
Let $U$ be an adapted continuous $H$-valued mild solution with
drift $\kappa_t b_t(U_t)$, $0\le\kappa_t\le1$, and diffusion $g(U_t)$,
on any probability space carrying cylindrical Brownian noise.
Under \eqref{eq:drift-diffusion-assumptions}--\eqref{eq:drift-linear-growth},
\begin{equation}
 \sup_{0\le t\le T}\mathbb E\|U_t\|_{L^2(\mathbb T)}^2
 \le C_T(1+\|u_0\|_{L^2(\mathbb T)}^2),
 \label{eq:drift-uniform-second-moment}
\end{equation}
where $C_T$ depends only on $T,L_b,L_g$, and is independent of
$\kappa$ and of the probability space.
\end{lemma}

\begin{proof}[Proof of Lemma~\ref{lem:drift-uniform-second-moment}]
\emph{Heat smoothing and localization.}
For $0<r\le T$, the heat kernel and the growth assumption give
\begin{equation}
 \begin{aligned}
 \|S_r\mathsf M_{g(v)}\|_{\mathrm{HS}}^2
 &=\|p_r\|_{L^2(\mathbb T)}^2
                 \int_{\mathbb T}|g(v(x))|^2\,dx\\
 &\le C_T(1+r^{-1/2})
              (1+\|v\|_{L^2(\mathbb T)}^2),\\
 \|b_t(v)\|_{L^2(\mathbb T)}^2
 &\le 2L_b^2(1+\|v\|_{L^2(\mathbb T)}^2).
 \end{aligned}
 \label{eq:drift-growth-heat-bound}
\end{equation}
Choose $R>\|u_0\|_{H}$ and set
$\sigma_R=\inf\{t:\|U_t\|_{H}\ge R\}$, with
$\inf\varnothing=\infty$. On the full deterministic interval define
\[
 U_t^{[R]}=S_tu_0+
 \int_0^t\mathbf1_{\{s\le\sigma_R\}}S_{t-s}\kappa_sb_s(U_s)\,ds
 +\int_0^t\mathbf1_{\{s\le\sigma_R\}}
                  S_{t-s}\mathsf M_{g(U_s)}\,dW_s.
\]
This stopped-integrand convolution agrees with $U_t$ before
$\sigma_R$, so
$\mathbf1_{\{s\le\sigma_R\}}\|U_s\|_{H}^2
\le\|U_s^{[R]}\|_{H}^2$.
The stopped integrands and \eqref{eq:drift-growth-heat-bound} imply
that $F_R(t):=\mathbb E\|U_t^{[R]}\|_{H}^2$ is finite and bounded
on $[0,T]$ for each $R$.

\emph{A bound independent of the localization level.}
Heat contraction, Cauchy--Schwarz for the drift integral and
stochastic isometry give
\begin{equation}
 F_R(t)\le C_T(1+\|u_0\|_{H}^2)
       +C_T\int_0^t[1+(t-s)^{-1/2}]F_R(s)\,ds.
 \label{eq:drift-moment-volterra}
\end{equation}
Multiply by $e^{-\omega t}$ and take the time supremum.
Choosing $C_T(\omega^{-1}+\sqrt\pi\,\omega^{-1/2})\le1/2$
absorbs the convolution and gives
$F_R(t)\le C_T(1+\|u_0\|_H^2)$ uniformly in $R,t$.
Since continuous $H$ paths eventually satisfy $\sigma_R>T$,
Fatou proves \eqref{eq:drift-uniform-second-moment}.
\end{proof}

\subsubsection{The Girsanov densities have expectation one}
\label{sec:drift-density-mass}
We first isolate the elementary passage from stopped densities to
the full density. Its hypothesis is a moment bound under the
stopped changed measures.

\begin{lemma}[An energy criterion for the full density]
\label{lem:drift-energy-mass-one}
Let $\theta$ be a predictable $H$-valued process, let
$A_t=\int_0^t\|\theta_s\|_{H}^2\,ds$, and assume $A_T<\infty$
almost surely. Define
\[
 R_t=\exp\!\left(\int_0^t\langle\theta_s,dW_s\rangle_{H}
                         -\tfrac12A_t\right),\qquad
 \tau_n=\inf\{t:A_t\ge n\},\quad n\ge1.
\]
The stopped densities have expectation one; set
$d\mathbb Q_n=R_{T\wedge\tau_n}\,d\mathbb P$.
If $\sup_n\mathbb E_{\mathbb Q_n}A_T\le C<\infty$, then
$\mathbb E_{\mathbb P}R_T=1$, and $R$ is a martingale on $[0,T]$.
\end{lemma}

\begin{proof}[Proof of Lemma~\ref{lem:drift-energy-mass-one}]
Continuity of $A$ gives $A_{T\wedge\tau_n}\le n$.
Novikov's criterion therefore proves
$\mathbb E R_{T\wedge\tau_n}=1$.
On $\{A_T<n\}$ the stopped and unstopped densities coincide.
Since a nonnegative local martingale is a supermartingale,
\begin{equation}
 \begin{aligned}
 1\ge\mathbb E_{\mathbb P}R_T
 &\ge\mathbb E_{\mathbb P}[R_T\mathbf1_{\{A_T<n\}}]\\
 &=\mathbb Q_n(A_T<n)
 \ge1-\frac{\mathbb E_{\mathbb Q_n}A_T}{n}
 \ge1-\frac Cn.
 \end{aligned}
 \label{eq:drift-density-no-mass-loss}
\end{equation}
Let $n\to\infty$. The supermartingale inequalities then give
$1=\mathbb E R_T\le\mathbb E R_t\le1$ for all $t\le T$.
For $s<t$, the nonnegative difference
$R_s-\mathbb E[R_t\mid\mathcal F_s]$ has expectation zero.
Thus the conditional supermartingale inequality is an equality,
which proves the martingale assertion.
\end{proof}

For any solution path $U$ define the drift-removing noise direction
and its energy by
\begin{equation}
 \begin{aligned}
 \theta_t(U)(x)&:=\frac{b(t,x,U_t(x))}{g(U_t(x))},\\
 A_t(U)&:=\int_0^t\|\theta_s(U)\|_{H}^2\,ds,\\
 \|\theta_t(U)\|_{H}^2
 &\le 2c_0^{-2}L_b^2(1+\|U_t\|_{H}^2).
 \end{aligned}
 \label{eq:drift-noise-direction}
\end{equation}
The denominator is nonzero because $|g|\ge c_0$.
In particular $A_T(U)<\infty$ on every
continuous $H$ path.

\begin{proof}[Proof of Proposition~\ref{prop:linear-growth-girsanov-transfer}
(existence and removal of the drift)]
\emph{Adding the drift.}
Take a drift-free solution $(X,W)$ with state law $\mu_0$.
Use the density in Lemma~\ref{lem:drift-energy-mass-one} with
$\theta=\theta(X)$ and $\tau_n=\inf\{t:A_t(X)\ge n\}$.
Under its stopped changed measure $\mathbb Q_n$,
\[
 W_t^{(n)}=W_t-\int_0^{t\wedge\tau_n}\theta_s(X)\,ds
\]
is cylindrical Brownian motion. Testing the equation against smooth
spatial functions and using $g(X_s)\theta_s(X)=b_s(X_s)$ shows that
$X$ has drift $\mathbf1_{\{s\le\tau_n\}}b_s(X_s)$ under
$\mathbb Q_n$. Variation of constants gives its mild equation.
Lemma~\ref{lem:drift-uniform-second-moment}, applied under
$\mathbb Q_n$, and \eqref{eq:drift-noise-direction} yield
\begin{equation}
 \sup_n\mathbb E_{\mathbb Q_n}A_T(X)
 \le C_T(1+\|u_0\|_{H}^2).
 \label{eq:drift-uniform-changed-energy}
\end{equation}
The energy criterion gives $\mathbb E R_T=1$.
With $d\mathbb Q=R_T\,d\mathbb P$, the full cylindrical Girsanov
theorem makes $W-\int_0^\cdot\theta_s(X)\,ds$ Brownian and gives
\eqref{eq:drift-mild-equation} for $X$ under $\mathbb Q$.
The density is finite and strictly positive, so the two probability
measures, and hence their state-path laws, are equivalent.

\emph{Removing the drift from any solution.}
Now let $(U,W)$ be any weak mild solution of
\eqref{eq:drift-mild-equation} under $\mathbb P$.
Form
\begin{equation}
 \widehat R_t
 =\exp\!\left(-\int_0^t\langle\theta_s(U),dW_s\rangle_{H}
                         -\tfrac12A_t(U)\right).
 \label{eq:drift-reverse-density}
\end{equation}
For the stopped density, Girsanov makes
$W+\int_0^{\cdot\wedge\tau_n}\theta_s(U)\,ds$ Brownian.
The remaining drift is
$\mathbf1_{\{s>\tau_n\}}b_s(U_s)$, whose multiplier also lies in
$[0,1]$. Lemma~\ref{lem:drift-uniform-second-moment} therefore
gives \eqref{eq:drift-uniform-changed-energy} for these reverse
changed measures as well. Apply
Lemma~\ref{lem:drift-energy-mass-one} with the direction
$-\theta(U)$ to conclude $\mathbb E\widehat R_T=1$.
Under $d\widehat{\mathbb Q}=\widehat R_T\,d\mathbb P$,
\[
 \widehat W_t=W_t+\int_0^t\theta_s(U)\,ds
\]
is Brownian and $U$ is drift-free. Its state-path law under
$\widehat{\mathbb Q}$ is therefore $\mu_0$.
The identification of its original law is proved next.
\end{proof}

\subsubsection{A canonical density and uniqueness of the state-path law}
\label{sec:drift-canonical-density}
Let $X_t(\omega)=\omega(t)$ be the coordinate process on
$\mathcal X_T$, equipped with $\mu_0$ and its natural filtration.
For a smooth test $\varphi:\mathbb T\to\mathbb R$, define
\begin{equation}
 \mathcal M_t(\varphi)
 :=\langle X_t,\varphi\rangle_{H}
   -\langle u_0,\varphi\rangle_{H}
   -\int_0^t\langle X_s,\Delta\varphi\rangle_{H}\,ds.
 \label{eq:drift-canonical-state-martingale}
\end{equation}
The drift-free equation implies that these are continuous square-integrable
martingales, with cross brackets
\begin{equation}
 \langle\mathcal M(\varphi),\mathcal M(\psi)\rangle_t
 =\int_0^t\int_{\mathbb T}
       \varphi(x)\mathfrak a(X_s(x))\psi(x)\,dx\,ds.
 \label{eq:drift-canonical-brackets}
\end{equation}
Square integrability follows from
Lemma~\ref{lem:drift-uniform-second-moment} with zero drift and
the growth bound on $g$.
Put
\[
 \psi_t^b(X)(x)
 :=\frac{b(t,x,X_t(x))}{g(X_t(x))^2}.
\]
Then $g(X_t)\psi_t^b(X)=\theta_t(X)$.

\begin{lemma}[The stochastic density is determined by the state law]
\label{lem:drift-canonical-integral}
There is a Borel function $J_T^b:\mathcal X_T\to\mathbb R$ such that,
on every realization $(X,W)$ of the drift-free equation with state
law $\mu_0$,
\begin{equation}
 J_T^b(X)=\int_0^T\langle\theta_t(X),dW_t\rangle_{H}
 \quad\hbox{almost surely}.
 \label{eq:drift-canonical-integral}
\end{equation}
It is obtained by integrating $\psi^b$ against the state martingales
\eqref{eq:drift-canonical-state-martingale}.
\end{lemma}

\begin{proof}[Proof of Lemma~\ref{lem:drift-canonical-integral}]
For a simple predictable spatial test
$v_t(x)=\sum_{j=1}^m\eta_j\mathbf1_{(r_j,t_j]}(t)\varphi_j(x)$,
where $\eta_j$ is a bounded Borel function of the canonical path
up to $r_j$ and $\varphi_j$ is smooth, define
\[
 I(v)=\sum_{j=1}^m\eta_j
       [\mathcal M_{t_j}(\varphi_j)-\mathcal M_{r_j}(\varphi_j)].
\]
This is a Borel function of the state path. The bracket identity
gives the isometry
\begin{equation}
 \mathbb E_{\mu_0}|I(v)-I(w)|^2
 =\mathbb E_{\mu_0}\int_0^T\int_{\mathbb T}
       |g(X_t(x))(v_t(x)-w_t(x))|^2\,dx\,dt.
 \label{eq:drift-canonical-isometry}
\end{equation}
On any drift-free realization $I(v)=\int\langle g(X)v,dW\rangle$.
The target $\psi^b$ has squared norm
$\mathbb E_{\mu_0}A_T(X)<\infty$ in this isometry. The measure
$g(X_t(x))^2d\mu_0\,dt\,dx$ is finite by the second-moment bound.
Density of simple predictable tests with smooth spatial factors
therefore gives $v^{(m)}$ of the displayed form such that
$g(X)v^{(m)}\to g(X)\psi^b=\theta(X)$ in $L^2(d\mu_0\,dt\,dx)$.

The isometry makes $I(v^{(m)})$ converge in $L^2(\mu_0)$.
Choose a subsequence converging $\mu_0$-almost surely and define
$J_T^b$ as its pointwise limit on the Borel set of convergence,
and zero elsewhere. On every realization with state law $\mu_0$,
the same subsequence converges almost surely to $J_T^b(X)$.
Stochastic isometry also makes it converge in $L^2$ to the integral
in \eqref{eq:drift-canonical-integral}. Uniqueness of the limit in
probability proves the claimed identity.
\end{proof}

\begin{proof}[Proof of Proposition~\ref{prop:linear-growth-girsanov-transfer}
(uniqueness)]
Define the Borel density on canonical state-path space by
\begin{equation}
 \mathcal D_T^b(X)
 :=\exp\!\left(J_T^b(X)-\tfrac12A_T(X)\right).
 \label{eq:drift-canonical-density}
\end{equation}
The existence part and Lemma~\ref{lem:drift-canonical-integral}
show that $\int\mathcal D_T^b\,d\mu_0=1$.
For any drifted solution, use its reverse change of measure
$\widehat{\mathbb Q}$ constructed above. Under this measure
$U$ is drift-free with driver $\widehat W$ and state law $\mu_0$.
Substituting $dW=d\widehat W-\theta(U)\,dt$ into
\eqref{eq:drift-reverse-density} gives
\[
 \widehat R_T^{-1}
 =\exp\!\left(\int_0^T\langle\theta_t(U),d\widehat W_t\rangle_{H}
                         -\tfrac12A_T(U)\right)
 =\mathcal D_T^b(U).
\]
The substitution is valid first up to each finite energy stop,
and then on $[0,T]$ since $A_T(U)<\infty$.
Consequently, for every bounded Borel $F$ on $\mathcal X_T$,
\begin{equation}
 \mathbb E_{\mathbb P}F(U)
 =\mathbb E_{\widehat{\mathbb Q}}
                [F(U)\widehat R_T^{-1}]
 =\int_{\mathcal X_T}F(X)\mathcal D_T^b(X)\,\mu_0(dX).
 \label{eq:drift-law-identification}
\end{equation}
This identifies every drifted state-path law with
$\mu_b=\mathcal D_T^b\mu_0$. The density is strictly positive and
finite $\mu_0$-almost surely, proving mutual absolute continuity.
\end{proof}

\begin{proof}[Proof of Theorem~\ref{thm:main-extensions}]
For (i), Theorem~\ref{thm:main-weak-uniqueness} gives the drift-free
input. For (ii), use
Proposition~\ref{thm:unbounded-coefficient-linf-localization}
when $\beta<1$, and the usual Lipschitz mild theory when $\beta=1$
\cite[Chapter~7]{DaPratoZabczyk2014}.
Global H\"older continuity gives at most linear growth, so in both
cases \eqref{eq:drift-diffusion-assumptions} holds.
On each finite interval apply
Proposition~\ref{prop:linear-growth-girsanov-transfer} with
$L_b=L_{b,T}$. It gives weak existence, uniqueness of the continuous
state-path law and mutual absolute continuity with the drift-free law.
\end{proof}

\clearpage
\appendix
\raggedbottom
\section{Gaussian derivative weights and heat-span geometry}
\label{app:gaussian-geometry}
This appendix proves the Gaussian input estimates from
Sections~\ref{sec:multipoint-score-package}
and~\ref{sec:conditional-gaussian-inputs}, and collects the supporting
heat-span and general leverage lemmas. The final subsection gives
the Gaussian product differentiation identity used in the comparison proofs.

\subsection{Moments of Gaussian derivative weights and one spatial input}
\label{app:score-moments}
Use the notation and hypotheses of
Lemma~\ref{lem:gaussian-one-input-scores}.
\nopagebreak[4]

\begin{proof}[Proof of Lemma~\ref{lem:gaussian-one-input-scores}]

\emph{Moments of Gaussian derivative weights.}
Differentiating the Gaussian density on its covariance range gives
the two displayed derivative weight formulas. After whitening on that range, each derivative weight is a polynomial in a
standard Gaussian vector. The fixed-order Gaussian polynomial
moment bounds \cite{Janson1997} give the stated $L^p$ estimates.
Ellipticity keeps the covariance range fixed.

\emph{Integrated leverage.}
The unweighted integral is $\operatorname{tr}((\mathsf G^Z)^+\mathsf G^Z)\le q$.
For the time weight, split at $s_Z/2$. On $u\le s_Z/2$, all terminals
are active and heat contraction gives
\[
 \mathsf H_u^Z:=\int_{\mathbb T}Q_u^Z(\xi)Q_u^Z(\xi)^*\,d\xi
   \le \mathsf H_{s_Z/2}^Z\le\frac2{s_Z}\int_{s_Z/2}^{s_Z}\mathsf H_t^Z\,dt
   \le\frac2{s_Z}\mathsf G^Z.
\] Taking the trace after multiplication by $(\mathsf G^Z)^+$ bounds the
spatial leverage mass at each $u\le s_Z/2$ by $2q/s_Z$.
The two weighted pieces are therefore bounded by
\[
 \frac{2q}{s_Z}\int_0^{s_Z/2}u^{\gamma_{\mathrm{wt}}}\,du
       \le C_{q,\gamma_{\mathrm{wt}}}s_Z^{\gamma_{\mathrm{wt}}},
 \qquad
 (s_Z/2)^{\gamma_{\mathrm{wt}}}\int_{s_Z/2}^{h}
             \int_{\mathbb T}\mathsf{Lev}^Z(u,\xi)\,d\xi\,du
       \le C_qs_Z^{\gamma_{\mathrm{wt}}}.
\]
The first uses $\gamma_{\mathrm{wt}}>-1$ and the second uses
$\gamma_{\mathrm{wt}}<0$ together with the unweighted mass bound.

\emph{Spatial integration of the input.}
Use the explicit average $I_f(x)$ from
\eqref{eq:heat-averaged-covariance-input}, with $x_1=x$ and
$Z=((T_1,x),(T_2,x_2),\ldots,(T_q,x_q))$.
Integrate the other centers against the mass-one heat-kernel
product; this average is absent for $q=1$. The resulting positive kernel from the input point
$\xi$ to the output point $x_1=x$ depends only on $x-\xi$:
translating every center and $\xi$ by the same amount leaves both
the heat-kernel product and $Q_u^Z(\xi)^*(\mathsf G^Z)^+Q_u^Z(\xi)$ unchanged.
Its equal row and column masses are bounded, using the
weighted leverage estimate, by
$C_{q,\gamma_{\mathrm{wt}}}s_Z^{\gamma_{\mathrm{wt}}}$;
at $\gamma_{\mathrm{wt}}=0$ use the unweighted estimate.
Young's inequality at each time and Minkowski's inequality in time,
using $\|f_u\|_{L^2(\mathbb T)}\le
\mathsf A_{\mathrm{in}}u^{\gamma_{\mathrm{wt}}}$, prove
\eqref{eq:multipoint-weighted-input-L2}.
The constant in part~(iii) is therefore independent of every
averaging time $r_i>0$.
\end{proof}

\subsection{Physical heat directions}
\label{app:physical-directions}
Use the observation covariance and heat-kernel averages of
Lemma~\ref{lem:physical-heat-score}.

\begin{proof}[Proof of Lemma~\ref{lem:physical-heat-score}]
\emph{The common method: realize the mean change by a heat forcing.}
The matrix $\mathsf G^Z$ defined before
Lemma~\ref{lem:gaussian-one-input-scores} is the covariance of unit
additive noise. If a deterministic forcing $f$ produces a function
$v$ satisfying $(\partial_t-\Delta)v=f$, $v_0=0$, and
$(v_{T_i}(x_i))_{i=1}^q=\mathbf b$, then the minimum-energy formula gives
\[
 \mathbf b^*(\Sigma^Z)^+\mathbf b
 \le C\mathbf b^*(\mathsf G^Z)^+\mathbf b
 \le C\int_0^{h}\|f_t\|_{L^2(\mathbb T)}^2\,dt.
\]
The observation map
$f\mapsto(\int_0^{T_i}S_{T_i-t}f_t(x_i)\,dt)_i$ has Gram
$\mathsf G^Z$ and minimum input energy
$\mathbf b^*(\mathsf G^Z)^+\mathbf b$.
Ellipticity gives the first inequality on the common range.

\emph{(i) A point input.}
Put $b_t(y)=p_t(y-\xi)$. For each observation choose a smooth cutoff
$0\le\eta_i(t,y)\le1$ equal to one at $(T_i,x_i)$, supported where
$T_i/2<t<3T_i/2$ and
$\operatorname{dist}_{\mathbb T}(y,x_i)\le C\sqrt{T_i}$.
Its time derivative and second spatial derivative are bounded by
$C/T_i$, and its first spatial derivative by $C/\sqrt{T_i}$.
Set $\eta=1-\prod_{i=1}^q(1-\eta_i)$ and $v=\eta b$.
Then $v_0=0$ and $v_{T_i}(x_i)=p_{T_i}(x_i-\xi)$.
Since $b$ solves the heat equation, its forcing is
\[
 f=(\partial_t\eta-\Delta\eta)b-2\partial_y\eta\,\partial_yb.
\]
The product rule and Cauchy--Schwarz give
\[
 |f|^2\le C_q\sum_i\mathbf1_{\operatorname{supp}\eta_i}
       (T_i^{-2}|b|^2+T_i^{-1}|\partial_yb|^2).
\]
On the $i$th support, the Gaussian heat-kernel estimates give
\[
 |b_t(y)|^2+T_i|\partial_yb_t(y)|^2
 \le C_qT_i^{-1}
 e^{-c_q\operatorname{dist}_{\mathbb T}(x_i,\xi)^2/T_i}.
\]
That support has time--space volume at most $CT_i^{3/2}$.
Consequently
\[
 \int_0^{h}\|f_t\|_{L^2(\mathbb T)}^2\,dt
 \le C_q\sum_iT_i^{-3/2}
 e^{-c_q\operatorname{dist}_{\mathbb T}(x_i,\xi)^2/T_i}.
\]
Take square roots in the minimum-energy bound and use
$T_i^{-3/4}e^{-c_q\operatorname{dist}_{\mathbb T}(x_i,\xi)^2/(2T_i)}
\le C_qT_i^{-1/4}p_{C_qT_i}(x_i-\xi)$ to obtain
\eqref{eq:multipoint-first-score-envelope}.

\emph{(ii) An $L^2$ initial input.}
Use cutoffs of a common spatial radius $C\sqrt{s_Z}$ and time width
$s_Z/2$, centered at $(T_i,x_i)$. They are supported in $t\ge s_Z/2$;
their derivative bounds are $C/s_Z$ and $C/\sqrt{s_Z}$.
With the same union cutoff $\eta$, set $v=\eta d_t$.
Now $(\partial_t-\Delta)d_t=-\Lambda d_t$, so
\[
 f=(\partial_t\eta-\Delta\eta)d_t
        -2\partial_y\eta\,\partial_yd_t-\Lambda\eta d_t.
\]
Average against the heat-kernel product in
\eqref{eq:multipoint-mean-energy}. Successive branch integration
gives uniform center marginals, so each spatial support has
averaged mass $C\sqrt{s_Z}$ and time length $Cs_Z$.
On these supports, heat contraction and the heat-gradient estimate give
\[
 \|d_t\|_{L^2(\mathbb T)}^2\le\|d\|_{L^2(\mathbb T)}^2,
 \qquad
 \|\partial_yd_t\|_{L^2(\mathbb T)}^2
       \le Cs_Z^{-1}\|d\|_{L^2(\mathbb T)}^2,
 \qquad
 \Lambda^2\|d_t\|_{L^2(\mathbb T)}^2
       \le Cs_Z^{-2}\|d\|_{L^2(\mathbb T)}^2.
\]
For the last inequality use $t\ge s_Z/2$ and
$\Lambda^2e^{-2\Lambda t}\le Cs_Z^{-2}$.
Thus each of the three forcing terms contributes at most
$C_q\sqrt{s_Z}\,s_Z\,s_Z^{-2}\|d\|_{L^2(\mathbb T)}^2
=C_qs_Z^{-1/2}\|d\|_{L^2(\mathbb T)}^2$ to the averaged energy.
The minimum-energy bound proves \eqref{eq:multipoint-mean-energy}.

\emph{(iii) A common constant shift.}
Use the cutoffs from (ii), with the constant heat solution $1$ in
place of $d_t$ and with no damping. Their union $v=\eta$ equals one
at every observation. For fixed centers its forcing has size at
most $C_q/s_Z$ and support volume at most $C_qs_Z^{3/2}$, so
\[
 \mathbf1^*(\Sigma^Z)^+\mathbf1\le C_qs_Z^{-1/2}.
\]
The construction also puts $\mathbf1$ in the covariance range.
For every $\mathbf a\in\mathbb R^q$, Cauchy--Schwarz on that range yields
$(\mathbf a^*\mathbf1)^2
\le(\mathbf a^*\Sigma^Z\mathbf a)(\mathbf1^*(\Sigma^Z)^+\mathbf1)$.
This is exactly \eqref{eq:multipoint-common-translation}.
\end{proof}

\subsection{Finite heat spans and collisions}
\label{app:heat-spans}
The following supporting lemma gives the heat-span geometry used in
Corollary~\ref{cor:heat-smoothed-observation-score} and in the
centered spatial heat-kernel bounds of Section~\ref{sec:centered-word-calculus}.
On $\mathbb R$, $p_t(z)=(4\pi t)^{-1/2}e^{-z^2/(4t)}$;
on $\mathbb T$, $p_t$ denotes its periodization.

\begin{lemma}[Finite heat spans]
\label{lem:finite-line-heat-span}
For each fixed $q$, arbitrary $t_i>0$, $x_i\in\mathbb R$, and
$f\in\operatorname{span}\{p_{t_i}(x_i-\cdot)\}_{i=1}^q$,
\begin{equation}
 |f(z)|^2\le C_q\|f\|_{L^2(\mathbb R)}^2\sum_{i=1}^q t_i^{-1/2}e^{-c|z-x_i|^2/t_i}.
 \label{eq:finite-line-span-envelope}
\end{equation}
\label{lem:finite-torus-heat-span}
For some $h_q^{\mathrm{span}}>0$, the periodic version holds whenever $t_i\le h_q^{\mathrm{span}}$:
\begin{equation}
 |f(z)|^2\le C_q\|f\|_{L^2(\mathbb T)}^2
       \sum_{i=1}^q t_i^{-1/2}e^{-c_q\operatorname{dist}(z,x_i)^2/t_i}.
 \label{eq:finite-torus-span-envelope}
\end{equation}
Both bounds include collisions and arbitrary scale ratios.
\end{lemma}

\begin{proof}[Proof of Lemma~\ref{lem:finite-line-heat-span}]

\emph{A normalized cluster.}
Scale the largest time to one and translate one center to zero.
First suppose the remaining times are bounded below and the centers
lie in a compact set. Up to scalar factors, the kernels are
$e^{-b_{2,i}z^2+b_{1,i}z}$, with
$b_{2,i}=1/(4t_i)\ge1/4$ and $b_{1,i}=x_i/(2t_i)$ in a compact set.
Their span satisfies a monic ODE of bounded order with coefficients
in a compact family, including at collisions. Indeed, an order-$m$
annihilator sends the next kernel to a polynomial of degree at most
$2m$ times that kernel; adjoining
$(\partial_z+2b_{2,i}z-b_{1,i})^{2m+1}$ annihilates it as well.
The order of this ODE is bounded in terms of $q$, and its
coefficients vary in a compact family. Its initial data at zero
are controlled by $\|f\|_{L^2(-1,1)}$: otherwise normalize those
data, take a convergent subsequence of equations and solutions,
and obtain a nonzero initial datum whose solution vanishes on
$(-1,1)$, contradicting uniqueness for the limiting ODE.

For completeness, each first-order factor can be solved by
\[
 (\partial_z+2bz-c_1)v=g
 \quad\Longrightarrow\quad
 v(z)=e^{-bz^2+c_1z}
       \left[v(0)+\int_0^ze^{by^2-c_1y}g(y)\,dy\right].
\]
Here $b\ge1/4$, and $b,c_1$ stay in bounded parameter sets.
Repeatedly using this formula for the finitely many factors
bounds the solution by a Gaussian with any fixed smaller decay
rate. Polynomial factors are absorbed by decreasing that rate.
Together with the initial-data bound this gives, for some fixed
$0<c<1/4$,
\[
 \sup_z e^{cz^2}|f(z)|\le C\|f\|_{L^2(-1,1)}
                         \le C\|f\|_{L^2(\mathbb R)}.
\]

\emph{Separated clusters.}
If the line estimate failed, choose a counterexample subsequence
along which all scale ratios and normalized center separations
converge. Group comparable scales with bounded normalized separation.
Each group has the preceding envelope. Groups separated in space
have vanishing overlap. For groups of scales $t\ll T$, the same
envelope gives, for functions $f_t,f_T$ in their respective spans,
\[
 |\langle f_t,f_T\rangle|
 \le\|f_t\|_1\|f_T\|_\infty
 \le C(t/T)^{1/4}\|f_t\|_{L^2(\mathbb R)}\|f_T\|_{L^2(\mathbb R)}.
\]
Thus the squared norm of the sum is at least half the sum of the
group squared norms for large indices. Each group envelope is
bounded by a constant times the Gaussian term of a largest-time
representative. Summing the at most $q$ group bounds contradicts
failure of \eqref{eq:finite-line-span-envelope}.

\emph{Periodization.}
A failing sequence of torus estimates with
$t_{\max}:=\max_i t_i\to0$ admits the same grouping. Each group
lifts to a short arc on the line. Its Gaussian envelope makes overlap
between distinct periodic images exponentially small in
$1/t_{\max}$, so its line and torus squared norms differ by a relative
$O(e^{-c/t_{\max}})$. Periodizing the envelopes and applying the same
cluster argument proves \eqref{eq:finite-torus-span-envelope}.
\end{proof}

\subsection{Total leverage mass}
\label{app:leverage}
The trace calculation in Corollary~\ref{cor:bounded-profile-score}
extends to arbitrary observation response rows. This general form
also covers stopped observations and reserved Gaussian fractions.

\begin{lemma}[Uniform total leverage mass]
\label{lem:total-leverage}
Fix $q\ge1$ observations and a time horizon $h_{\mathrm{obs}}>0$.
For a noise input at time $u$ and position $\xi$, let
$r(u,\xi)=(r_i(u,\xi))_{i=1}^q\in\mathbb R^q$ be its vector of responses
at the observations. In particular, for the heat-smoothed
observations of Corollary~\ref{cor:heat-smoothed-observation-score},
started at time zero with $\sigma_i\le h_{\mathrm{obs}}$,
\[
 r_i(u,\xi)=
 \begin{cases}
 p_{T_i-u}(x_i-\xi),&0\le u<\sigma_i,\\
 0,&u\ge\sigma_i.
 \end{cases}
\]
Thus observation $i$ receives no further noise after $\sigma_i$.
The estimate below applies to any measurable response vector $r$
with finite covariance.

Let $\kappa(u,\xi)$ be the local Gaussian noise variance, with a
uniform lower bound $\kappa(u,\xi)\ge\kappa_->0$.
For the references in Definition~\ref{def:coefficient-band-data},
this is $\kappa(u,\xi)=\bar{\mathfrak a}_u(\xi)=m_0^2+\bar a_u(\xi)$.
The resulting covariance matrix of the observation vector is
\[
 \Sigma_{\mathrm{obs}}:=\int_0^{h_{\mathrm{obs}}}\int_{\mathbb T}
              \kappa(u,\xi)r(u,\xi)r(u,\xi)^*\,d\xi\,du.
\]
Assume $\Sigma_{\mathrm{obs}}$ is finite, and let $\Sigma_{\mathrm{obs}}^+$ be its inverse on its range,
extended by zero on its null space. The quantity
$r(u,\xi)^*\Sigma_{\mathrm{obs}}^+r(u,\xi)$ is the covariance-normalized squared size of
the noise response, called its leverage. Its total mass satisfies
\begin{equation}
 \int_0^{h_{\mathrm{obs}}}\int_{\mathbb T}
                r(u,\xi)^*\Sigma_{\mathrm{obs}}^+r(u,\xi)\,d\xi\,du\le q/\kappa_-.
 \label{eq:total-leverage}
\end{equation}
In particular, for any bounded measurable covariance-profile
error $f=f(u,\xi)$,
\begin{equation}
 \int_0^{h_{\mathrm{obs}}}\int_{\mathbb T}
       |f(u,\xi)|\,r(u,\xi)^*\Sigma_{\mathrm{obs}}^+r(u,\xi)\,d\xi\,du
       \le\frac q{\kappa_-}\|f\|_{L^\infty([0,h_{\mathrm{obs}}]\times\mathbb T)}.
 \label{eq:bounded-input-leverage-general}
\end{equation}
\end{lemma}
\nopagebreak[4]

\begin{proof}[Proof of Lemma~\ref{lem:total-leverage}]
The covariance definition gives the exact identity
\[
 \int_0^{h_{\mathrm{obs}}}\int_{\mathbb T}
 \kappa(u,\xi)r(u,\xi)^*\Sigma_{\mathrm{obs}}^+r(u,\xi)\,d\xi\,du
 =\operatorname{tr}(\Sigma_{\mathrm{obs}}^+\Sigma_{\mathrm{obs}})
 =\operatorname{rank}\Sigma_{\mathrm{obs}}\le q.
\]
The integrand is nonnegative and $\kappa\ge\kappa_-$, proving
\eqref{eq:total-leverage}. Multiplying by $\|f\|_\infty$ proves
\eqref{eq:bounded-input-leverage-general}.
\end{proof}

\subsection{Differentiating Gaussian products}
\label{sec:full-product-identity}
The following supporting identity separates changes in the Gaussian
mean, covariance and the subtracted scalar functions. It is applied in
the comparison of Gaussian average terms in
Proposition~\ref{prop:all-order-predictable-pair} and in the covariance
interpolation for Proposition~\ref{prop:joint-first-pair}.

Fix an integer $n\ge0$, a linear subspace $E\subset\mathbb R^{n+1}$,
and a parameter $\theta$ in a real interval. Let $X_\theta$ be a
Gaussian vector supported on a fixed affine space parallel to $E$,
with absolutely continuous mean $\mu_\theta$ and covariance
$\Sigma_\theta$, positive on $E$. Assume $\mu_\theta^{\prime}\in E$
at almost every parameter. For bounded Borel functions
$f_i:\mathbb R\to\mathbb R$, $0\le i\le n$, and absolutely continuous scalar functions $b_i(\theta)$, $1\le i\le n$,
put
\[
 H_\theta(X)=f_0(X_0)\prod_{i=1}^n(f_i(X_i)-b_i(\theta)),
 \qquad I(\theta)=\mathbb EH_\theta(X_\theta).
\]
Here $X_i$ denotes the $i$th coordinate of $X$ and a prime denotes
 differentiation with respect to $\theta$. Define the mean and
quadratic Gaussian weights, with $G=X_\theta-\mu_\theta$ and inverses on $E$, by
\begin{align*}
 S_\mu&=\langle\Sigma_\theta^+G,\mu_\theta'\rangle,\\
 S_\Sigma&=\tfrac12\left(G^*\Sigma_\theta^+
       \Sigma_\theta'\Sigma_\theta^+G
                    -\operatorname{tr}(\Sigma_\theta^+\Sigma_\theta')\right).
\end{align*}

\begin{samepage}
\begin{lemma}[Three-leg Gaussian product identity]
\label{lem:three-leg-product}
The variation of $I(\theta)$ is the sum of a mean-derivative term, a
covariance-derivative term, and the derivatives $b_i^{\prime}(\theta)$
of the subtracted scalar functions.
At almost every parameter where the derivatives exist,
\begin{align}
 I'(\theta)={}&\mathbb E[H_\theta(X_\theta)S_\mu]
             +\mathbb E[H_\theta(X_\theta)S_\Sigma]\nonumber\\
 &-\sum_{i=1}^n b_i'(\theta)
       \mathbb E\left[f_0(X_{\theta,0})
               \prod_{\substack{1\le k\le n\\k\ne i}}(f_k(X_{\theta,k})-b_k(\theta))\right],
 \label{eq:three-leg-product}
\end{align}
If these three
legs have an integrable common majorant, the endpoint difference
is the integral of \eqref{eq:three-leg-product}.
\end{lemma}
\end{samepage}
\begin{proof}[Proof of Lemma~\ref{lem:three-leg-product}]
On the covariance range $E$, the Gaussian density $\varrho_\theta$ satisfies
$\partial_\theta\varrho_\theta=\varrho_\theta(S_\mu+S_\Sigma)$.
On compact parameter intervals its covariance is bounded below,
so differentiation under the integral is justified by a Gaussian
envelope times the integrable coefficient derivatives.
Apply the product rule to $H_\theta\varrho_\theta$: differentiating
$\varrho_\theta$ gives the first two terms of \eqref{eq:three-leg-product},
and differentiating each $-b_i(\theta)$ gives the corresponding
omitted-factor term. No state derivative of $f_i$ is taken, so
bounded Borel functions are allowed. Integration under the stated
common majorant gives the endpoint formula.
\end{proof}

\section{Elementary estimates for insertion times}
\label{app:time-kernels}
The root weight $k_\lambda$ and its mass $L_\lambda$ are defined
in Section~\ref{sec:gaussian-input-estimates}.

\subsection{Simplex and cube bounds}
\label{app:positive-gaps}
\begin{proof}[Proof of Lemma~\ref{lem:insertion-gap-integration}]
The geometric-mean inequality
$\min_i\delta_i\le\prod_i\delta_i^{1/m}$ gives the pointwise
bound in \eqref{eq:minimum-gap-cube-bound}. Integrating its product
over the cube gives $(m/\eta)^m h^\eta$.

For \eqref{eq:insertion-gap-simplex-bound}, the $n+1$ gaps
$\delta_1,\ldots,\delta_n,s$ sum to $h$; partition according to
their largest member, which is at least $h/(n+1)$.
If that member is $s$, bound $s^\sigma$ by $C_{n,\sigma}h^\sigma$
and integrate every other factor over $(0,h)$.
If it is a power gap $\delta_i$, use $s$ instead of $\delta_i$ as
an integration variable, with Jacobian one. Bound
$\delta_i^{\alpha_i-1}$ by $C_{n,\alpha_i}h^{\alpha_i-1}$ and
integrate $s^\sigma$ and the remaining factors over $(0,h)$.
Both cases give $Ch^{\gamma+\sigma}L_\lambda(h)^r$.
If the largest gap is a root gap $\delta_i$, the same substitution
uses
\[
 k_\lambda(\delta_i)\le(n+1)k_\lambda(h)
                      \le(n+1)h^{-1}L_\lambda(h).
\]
Here $t k_\lambda(t)$ is increasing and $k_\lambda(t)$ is decreasing.
The remaining root gaps supply $L_\lambda(h)^{r-1}$, the power
gaps supply $Ch^\gamma$, and $\int_0^h s^\sigma ds$
supplies $h^{\sigma+1}/(\sigma+1)$. This gives the same bound.
All integrations are finite because $\alpha_j>0$ and $\sigma>-1$.
The case $n=0$ is immediate.
\end{proof}

\subsection{Squared time integrals}
\label{app:clipped-convolutions}
\begin{proof}[Proof of Lemma~\ref{lem:outer-square-integrals}]
Monotonicity of $L_\lambda$ proves
\eqref{eq:positive-square-integral}.
Below $\lambda^{-4}$, direct integration of
$L_\lambda(w)^{2r}=(4\lambda)^{2r}w^{r/2}$ bounds the critical
integral by $2\cdot4^{2r}/r$. Above that scale,
$dL_\lambda(w)=dw/w$, so integration of $L_\lambda^{2r}dL_\lambda$
proves \eqref{eq:root-critical-square-integral}.

For \eqref{eq:current-mean-square-integral}, put
$a_*=(\chi-2\rho)_+<1$.
It suffices to prove the bound for $0<\varepsilon<1-a_*$.
If $r\ge1$, put $\eta=(a_*+\varepsilon)/r\in(0,1)$.
The elementary bound $\min(1,z)\le z^\eta$ gives
\[
 k_\lambda(t)\le\lambda^\eta t^{\eta/4-1},\qquad
 L_\lambda(w)^{2r}\le C_{r,\eta}
          \lambda^{2(a_*+\varepsilon)}w^{(a_*+\varepsilon)/2}.
\]
The exponent in the resulting integral is
$b\beta/2+(2\rho-\chi)_+/2+\varepsilon/2-1>-1$.
Integrating and taking the square root gives the asserted bound,
with an extra factor $h^{\varepsilon/4}\le1$.
If $r=0$, then $b\ge1$ and
$b\beta+2\rho-\chi>0$. Direct power integration gives
$Ch^{b\beta/4+(2\rho-\chi)/4}$; when $\chi>2\rho$, use
$h^{-(\chi-2\rho)/4}\le\lambda^{\chi-2\rho}$.
This proves the same estimate.
For later use on initial windows, if $h\le\lambda^{-4}$ the
integral inside the square root in
\eqref{eq:current-mean-square-integral} is exactly
\[
 \frac{(4\lambda)^{2r}
 h^{b\beta/2+\rho-\chi/2+r/2}}
      {b\beta/2+\rho-\chi/2+r/2}.
\]
Its denominator is positive: use $\chi<1$ when $r\ge1$, and
$\chi<\beta+2\rho$ when $r=0$, $b\ge1$.
\end{proof}

\section{Auxiliary calculations for finite-word comparisons}
\label{app:finite-word-details}
Section~5 states the working estimates before its four main comparison
proofs. This appendix supplies the complete first-two-order expansion
examples, verifies the common centered-word kernel and proves the
physical-input, moment and truncation estimates used by those proofs.

\begin{samepage}
\subsection{The first two coefficient-expansion orders}
\label{sec:first-two-expansion-orders}
Fix $x\in\mathbb T$, a block $a_\lambda$ and a window $[t_0,t_0+h]$.
The physical process is $U_{t_0+q}$ at elapsed time $q$.
Use the $\theta=1$ reference of
Definition~\ref{def:predictable-mean-functional}, with mean
$r\mapsto S_rU_{t_0}$, profile $\bar a=\bar a^1$ and Gaussian convolution
$\mathcal P^Z=\mathcal P^{Z,1}$. We compute the initial average,
stochastic term and remainder at each of the first two orders.

\end{samepage}
\paragraph{\textbf{First order: the initial average, fluctuation and first remainder.}}
For the single target $Z_0=((h,x))$, the backward observation is
$Y_q^{Z_0,1}=(S_{h-q}U_{t_0+q})(x)$, with
$Y_0^{Z_0,1}=(S_hU_{t_0})(x)$ and $Y_h^{Z_0,1}=U_{t_0+h}(x)$.
The initial average with no covariance insertion is
\begin{equation}
 \mathsf P_{1,\lambda}
 =\bigl(\mathcal P^{Z_0}_{0,h}a_\lambda\bigr)((S_hU_{t_0})(x)).
 \label{eq:explicit-first-predictable}
\end{equation}
The stochastic term is $N_{1,\lambda}=N_{1,\lambda}^U(h,x)$
from \eqref{eq:explicit-first-martingale}. Indeed,
$F_q^U(x,z)=(\mathcal P^{Z_0}_{q,h}a_\lambda)(z)$, and the
noise in $Y^{Z_0,1}$ is
$(m_0S_{h-q}\,dW_q^{0,t_0}+S_{h-q}L_U(q)\,dB_q^{t_0})(x)$.
Thus $\mathsf P_{1,\lambda}$ averages the terminal coefficient band using the
reference Gaussian law fixed by the initial state, while $N_{1,\lambda}$
records the fluctuation along the actual backward observation.
\begin{samepage}
\noindent The first covariance correction is
\begin{equation}
 \begin{aligned}
 R_{1,\lambda}={}&\frac12\int_0^h\int_{\mathbb T}
 [a(U_{t_0+s_1}(\xi))-\bar a_{s_1}(\xi)]p_{h-s_1}(x-\xi)^2\\
 &\hspace{2em}\times
 \bigl(\partial_z^2\mathcal P^{Z_0}_{s_1,h}a_\lambda\bigr)
       ((S_{h-s_1}U_{t_0+s_1})(x))\,d\xi\,ds_1.
 \end{aligned}
 \label{eq:explicit-first-remainder}
\end{equation}
\end{samepage}
The first decomposition is therefore
\begin{equation}
 R_{0,\lambda}=a_\lambda(U_{t_0+h}(x))=\mathsf P_{1,\lambda}+N_{1,\lambda}+R_{1,\lambda}.
 \label{eq:explicit-first-decomposition}
\end{equation}

\paragraph{\textbf{Second order: decomposing the first remainder.}}
For each insertion point $(s_1,\xi)$, $0<s_1<h$, use the two targets
$Z_1=((h,x),(s_1,\xi))$ and the explicit two-variable test
\begin{equation}
 F_1^{s_1,\xi}(z,w):=
 \frac12[a(w)-\bar a_{s_1}(\xi)]p_{h-s_1}(x-\xi)^2
 \bigl(\partial_z^2\mathcal P^{Z_0}_{s_1,h}a_\lambda\bigr)(z).
 \label{eq:explicit-one-insertion-test}
\end{equation}
The superscript records the insertion point $(s_1,\xi)$;
this is $F_1^{U,\lambda}$ in \eqref{eq:actual-coupled-inserted-tests}
for the displayed history, so its block dependence is through $a_\lambda$.
Its observation vector is
$\mathbf Y_q^{Z_1}=((S_{h-q}U_{t_0+q})(x),(S_{s_1-q}U_{t_0+q})(\xi))$,
$0\le q\le s_1$, so its second coordinate at $q=s_1$ is $U_{t_0+s_1}(\xi)$.
In particular, $R_{1,\lambda}=\int_0^h\int_{\mathbb T}
F_1^{s_1,\xi}(\mathbf Y_{s_1}^{Z_1})\,d\xi\,ds_1$.
Applying the same decomposition to this two-coordinate observation gives
\begin{align}
 \mathsf P_{2,\lambda}
 &=\int_0^h\int_{\mathbb T}
 \bigl(\mathcal P^{Z_1}_{0,s_1}F_1^{s_1,\xi}\bigr)
                    (\mathbf Y_0^{Z_1})\,d\xi\,ds_1,
 \label{eq:explicit-second-predictable}\\
 N_{2,\lambda}
 &=\int_0^h\int_{\mathbb T}
 \left[\int_0^{s_1}
 \nabla\bigl(\mathcal P^{Z_1}_{q,s_1}F_1^{s_1,\xi}\bigr)
       (\mathbf Y_q^{Z_1})\cdot d\mathbf Y_q^{Z_1}\right]d\xi\,ds_1.
 \label{eq:explicit-second-martingale}
\end{align}
Here $Z_1$ depends on the outer variables $(s_1,\xi)$, and
$\mathcal P^{Z_1}$ smooths both state coordinates jointly, preserving
their reference covariance. Thus $\mathsf P_{2,\lambda}$ is the initial-past
Gaussian average of the one-insertion test, and $N_{2,\lambda}$ is its
integrated, conditionally centered fluctuation.

The second remainder has two covariance insertions. Writing the
second insertion as $(s_2,\eta)$, with $0<s_2<s_1<h$, it is explicitly
\begin{equation}
 \begin{aligned}
 R_{2,\lambda}={}&\frac12\int_0^h\int_{\mathbb T}\int_0^{s_1}\int_{\mathbb T}
 [a(U_{t_0+s_2}(\eta))-\bar a_{s_2}(\eta)]\\
 &\quad\times\left[
 \bigl(p_{h-s_2}(x-\eta)\partial_z
       +p_{s_1-s_2}(\xi-\eta)\partial_w\bigr)^2
 \bigl(\mathcal P^{Z_1}_{s_2,s_1}F_1^{s_1,\xi}\bigr)
 \right](\mathbf Y_{s_2}^{Z_1})\,d\eta\,ds_2\,d\xi\,ds_1.
 \end{aligned}
 \label{eq:explicit-second-remainder}
\end{equation}
Here the first insertion used the one-coordinate direction
$p_{h-s_1}(x-\xi)\partial_z$. The second insertion acts on the jointly
smoothed two-coordinate test, using both heat rows. Its derivative
operator expands as
\[
 \begin{aligned}
 &\bigl(p_{h-s_2}(x-\eta)\partial_z
       +p_{s_1-s_2}(\xi-\eta)\partial_w\bigr)^2\\
 &\quad=p_{h-s_2}(x-\eta)^2\partial_z^2
 +2p_{h-s_2}(x-\eta)p_{s_1-s_2}(\xi-\eta)\partial_z\partial_w
 +p_{s_1-s_2}(\xi-\eta)^2\partial_w^2.
 \end{aligned}
\]
The mixed derivative records the covariance between the two old
observations. The $w$ derivatives act on the Gaussian-smoothed
first leaf; the newly inserted discrepancy
$a(U_{t_0+s_2}(\eta))-\bar a_{s_2}(\eta)$ is multiplied afterwards and stays
outside those derivatives. Consequently the first two orders read
\begin{equation}
 \boxed{\begin{aligned}
 R_{1,\lambda}&=\mathsf P_{2,\lambda}+N_{2,\lambda}+R_{2,\lambda},\\
 a_\lambda(U_{t_0+h}(x))&=\mathsf P_{1,\lambda}+N_{1,\lambda}
                      +\mathsf P_{2,\lambda}+N_{2,\lambda}+R_{2,\lambda}.
 \end{aligned}}
 \label{eq:explicit-second-decomposition}
\end{equation}
Both $\mathsf P_{2,\lambda}$ and $N_{2,\lambda}$ contain one insertion;
$R_{2,\lambda}$ contains two and is left unexpanded.

\subsection{The common time kernel after centering}
\label{app:centered-word-kernel}
Use the histories and assignments of
Definition~\ref{def:centered-word-assignments} and the root
mixed-difference bound \eqref{eq:capped-root-mixed-difference}.
Keep the old quadratic weights on the constant Gaussian part of
\eqref{eq:interval-gaussian-derivative-weight}, with fixed variance
density $\kappa_*>0$. On interval $j$, $Q_{s_j}(\xi)$ denotes the
column of active heat rows $p_{T_i-s_j}(x_i-\xi)$, $0\le i<j$.

\begin{corollary}[Spatial envelopes for interval derivatives]
\label{cor:interval-derivative-spatial-envelopes}
Write $(T_i,x_i)=(s_i,\xi_i)$ for the observation targets, with
$(T_0,x_0)=(h,x)$, as in
Definition~\ref{def:centered-word-assignments}.
For $\ell_{\mathrm{sc}}$ successive Gaussian derivatives, use
$\delta_j=s_{j-1}-s_j$, $j=1,\ldots,\ell_{\mathrm{sc}}$.
Observation $0$ is the terminal root and observations
$1,\ldots,j-1$ are the inserted points available at step $j$.
Assume $0<h\le h_{\ell_{\mathrm{sc}}-1}$ with the threshold
in \eqref{eq:finite-order-window-threshold}.
For each fixed Gaussian moment exponent, the quadratic derivative
weight on interval $j$ and the final linear derivative weight at
$q=s_{\ell_{\mathrm{sc}}}$ have the respective spatial envelopes
\begin{align}
 &\frac{C_{\ell_{\mathrm{sc}}}}{\delta_j}\sum_{i=0}^{j-1}
                p_{C_{\ell_{\mathrm{sc}}}(T_i-s_j)}(x_i-x_j),
 \label{eq:all-order-second-score-envelope}\\
 &\frac{C_{\ell_{\mathrm{sc}}}}{\sqrt{\delta_{\ell_{\mathrm{sc}}}}}\sum_{i=0}^{\ell_{\mathrm{sc}}-1}
       (T_i-q)^{1/4}p_{C_{\ell_{\mathrm{sc}}}(T_i-q)}(x_i-\xi),\qquad q=s_{\ell_{\mathrm{sc}}}.
 \label{eq:all-order-first-score-envelope}
\end{align}
In the second line, $q=s_{\ell_{\mathrm{sc}}}$ is the stochastic
integration time and $\xi\in\mathbb T$ is the spatial noise-input
point; the old targets are $(T_i,x_i)=(s_i,\xi_i)$.
\end{corollary}
\begin{proof}[Proof of Corollary~\ref{cor:interval-derivative-spatial-envelopes}]
Apply Corollary~\ref{cor:heat-smoothed-observation-score} on each interval
$(s_j,s_{j-1})$, whose length is $\delta_j$, to the active
observation list $0,\ldots,j-1$. Together with the linear Gaussian
weight in \eqref{eq:gaussian-score-definitions}, it gives
\eqref{eq:all-order-first-score-envelope} on the final interval.
For a quadratic weight, its fixed Gaussian moment is bounded by
the square of the corresponding linear weight, by
\eqref{eq:multipoint-score-products}. Use
$(\sum_i b_i)^2\le j\sum_i b_i^2$ and the heat-kernel bound
$(t^{1/4}p_{Ct}(y))^2\le C'p_{C't}(y)$ for $0<t\le1$.
This proves \eqref{eq:all-order-second-score-envelope}.
\end{proof}

\begin{proof}[Proof of Lemma~\ref{lem:centered-word-time-kernel}]

\emph{Centering and spatial integration.}
Lemma~\ref{lem:consistent-current-interpolation} and
\eqref{eq:hybrid-controlled-future-mean} give conditional displacement
moments $C_ph^{1/4}$, also for partially Gaussianized values. Hence
an unassigned coefficient factor has moment $C_ph^{\beta/4}$.

Apply Gaussian integration by parts to the independent physical
increments on $\delta_j=s_{j-1}-s_j$. Each second derivative uses
the quadratic weight \eqref{eq:interval-gaussian-derivative-weight};
a terminal first derivative uses the linear weight in
\eqref{eq:gaussian-score-definitions} on the last interval.
The weights are fixed by the profile and history.
Their mean-zero identity \eqref{eq:procedure-score-centering}
and the discrete product rule assign each interval difference
to the root or an affected rough factor.

For the root use \eqref{eq:capped-root-mixed-difference}.

For each rough coefficient factor, collect the centered interval indices whose
differences are assigned to that factor into a block $I$.
A nonempty block $I$ has bound
$C_n\min_{j\in I}\delta_j^{\beta/4}$, since any mixed difference
is at most $2^{|I|-1}$ times any constituent first difference.
An unassigned factor contributes $h^{\beta/4}$, also after any subset of
increments is suppressed. Fixed sufficiently high moments and
H\"older control the products of coefficient factors.

Use the spatial heat-kernel bounds
\eqref{eq:all-order-second-score-envelope} and
\eqref{eq:all-order-first-score-envelope}.
Indeed the interval covariance dominates
$\kappa_*\delta_j\int_{\mathbb T}Q_{s_j}(\xi)Q_{s_j}(\xi)^*\,d\xi$
on the fixed Gaussian fraction; square the first-derivative heat-kernel bound and
use $t^{1/2}p_{Ct}^2\le Cp_{C't}$ for the quadratic Gaussian weight. Choosing a
parent $i<j$ in every summand creates a spatial heat tree;
these parent choices need not equal those of the centering tree.
Integrate branches off its root-to-final-input path using unit heat
mass; path times telescope to $h-q$. With a final linear Gaussian weight the
remaining envelope is
$C_{\ell_{\mathrm{sc}}}\prod_{j=1}^{\ell_{\mathrm{sc}}-1}\delta_j^{-1}\delta_{\ell_{\mathrm{sc}}}^{-1/2}(h-q)^{1/4}
 p_{C_{\ell_{\mathrm{sc}}}(h-q)}(x_0-\xi)$.
Its negative-norm Hilbert--Schmidt factor
$(h-q)^{\rho/2-1/4}$ cancels the quarter power. Actual noise is
bounded on the right; random coefficient moments are taken before
this positive-kernel estimate.

\emph{The resulting time kernels.}
For each root-assigned interval $j\in I_0$, Gaussian increment
moments replace its translation minimum by
$C_n\min(1,\lambda\delta_j^{1/4})$. Fixed higher moments of the
normalized derivative weights and H\"older's inequality combine
these bounds with the coefficient differences.

Multiplying by $\delta_j^{-1}$ gives $k_\lambda(\delta_j)$ for each
$j\in I_0$, with a single root amplitude $\lambda^{-\beta}$. For a rough
block $I$, use \eqref{eq:minimum-gap-cube-bound} pointwise:
\[
 \left(\min_{j\in I}\delta_j\right)^{\beta/4}
 \prod_{j\in I}\delta_j^{-1}
 \le\prod_{j\in I}\delta_j^{\beta/(4|I|)-1}.
\]
Thus the blocks have total degree $n_{\mathrm b}\beta/4$,
the unassigned factors give $h^{u\beta/4}$, and the newest
factor gives $s^{\beta/4}$. These prove
\eqref{eq:all-order-predictable-common-ledger}.
Lemma~\ref{lem:insertion-gap-integration}, with $\sigma=\beta/4$,
then gives
\begin{equation}
 C_n\lambda^{-\beta}h^{u\beta/4}
 h^{(n_{\mathrm b}+1)\beta/4}L_\lambda(h)^{n_{\mathrm r}}
 =C_n\lambda^{-\beta}h^{n\beta/4}L_\lambda(h)^{n_{\mathrm r}}.
 \label{eq:procedure-final-predictable-integral}
\end{equation}
This argument includes $n_{\mathrm b}=0$, with no power gaps.
Summing over the finitely many assignments and using
$L_\lambda(h)\le C\mathcal L(\lambda,h)$ proves
\eqref{eq:centered-Gaussian-average-size}.
Every gap exponent $\alpha_j=\beta/(4|I_{\pi(j)}|)$ is positive;
at fixed order, later small losses are chosen below these finitely
many margins.
\end{proof}

\subsection{Physical inputs, moments and truncation}
\label{app:physical-input-details}
Fix an order $j\ge2$, a stochastic integration time $q<h$,
the outer age $w=h-q$ and the gap of the last Gaussian derivative $0<\delta\le w$.
Use the initial difference $d=U_{t_0}-V_{t_0}$ and
$Q_\Lambda=(\Lambda-\Delta)^{-\rho}$.

\begin{proof}[Proof of Lemma~\ref{lem:comparison-spatial-inputs}]
\emph{Part~(i): the coefficient input.}
The estimate \eqref{eq:native-integrator-heat-l2} is the Frobenius
prefix estimate followed by heat summation.

\emph{Part~(ii): the solution input.}
We prove \eqref{eq:one-weight-heat-tree} by Plancherel before
integrating any additional heat branches. Multiplication by $d$
convolves Fourier coefficients, so reindexing by $\ell=k-j$ gives
\[
 \begin{aligned}
 \|Q_\Lambda^{1/2}S_t\mathsf M_dS_s\|_{\mathrm{HS}}^2
 =\sum_{\ell\in\mathbb Z}|\widehat d(\ell)|^2
       \sum_{k\in\mathbb Z}
       (\Lambda+4\pi^2k^2)^{-\rho}
       e^{-8\pi^2tk^2-8\pi^2s(k-\ell)^2}.
 \end{aligned}
\]
The inner sum is maximal at $\ell=0$: the two nonnegative
symmetric decreasing sequences have centered interval level
sets, whose overlap is maximal there. Hence it is at most
\[
 \sum_{k\in\mathbb Z}(\Lambda+4\pi^2k^2)^{-\rho}
                   e^{-8\pi^2(t+s)k^2}
 \le C_\rho(t+s)^{\rho-1/2},\qquad t+s\le1.
\]
Separate $k=0$ and compare the remaining sum with its finite
Gaussian integral, using $\rho<1/2$. Plancherel proves
\eqref{eq:one-weight-heat-tree}; additional heat branches have
unit mass and are integrated by Schur and Cauchy--Schwarz.

Taking the square root in \eqref{eq:multipoint-mean-energy},
with earliest remaining age at least $\delta$, gives the
$C_j\delta^{-1/4}$ mean-derivative cost after heat averaging.
The argument applies conditionally to the current input
$D_{t_0+q}$; heat contraction and time translation keep its global
norm at most $\|U-V\|_{\mathscr L_0^2(H)}$.

\emph{The compensated old mean.}
Use the controlled mean identity
\eqref{eq:hybrid-controlled-future-mean}.
The old remaining-mean direction
$(e^{-\Lambda T}-e^{-\Lambda q})S_TD_{t_0}$ is also kept. Start
its localized control after $q+\delta/2$. Its cutoff derivative and
damping derivative are bounded by $C\delta^{-1}e^{-\Lambda q}$,
using $\Lambda e^{-\Lambda(t-q)}\le C/\delta$ there. Its physical
derivative weight therefore has the same $C\delta^{-1/4}\|D_{t_0}\|_{L^2(\mathbb T)}$ cost,
uniformly in $\Lambda$.

\end{proof}

\begin{remark}[Moment interpolation for actual inputs]
\label{rem:actual-input-moment-interpolation}
We use the following estimate to combine an actual input with rough amplitudes. For a power $\|U-V\|_{\mathscr L_0^2(H)}^{\chi-\varepsilon}$, choose a small moment increment $\eta_{\mathrm{mom}}>0$ and
interpolate the global $L^{2+\eta_{\mathrm{mom}}}$ norm of $D$ between its $L^2$ norm and a bounded higher
moment, then apply H\"older. Such higher moments follow from bounded
raw noise and $D_0=0$.
\end{remark}

\begin{samepage}
\begin{lemma}[Truncation of a normalized last kernel]
Write $B_\Delta=L(a(U))-L(a_{n_{\mathrm{sm}}}(V))$ for the actual
residual-noise difference.\label{lem:normalized-kernel-truncation}
Use the right heat factor and kernel envelope
\eqref{eq:native-integrator-exact-factor} from Step~2 of the
higher-martingale proof. For a signed covariance input and a
linear $\mathsf d_{\mathrm{coef}}(U,V)$ row, use a normalized kernel
of the form
$\widetilde{\mathscr K}(x,y)=\mathsf A_{\mathrm{ker}} \delta^{-1/2}w^{1/4}p_{Cw}(x-y)\Xi(x,y)$,
where $\mathsf A_{\mathrm{ker}}\ge0$ is the deterministic product
of old-gap powers and root clips, and $\Xi(x,y)$ is the remaining
normalized random kernel amplitude.
The finite moments of the derivative weights and conditional rough-leaf moments, followed
by Gaussian Jensen/H\"older, give
$\mathbb E_{t_0}|\Xi_q(x,y)|^p\le C_{j,p}$ uniformly in $q,x,y$ after
the displayed normalization. The kernel at $q$ is predictable at
$t_0+q$; its moment bound is conditional at the older past $t_0$.
For every $\varepsilon>0$ and prescribed tail power
$p_{\mathrm{tail}}>0$, split at $|\Xi|=\lambda^\varepsilon$.
The bounded part obeys the pathwise signed trace with the additional
factor $\lambda^\varepsilon$. The remaining part contributes an
arbitrarily small deterministic $\lambda^{-p_{\mathrm{tail}}}$ tail
in the integrated comparison norm, at each fixed order $j\ge2$.
\end{lemma}
\end{samepage}
\begin{proof}[Proof of Lemma~\ref{lem:normalized-kernel-truncation}]
For any $a>0$,
\[
 \mathbb E_{t_0}\bigl[|\Xi|^2\mathbf1_{|\Xi|>\lambda^\varepsilon}\bigr]
 \le\lambda^{-a\varepsilon}\mathbb E_{t_0}|\Xi|^{2+a}.
\]
Choose the fixed moment order $2+a$ large enough for the prescribed
tail exponent. Split at $|\Xi|=\lambda^\varepsilon$, keeping
$S_{\delta/2}$ on the right of both pieces. The good kernel satisfies
the pathwise signed trace. For the bad kernel, discard that gain and use
$\|B_\Delta\|_{\rm op}\le C$, $Q_\Lambda\le I$ and
\[
 \mathbb E\|Q_\Lambda^{1/2}\widetilde{\mathscr K}_{\rm bad}
                   S_{\delta/2}B_\Delta\|_{\rm HS}^2
 \le C\mathbb E\|\widetilde{\mathscr K}_{\rm bad}\|_{\rm HS}^2
 \le C_{j,p_{\mathrm{tail}} }\mathsf A_{\mathrm{ker}} ^2\delta^{-1}\lambda^{-p_{\mathrm{tail}} }.
\]
Here $\|p_{Cw}\|_{L^2(\mathbb T)}^2\lesssim w^{-1/2}$
cancels $w^{1/2}$; use a sufficiently high moment in Markov's inequality.
The last-gap density is now $\delta^{-1/2}$, with positive degree
$1/2$. For $j\ge2$, old nonempty rough blocks or root clips give
the integrable square corners checked in
Section~\ref{sec:higher-martingale-main-proof}. Thus the bad part yields
an arbitrary deterministic $\lambda^{-p_{\mathrm{tail}} }$ tail after increasing the
chosen moment.  This tail is separate from the
coefficient-approximation error.
These procedures are applied after the positive one-input estimate,
and all old derivative weights preserve independent reserved fractions and centering.
\end{proof}

\section{Coefficient-error insertions and approximation limits}
\label{app:approximation-details}
We first define the error insertion and its time-gap bounds. The
following subsections prove those bounds, identify every replacement
word and justify the exact brackets and approximation limits used in
Section~5.

\begin{samepage}
\subsection{Coefficient errors near colliding observation times}
\label{sec:uniform-coefficient-diagonal}
Fix a finite depth $N_\beta$, a damping strength $\Lambda$ and
one physical solution $U$, relabelled on a window so that $U_t$
means $U_{t_0+t}$. The index $n_{\mathrm{sm}}$ changes
the mollified coefficient $a_{n_{\mathrm{sm}}}$, while $U$ stays fixed.
The goal is to show that inserting $a-a_{n_{\mathrm{sm}}}$ has
vanishing total contribution even when new and old times collide.

\end{samepage}
At a coefficient-insertion time $q$, the finite expansion on the
fixed physical solution $U$ with the smooth residual
$a_{n_{\mathrm{sm}}}$ gives the bracket decomposition
\begin{equation}
 a(U_q)-\bar a_q
   =[a_{n_{\mathrm{sm}}}(U_q)-\bar a_q]
      +[a(U_q)-a_{n_{\mathrm{sm}}}(U_q)].
 \label{eq:diagonal-exact-split}
\end{equation}
The second bracket is the coefficient-error insertion. Its amplitude
is at most $\epsilon_{n_{\mathrm{sm}}}:=
\|a-a_{n_{\mathrm{sm}}}\|_\infty$.
We need both its small amplitude and an integrable bound when the
new insertion approaches an old observation. The argument combines
a bound that preserves the coefficient error with a second bound
that uses centering to gain time regularity.

Fix a dyadic upper cutoff $\lambda_{\mathrm{cut}}\in\mathscr S_{\mathrm z}$:
in this approximation argument we use only the coefficient blocks
$a_\lambda$ with state-frequency labels $\lambda\le\lambda_{\mathrm{cut}}$.
Fix also the depth $N_\beta$ and the damping $\Lambda$.
For each fixed $\beta$, $\Lambda$ and $N_\beta$, choose
$\lambda_{\mathrm{cut}}$ large enough for the required tail accuracy
using \eqref{eq:bands-bernstein-tail}, and hold it fixed as
$n_{\mathrm{sm}}\to\infty$.
Constants in this argument depend on these fixed choices and the
common boundedness, ellipticity and H\"older bounds.
For an old history $h=s_0>s_1>\cdots>s_l>0$, write
\[
 \boldsymbol\delta=(\delta_1,\ldots,\delta_l),\qquad
 \delta_i=s_{i-1}-s_i>0,\quad i=1,\ldots,l.
\]
Thus $\delta_i$ is the length of the $i$th old time interval
$(s_i,s_{i-1})$. For $l=0$, the vector is empty and $s_l=s_0=h$.
For a new insertion at $0\le q<s_l$, write $\tau=s_l-q>0$
for its gap to the earliest old observation.
The new derivative weight contributes $\tau^{-1}$, and the old derivative weights
contribute $\prod_{i=1}^l\delta_i^{-1}$.

For one centered assignment as in
Definition~\ref{def:centered-word-assignments}, let $I_0$ contain
the root-assigned intervals and $\mathscr B_{\mathrm{rough}}$ the
nonempty interval sets assigned to inserted coefficient factors. Define
\[
 B(\boldsymbol\delta)
  =\left[\prod_{I\in\mathscr B_{\mathrm{rough}}}
                         \left(\min_{i\in I}\delta_i\right)^{\beta/4}\right]
    \begin{cases}
     \min\left(1,\lambda_{\mathrm{cut}}
                    (\min_{i\in I_0}\delta_i)^{1/4}\right),&I_0\ne\varnothing,\\
     1,&I_0=\varnothing.
    \end{cases}
\]
The product has one factor for each inserted coefficient receiving at least
one centered difference. Its value is one if
$\mathscr B_{\mathrm{rough}}=\varnothing$.
For $I_0\ne\varnothing$, the last factor records the centering gain
of the root: a Gaussian increment over a gap $\delta_i$ has size
$\delta_i^{1/4}$, and Bernstein's inequality supplies the frequency
factor $\lambda\le\lambda_{\mathrm{cut}}$.
Bounding a mixed difference through any one of its increments gives
the minimum over $i\in I_0$; the cap at one uses the bounded size
of the root.

\begin{lemma}[Fresh centering with half of the old gains]
\label{lem:uniform-fresh-error}
For roots $a_\lambda$ with $\lambda\le\lambda_{\mathrm{cut}}$ and each
of their finitely many centered assignments, the RMS integrand of
an insertion of $(a-a_{n_{\mathrm{sm}}})(U_q)$ is bounded by
\begin{equation}
 C_{\Lambda,\lambda_{\mathrm{cut}} ,N_\beta}\epsilon_{n_{\mathrm{sm}}}^{1/2}
       \tau^{\beta/8-1}
       \prod_{i=1}^l \delta_i^{-1}\,B(\boldsymbol\delta)^{1/2}.
 \label{eq:diagonal-fresh-cap}
\end{equation}
RMS means the $L^2(\Omega;L^2(\mathbb T))$ norm after integrating
the other spatial centers against the positive heat kernels;
the remaining spatial variable is the output point $x$. The assertion remains valid for the
varying companions $V^{n_{\mathrm{sm}},\Lambda}$.
\end{lemma}
The proof is given in Appendix~\ref{app:fresh-centering}.

The bound in \eqref{eq:diagonal-fresh-cap} has finite mass
at simultaneous old and fresh time collisions. The calculation,
including the minimum-gap factors, is
Lemma~\ref{lem:coefficient-error-gap-integrability}.

\begin{corollary}[Integrated coefficient-error insertions vanish]
\label{cor:integrated-coefficient-error}
Under Lemma~\ref{lem:uniform-fresh-error}, fix the frequency
cutoff, depth and damping. For an error insertion with $l$ old gaps,
let $I_{n_{\mathrm{sm}}}(q,\boldsymbol\delta)$ be its spatially
integrated, finite-band-summed integrand. Its integral over the
ordered time domain, including initial triangles, has
$L^2(\Omega;H)$ norm at most
$C_{\Lambda,\lambda_{\mathrm{cut}},N_\beta}
\epsilon_{n_{\mathrm{sm}}}^{1/2}$.
The same conclusion holds with negative Sobolev output norm.
\end{corollary}
\begin{proof}[Proof of Corollary~\ref{cor:integrated-coefficient-error}]
For a finite-variation insertion one must integrate this RMS
density before squaring. With $l$ old gaps, use the time domain
\[
 \mathcal H_l^{\mathrm{gap}}(h)
  =\bigl\{(q,\boldsymbol\delta)\in(0,h)^{l+1}:
                               q+\textstyle\sum_{i=1}^l \delta_i<h\bigr\},
 \qquad d\boldsymbol\delta=\prod_{i=1}^l d\delta_i.
\]
Here $q$ is the current insertion time, the youngest old observation
is at $h-\sum_{i=1}^l \delta_i$, and the fresh gap is
$\tau=h-q-\sum_{i=1}^l \delta_i>0$. Positive-gap cutoffs restrict this domain;
the collision estimates of Appendix~\ref{app:collision-integrals}
permit their removal. Let
$I_{n_{\mathrm{sm}}}(q,\boldsymbol\delta)$ denote an error integrand after
its spatial integrations and finite root-band sum. Its sample/root-space
norm is $L^2(\Omega;H)$. Minkowski gives
\begin{equation}
 \begin{aligned}
 \left\|\int_{\mathcal H_l^{\mathrm{gap}}(h)}
       I_{n_{\mathrm{sm}}}(q,\boldsymbol\delta)\,dq\,d\boldsymbol\delta
                    \right\|_{L^2(\Omega;H)}
 &\le\int_{\mathcal H_l^{\mathrm{gap}}(h)}
       \|I_{n_{\mathrm{sm}}}(q,\boldsymbol\delta)\|_{L^2(\Omega;H)}
                                                 \,dq\,d\boldsymbol\delta\\
 &\le C_{\Lambda,\lambda_{\mathrm{cut}},N_\beta}
                                      \epsilon_{n_{\mathrm{sm}}}^{1/2}.
 \end{aligned}
 \label{eq:diagonal-rms-first}
\end{equation}
Negative Sobolev output norms only improve
\eqref{eq:diagonal-rms-first}. The same positive convolution proof
works on initial triangles by extending the input by zero before
time zero.
The integrability needed here is verified in
Lemma~\ref{lem:coefficient-error-gap-integrability} below.
\end{proof}

Thus the full finite collection of $U$-side coefficient-error
insertions tends to zero in $L^2(\Omega;H)$. This
is the derivative-independent part of the approximation diagonal.
Finite-band root errors and profile errors at the same actual means
are treated by uniform convergence and Gaussian interpolation.
High root frequencies are removed separately by the uniform
$\lambda_{\mathrm{cut}} ^{-\beta}$ tail. The resulting bounds involve
$\|a_{n_{\mathrm{sm}}}-a\|_\infty$ and the common rough coefficient constants.

\subsection{Fresh centering and the two complementary bounds}
\label{app:fresh-centering}
Use the notation and hypotheses of
Lemma~\ref{lem:uniform-fresh-error}.
\nopagebreak[4]

\begin{proof}[Proof of Lemma~\ref{lem:uniform-fresh-error}]

Condition at the actual current sigma field. The value $(a-a_{n_{\mathrm{sm}}})(U_q)$
is fixed there. Reserve a fresh independent Gaussian fraction for
the new quadratic Gaussian weight, leaving every old derivative weight on its previously
fixed fraction. Its normalized moments are bounded; its physical
normalization costs $C\tau^{-1}$.

\emph{First bound: preserve the small coefficient error.}
Keep the old centerings and use
$\|a-a_{n_{\mathrm{sm}}}\|_\infty=\epsilon_{n_{\mathrm{sm}}}$.
The centered block calculus gives
\[
 C\epsilon_{n_{\mathrm{sm}}}
 \underbrace{\tau^{-1}}_{\text{fresh derivative weight}}
 \underbrace{\prod_{i=1}^l\delta_i^{-1}}_{\text{old derivative weights}}
 B(\boldsymbol\delta).
\]

\emph{Second bound: gain integrability in the fresh gap.}
Use uniform boundedness of the coefficient error and the mean-zero
fresh quadratic Gaussian weight. Subtract from its test the value with the fresh
Gaussian increment set to zero.
Every coefficient leaf has a common $C^\beta$ bound, every subtracted value
is bounded, and the finitely many root bands have a fixed
$C^\beta$ bound depending on $\lambda_{\mathrm{cut}} $. Telescoping the finite product
therefore bounds this increment by $C_{\lambda_{\mathrm{cut}} ,N_\beta}|G_\tau|^\beta$,
where $G_\tau$ has moments $C_p\tau^{1/4}$. Gaussian
H\"older with the fixed moments of the normalized derivative weights gives
\[
 C\tau^{\beta/4-1}\prod_{i=1}^l \delta_i^{-1}.
\]
This uses uniform companion moments instead of $B(\boldsymbol\delta)$;
both bounds are uniform in the actual means.

\emph{Combine the two bounds.}
The inequality $\min\{u,v\}\le\sqrt{uv}$ for $u,v\ge0$ gives
\eqref{eq:diagonal-fresh-cap}. The geometric mean preserves the
small factor $\epsilon_{n_{\mathrm{sm}}}^{1/2}$, gains the integrable
fresh-gap power $\tau^{\beta/8-1}$, and preserves half the old
gains through $B(\boldsymbol\delta)^{1/2}$. Positive spatial heat-row integration
and fixed higher moments justify the same bounds in RMS before
this interpolation. Bounded raw-noise moments and the common
H\"older constants are uniform in $n_{\mathrm{sm}}$, so the proof applies to
each fixed actual process $V^{n_{\mathrm{sm}},\Lambda}$.
\end{proof}

\subsection{Integrating old and fresh time gaps}
\label{app:collision-integrals}
\begin{lemma}[Integrability at simultaneous old and fresh collisions]
\label{lem:coefficient-error-gap-integrability}
At fixed $\lambda_{\mathrm{cut}},\Lambda,N_\beta$ and $h\le1$,
the time-gap majorant in \eqref{eq:diagonal-fresh-cap} is integrable
over the ordered insertion domain, with total mass at most
$C_{\Lambda,\lambda_{\mathrm{cut}},N_\beta}
\epsilon_{n_{\mathrm{sm}}}^{1/2}$.
\end{lemma}
\begin{proof}[Proof of Lemma~\ref{lem:coefficient-error-gap-integrability}]

We first integrate the old-gap factors, then the fresh-gap factor. For a rough
block of size $r$ and any $\gamma>0$,
\begin{equation}
 \int_{(0,h)^r}\prod_{i=1}^r\frac{d\delta_i}{\delta_i}
                      \min_i \delta_i^\gamma
       \le(r/\gamma)^r h^\gamma.
 \label{eq:diagonal-block-cube}
\end{equation}
This is \eqref{eq:minimum-gap-cube-bound}.
In \eqref{eq:diagonal-fresh-cap} the rough blocks have
$\gamma=\beta/8>0$. For a root block, its minimum-gap clip to the
power $1/2$ is at most
$\lambda_{\mathrm{cut}}^{1/2}(\min_i\delta_i)^{1/8}$, so the same
bound applies with $\gamma=1/8$. Here $\lambda_{\mathrm{cut}}$ is fixed.
Finally
\[
 \int_0^h\tau^{\beta/8-1}d\tau
                         =\frac8\beta h^{\beta/8}<\infty.
\]
Enlarging the actual ordered simplex to this product domain proves
an integrable bound at every simultaneous old/fresh collision.
\end{proof}

\subsection{Explicit coefficient-replacement words}
\label{app:coefficient-replacement-words}
We give the exact finite-sum definition of
$\operatorname{Err}_{n_{\mathrm{sm}},\lambda}$ used in Section~5.
\begin{definition}[Replacement tests and outputs at a fixed physical pair]
\label{def:replacement-tests-outputs}
Fix an integer depth $N_\beta\ge1$, the block $\lambda$, the
window $[t_0,t_0+h]$ and the physical pair
$(U,V^{n_{\mathrm{sm}},\Lambda})$ throughout the replacements.
Put $a^{[0]}=a$, $a^{[1]}=a_{n_{\mathrm{sm}}}$.
For $i\in\{U,V\}$, keep the mean $\phi^i$ of
Definition~\ref{def:coupled-backward-martingales} fixed, and let
$\kappa^{i,b}$ be the self-consistent profile constructed from
$a^{[b]}$ at this mean, $b\in\{0,1\}$.
Thus $\kappa^{i,0}=\bar a^i$.
All word integrals use the finite frequency and positive-gap
restrictions of Lemma~\ref{lem:exact-bracket-cutoff-limits}.

Fix a history $h=s_0>s_1>\cdots>s_q>0$, $\xi_0=x$ and
$Z_j=((s_k,\xi_k))_{k=0}^j$.
For $\vec b=(b_0,\ldots,b_q)\in\{0,1\}^{q+1}$ and
$b\in\{0,1\}$, set $F_0^{i,\vec b,b}=(a^{[b_0]})_\lambda$ and
\begin{equation}
 \begin{aligned}
 F_j^{i,\vec b,b}:={}&\frac12[a^{[b_j]}(z_j)-\kappa_{s_j}^{i,b}(\xi_j)]\\
 &\times
 \left(\sum_{k=0}^{j-1}p_{s_k-s_j}(\xi_k-\xi_j)\partial_{z_k}\right)^2
 \mathcal P^{Z_{j-1},\kappa^{i,b}}_{s_j,s_{j-1}}
 F_{j-1}^{i,\vec b,b},\qquad 1\le j\le q.
 \end{aligned}
 \label{eq:approximation-test-recursion}
\end{equation}
Here $\mathcal P^{Z,\kappa}_{r,s}$ has the covariance in
\eqref{eq:predictable-mean-covariance}, with $\bar a^\theta$
replaced by $\kappa$. Thus $b_0$ chooses the root, $b_j$ the
$j$th leaf coefficient, and $b$ all subtracted values and Gaussian covariances.
Define the three word outputs by
\begin{equation}
 \begin{aligned}
 \mathcal O^i_{\mathsf P,q}[\vec b,b]
 &:=\int_{\mathcal H_q(h)}
 (\mathcal P^{Z_q,\kappa^{i,b}}_{0,s_q}F_q^{i,\vec b,b})
                                  (\mathbf Y_0^{i,Z_q})\,d\mathbf s\,d\boldsymbol\xi,\\
 \mathcal O^i_{N,q}[\vec b,b]
 &:=\int_{\mathcal H_q(h)}\int_0^{s_q}
 \nabla(\mathcal P^{Z_q,\kappa^{i,b}}_{r,s_q}F_q^{i,\vec b,b})
       (\mathbf Y_r^{i,Z_q})\cdot d\mathbf M_r^{i,Z_q}
                                         \,d\mathbf s\,d\boldsymbol\xi,\\
 \mathcal O^i_{R,q}[\vec b,b]
 &:=\int_{\mathcal H_q(h)}F_q^{i,\vec b,b}
                         (\mathbf Y_{s_q}^{i,Z_q})\,d\mathbf s\,d\boldsymbol\xi.
 \end{aligned}
 \label{eq:approximation-word-outputs}
\end{equation}
The domain $\mathcal H_q(h)$ is the ordered time--space domain of
Definition~\ref{def:finite-coefficient-band-words}; at $q=0$ the
outer integral is absent and $s_0=h$.
For a target $(T,\xi)$ in $Z_q$, the corresponding component of
$\mathbf Y_r^{i,Z_q}$ is the backward field in
Definition~\ref{def:coupled-backward-martingales}(i), with $h$ replaced
by $T$ and evaluated at $\xi$. Its stochastic part is
$\mathbf M^{i,Z_q}$: the same component has differential
$[m_0S_{T-r}\,dW_r^{0,t_0}+S_{T-r}L_i(r)\,dB_r^{t_0}](\xi)$.
The physical fields and their noise operators stay fixed.
For the complete correction, use the common-noise reference
\[
 X_s^{i,b}:=\phi_s^i+m_0\int_0^sS_{s-r}\,dW_r^{0,t_0}
             +\int_0^sS_{s-r}L(\kappa_r^{i,b})\,dB_r^{t_0}.
\]
At $q=0$, with $\vec b=(b_0)$, set
\[
 \mathcal O^i_{\Gamma,0}[\vec b,b]:=\mathcal O^i_{N,0}[\vec b,b]
 -\left[(a^{[b_0]})_\lambda(X_h^{i,b})
       -\mathbb E_{t_0}[(a^{[b_0]})_\lambda(X_h^{i,b})]\right].
\]
All these outputs use the fixed block $\lambda$.

In the notation of \eqref{eq:actual-higher-martingale-formula},
the actual higher stochastic words are
$N_{j,\lambda}^U=\mathcal O^U_{N,j-1}[\mathbf0,0]$ and
$N_{j,\lambda}^V=\mathcal O^V_{N,j-1}[(0,1,\ldots,1),0]$,
where $\mathbf0$ denotes the zero tuple. Both use the root
$a_\lambda$ and the profile $\bar a^i$; the $V$ leaves use
$a_{n_{\mathrm{sm}}}$.

\end{definition}

\begin{definition}[Finite replacement-error sum]
\label{def:finite-replacement-error-sum}
The word types and depths included in the sum are exactly
\[
 \mathscr W_{N_\beta}
 :=\{(\mathsf P,q),(N,q):0\le q<N_\beta\}
   \cup\{(R,q):0\le q\le N_\beta\}\cup\{(\Gamma,0)\}.
\]
For $0\le\ell\le q$, write $\vec b^{\ell\leftarrow c}$ for
$\vec b$ with entry $\ell$ replaced by $c\in\{0,1\}$.
Use $\|\cdot\|=\|\cdot\|_{\mathscr L_h^2(H_\Lambda^{-\rho})}$.
Set
\begin{equation}
 \begin{aligned}
 \operatorname{Err}_{n_{\mathrm{sm}},\lambda}
 :={}&\sum_{i\in\{U,V\}}\sum_{(W,q)\in\mathscr W_{N_\beta}}
 \sum_{\vec b\in\{0,1\}^{q+1}}\\
 &\quad\times\Bigg[
 \sum_{b=0}^1\sum_{\ell=0}^q
 \big\|\mathcal O^i_{W,q}[\vec b^{\ell\leftarrow1},b]
               -\mathcal O^i_{W,q}[\vec b^{\ell\leftarrow0},b]\big\|\\
 &\qquad+\big\|\mathcal O^i_{W,q}[\vec b,1]-\mathcal O^i_{W,q}[\vec b,0]\big\|\Bigg].
 \end{aligned}
 \label{eq:bandwise-approximation-error}
\end{equation}
The inner $\ell=0$ term replaces the root by
$(a_{n_{\mathrm{sm}}}-a)_\lambda$; for $\ell\ge1$ it replaces
one leaf by $a_{n_{\mathrm{sm}}}-a$, keeping the subtracted value and
Gaussian covariance fixed. In particular, a fresh $U$-side input
$(a-a_{n_{\mathrm{sm}}})(U)$ is the negative of such a last-leaf
correction in $\mathcal O^U_{R,q}$. The last norm changes the profile,
hence its subtracted values and covariance, at the same mean.

\end{definition}

The dependence of this error on the fixed
$T,\rho,\Lambda,h,N_\beta$ and cutoffs is suppressed.
Lemma~\ref{lem:exact-bracket-cutoff-limits} proves its convergence
to zero as $\epsilon_{n_{\mathrm{sm}}}\to0$ at fixed cutoffs,
and specifies how the cutoff tails are removed.
The amplitude-truncation tail $\lambda^{-p_{\mathrm{tail}}}$ is
kept separately.

For clarity, set
$\vec b^{[\ell]}=(\underbrace{1,\ldots,1}_{\ell},0,\ldots,0)$,
$0\le\ell\le q+1$, for tuples of length $q+1$.
The exact telescoping identity is
\[
 \begin{aligned}
 \mathcal O^i_{W,q}[\mathbf1,1]-\mathcal O^i_{W,q}[\mathbf0,0]
 ={}&\sum_{\ell=0}^q
 \bigl(\mathcal O^i_{W,q}[\vec b^{[\ell+1]},0]
             -\mathcal O^i_{W,q}[\vec b^{[\ell]},0]\bigr)\\
 &+\mathcal O^i_{W,q}[\mathbf1,1]-\mathcal O^i_{W,q}[\mathbf1,0].
 \end{aligned}
\]
Here $\mathbf0,\mathbf1$ are the constant tuples.
Every difference on the right is a summand of
\eqref{eq:bandwise-approximation-error} before taking its norm.

\begin{definition}[Replacement errors for the joint expansion]
\label{def:joint-replacement-errors}
For the joint comparisons in Section~\ref{sec:single-species-joint-flow-recursion}, use
\eqref{eq:bandwise-approximation-error} with the explicitly enlarged
word index set
\begin{equation}
 \begin{aligned}
 \mathscr W_{N_\beta}^{\mathrm{joint}}
 :=\mathscr W_{N_\beta}
 &\cup\{(\mathsf P^{\mathrm{joint}},q),(N^{\mathrm{joint}},q):0\le q<N_\beta\}\\
 &\cup\{(R^{\mathrm{joint}},q):0\le q\le N_\beta\}
 \cup\{(\Gamma^{\mathrm{joint}},0)\}
 \cup\{(\mathfrak J,q):0\le q<N_\beta\}.
 \end{aligned}
 \label{eq:joint-approximation-word-index-set}
\end{equation}
The extra types are the joint words in
\eqref{eq:joint-V-finite-expansion}, their complete first correction,
and the marks $\mathfrak J_{q,\lambda,h}$.
On the $U$ side these joint words equal the ordinary $U$ words,
and the marks are zero. The complete correction is formed from
the joint first stochastic term by subtracting the centered
reference terminal test, exactly as in
Definition~\ref{def:coupled-backward-martingales}(iii).
For them the root, leaf and profile choices are precisely those of
\eqref{eq:approximation-test-recursion}; their joint observation
fields and the physical bracket defect
\eqref{eq:single-species-direct-defect} stay fixed.
The error sum is enlarged by nonnegative terms, so all Section~5
bounds remain valid with this same larger error.
\end{definition}

\subsection{Exact brackets and removal of the approximations}
\label{app:approximation-limits}
Use the fixed coupled pair and common coefficient bounds in
Lemma~\ref{lem:exact-bracket-cutoff-limits}.

\begin{proof}[Proof of Lemma~\ref{lem:exact-bracket-cutoff-limits}]

Extend each observation's noise row $\mathsf h_s$ by zero after
its stopping time. Bounded covariance and
$\|p_t\|_{L^2(\mathbb T)}^2\lesssim t^{-1/2}$ give
$\mathbb E\int_0^h\|\mathsf h_s\|_{\mathrm{HS}}^2ds<\infty$.
The noise-series construction of the It\^o integral
\cite[Chapter~4]{DaPratoZabczyk2014} therefore gives convergence
of the finite-coordinate integrals in $L^2$. Their brackets converge
in $L^1$ as well, since for orthogonal projections
$\mathsf{Proj}_d\to I$,
\[
 \|\mathsf h_s\mathsf h_s^*
      -\mathsf h_s\mathsf{Proj}_d\mathsf h_s^*\|_{\mathrm F}
 \le\|\mathsf h_s(I-\mathsf{Proj}_d)\|_{\mathrm{HS}}^2,
 \qquad
 \mathbb E\int_0^h\|\mathsf h_s(I-\mathsf{Proj}_d)\|_{\mathrm{HS}}^2ds\to0.
\]

Fix the gap cutoff $\eta>0$ and the depth. Regularize rows
and leaf tests at scale $\zeta$, with the physical brackets fixed.
Induction uses Gaussian smoothing for derivatives, fixed-order
moments for products and isometry for stochastic integrals.
The range and heat-span bounds cover coincident observations.
Send $\zeta\to0$ first. The gap and squared-integral bounds of
Lemmas~\ref{lem:insertion-gap-integration}--\ref{lem:outer-square-integrals}
then allow $\eta\downarrow0$ and justify stochastic Fubini,
including the initial time strips. All these bounds use only the
common boundedness, ellipticity and H\"older constants.

For coefficient replacement, put
$\mathfrak f_{n_{\mathrm{sm}}}=a(U)-a_{n_{\mathrm{sm}}}(V)$ and
$\mathfrak f_{n_{\mathrm{sm}},\lambda_{\mathrm{cut}}}
=a_{\le\lambda_{\mathrm{cut}}}(U)-(a_{n_{\mathrm{sm}}})_{\le\lambda_{\mathrm{cut}}}(V)$.
Uniform state-band tails give
\[
 \|\mathfrak f_{n_{\mathrm{sm}}}-\mathfrak f_{n_{\mathrm{sm}},\lambda_{\mathrm{cut}}}\|_{\mathscr L_0^2(H_\Lambda^{-\rho})}
 \le C_T\lambda_{\mathrm{cut}}^{-\beta},\qquad
 \mathsf d_{\mathrm{coef}}(U,V^{n_{\mathrm{sm}},\Lambda})\le
 \|\mathfrak f_{n_{\mathrm{sm}},\lambda_{\mathrm{cut}}}\|_{\mathscr L_0^2(H_\Lambda^{-\rho})}
 +C_T\lambda_{\mathrm{cut}}^{-\beta}+\sqrt{\varepsilon_{\mathrm{app}}},
\]
where $\mathsf d_{\mathrm{coef}}(U,V^{n_{\mathrm{sm}},\Lambda})$ is the
coefficient discrepancy \eqref{eq:coefficient-discrepancy-norm}
for the displayed companion.
At fixed $\lambda_{\mathrm{cut}}$ and $\eta$, Bernstein's estimate gives
$\|\partial_z^j(a_{n_{\mathrm{sm}}}-a)_\lambda\|_\infty
\le C_{\lambda_{\mathrm{cut}},j}\|a_{n_{\mathrm{sm}}}-a\|_\infty$.
For the fixed means and their interpolation,
profile contraction gives
$\|\kappa^{i,1}-\kappa^{i,0}\|_\infty
\le C\|a_{n_{\mathrm{sm}}}-a\|_\infty$.
Gaussian covariance interpolation and
\eqref{eq:bounded-input-leverage} control each covariance and
subtracted-value change by this error.
At fixed profile index, multilinearity of
\eqref{eq:approximation-test-recursion} telescopes each test
change into one root or leaf replaced by $a_{n_{\mathrm{sm}}}-a$.
The Bernstein bound controls the root, uniform convergence the
leaves at positive gaps, and the heat one-input estimate the
stochastic integrals. The fresh $U$ insertion uses
Lemma~\ref{lem:uniform-fresh-error} and
\eqref{eq:diagonal-rms-first}. All bounds are uniform in the
actual means and common coefficient class.
For the reference terminal term in $\mathcal O^i_{\Gamma,0}$, couple its
two references with the same noises. It\^o isometry and
\eqref{eq:haar-heat-signed} at $\rho=0$ give
\[
 \mathbb E_{t_0}\|X_h^{i,1}-X_h^{i,0}\|_{L^2(\mathbb T)}^2
 \le C\sqrt h\,\|\kappa^{i,1}-\kappa^{i,0}\|_\infty^2.
\]
The fixed-band derivative bound and conditional Jensen then
control both its terminal test and its conditional expectation.
The same bounds cover the finite joint corrections with
integrators and defects fixed. Hence the error sum has a common
vanishing modulus at fixed cutoffs.

For any prescribed tolerance, choose $\lambda_{\mathrm{cut}}$ large and
$\eta$ small using the uniform tails, then choose a finite
$n_{\mathrm{sm}}$ for the remaining replacement errors. This gives
the stated modulus and, with tolerance $\Lambda^{-p_{\mathrm{tail}}}$,
the deterministic diagonal used in Section 7. The amplitude-truncation
tail is estimated separately.
\end{proof}

\section{Stopped heat histories and Gaussian derivative estimates}
\label{app:stopped-heat-details}
Section~\ref{sec:stopped-heat-contraction-joint-defect} states the
heat-history and Gaussian derivative bounds used by the Section~6 comparisons.
This appendix proves the three supporting estimates. The definitions
of stopping and target ages, heat rows, and age groups remain with
their statements in the main text.

\subsection{Damping stopped heat histories}
\label{app:stopped-damping}
Use the notation and hypotheses of Lemma~\ref{lem:stopped-history-contraction}.

\begin{proof}[Proof of Lemma~\ref{lem:stopped-history-contraction}]

\emph{Damping the full history.}
Extend $F$ by zero before $q$. The semigroup identity gives
$G=(I-\Lambda\mathsf{Res}_\Lambda)F$ on $u\ge q$, where
$\mathsf{Res}_\Lambda F(u)=\int_u^\infty e^{-\Lambda(s-u)}S_{s-u}F(s)ds$.
With time frequency $\omega\in\mathbb R$ and spatial frequency
$k\in\mathbb Z$, its full-line time/spatial Fourier multiplier is
$(4\pi^2 k^2-i\omega)/(\Lambda+4\pi^2 k^2-i\omega)$, of modulus at most
one. Plancherel and restriction to $u\ge q$ prove
\eqref{eq:stopped-heat-damping-contraction}. Approximation in
time--space $L^2$ is valid since the heat-square singularity
$t^{-1/2}$ is integrable.

\emph{Recovering the current row.}
Before the first stop, the semigroup identity gives
$G(q)=e^{-\Lambda(t-q)}S_{t-q}G(t)$.
Heat contraction and $e^{-\Lambda(t-q)}\le1$ imply
$\|G(q)\|_{L^2(\mathbb T)}\le\|G(t)\|_{L^2(\mathbb T)}$.
Average over an interval of length $\delta/2$:
\[
 \|G(q)\|_{L^2(\mathbb T)}^2
 \le\frac2\delta\int_q^{q+\delta/2}\|G(t)\|_{L^2(\mathbb T)}^2\,dt
 \le\frac2\delta\|F\|_{L^2([q,\infty)\times\mathbb T)}^2.
\]
The same calculation with $\Lambda=0$ treats the undamped row.
For comparable ages, $\delta$ is comparable to $w$. Since $\Sigma\ge c\mathsf G$, these inequalities
compare the current rows with the observation covariance.

\emph{The extra gain in the defect row.}
For $\Lambda w\le1$, the resolvent identity also gives
$\|Q_{\mathrm{def}} ^*z\|_{L^2(\mathbb T)}\le\Lambda\int_q^{q+Cw}\|F(u)\|_{L^2(\mathbb T)}du
\le C\Lambda\sqrt w\|F\|_{L^2}$, proving the extra $r_w$.
For $\Lambda w\ge1$, subtract the two already bounded rows.
The joint row $FL_U-G(L_U-L_V)$ therefore has covariance at
most $C\mathsf G$; the additive channel gives the lower bound
$m_0^2\mathsf G$.
\end{proof}

\subsection{Fourier and physical localization}
\label{app:stopped-localization}
Use the notation and hypotheses of Lemma~\ref{lem:comparable-heat-span-Fourier}.

\begin{proof}[Proof of Lemma~\ref{lem:comparable-heat-span-Fourier}]

\emph{Fourier localization of the span.}
For sufficiently small $w$, an orthonormal basis $(\psi_l)_{l=1}^{\dim E}$ satisfies
\begin{equation}
 \sum_{l=1}^{\dim E}|\widehat\psi_l(\xi)|^2\le C_n\sqrt w\,e^{-c_nw\xi^2}.
 \label{eq:comparable-span-Fourier}
\end{equation}
On the torus the Fourier variable in the first bound is
$\xi\in\mathbb Z$.
For $w=1$ on the line and bounded centers, apply the normalized-cluster
argument in the proof of Lemma~\ref{lem:finite-line-heat-span} to
$\exp(-T_i\xi^2-i x_i\xi)$. The first-order factors are now
$\partial_\xi+2T_i\xi+i x_i$; their imaginary coefficients contribute
only phases to the integrating factors. The same Gaussian-weight
bound therefore gives
$|\widehat f(\xi)|^2\le Ce^{-c\xi^2}\|f\|_{L^2(\mathbb R)}^2$, uniformly at collisions. Unbounded-center
sequences split into bounded-separation clusters. Their physical
Gaussian envelopes make distinct clusters asymptotically orthogonal,
and translating a cluster changes only a Fourier phase. The
Christoffel variational identity and scaling prove
\eqref{eq:comparable-span-Fourier}.
For the torus, lift each small-time cluster to the line. Gaussian
tails make periodization an asymptotic isometry, and torus Fourier
coefficients sample the line transforms, giving the same bound.
To see how this gives the translated projection bound, put
$e_k(z)=e^{2\pi ikz}$ and use the complex Fourier representation of
$H$. Orthogonality in the common translation variable gives
\[
 \begin{aligned}
 \int_{\mathbb T}\|\mathsf{Proj}_{E,x}B\|_{\mathrm{HS}}^2\,dx
 &=\sum_{k\in\mathbb Z}\sum_l|\widehat\psi_l(k)|^2
                                      \|B^*e_k\|_{L^2(\mathbb T)}^2\\
 &\le C_n\sqrt w\sum_k e^{-c_nwk^2}
                                      \|B^*e_k\|_{L^2(\mathbb T)}^2.
 \end{aligned}
\]
The last sum is $\|S_{c_nw}B\|_{\mathrm{HS}}^2$ after
adjusting $c_n$, proving \eqref{eq:comparable-span-input-heat}
with the signed operator inside the norm.

\emph{The physical Gaussian derivative response row.}
To obtain a positive spatial envelope, factor $\Sigma^{-1/2}\mathsf{Damp}_\Lambda(q)Q=T\mathsf{Proj}_E$,
where $T:E\to\mathbb R^n$ has norm at most $C_n\delta^{-1/2}$ by the preceding
range estimate and the physical Christoffel envelope of the same
finite-span lemma. Finally
\[
 e^{-\Lambda(\sigma_i-q)}p_{T_i-q}
 =\big[e^{-\Lambda\delta/2}
 e^{-\Lambda(\sigma_i-q-\delta/2)}p_{T_i-q-\delta/2}\big]*p_{\delta/2}.
\]
Normalize the bracketed row using a reserved future covariance
fraction after $q+\delta/2$. Its remaining stops are at least
$\delta/2$ and its targets are comparable to the original ones.
The preceding range and physical bounds apply, proving both assertions.
\end{proof}

\subsection{Combining observation ages}
\label{app:stopped-age-groups}
Use the notation and hypotheses of Lemma~\ref{lem:all-age-score-frame}.

\begin{proof}[Proof of Lemma~\ref{lem:all-age-score-frame}]

\emph{Separation of the age groups.}
For distinct groups $g\ne g^{\prime}$, their histories satisfy
\begin{equation}
 |\langle F_g,F_{g'}\rangle|\le(\mathfrak r_{\mathrm{sep}} -1)^{-1/2}\|F_g\|_{L^2((0,\infty)\times\mathbb T)}\|F_{g'}\|_{L^2((0,\infty)\times\mathbb T)},
 \label{eq:age-cluster-canonical-correlation}
\end{equation}
\begin{equation}
                  \mathsf G\ge\tfrac12\operatorname{diag}(\mathsf G_g).
 \label{eq:age-cluster-covariance-frame}
\end{equation}

For successive groups let $\tau_{\mathrm{young}} =\max t_i$ in the younger and
$\tau =\min t_i$ in the older. The older history's spatial norm is
nondecreasing before $\tau $, by heat contraction. Hence
\[
 \int_0^{\tau_{\mathrm{young}}} \|F_{g'}(u)\|_{L^2(\mathbb T)}^2du
 \le\frac{\tau_{\mathrm{young}}}{\tau -\tau_{\mathrm{young}} }\int_{\tau_{\mathrm{young}}} ^{\tau} \|F_{g'}(u)\|_{L^2(\mathbb T)}^2du
 \le(\mathfrak r_{\mathrm{sep}} -1)^{-1}\|F_{g'}\|_{L^2((0,\infty)\times\mathbb T)}^2.
\]
Cauchy--Schwarz gives \eqref{eq:age-cluster-canonical-correlation};
summing cross terms yields
\[
 \left\|\sum_gF_g\right\|_{L^2}^2
 \ge\left[1-\frac{n-1}{\sqrt{\mathfrak r_{\mathrm{sep}}-1}}\right]
             \sum_g\|F_g\|_{L^2}^2
 \ge\frac12\sum_g\|F_g\|_{L^2}^2.
\]
The norms are over time and space. Within each group the age
ratio is at most $\mathfrak r_{\mathrm{sep}}^{n-1}$, so the
comparable-age estimates apply with constants depending on $n$.

\emph{The signed covariance input.}
Reserve the block-diagonal fraction for the new derivative weight. Factor its
signed row as $\Sigma_g^{-1/2}Q_{\mathrm{def},g}=T_g\mathsf{Proj}_{E_g}$ with
$\|T_g\|\le Cr(\tau _g)\tau _g^{-1/2}$ and apply the signed Fourier
projection and Haar heat estimate. The projection bound contributes $\tau_g^{1/4}$ after taking
square roots. The signed heat trace contributes
$\Lambda^{\rho/2}\tau_g^{-1/4}+\tau_g^{-1/4-\rho/2}$.
Multiplying these with the range cost $r(\tau_g)\tau_g^{-1/2}$,
the companion cost $C\tau_{g'}^{-1/2}$, and the output embedding
$\Lambda^{-\rho/2}$ gives exactly
\[
 C_nr(\tau_g)\tau_g^{-1/2}\tau_{g'}^{-1/2}
            [1+\Lambda^{-\rho/2}\tau_g^{-\rho/2}].
\]
The global signed input has norm $\mathsf d_{\mathrm{coef}}(U,V)$, and the remaining factors
have bound $\mathsf A_0$, proving \eqref{eq:all-age-one-input-score}.
Setting both group ages equal to $w$ gives
\eqref{eq:comparable-one-input-covariance-score}.
A second bounded $\mathcal B$ is an admissible companion.

\emph{The mean input.}
For the squared mean estimate, the range factorization costs
$C\tau^{-1}$ after squaring. Apply the translated projection bound
to the rank-one input with image $d$; its factor $\sqrt\tau$
leaves $C\tau^{-1/2}\|S_{c_n\tau}d\|_{L^2(\mathbb T)}^2$,
which is \eqref{eq:joint-native-mean-one-input}.
For the positive row, integrate the physical spatial heat-kernel bound against
$|d|$. Each heat kernel has mass one, so Jensen gives
\[
 \bigl[(p_t*|d|)(x)\bigr]^2\le(p_t*|d|^2)(x).
\]
Thus squaring preserves the spatial $L^2$ input inside the original
heat-tree integrals. Independent Gaussian fractions for old derivative weights, this new
derivative weight and annular extraction preserve all stated amplitudes.
\end{proof}

\end{document}